\documentclass[11pt, twoside]{amsart}

\usepackage[T1]{fontenc}
\usepackage{amssymb,amsmath,amstext}

\newtheorem{theor}{Theorem}[section]
\newtheorem{lem}[theor]{Lemma} 
\newtheorem{defin}[theor]{Definition}

\newtheorem{prop}[theor]{Proposition} 

\newtheorem{notation}[theor]{Notation}
\newtheorem{exam}[theor]{Example}

\newtheorem{cor}[theor]{Corollary}
\newtheorem{rem}[theor]{Remark}

\newtheorem{assump}[theor]{Assumption}
\newtheorem{observation}[theor]{Observation}
\newtheorem{termin}[theor]{Terminology}
\newtheorem{assumptionsandterminology}[theor]{Assumptions and terminology}

\newcommand{\cl}{\mathrm{cl}}

\newcommand{\mbb}{\mathbb}

\newcommand{\es}{\emptyset}

\newcommand{\uhr}{\upharpoonright}

\newcommand{\nts}{\negthickspace}
\newcommand{\uhrc}{\nts \upharpoonright \nts}

\newcommand{\trir}{\triangleright}
\newcommand{\mult}{\mathrm{mult}}
\newcommand{\smult}{\mathrm{smult}}
\newcommand{\perm}{\mathrm{perm}}

\newcommand{\mcA}{\mathcal{A}}
\newcommand{\mcB}{\mathcal{B}}
\newcommand{\mcC}{\mathcal{C}}
\newcommand{\mcD}{\mathcal{D}}
\newcommand{\mcE}{\mathcal{E}}
\newcommand{\mcF}{\mathcal{F}}
\newcommand{\mcG}{\mathcal{G}}
\newcommand{\mcH}{\mathcal{H}}

\newcommand{\mcK}{\mathcal{K}}

\newcommand{\mcM}{\mathcal{M}}
\newcommand{\mcN}{\mathcal{N}}

\newcommand{\mcP}{\mathcal{P}}

\newcommand{\mcS}{\mathcal{S}}

\newcommand{\mcU}{\mathcal{U}}
\newcommand{\mcV}{\mathcal{V}}

\newcommand{\mbC}{\mathbf{C}}

\newcommand{\mbE}{\mathbf{E}}
\newcommand{\mbF}{\mathbf{F}}
\newcommand{\mbG}{\mathbf{G}}

\newcommand{\mbP}{\mathbf{P}}

\newcommand{\mbK}{\mathbf{K}}
\newcommand{\mbS}{\mathbf{S}}

\newcommand{\mbX}{\mathbf{X}}
\newcommand{\mbY}{\mathbf{Y}}
\newcommand{\mbZ}{\mathbf{Z}}

\newcommand{\mbSK}{\mathbf{SK}}

\newcommand{\mbbP}{\mathbb{P}}
\newcommand{\mbbN}{\mathbb{N}}
\newcommand{\mbbR}{\mathbb{R}}

\newcommand{\rng}{\mathrm{rng}}

\title[Asymptotic probabilities of extension properties]
{Asymptotic probabilities of extension properties\\
and random $l$-colourable structures}

\author{Vera Koponen}

\address{Vera Koponen, Department of Mathematics, Uppsala University, Box 480,
75106 Uppsala, Sweden.}

\email{vera@math.uu.se}

\begin{document}

\begin{abstract} 
\noindent 
We consider a set $\mbK = \bigcup_{n \in \mbbN}\mbK_n$ of {\em finite} structures such that
all members of $\mbK_n$ have the same universe, the cardinality of which approaches
$\infty$ as $n\to\infty$. Each structure in $\mbK$ may
have a nontrivial underlying pregeometry and on each $\mbK_n$ we consider
a probability measure, either the uniform measure, 
or what we call the {\em dimension conditional measure}.
The main questions are: What conditions imply that for every extension axiom
$\varphi$, compatible with the defining properties of $\mbK$, the probability
that $\varphi$ is true in a member of $\mbK_n$ approaches 1 as $n \to \infty$?
And what conditions imply that this is not the case,
possibly in the strong sense that the mentioned probability approaches 0 for some $\varphi$?

If each $\mbK_n$ is the set of structures with universe $\{1, \ldots, n\}$, 
in a fixed relational language, 
in which certain ``forbidden'' structures cannot be weakly embedded and $\mbK$ 
has the disjoint amalgamation property, then there is a condition 
(concerning the set of forbidden structures) which,
if we consider the uniform measure,
gives a dichotomy; i.e. the condition holds if and only if the answer to the first question is `yes'.
In general, we do not obtain a dichotomy, but we do obtain a condition guaranteeing
that the answer is `yes' for the first question, as well as a condition guaranteeing that
the answer is `no'; and we give examples showing that in the gap between these conditions
the answer may be either `yes' or `no'. 
This analysis is made for both the uniform measure and for the dimension conditional measure.
The later measure has closer relation to random generation of structures and
is more ``generous'' with respect to satisfiability of extension axioms.

Random $l$-colour{\em ed} structures fall naturally into the framework
discussed so far, but random $l$-colour{\em able} structures need further
considerations. It is not the case that every extension axiom compatible with
the class of $l$-colourable structures almost surely holds in an $l$-colourable structure.
But a more restricted set of extension axioms turns out to hold almost surely,
which allows us to prove a zero-one law for random $l$-colourable structures,
using a probability measure which is derived from the dimension conditional measure,
and, after further combinatorial considerations, also for the uniform probability measure.
\\
{\Small {\em Keywords:} Model theory, finite structure, asymptotic probability, 
extension axiom, zero-one law, colouring.}
\end{abstract}

\maketitle
\let\thefootnote\relax\footnote{This work was carried out in part while the author was 
a visiting researcher at Institut Mittag-Leffler during the autumn 2009.}

\begin{center}{\large Contents}\end{center}
{\small
\contentsline{section}{\ref{introduction}. Introduction}{\pageref{introduction}}
\contentsline{section}{\ref{preliminaries}. Preliminaries}{\pageref{preliminaries}}
\contentsline{section}{\ref{forbidden structures}. 
Permitted structures and substitutions}{\pageref{forbidden structures}}
\contentsline{section}{\ref{examples concerning notion of admittance}. 
Examples}{\pageref{examples concerning notion of admittance}}
\contentsline{section}{\ref{proving upper limit on the probability of an extension axiom}.
Proof of Theorem~\ref{upper limit on the probability of an extension axiom}}
{\pageref{proving upper limit on the probability of an extension axiom}}
\contentsline{section}{\ref{conditional measures}. 
Conditional probability measures}{\pageref{conditional measures}}
\contentsline{section}{\ref{pregeometries}. 
Underlying pregeometries}{\pageref{pregeometries}}
\contentsline{section}{\ref{proofs of main results}.
Proofs of Theorems~\ref{0-1 law for pregeometries}, 
\ref{third part of 0-1 law for pregeometries} and~\ref{0-law for extension axioms in pregeometries}}
{\pageref{proofs of main results}}
\contentsline{section}{\ref{l-colourable structures}.
Random $l$-colourable structures}{\pageref{l-colourable structures}}
\contentsline{section}{\ref{the uniform measure and l-colourable structures}.
The uniform probability measure and the typical distribution of colours}
{\pageref{the uniform measure and l-colourable structures}}
\contentsline{section}{References}{\pageref{the bibliography}}
}

\section{Introduction}\label{introduction}

\noindent
Extension axioms have been used as a technical tool for proving
zero-one laws \cite{EF, Fag, Hod, Gleb, KPR},
but they also have other implications which will be 
explained below. 
Extension axioms, by their definition, express 
possibilities of extending a structure that are compatible (or ``consistent'') with 
the definition of a given 
class of structures under consideration.
So given a structure $\mcM$ from this class, 
the set of extension axioms which are satisfied in $\mcM$ tells which
possibilities of extending substructures of $\mcM$,
in ways compatible with the context,
are actually realized in the particular structure $\mcM$.
Thus, extension axioms have a combinatorial interest of their own.

If we consider the class of all finite $L$-structures, where $L$ is a language
with finite relational vocabulary, then it follows from the proof of
the zero-one law (as presented in \cite{EF, Fag, Hod}) that, 
for every extension axiom, almost all sufficiently large finite
$L$-structures satisfy it.
Hence the interesting case to study is the case when there are some
restrictions on the structures under consideration. 
For example, we could restrict ourselves to the class of finite structures in which
some particular structure cannot be (weakly) embedded; for instance,
the class of triangle-free graphs. 
Specific classes of this kind have been studied extensively.
An overview with emphasis on graphs and partial orders is found in \cite{PST};
see also \cite{KPR, PS} and recent results \cite{ABBM, BBS08, BBS10}.
An overview with focus on zero-one laws is found in \cite{Win};
it takes up, among other things, the number theoretic approach to zero-one laws which was
first developed by K. Compton, and which is the subject of a book by S. Burris \cite{Bur-book}.
However, none of the previously published research focuses specifically on searching for ``dividing lines'' for
asymptotic probabilities of extension properties in a general model theoretic setting.
That is the purpose of this article, as well as deriving consequences such as zero-one laws
and, finally, studying random $l$-colourable structures.

The general framework of this article is the following.
For some language $L$, $\mbK = \bigcup_{n \in \mbbN}\mbK_n$ is a set of 
{\em finite} $L$-structures such that
all members of $\mbK_n$ have the same universe; often an initial 
segment of $\{1, 2, 3, \ldots \}$.
In addition, each $\mcM \in \mbK$ may have a nontrivial closure operator which
makes it into a pregeometry; in this case, the closure operator is uniformly
definable on all members of $\mbK$ in the sense described in 
Definition~\ref{definition of a structure being a pregeometry}.
An important special case is when the closure (and pregeometry) is {\em trivial}, by which we mean
that every subset of any structure from $\mbK$ is closed.
If $P$ is a property, then the expression that `{\em a member of $\mbK$ almost surely
has property $P$}' is shorthand for saying that, with respect to some probability measure
$\mu_n$ on $\mbK_n$, the probability that $\mcM \in \mbK_n$ has $P$ approaches 1 as $n \to \infty$.
If $\mu_n(\mcM) = 1/|\mbK_n|$ for all $n$ and all $\mcM \in \mbK_n$ (the uniform probability measure), then
we may instead say that `{\em almost all $\mcM \in \mbK$ have $P$}'.
By a {\em zero-one law for $\mbK$} we mean that for every $L$-sentence $\varphi$, either it or its negation
almost surely holds in $\mbK$.

Suppose that $\mcA \subseteq \mcB \subseteq \mcM \in \mbK$ and that
$A$ and $B$ are closed subsets of $M$. 
For a structure $\mcN$, the 
{\em $\mcB/\mcA$-extension axiom}
holds for $\mcN$ if for every embedding
$\tau$ of $\mcA$ into $\mcN$ there
is an embedding $\pi$ of $\mcB$ into $\mcN$ which
extends $\tau$.
If the dimension of $B$ is at most $k+1$, then we call it a {\em $k$-extension axiom} of $\mbK$.
If the closure is trivial then dimension is the same as cardinality.

If $L$ has no constant symbols we allow the universe of 
$\mcA$ to be empty, and in this case the $\mcB/\mcA$-extension axiom expresses
that there exists a copy of $\mcB$ in the ambient structure.
Hence, if $\mcM$ satisfies all $k$-extension axioms, then every
$\mcB \in \mbK$ of dimension at most $k+1$ can be embedded into $\mcM$;
in this case one may say that $\mcM$ is `$(k+1)$-universal for $\mbK$'.
By involving pebble games \cite{Imm, Poi} it follows that 
if $L$ is relational and $\mcM \in \mbK$ satisfies all $k$-extension axioms of $\mbK$,
then $\mcM$ has the following `homogeneity property, up to $k$-variable expressibility':
Whenever $\bar{a}, \bar{a}'$ are tuples of elements and there is an
isomorphism from the closure of $\bar{a}$ to the closure of $\bar{a}'$
which sends $a_i$ to $a'_i$, then $\bar{a}$ and $\bar{a}'$ satisfy exactly
the same formulas in which at most $k$ distinct variables occur.

If the class $\mbK^*$ of all structures which can be embedded into some member of
$\mbK$ has (up to taking isomorphic copies) the joint embedding property and
the amalgamation property, then a structure $\mcM$ exists which satisfies
all $k$-extension axioms of $\mbK$ for every $k \in \mbbN$;
because we can let $\mcM$ be the so-called Fra\"{i}ss\'{e} limit of $\mbK^*$.
However, if $\mbK$ contains arbitrarily large (finite) structures,
then the Fra\"{i}ss\'{e} limit of $\mbK^*$ is infinite.
The question whether, for every $k \in \mbbN$, there exists a {\em finite} $\mcM \in \mbK$ which
satisfies every $k$-extension axiom of $\mbK$ may be hard.
For instance, the problem \cite{Che92} whether there is a finite triangle-free graph
which satisfies every 4-extension axiom is still open.
By using the fact that the proportion of triangle-free graphs with vertices $1, \ldots, n$
which are bipartite approaches 1 as $n$ approaches infinity \cite{EKR, KPR},
it is straightforward to derive that the proportion of all triangle-free graphs
with vertices $1, \ldots, n$ which satisfy all 3-extension axioms approaches 0 as 
$n$ approaches infinity.
The main results in Sections~\ref{forbidden structures}~--~\ref{pregeometries} 
are concerned with the question of
when, for some $k$ and large enough $n$, it is usual (or unusual), in senses to be made precise, 
that structures in $\mbK_n$ satisfy all $k$-extension axioms.

For the moment, assume that, for each $n$, $\mu_n$ is a probability 
measure on $\mbK_n$.
Let $Th_{\mu}(\mbK)$ be the set of sentences $\varphi$ such that
the $\mu_n$-probability that $\varphi$ is true in
a member of $\mbK_n$ approaches 1 as $n$ approaches infinity.
Also assume that $\mbK^*$, as defined above, satisfies the joint embedding and
amalgamation properties and let $Th_{\mbb{F}}(\mbK)$ be the complete theory
of the Fra\"{i}ss\'{e} limit of $\mbK^*$. 
If, moreover, the closure is trivial on all members of $\mbK$,
it is straightforward to see that $Th_{\mu}(\mbK) = Th_{\mbb{F}}(\mbK)$
if and only if $Th_{\mu}(\mbK)$ contains all extension axioms of $\mbK$.
(We can get rid of the assumption that the closure is trivial
if we assume that it is ``well-behaved'', as in Section~\ref{pregeometries}; and then we argue
like in Section~\ref{proof of third part of 0-1 law for pregeometries}.)

The rest of the introduction is devoted to explaining, roughly, 
the results of this article. We try to appeal to the reader's intuition
rather than giving the full definitions of notions involved;
but sometimes references to these definitions are given.

We start, in Sections~\ref{forbidden structures} 
--~\ref{proving upper limit on the probability of an extension axiom}
, by considering $\mbK$
such that all $\mcM \in \mbK$ have trivial closure, so dimension is the same
as cardinality. Also, until Section~\ref{conditional measures} we consider only the uniform measure.
The first result, 
Theorem~\ref{dichotomy for forbidden weak substructures}, 
gives a dichotomy for the special
case when, for a fixed language $L$, with finite relational vocabulary,
and set $\mbF$ of ``forbidden'' $L$-structures,
$\mbK_n$ is defined to be the set of all $L$-structures $\mcM$
with universe $\{1, \ldots, n\}$ such that no $\mcF \in \mbF$ can
be {\em weakly embedded} into $\mcM$ (see Section~\ref{languages, structures and embeddings}).
If every $\mcF \in \mbF$ is ``simple'' in a sense which is made precise in
Theorem~\ref{dichotomy for forbidden weak substructures}, 
then for every extension axiom $\varphi$ of $\mbK$, 
the proportion of $\mcM \in \mbK_n$ which satisfy $\varphi$
approaches 1 as $n$ approaches infinity; and $\mbK$ has a zero-one law.
On the other hand, if there is at least one ``non-simple'' $\mcF \in \mbF$,
then for some $0 \leq c < 1$ and $2|F|$-extension axiom $\varphi$, 
the proportion of $\mcM \in \mbK_n$ in which $\varphi$ is true never
exceeds $c$; if the language has no unary relation symbols, then this
proportion approaches 0 as $n$ approaches infinity.
It may nevertheless be the case that $\mbK$ has a zero-one law,
as in the example of triangle-free graphs \cite{KPR}.

Theorem~\ref{dichotomy for forbidden weak substructures}, 
just described, is proved by using the more general 
Theorems~\ref{multiplicity when substitutions are admitted} 
and~\ref{upper limit on the probability of an extension axiom}. 
In Theorems~\ref{multiplicity when substitutions are admitted} 
and~\ref{upper limit on the probability of an extension axiom} 
we have no assumptions about how $\mbK$ is defined.
We will call a structure $\mcA$ {\em permitted} if it can be embedded into
some structure in $\mbK$. 
For the sake of simplifying this introductory description of the results,
let's assume that every permitted structure is isomorphic to some structure in
$\mbK$; in other words, we assume that $\mbK$ is, up to taking isomorphic copies,
closed under substructures (the `hereditary property').
The key concept will be that of {\em substitutions of permitted structures}
in a permitted (super)structure $\mcM$,
that is, the act of replacing, in $\mcM$, the interpretations (of relation symbols)
on the universe of $\mcA \subseteq \mcM$ by the interpretations in
another permitted structure $\mcA'$ with the same universe as $\mcA$.
If whenever $\mcA$, $\mcA'$, $\mcM$ are permitted, $\mcA \subseteq \mcM$ and
$\mcA$ and $\mcA'$ have the same universe, the result of  ``replacing
$\mcA$ by $\mcA'$ in $\mcM$'', denoted $\mcM[\mcA \trir \mcA']$, is
a permitted structure, then, for every extension axiom of $\mbK$, 
the proportion of structures in $\mbK_n$ in which it is true approaches 1 as $n$ approaches infinity.
This statement is a consequence of 
Theorem~\ref{multiplicity when substitutions are admitted} which, essentially,
is a reformulation, with the terminology used here, of known results -- although
this may not be obvious at first sight.

If, however, there exist permitted $\mcA$, $\mcA'$, $\mcM$ such that
$\mcM[\mcA \trir \mcA']$ is not permitted -- we say ``forbidden'' --
but the reverse substitution, that is, the replacement of $\mcA'$ by $\mcA$,
{\em never} produces a forbidden structure from a permitted one,
then one of the following holds:
(a) $\mbK$ fails to satisfy the disjoint amalgamation property, or
(b) there is an extension axiom $\varphi$ of $\mbK$ and $0 \leq c < 1$ such
that the proportion of $\mcM \in \mbK_n$ which satisfy $\varphi$ never exceeds $c$;
and if there are no unary relation symbols, then this proportion approaches 0.
Consider the example when $\mbK$ is the set of triangle-free graphs
and $\mcA$ and $\mcA'$ are graphs with vertex set $\{i,j\}$ where $i$ and $j$  
are adjacent in $\mcA'$ but not in $\mcA$.
Then we can find $\mcM \in \mbK$ such that $\mcM[\mcA \trir \mcA']$ 
is forbidden, but since the removal of an edge from a triangle-free graph
never produces a triangle and the class of triangle-free graphs has the
disjoint amalgamation property we are in case (b).
The statement before this example is a consequence of 
Theorem~\ref{upper limit on the probability of an extension axiom} and its corollary.
From these results we also get information, in case (a), 
about an instance of disjoint amalgamation which fails, and in case (b),
about the extension axiom $\varphi$.
Theorem~\ref{upper limit on the probability of an extension axiom} is proved by a counting argument.
One proves, under the assumption that
$\mbK$ has the disjoint amalgamation property, that for a properly chosen 
extension axiom $\varphi$ it is the case that for every $\mcM \in \mbK_n$ which satisfies 
$\varphi$, there are sufficiently many $\mcN \in \mbK_n$
which do {\em not} satisfy $\varphi$.

There is a third possibility, other than those considered in
the previous two paragraphs.
It is possible that there are permitted $\mcA$ and $\mcA'$ with the same universe
such that the substitution of $\mcA'$ for $\mcA$ in some permitted (super)structure $\mcM$
may produce a forbidden (not permitted) structure, but whenever this happens
then the reverse substitution of $\mcA$ for $\mcA'$ in some permitted $\mcN$, say,
may also produce a forbidden structure.
In this case it is possible that for every extension axiom $\varphi$ of $\mbK$,
the proportion of structures in $\mbK_n$ in which $\varphi$ is true approaches 1
as $n$ approaches infinity. But it is also possible that for some
extension axiom $\varphi$ of $\mbK$, the proportion of structures in $\mbK$
in which $\varphi$ is true approaches 0 as $n$ approaches infinity.
Section~\ref{examples concerning notion of admittance} gives examples showing this.
The same section also gives examples for which
Theorem~\ref{upper limit on the probability of an extension axiom} applies.
These examples show how the rather technical 
Theorem~\ref{upper limit on the probability of an extension axiom} and its 
(less technical) corollary can be used. 
Some examples in Section~\ref{examples concerning notion of admittance} 
also serve the purpose of 
illustrating differences between the uniform probability measure and 
conditional probability measures, 
which are introduced in Section~\ref{conditional measures};
these examples will be re-examined in Section~\ref{conditional measures}.
Section~\ref{proving upper limit on the probability of an extension axiom} 
is devoted to the proof of Theorem~\ref{upper limit on the probability of an extension axiom}.

In Section~\ref{conditional measures} conditional probability measures (on $\mbK_n$) 
are introduced,
motivated and 
illustrated with examples 
(that we have already met in Section~\ref{examples concerning notion of admittance}). 
One reason for introducing these are that the conditions which,
according to Theorem~\ref{multiplicity when substitutions are admitted}, guarantee
that for every extension axiom of $\mbK = \bigcup_{n \in \mbbN}\mbK_n$, 
the proportion of structures in
$\mbK_n$ which satisfy it approaches 1 as $n$ approaches infinity,
are rather restrictive.
The conditional measures that we consider -- or the {\em dimension conditional measures}, to be precise --
are more permissive with respect to satisfiability of extension axioms.
This is made precise by Lemma~\ref{admittance of substitutions implies acceptance} 
and 
Example~\ref{graphs with a unary predicate}, for instance.
Another motivation for considering conditional measures is that they 
are more closely related to random generation of finite structures.
While the uniform measure is conceptually simple it may,
for some $\mbK = \bigcup_{n \in \mbbN}\mbK_n$, be unclear what type of 
random generation procedure will, for any $\mcM \in \mbK_n$, generate $\mcM$ with
probability exactly $1 / |\mbK_n|$. 
Often, as in the case of $l$-coloured, or $l$-colourable, structures (graphs, for example),
the most obvious generation procedure -- first 
randomly assign colours, then randomly assign
relationships (e.g. edges) so that the colouring is respected --
corresponds to conditional measures, in the sense of this article.
A third reason for considering conditional measures is simply that they
may, in some situations, offer a simpler analysis of asymptotic problems than does the
uniform measure, while they are still natural in the sense of being
related to random generation of finite structures.
Finally we note that in some cases, as that of random $l$-colourable structures,
the conditional measure considered here coincides with the uniform probability measure
on properties which are first-order definable.
This follows from the proofs of the main theorems in 
Sections~\ref{l-colourable structures} 
and~\ref{the uniform measure and l-colourable structures} 

In Section~\ref{pregeometries} we start working in a context where the structures that
we consider have underlying (possibly nontrivial) pregeometries, and
`dimension' takes over the role of `cardinality'.
By a pregeometry on a structure $\mcM$ we mean a closure operator $\cl_{\mcM}$
which operates on subsets of the universe of $\mcM$ and satisfies certain
conditions \cite{Aig, Hod}; moreover we require that $\cl_{\mcM}$ is uniformly
definable in all structures considered 
(Definition~\ref{definition of a structure being a pregeometry} and 
Assumption~\ref{assumptions about pregeometries}).
The context considered previously is a special case of the 
framework of Section~\ref{pregeometries}.
The main results of this section,
Theorems~\ref{0-1 law for pregeometries}, 
\ref{third part of 0-1 law for pregeometries} 
and~\ref{0-law for extension axioms in pregeometries}, 
apply to the dimension conditional measure,
which is a conditional measure that ``considers'' closed sets of dimension 0 first,
then closed subsets of dimension 1, then of dimension 2, and so on. 
These theorems are related to 
Theorems~\ref{multiplicity when substitutions are admitted} 
and~\ref{upper limit on the probability of an extension axiom}.
Theorems~\ref{0-1 law for pregeometries} 
and~\ref{third part of 0-1 law for pregeometries} 
represent the ``positive'' side of things, like
Theorem~\ref{multiplicity when substitutions are admitted}, showing
that if certain conditions are satisfied, then for every extension axiom
of $\mbK$ the probability (with the dimension conditional measure) that it
holds in a member of $\mbK_n$ approaches 1 as $n$ approaches infinity;
and from this a zero-one law is derived.
The conditions in question require, as in Section~\ref{forbidden structures},
that whenever $\mcM$ is permitted, then certain ``substitutions'',
or ``replacements'', of interpretations can be made in $\mcM$
without producing a forbidden (not permitted) structure.
Also, there is a requirement that the underlying pregeometry,
and possibly some other structure which is never changed,
is {\em polynomially $k$-saturated}. This roughly means that
for every $k \in \mbbN$ and all sufficiently large $n$ and every
$\mcM \in \mbK_n$, the reduct of $\mcM$ to the sublanguage which defines the pregeometry
satisfies every $k$-extension axiom (with respect to the set of such reducts);
and moreover, the truth of a $k$-extension axiom
has many different witnesses compared to the size of the universe.

The last result of Section~\ref{pregeometries}, 
Theorem~\ref{0-law for extension axioms in pregeometries},
is a ``cousin'' of Theorem~\ref{upper limit on the probability of an extension axiom} 
and its corollary,
and tells that if there are
permitted $\mcA$ and $\mcA'$ such that $\mbK$ {\em accepts} 
(Definition~\ref{k-independence hypothesis for classes of structures}) the
substitution $[\mcA \trir \mcA']$ but {\em not} the reverse subsitution
$[\mcA' \trir \mcA]$, then either $\mbK$ fails to have the {\em independent amalgamation
property}, or for some extension axioms $\varphi$ and $\psi$, the probability,
with the dimension conditional measure, that $\varphi \wedge \psi$ holds in a member of $\mbK_n$
approaches 0 as $n$ approaches infinity.
The analogue of Theorem~\ref{dichotomy for forbidden weak substructures} 
in the setting of underlying pregeometries 
and the dimension conditional measure is given by the corollary in 
Example~\ref{Henson examples}.
The proofs of Theorems~\ref{0-1 law for pregeometries}, 
\ref{third part of 0-1 law for pregeometries} 
and~\ref{0-law for extension axioms in pregeometries} appear in 
Section~\ref{proofs of main results}.

Section~\ref{proofs of main results} gives the proofs of the main 
theorems of Section~\ref{pregeometries}.
The definitions appearing in 
Sections~\ref{proving upper limit on the probability of an extension axiom}
and~\ref{proofs of main results}
are only used within those sections.

Sections~\ref{l-colourable structures} 
and~\ref{the uniform measure and l-colourable structures} 
study asymptotic properties of random {\em $l$-colourable}, as well as {\em strongly}
$l$-colourable, structures
in a fixed (but arbitrary) relational language in which the arity of each symbol is at least 2.
Examples~\ref{example of coloured structures of first kind} 
and~\ref{example of coloured structures of second kind} 
show that {\em $l$-coloured} structures
can be treated within the context developed in Section~\ref{pregeometries}.
Theorem~\ref{third part of 0-1 law for pregeometries} implies that $l$-coloured structures
satisfy a zero-one law with respect to the dimension conditional measure.
Since $l$-colourable structures can be viewed as reducts of $l$-coloured
structures we will also consider a ``reduct version'' of the dimension
conditional measure. With this probability measure it is {\em not} true that
all extension axioms of $l$-colourable structures hold almost surely; but 
we can show that all extension axioms of a certain kind, called
the {\em $l$-colour compatible extension axioms}, hold almost surely in
sufficiently large structures; and this is enough for subsequently deriving
a zero-one law for $l$-colourable structures, when using the 
probability measure derived from the dimension conditional measure
(Theorem~\ref{0-1 law for l-colourable structures}).
We also prove a result saying that if almost all $l$-colourable structures
have an $l$-colouring with sufficiently even distribution of colours,
then, with the uniform probability measure, 
every $l$-colour compatible extension axiom holds almost surely, 
almost every $l$-colourable
structure has a unique $l$-colouring (up to permutation of the colours), 
and a zero-one law holds with
the uniform probability measure as well
(Theorem~\ref{transfer of 0-1 laws to uniform measure}
and Proposition~\ref{definability of colourings under richness condition}).

In Section~\ref{the uniform measure and l-colourable structures}
we prove, by combinatorial arguments, that, indeed, almost all
$l$-colourable structures have an $l$-colouring 
with sufficiently even distribution of colours
(Theorem~\ref{new theorem}). 
Thereby we confirm that almost all
$l$-colourable structures have a unique $l$-colouring and that, also with the 
uniform probability measure, a zero-one law holds for $l$-colourable structures
(Theorems~\ref{main theorem} and~\ref{main theorem for strong colourings}). 

All results of the article hold also if one restricts
attention to structures in which certain relation symbols
(of arity at least 2) are always interpreted
as irreflexive and symmetric relations.
All arguments, except those in Section~\ref{the uniform measure and l-colourable structures},
work out in the same way under this assumption.

The results in Sections~\ref{l-colourable structures} 
and~\ref{the uniform measure and l-colourable structures} 
may be useful in contexts which do not directly speak about colourings. 
Suppose that for some $\mbK = \bigcup_{n \in \mbbN}\mbK_n$ and
probability measure $\mu_n$ on $\mbK_n$ there is $l \in \mbbN$ such that,
for $n$ large enough,
$\mcM \in \mbK_n$ is almost surely $l$-colourable.
If we know that every $l$-colourable structure (with universe an
initial segment of $\{1, 2, \ldots \}$) belongs to $\mbK$ and that
the set of $L$-structures which are $l$-colourable
has a zero-one law for the measures $\mu_n$, then also
$\mbK$ has a zero-one law for the same measures.
This approach was used in \cite{KPR} when proving that
if $\mbK_n$ is the set of $(l+1)$-clique-free graphs (or $\mcK_{l+1}$-free graphs)
with universe $\{1, \ldots,n\}$,
then $\mbK$ has a zero-one law for the uniform probability measure.
The authors of \cite{KPR} first proved that almost all $(l+1)$-clique-free graphs are $l$-colourable,
with a relatively even distribution of colours,
and then that the $l$-colourable graphs have a zero-one law.

The notions of `polynomial $k$-saturation' and `acceptance of substitutions'
in Section~\ref{pregeometries}
are versions, adapted to the context of this article,
of the notions `polynomial $k$-saturation' and `$k$-independence hypothesis'
in \cite{Dj06a}.
This is sufficiently clear for polynomial $k$-saturation, but it is perhaps
harder to see the relationship between admittance of ($k$-)substitutions and
the $k$-independence hypothesis. 
However, in both cases the essential difference between Section~\ref{pregeometries}
of this article and
\cite{Dj06a} is that in \cite{Dj06a} complete types of an infinite structure are
considered, while here we consider types with only quantifier-free formulas 
of tuples enumerating the universe
of a closed substructure of some permitted structure. 
But in this article we avoid speaking about such types since it is equally convenient
to speak about (sub)structures and formulas describing them up to isomorphism.
Lemmas~\ref{changes of relations in one step} 
--~\ref{induction step for pregeometries},
as well as their proofs, are adaptations to the context of this article of 
Lemmas~2.16 --~2.22 in \cite{Dj06a}.
The results of this article have their beginnings in considerations from two 
directions. On the one hand, trying to understand asymptotic satisfiability
of extension axioms -- conditions implying that they almost surely hold, and
conditions implying that some almost surely fail -- and
on the other hand, trying to understand if some zero-one laws for finite structures were hidden
in the probabilistic arguments used in \cite{Dj06a}.
\\

\noindent
{\bf Acknowledgements.}
I thank Svante Janson for helpful suggestions concerning
Section~\ref{the uniform measure and l-colourable structures}, 
which shortened some proofs there.
I also thank the anonymous referee for having read the article so carefully,
for valuable suggestions and for pointing out some errors, now corrected.

\section{Preliminaries}\label{preliminaries}

\subsection{Languages, structures and embeddings}
\label{languages, structures and embeddings}
For basic notions not explained here the reader is refered to \cite{Hod, EF}.
By a {\bf \em language} $L$ we mean the 
set of (first-order) formulas that can be built up
from a {\bf \em vocabulary} (also called {\bf \em signature})
which is a set of relation, constant and/or function symbols.
We consider the identity symbol `='  as a logical symbol which
we may always use, together with connectives and (first-order) 
quantifiers, to build formulas; so `=' is never mentioned when
we describe the symbols of a vocabulary.
If the vocabulary has no constant or function symbols,
then we call it {\bf \em relational}.

Structures will be denoted by ``calligraphic'', letters:
$\mcA$, $\mcB$, $\ldots$, $\mcM$, $\mcN$, $\ldots$.
Their universes
will be denoted by the corresponding non-calligraphic letter
$A$, $B$, $\ldots$, $M$, $N$, $\ldots$, or with bars around the letter;
for instance, $|\mcM|$ as well as $M$ denote the universe of $\mcM$.
The cardinality of a set $X$ is denoted by $|X|$;
and the cardinality of the (universe of) the structure $\mcM$ is
denoted by $\left\|\mcM\right\|$, or by $|M|$.
Boldface letters always denote classes, usually sets, of structures.
Sequences, or tuples, of elements are denoted by $\bar{a}, \bar{b}, \ldots$;
and $|\bar{a}|$ denotes the length of the sequence $\bar{a}$.
By `$\bar{a} \in M$' we mean that $\bar{a}$ is a sequence such that all of
its elements belong to the set $M$. Sometimes we write $\bar{a} \in M^n$
to show that $\bar{a}$ has length $n$.
By $\rng(\bar{a})$, the {\bf \em range} (or {\bf \em image}) of $\bar{a}$,
we denote the set of all elements that occur in $\bar{a}$.
In the last section we often use the abbreviation $[n] = \{1, \ldots, n\}$
if $n$ is a positive integer.
For $\alpha \in \mbbR$, $\lfloor \alpha \rfloor$ denotes the largest integer $m$ such
that $m \leq \alpha$.
If $f : A \to B$ and $\bar{a} = (a_1, \ldots, a_r) \in A^r$, then
$f(\bar{a})$ denotes the sequence $(f(a_1), \ldots, f(a_r))$.
If $L$ has no constant symbols, then we allow an $L$-structure to have an
empty universe.

Suppose that $\mcM$ and $\mcN$ are $L$-structures, where $L$ is, as usual, a language.
A function $f : M \to N$ is called a {\bf \em weak embedding} of $\mcM$ (in)to $\mcN$ if
$f$ is {\em injective} and:
\begin{itemize}
	\item[(1)] For every constant symbol $c$, $f(c^{\mcM}) = c^{\mcN}$.
	\item[(2)] For every function symbol $g$, of arity $r$, say, and every
	$\bar{a} \in M^r$, $f(g^{\mcM}(\bar{a})) = g^{\mcN}(f(\bar{a}))$.
	\item[(3)] for every relation symbol, $R$, of arity $r$, say, 
	if $\bar{a} \in R^{\mcM}$ then $f(\bar{a}) \in R^{\mcN}$.
\end{itemize}
We say that $f$ is an {\bf \em embedding} if $f$ is {\em injective} and (1), (2) and 
the following hold:
\begin{itemize}
	\item[(3')] for every relation symbol, $R$, of arity $r$, say, 
	$\bar{a} \in R^{\mcM}$ if and only if $f(\bar{a}) \in R^{\mcN}$.
\end{itemize}
Thus, embeddings are injective and a bijective embedding is the same as an isomorphism.
We say that $\mcM$ is {\bf \em (weakly) embeddable} into $\mcN$ if there
exists a (weak) embedding from $\mcM$ to $\mcN$.
We say that $\mcM$ is a {\bf \em weak substructure} of $\mcN$, denoted $\mcM \subseteq_w \mcN$,
if $M \subseteq N$ and the identity mapping $id : M \to N$ is a weak embedding.
We call $\mcM$ a {\bf \em substructure} of $\mcM$, denoted $\mcM \subseteq \mcN$,
if $M \subseteq N$ and the identity mapping $id : M \to N$ is an embedding.
$\mcA$ is a {\em proper} (weak) substructure of $\mcM$ if $\mcA$ is
a (weak) substructure of $\mcM$ and $\mcA \neq \mcM$. 
The symbol `$\subset$' means `proper subset' or `proper substructure'.

If $\mcM$ is a structure and $A \subseteq M$, then $\mcM \uhrc A$ denotes the
substructure of $\mcM$ which is generated by $A$
(the smallest substructure $\mcN$ of $\mcM$ such that $A \subseteq N$); 
so if the vocabulary is relational,
then $|\mcM \uhrc A| =A$.
If $L_0$ is a language such that $L_0 \subseteq L$ and $\mcM$ is an $L$-structure,
then $\mcM \uhrc L_0$ denotes the reduct of $\mcM$ to $L_0$.
Observe that if all constant and function symbols of $L$ belong to the vocabulary
of $L_0$, then the reduct $\mcM \uhrc L_0$ can also be viewed as an $L$-structure
in which the interpretation of every relation symbol which belongs to the vocabulary
of $L$ but not to the vocabulary of $L_0$ is empty. 
So provided that the smaller language $L_0$ contains all constant and function symbols we have 
$\mcM \uhrc L_0 \subseteq_w \mcM$, from which it is apparent that 
the notion of weak substructure generalizes the notion of reduct, as well as the notion of substructure.

Since we will several times speak about graphs, we note that, with graph
theoretic terminology,
if $\mcM$ and $\mcN$ are graphs, then $\mcM$ is a {\bf \em subgraph} of $\mcN$
if and only if $\mcM$ is a weak substructure of $\mcN$;
and $\mcM$ is an {\bf \em induced subgraph} of $\mcN$ if and only if
$\mcM$ is a substructure of $\mcN$.

Suppose that $R$ is a relation symbol from the vocabulary of 
the language of $\mcM$. Then a tuple $\bar{a}$ of elements from
$M$ is called an {\bf \em $R$-relationship} of $M$ if
$\bar{a} \in R^{\mcM}$ (or equivalently, if $\mcM \models R(\bar{a})$). 
If the symbol `$R$' is clear from the context, or if it does not
matter which $R$ we refer to, then we may just call an $R$-relationship
a {\bf \em relationship}.
Sometimes we consider only structures $\mcM$ in which certain relation symbols $R_1, \ldots, R_k$
are interpreted as irreflexive and symmetric relations 
(see Remark~\ref{remark that the theorems generalize to symmetric structures}).
In this case an {\em $R_i$-relationship} of $\mcM$ (for $i = 1, \ldots, k$) is a set $\rng(\bar{a})$
such that $\bar{a} \in (R_i)^{\mcM}$.
So for graphs in general, a relationship is the same as a directed edge;
and if we consider only undirected graphs, a relationship is the same as
an (undirected) edge.

\begin{rem}\label{remark that the theorems generalize to symmetric structures}{\rm
Suppose that $R$ is an $n$-ary relation on a set $A$.
Then $R$ is called {\bf \em irreflexive} if $(a_1, \ldots, a_n) \in R$ implies that
$a_i \neq a_j$ if $i \neq j$.
If $(a_1, \ldots, a_n) \in R$ implies that $(a_{\pi(1)}, \ldots, a_{\pi(n)}) \in R$ for 
every permutation $\pi$ of $\{1, \ldots, n\}$, then we say that $R$ is {\bf \em symmetric}.
All results in this article hold also
if we assume that the interpretations of certain relation symbols are always
irreflexive and symmetric.
The proofs in this case are either the same as, or obvious modifications of, the given proofs.
}\end{rem}

\subsection{Amalgamation}
\label{amalgamation and fraisse limits}

Let $\mbK$ be a class of finitely generated
$L$-structures, where $L$ has a countable vocabulary,
and let $\widehat{\mbK}$ be the class consisting of all
$L$-structures $\mcM$ such that $\mcM$ is isomorphic to a member of $\mbK$;
so $\widehat{\mbK}$ is ``closed under isomorphism''.
See \cite{Hod} (Chapter 7), for example, for definitions of the following notions:
{\bf \em hereditary property}, or being {\bf \em closed under substructures} as we
sometimes say here, {\bf \em amalgamation property} and
{\bf \em joint embedding property}.
We say that $\mbK$ has any of these properties if $\widehat{\mbK}$ has it.
If the vocabulary of $L$ has only relation symbols then the
amalgamation property implies the joint embedding property;
but in general the later property is not implied by the first.

If $\widehat{\mbK}$ has all three properties, then the so-called
{\bf \em Fra\"{i}ss\'{e} limit} $\mcM_{\mbK}$ of $\widehat{\mbK}$ exists \cite{Hod}.
$\mcM_{\mbK}$ has the following properties: $\mcM_{\mbK}$ is countable,
every finitely generated $\mcA \subseteq \mcM_{\mbK}$ belongs to $\widehat{\mbK}$;
every $\mcA \in \mbK$ can be embedded into $\mcM_{\mbK}$, and
if $\mcA \subseteq \mcM_{\mbK}$ is finitely generated and $\mcA \subseteq \mcB \in \widehat{\mbK}$,
then there is an embedding $f : \mcB \to \mcM_{\mbK}$ such that $f \uhrc A$ is the 
identity function \cite{Hod}. 
The Fra\"{i}ss\'{e} limit $\mcM_{\mbK}$ of $\widehat{\mbK}$, if it exists, is
also called the Fra\"{i}ss\'{e} limit of $\mbK$.

We will consider the following (stronger) variant of the amalgamation property:
We say that $\widehat{\mbK}$ (and $\mbK$) has the 
{\bf \em disjoint amalgamation property} 
if whenever $\mcA, \mcB, \mcC \in \widehat{\mbK}$, 
$\mcA \subseteq \mcB$, $\mcA \subseteq \mcC$ and $B \cap C = A$,
then there is $\mcD \in \widehat{\mbK}$ such that
$\mcB \subseteq \mcD$, $\mcC \subseteq \mcD$.

\subsection{Pregeometries}
\label{introduction to pregeometries}

The notion of {\bf \em (combinatorial) pregeometry}, also called 
{\bf \em matroid}, will play a role in sections~\ref{pregeometries}
and~\ref{proofs of main results}.
See \cite{Hod} (Chapter 4.6), or \cite{Aig} (Chapter II.3), for a definition.
We use the following terminology when $(A, \cl)$ is a pregeometry, with {\bf \em
closure operator} $\cl$ which maps every $X \subseteq A$ to some closed $Y \subseteq A$.
For $X, Y, Z \subseteq A$, {\bf \em $X$ is independent from $Y$ over $Z$}
if for every $a \in X$, 
$a \in \cl(Y \cup Z)$ $\Longleftrightarrow$ $a \in \cl(Z)$.
In the special case that $Z = \es$ we say that {\bf \em $X$ is independent from $Y$}.
Because of the `exchange property' of pregeometries, independence is symmetric
with respect to $X$ and $Y$.
We say that $a \in A$ is independent from $Y \subseteq A$ over $Z \subseteq A$
if $\{a\}$ is independent from $Y$ over $Z$.
A set $X \subseteq A$ is called {\bf \em independent} if 
for every $a \in X$, $a$ is independent from $X - \{a\}$ (over $\es$).
The {\bf \em dimension} of $X \subseteq A$ is the supremum of the cardinalities of independent
subsets of $X$.
A set $X \subseteq A$ is called {\bf \em closed} if $\cl(X) = X$.
If $\cl(X) = X$ for every $X \subseteq A$ then we call
$(A, \cl)$ the {\bf \em trivial pregeometry} on $A$.

\subsection{Zero-one laws.}\label{preliminaries about zero-one laws}

Suppose that, for $n \in \mbbN$, $\mbK_n$ is a set of $L$-structures
and that $\mu_n$ is a probability measure on $\mbK_n$.
If $\mu_n(\mcM) = 1/|\mbK_n|$ for all $\mcM \in \mbK_n$,
then we call $\mu_n$ the {\bf \em uniform probability measure} on $\mbK_n$.
We say that $\mbK = \bigcup_{n \in \mbbN}\mbK_n$
has a {\bf \em zero-one law} if for every $L$-sentence $\varphi$,
$\lim_{n\to\infty}\mu_n\big(\{\mcM \in \mbK_n : \mcM \models \varphi\}\big)$
exists and is 0 or 1.
When saying ``$\varphi$ is almost surely true (or false)''  we mean that
the limit is 1 (or 0).
If $\mu_n$ is the uniform probability measure for all $n$ we may instead say that
``almost all $\mcM \in \mbK$ satisfy $\varphi$'' if the limit is 1.
By the {\bf \em almost sure theory of $\mbK$} (with respect to
the measures $\mu_n$), we mean the set of sentences $\varphi$ such that
the probability that $\varphi$ is true in $\mbK_n$ approaches 1 as $n \to \infty$.

\section{Permitted structures and substitutions}\label{forbidden structures}

\noindent
From this section and until Section~\ref{pregeometries} we work within the
following framework:

\begin{assumptionsandterminology}\label{basic asumptions}{\rm
Fix a first-order language $L$ with finite relational vocabulary.
Let $(m_n : n \in \mbbN)$ be a sequence of positive integers
such that $\lim_{n\to\infty} m_n =~\infty$.
For every $n \in \mbbN$ let $\mbK_n$ be a set of $L$-structures 
with universe $\{1, \ldots, m_n\}$; 
and let $\mbK = \bigcup_{n \in \mbbN}\mbK_n$.
A structure $\mcM$ is called {\bf \em represented (with respect to $\mbK$)} if it is isomorphic to
a structure in $\mbK$. 
A structure $\mcM$ is called {\bf \em permitted (with respect to $\mbK$}) if it is embeddable
into a structure in $\mbK$. 
A structure which is not permitted is called {\bf \em forbidden}.
Since we fix $\mbK$ for rest of the section we sometimes omit the phrase 
``with respect to $\mbK$''.
}\end{assumptionsandterminology}

\noindent
Observe that if $\mbK$ has the hereditary property,
then a structure is permitted if and only if it is represented. 
In this section and the next, all examples of $\mbK$ which are considered
in some detail have the hereditary property.
However, since the results do not depend on this we do not assume it.
(One example of $\mbK$ which is {\em not} closed under substructures is
given by letting $\mbK_n$ be the set of triangle-free graphs with universe $\{1, \ldots, n\}$
and diameter 2.)

\begin{defin}\label{definition of extension axiom}{\rm
Suppose that $\mcA$ and $\mcB$ are permitted structures and that 
$\mcA$ is a proper substructure of $\mcB$.\\
(i) The {\bf \em $\mcB / \mcA$-extension axiom} (or the {\bf \em $\mcB$-extension axiom over $\mcA$})
holds, by definition, in $\mcM$ if the following is true:
\begin{itemize}
\item[] For every embedding $\tau$ of $\mcA$ into $\mcM$ there exists an embedding
$\pi$ of $\mcB$ into $\mcM$ which extends $\tau$ (i.e. $\pi(a) = \tau(a)$ whenever
$a \in A$).
\end{itemize}
The $\mcB/\mcA$-extension axiom can be expressed by a first-order sentence 
of the form
$$\forall x_1, \ldots, x_n \exists y_1, \ldots, y_m \big(\varphi(x_1, \ldots, x_n) \longrightarrow
\psi(x_1, \ldots, x_n, y_1, \ldots, y_m)\big),$$
where $\varphi$ and $\psi$ are quantifier-free.
If the language has no constant symbols, 
then we allow the possibility that the universe of $\mcA$ is empty, in which
case the $\mcB / \mcA$-extension axiom is called the $\mcB / \es$-extension axiom.
It is then expressed by an existential formula
$$\exists y_1, \ldots, y_m \psi(y_1, \ldots, y_m).$$
(ii) If $\left| \mcB \right| \leq k+1$, then the $\mcB/\mcA$-extension axiom
is called a {\bf \em $k$-extension axiom of $\mbK$};
or, if we do not care about $k$, just an {\bf \em extension axiom of $\mbK$}. 
If $\mbK$ is clear from the context we may omit saying ``of $\mbK$''.
}\end{defin}

\begin{rem}\label{remark about extension axioms and 0-1 laws}{\rm
{\em If there are probability measures $\mu_n$ on $\mbK_n$, for $n \in \mbbN$,
such that for every extension axiom $\varphi$ of $\mbK$, the 
$\mu_n$-probability that $\mcM \in \mbK_n$ satisfies $\varphi$
approaches 1 as $n$ approaches $\infty$,
then $\mbK$ has a zero-one law for the measures $\mu_n$.}
The usual proof of this statement does not depend on the measures $\mu_n$.
It is proved in \cite{Fag, EF, Hod, Win} (for example) 
by collecting into a theory $T_{\mbK}$ all
extension axioms, together with sentences expressing the possible 
isomorphism types of substructures of members of $\mbK$.
The general idea of the argument is as follows.
By the assumptions in the above statement and compactness, $T_{\mbK}$ is consistent.
By a back-and-forth argument one then proves that $T_{\mbK}$ is countably
categorical and therefore complete.
The completeness of $T_{\mbK}$ (and compactness) implies that
$\mbK$ has a zero-one law.
}\end{rem}

\noindent
If we define $\mbK$ by forbidding certain weak substructures, and 
the thus obtained $\mbK$ has the disjoint amalgamation property, then
we have the following ``dichotomy''.

\begin{theor}\label{dichotomy for forbidden weak substructures}
Let $\mbF$ be a set of finite $L$-structures and, for every $n \in \mbbN$,
let $\mbK_n$ consist of exactly those $L$-structures $\mcM$ with universe 
$\{1, \ldots, n\}$ such that {\rm no} $\mcF \in \mbF$ is 
{\rm weakly} embeddable into $\mcM$
(so in particular, every member of $\mbF$ is forbidden).
Assume that $\mbK_n \neq \es$ for all sufficiently large $n$ and that
$\mbK$ has the disjoint amalgamation property.
Consider the following condition:
\begin{itemize}
\item[($\ast$)] There are $\mcF \in \mbF$, a relation symbol $R$ of arity $r$, say,
and $\bar{a} \in F^r$ such that $\rng(\bar{a})$ is a {\rm proper} subset of $F$,
$\bar{a} \in R^{\mcF}$, and if $\mcP$ is constructed by removing
the $R$-relationship $\bar{a}$, but making no other changes in $\mcF$,
then $\mcP$ is permitted.
\end{itemize}
If ($\ast$) is false, then, for every $k \in \mbbN$, 
the proportion of $\mcM \in \mbK_n$ which satisfy all $k$-extension axioms of $\mbK$ approaches 1
as $n \to \infty$.
If ($\ast$) is true, then letting $\mcF \in \mbF$, $R$ and $\bar{a}$ be any witnesses of 
property ($\ast$) and
letting $\alpha$ be the number of permitted structures with 
universe $\{1, \ldots, |\rng(\bar{a})|\}$, 
the proportion of $\mcM \in \mbK_n$ which satisfy all
$(2|F| - |\rng(\bar{a})| - 1)$-extension axioms of $\mbK$ never exceeds $1 - 1/(1 + \alpha)$.
Moreover, if $L$ has no unary relation symbols, then the proportion of
$\mcM \in \mbK_n$ which satisfy all $(2|F| - |\rng(\bar{a})| - 1)$-extension axioms
approaches 0 as $n \to \infty$.
\end{theor}

\noindent
Theorem~\ref{dichotomy for forbidden weak substructures} is a consequence of 
Theorems~\ref{multiplicity when substitutions are admitted}
and~\ref{upper limit on the probability of an extension axiom}.
Since one may see it as an application of these theorems,
we give the proof of Theorem~\ref{dichotomy for forbidden weak substructures}
as Example~\ref{example of forbidden weak substructures} in 
Section~\ref{examples concerning notion of admittance}.
The argument in Example~\ref{example of forbidden weak substructures}
gives some information about what happens if $\mbK$ does not have the disjoint amalgamation property.

\begin{rem}\label{remark about unary relation symbols}{\rm
(i) Suppose, as in Theorem~\ref{dichotomy for forbidden weak substructures}, 
that $\mbF$ is a set of finite $L$-structures and, for every $n \in \mbbN$,
let $\mbK_n$ consist of exactly those $L$-structures $\mcM$ with universe 
$\{1, \ldots, n\}$ such that {\rm no} $\mcF \in \mbF$ is 
{\rm weakly} embeddable into $\mcM$. Here is a condition on $\mbF$ which 
implies that $\mbK$ has the disjoint amalgamation property.
As in \cite{Hen72}, let us call an $L$-structure $\mcM$ {\em decomposable}
if there are different $L$-structures $\mcA$ and $\mcB$ such that $M = A \cup B$,
$\mcA \uhrc A \cap B = \mcB \uhrc A \cap B$
and for every relation symbol $R$, $R^{\mcM} = R^{\mcA} \cup R^{\mcB}$.
Otherwise we call $\mcM$ {\em indecomposable}.
It is now straightforward to show that if all structures in $\mbF$ are indecomposable,
then $\mbK$ has the disjoint amalgamation property.
(This statement is analogous to Theorem~1.2~(i) in \cite{Hen72}.)\\
(ii) One may ask if the assumption that there are no unary relation symbols is
necessary for the last statement of 
Theorem~\ref{dichotomy for forbidden weak substructures}.
The author does not have an example showing that this statement fails without
the assumption that there are no unary relation symbols, {\em if} we
assume, as in Theorem~\ref{dichotomy for forbidden weak substructures},
that $\mbK$ has the disjoint amalgamation property.
But Example~\ref{example of necessity of no unary symbols}
shows that when it is assumed that there are no unary relation symbols in
Theorem~\ref{upper limit on the probability of an extension axiom},
then this assumption is necessary.
}\end{rem}

\noindent
Two examples follow, one for which ($\ast$) in Theorem~\ref{dichotomy for forbidden weak substructures}
does not hold, and one for which ($\ast$) holds.

\begin{exam}{\rm
Suppose that $L$ has only one binary relation symbol $R$ and that
$\mbF = \{\mcA, \mcB\}$, where $A = \{1\}$, $R^{\mcA} = \{(1,1)\}$,
$B = \{1,2\}$ and $R^{\mcB} = \{(1,2), (2,1)\}$.
If $\mbK_n$ and $\mbK$ are defined as in Theorem~\ref{dichotomy for forbidden weak substructures},
then an $L$-structure is permitted if and only if it is an irreflexive and antisymmetric
directed graph. Moreover, the property ($\ast$) fails for $\mbF$.
}\end{exam}

\begin{exam}\label{example of l-clique-free graphs}{\rm ({\bf $\mcK_l$-free graphs})
It is not difficult to define $\mbF$ for which the property ($\ast$) holds,
but let us mention an example which has been studied in some detail \cite{KPR}.
Let $L$ have only one binary relation symbol $R$ and consider only structures
in which $R$ is interpreted as an irreflexive and symmetric relation, that is,
an undirected graph without loops. 
Let $l \geq 3$ and let $\mcK_l$ be the complete (undirected) graph
with vertices $1, \ldots, l$. If $\mbF  = \{\mcK_l\}$ then condition ($\ast$) holds,
since the removal of one edge from $\mcK_l$ creates a permitted graph.
It is easy to see that $\mbK$ has the disjoint amalgamation property.
By Theorem~\ref{dichotomy for forbidden weak substructures}, the proportion of $\mcM \in \mbK_l$
which satisfy all $(2l - 3)$--extension axioms of $\mbK$ approaches 0 as $n \to \infty$.
For $l = 3$ at least, this conclusion is not new.
Because the proportion of $\mcK_3$-free graphs ({\em triangle-free} graphs)
which are bipartite approaches 1 as $n \to \infty$ \cite{EKR, KPR};
and a graph is bipartite if and only if it has no cycle of odd length;
moreover, it is easy to see that a 5-cycle or 3-cycle exists in every 
$\mcK_3$-free graph which satisfies all 3-extension axioms.
}\end{exam}

\begin{rem}\label{remark about failure of extension axioms and 0-1 laws}{\rm
Even if, for some extension axiom $\varphi$ of $\mbK$, the proportion of $\mcM \in \mbK_n$ which satisfy
$\varphi$ does not approach 1, 
$\mbK$ may nevertheless have a zero-one law with respect to the uniform probability measure.
For example, it has been shown \cite{KPR} that, 
for every $l \geq 3$, if $\mbF = \{\mcK_l\}$
where $\mcK_l$ is the complete graph with $l$ vertices,
and $\mbK_n$ and $\mbK$ are defined as in Example~\ref{example of l-clique-free graphs}, 
then $\mbK$ has a zero-one law for the uniform probability measure.
}\end{rem}

\begin{defin}\label{definition of multiplicity}{\rm
Let $\mcM$, $\mcA$ and $\mcB$ be structures and suppose that $\mcA$ is a proper substructure of $\mcB$.\\
(i) We say that the {\bf \em $\mcB/\mcA$-multiplicity of $\mcM$ is at least $m$} 
(or that the {\bf \em $\mcB$-multiplicity over $\mcA$ in $\mcM$ is at least $m$}) if  the following holds:
\begin{itemize}
\item[] Whenever $\sigma$ is an embedding of $\mcA$ into $\mcM$, then there are
embeddings $\sigma_i$ of $\mcB$ into $\mcM$, for $i = 1, \ldots, m$,
such that each $\sigma_i$ extends $\sigma$ and if $i \neq j$ then
$\sigma_i(B) \cap \sigma_j(B) = \sigma(A)$.
\end{itemize}
The {\bf \em $\mcB/\mcA$-multiplicity is $m$} if it is at least $m$ but not at least $m+1$.\\
(ii) We say that {\bf \em $\mcM$ has (at least) $n$ copies of $\mcA$} if there are
(at least) $n$ different substructures $\mcA'_1, \ldots, \mcA'_n$ of $\mcM$ such that each
$\mcA'_i$ is isomorphic to $\mcA$. 
}\end{defin}

\begin{rem}\label{remarks on extension axioms and multiplicity}{\rm
Observe the following relationships between extension axioms and multiplicity,
where we assume that $\mcA \subset \mcB$.\\
(i) $\mcM$ satisfies the $\mcB / \mcA$-extension axiom if and only if
the $\mcB / \mcA$-multiplicity of $\mcM$ is at least 1.\\
(ii) Suppose that there are a structure $\mcC$ and embeddings
$\sigma_1 : \mcB \to \mcC$ and $\sigma_2 : \mcB \to \mcC$ such that
$\sigma_1 \uhrc A = \sigma_2 \uhrc A$ and
$\sigma_1(B) \cap \sigma_2(B) = \sigma_1(A)$.
If $\mcM$ satisfies the $\mcC/\mcA$-extension axiom then the $\mcB/\mcA$-multiplicity
of $\mcM$ is at least 2.
}\end{rem}

\begin{defin}\label{definition of substitutions}{\rm
Suppose that the vocabulary of $L$ does not contain any constant symbol.
Let $\mcA$, $\mcB$ and $\mcM$ be $L$-structures such that
$\mcA \subseteq \mcM$ and $|\mcA| = |\mcB|$.
We define $\mcM[\mcA \trir \mcB]$ to be the structure obtained by ``replacing $\mcA$
by $\mcB$ inside $\mcM$'', or more precisely, $\mcM[\mcA \trir \mcB]$ 
is defined to be the structure with the same universe as $\mcM$ which satisfies the following conditions:
For every $n$ and every relation symbol $R$ of arity $n$,
	\begin{itemize}
		\item[(1)] \ if $(a_1, \ldots, a_n) \in A^n$, then 
	$(a_1, \ldots, a_n) \in R^{\mcM[\mcA \trir \mcB]} \Longleftrightarrow 
	(a_1, \ldots, a_n) \in R^{\mcB}$, and
	\item[(2)] \ if $(a_1, \ldots, a_n) \in M^n - A^n$,
	then $(a_1, \ldots, a_n) \in R^{\mcM[\mcA \trir \mcB]} \Longleftrightarrow
	(a_1, \ldots, a_n) \in R^{\mcM}$.
	\end{itemize}
The notation $\mcM[\mcA \trir \mcB]$ may be read as {\bf \em $\mcM$ with $\mcA$ replaced by $\mcB$},
or {\bf \em $\mcM$ with $\mcB$ substituted for $\mcA$}.
}\end{defin}

\begin{defin}\label{definition of admitting substitutions}{\rm
Let $\mcA$ and $\mcB$ be permitted structures (with respect to $\mbK$) with the same universe. \\
(i) We say that  {\bf \em $\mbK$ admits the substitution $[\mcA \trir \mcB]$}
if for every represented $\mcM$ such that $\mcA \subseteq \mcM$, 
$\mcM[\mcA \trir \mcB]$ is a represented structure.\\
(ii) We say that  {\bf \em $\mbK$ weakly admits the substitution $[\mcA \trir \mcB]$}
if for every represented $\mcM$ such that $\mcA \subseteq \mcM$, 
$\mcM[\mcA \trir \mcB]$ is a permitted structure.\\
(iii) If $|\mcA| = |\mcB|$ and $\left\|\mcA\right\| \leq k$, then 
we call $[\mcA \trir \mcB]$ a {\bf \em $k$-substitution}.\\
(iv) If $\mbK$ (weakly) admits every $k$-substitution $[\mcA \trir \mcB]$, where 
$\mcA$ and $\mcB$ are permitted structures (with the same universe), 
then we say that $\mbK$ {\bf \em (weakly) admits $k$-substitutions}. 
}\end{defin}

\noindent
When speaking of a substitution $[\mcA \trir \mcB]$ we always assume that
$\mcA$ and $\mcB$ have the same universe.
Note that if every permitted structure is represented, 
which is the case if $\mbK$ has the hereditary property,
then $\mbK$ admits a substitution $[\mcA \trir \mcB]$ if and only if
$\mbK$ weakly admits $[\mcA \trir \mcB]$.

\begin{lem}\label{substitutions in permitted structures}
Suppose that $\mcA$ and $\mcB$ are permitted structures, with respect to $\mbK$,
with the same universe
such that the substitution $[\mcA \trir  \mcB]$ is weakly admitted with respect to $\mbK$.
Then for every permitted $\mcM$ such that $\mcA \subseteq \mcM$,
$\mcM[\mcA \trir \mcB]$ is permitted with respect to $\mbK$.
\end{lem}

\noindent
{\em Proof.}
Suppose that $\mcA$,  $\mcB$, $\mcM$ satisfy the premises of the lemma and that
$[\mcA \trir \mcM]$ is weakly admitted.
Since $\mcM$ is permitted there is a represented structure  $\mcN$ such that
$\mcM \subseteq \mcN$ (recall that the class of represented structures is closed under isomorphism).
Since the substitution $[\mcA \trir \mcB]$ is weakly admitted,
$\mcN[\mcA \trir \mcB]$ is permitted, so there is a represented $\mcN'$
such that $\mcN[\mcA \trir \mcB] \subseteq \mcN'$.
By assumption $\mcA \subseteq \mcM \subseteq \mcN$, so $\mcM[\mcA \trir \mcB] \subseteq \mcN[\mcA \trir \mcB]
\subseteq \mcN'$, which means that $\mcM[\mcA \trir \mcB]$ is permitted (because 
$\mcN'$ is represented).
\hfill $\square$

\begin{rem}\label{remark about failure of weak admittance of substitutions}{\rm
Suppose that $\rho$ is the supremum of the arities of relation symbols
in the vocabulary of $L$.\\
(i) It is straighforward to see that if $\mbK$ admits $\rho$-substitutions,
then $\mbK$ admits $k$-substitutions for every $k \in \mbbN$;
because every $k$-substitution can be achieved by performing, in sequence,
finitely many $\rho$-substitutions.\\
(ii) By using Lemma~\ref{substitutions in permitted structures}, it follows, much
as in (i), that if $\mbK$ weakly admits $\rho$-substitutions,
then $\mbK$ weakly admits $k$-substitutions for every $k$.
}\end{rem}

\noindent
Remember that a structure $\mcM$ satisfies the $\mcB/\mcA$-extension axiom
if and only if the $\mcB/\mcA$-multiplicity of $\mcM$ is at least $1$.
The next theorem is essentially a rephrasing, with the terminology
of this article, of a result of which has been used to prove
that every nontrivial parametric class of $L$-structures has a labeled zero-one law 
(\cite{Ober}, \cite{EF} Theorem 4.2.3).
A class $\mbC$ of finite $L$-structures is {\em nontrivial and parametric}, in the
sense of \cite{EF, Ober}, if and only if $\mbC$ is the class of represented structures
with respect to some $\mbK = \bigcup_{n \in \mbbN} \mbK_n$ which admits
$k$-substitutions for every $k$, and $\mbK_n$ is a nonempty set 
of $L$-structures with universe $\{1, \ldots, n\}$.
The result refered to in \cite{EF, Ober} is a generalization of
the well-known zero-one law for `random structures' \cite{Fag, Gleb}, which in the present
context amounts to considering the uniform probability measure on the 
set $\mbK_n$ of {\em all} $L$-structures with universe $\{1, \ldots, n\}$
(assuming that the arity of at least one relation symbol, besides `=', is greater than 1).

\begin{theor}\label{multiplicity when substitutions are admitted}
\cite{EF, Fag, Ober}
Let $\mbK = \bigcup_{n \in \mbb{N}} \mbK_n$ where each $\mbK_n$ is a nonempty set
of $L$-structures with universe $\{1, \ldots, m_n\}$ and $\lim_{n \to \infty} m_n = \infty$.
Suppose that $\mbK$ admits $k$-substitutions and let $p$ be any positive integer.
Whenever $\mcA \subset \mcB$ are permitted and $\left| \mcB \right| \leq k$, then
the proportion of structures $\mcM \in \mbK_n$
such that the $\mcB/\mcA$-multiplicity of $\mcM$ is at least $p$ approaches 1
as $n$ approaches $\infty$.
\end{theor}

\noindent
(The proof of the zero-one law by Glebski et al. \cite{Gleb}
does not use extension axioms, but a form of quantifier elimination,
which is why their article is not cited in Theorem~\ref{multiplicity when substitutions are admitted}.)
Since Theorem~\ref{multiplicity when substitutions are admitted} is not quite the
same as similar results refered to \cite{EF, Fag, Ober}, 
we give a sketch of its proof.
\\

\noindent
{\em Proof sketch.}
For simplicity, consider the  case when
$\left\| \mcB \right\| = \left\| \mcA \right\| + 1 \leq k$
and $p = 2$.
For every positive $d \in \mbbN$, let $\alpha_d$ denote the
number of different permitted structures with universe $\{1, \ldots, d\}$.
Let $\mcM \in \mbK_n$.
Since $\mbK$ admits $k$-substitutions it follows that
for every $d \leq k$, every permitted structure $\mcP$ with universe $\{1, \ldots, d\}$, 
all distinct $i_1, \ldots, i_d \in \{1, \ldots, m_n\}$, and $\mcM \in \mbK_n$,
the probability that $j \mapsto i_j$ is an embedding of $\mcP$ into $\mcM$
is $1/\alpha_d$, with the uniform probability measure. 
Suppose that $\mcA'$ is a copy of $\mcA$ with universe 
$A' = \{i_1, \ldots, i_d\} \subset M = \{1, \ldots, m_n\}$, so $d < k$.
For every $j \in \{1, \ldots, \lfloor m_n/2 \rfloor\} - A'$, the probability
that $\mcM \uhrc \{i_1, \ldots, i_d, j\}$
is a copy of  $\mcB$ is at least $1/\alpha_{d+1}$.
Therefore the probability that there is no 
$j \in \{1, \ldots, \lfloor m_n/2 \rfloor\} - A'$
such that  this holds is at most $\big(1 - 1/\alpha_{d+1}\big)^{\lfloor m_n/ 2 \rfloor - d}$.
There are at most $\binom{m_n}{d}$ copies of $\mcA$ in $\mcM$
and therefore the probability that some copy $\mcA' \subseteq \mcM$ of $\mcA$ cannot be extended
to a copy of $\mcB$ by adding an element from 
$\{1, \ldots, \lfloor m_n/2 \rfloor\} - A'$ is at most
$\binom{m_n}{d}\big(1 - 1/\alpha_{d+1}\big)^{\lfloor m_n / 2 \rfloor - d}$
which approaches 0 as $n$ approaches $\infty$ 
(because we assume that $\lim_{n \to \infty} m_n = \infty$).
In the same way, the probability that some copy $\mcA' \subseteq \mcM$ of $\mcA$ cannot be extended
to a copy of $\mcB$ by adding an element from 
$\{\lfloor m_n/2 \rfloor +1, \ldots, m_n\} - A'$
approaches 0 as $n \to \infty$. 
It follows that the probability that the $\mcA/\mcB$-multiplicity of $\mcM \in \mbK_n$ 
is less than 2 approaches 0 as $n \to \infty$.
\hfill $\square$
\\

\noindent
With Theorem~\ref{multiplicity when substitutions are admitted} 
at hand it remains to study what
happens, asymptotically, with extension axioms and multiplicities when
there are permitted $\mcA$ and $\mcB$ (with the same universe) such that
the substitution $[\mcA \trir \mcB]$ is not admitted with respect to $\mbK$.
The assumption that, for some permitted 
$\mcA$ and $\mcB$, the substitution $[\mcA \trir \mcB]$ is not admitted
is not enough, even if we assume that $\mbK$ has the
hereditary property and disjoint amalgamation property, to produce an extension
axiom $\varphi$ of $\mbK$ such that the proportion of structures in $\mbK_n$ satisfying $\varphi$
does {\em not} approach 1 as $n \to \infty$.
In this context it may, or may not, be the case that for every extension axiom,
the proportion of structures in $\mbK_n$ in which it is true approaches 1.
Examples~\ref{complete bipartite graphs} and~\ref{equivalence relations} show this.

But if $\mbK$ has the hereditary property and disjoint amalgamation property and
there are permitted $\mcA$ and $\mcB$ with the same universe such
that $[\mcA \trir \mcB]$ is admitted, and permitted $\mcM$ such
that  $\mcM[\mcB \trir \mcA]$ is forbidden,
then (by Corollary~\ref{simplification of upper limits theorem}) the proportion of
structures in $\mbK_n$ which satisfy all $(2|M| - 1)$-extension axioms never exceeds some $c < 1$;
and if there are no unary relation symbols, then this proportion approaches 0 as $n \to \infty$. 
If we do not assume that $\mbK$ has the hereditary and disjoint amalgamation properties,
then we can still obtain a related result 
(Theorem~\ref{upper limit on the probability of an extension axiom})
if we add another assumption on $\mcA$ and $\mcB$.
In the case that $\mbK$ has the hereditary property and disjoint amalgamation property,
Lemma~\ref{minimal substitutions}, below, implies that we can find
permitted $\mcA$ and $\mcB$ which satisfy this added assumption.

Recall that if $\mbK$ has the hereditary property,
then the notions `permitted structure' and `represented structure' coincide,
and therefore the notions `admit' (some substitution) and `weakly admit' (the same substitution)
coincide.

\begin{lem}\label{minimal substitutions}
Suppose that $\mbK$ has the hereditary property and the disjoint amalgamation property.
Assume that $\mcA$ and $\mcB$ are permitted structures with the same universe and
that the substitution $[\mcA \trir \mcB]$ is admitted, but $[\mcB \trir \mcA]$ is not admitted.
Then there are permitted $\mcA'$ and $\mcB'$ such that
\begin{itemize}
	\item[(1)] $A' = B' \subseteq A$,
	\item[(2)] the substitution $[\mcA' \trir \mcB']$ is admitted but $[\mcB' \trir \mcA']$ is not admitted, and
	\item[(3)] for every proper subset $U \subset A'$, $\mcA' \uhrc U = \mcB' \uhrc U$.
\end{itemize}
Moreover, if $\mcM$ is permitted and $\mcM[\mcB \trir \mcA]$ is forbidden,
then there is permitted $\mcM'$ with $M' = M$ such that $\mcM'[\mcB' \trir \mcA']$
is forbidden.
\end{lem}

\noindent
{\em Proof.}
With `$\subset$' we mean `proper subset' or `proper substructure'.
It suffices to prove that if $\mcA' = \mcA$ and $\mcB' = \mcB$ do not satisfy
(1) -- (3), then there is $U \subset A$ such that 
if $\mcA' = \mcA \uhrc U$ and $\mcB' = \mcB \uhrc U$, then
$[\mcA' \trir \mcB']$ is admitted but $[\mcB' \trir \mcA']$ is not admitted.
(Because if $A' = B'$ is a singleton set, then (3) trivially holds.)
The last statement of the lemma will follow from the proof that there exist
$\mcA'$ and $\mcB'$ satisfying (1) -- (3).

First we prove the following:
\\

\noindent
{\em Claim.} If $U \subset A$, $\mcU = \mcA \uhrc U$ and $\mcV = \mcB \uhrc U$,
then the substitution $[\mcU \trir \mcV]$ is admitted.
\\

\noindent
{\em Proof of Claim.}
Let $\mcM$ be any permitted structure and suppose that $\mcU \subseteq \mcM$.
We need to show that $\mcM[\mcU \trir \mcV]$ is permitted.
By the disjoint amalgamation property there is a permitted $\mcC$ such that
$\mcA \subseteq \mcC$, $\mcM \subseteq \mcC$ and $A \cap M = U$.
Since $[\mcA \trir \mcB]$ is admitted, $\mcC[\mcA \trir \mcB]$ 
is permitted. From $\mcM \subseteq \mcC$ and $A \cap M = U$ we  get
$\mcM[\mcU \trir \mcV] \subseteq \mcC[\mcA \trir \mcB]$,
so $\mcM[\mcU \trir \mcV]$ is permitted (and hence represented).
\hfill $\square$

Suppose that for some $U \subset A$, $\mcA \uhrc U \neq \mcB \uhrc U$.
(Otherwise $\mcA' = \mcA$, $\mcB' = \mcB$ satisfy (1) -- (3).)
Let $U_1, \ldots, U_l$ be an enumeration of all proper subsets $U_i \subset A = B$ such
that $\mcA \uhrc U_i \neq \mcB \uhrc U_i$.
By assumption there is a permitted $\mcM$ such that $\mcB \subset \mcM$
and $\mcN = \mcM[\mcB \trir \mcA]$ is forbidden.
For $i = 1, \ldots, l$, let $\mcU_i = \mcA \uhrc U_i$ and, by induction, 
define $\mcN_0 = \mcM$, $\mcV_1 = \mcM \uhrc U_1$,
$\mcN_{i+1} = \mcN_i[\mcV_{i+1} \trir \mcU_{i+1}]$ and
$\mcV_{i+1} = \mcN_i \uhrc U_{i+1}$.
Let $\mcA' = \mcN_l \uhrc A$.
Then $\mcN = \mcM[\mcB \trir \mcA] = \mcN_l[\mcA' \trir \mcA]$.

If every one of the substitutions
$[\mcV_1 \trir \mcU_1], \ldots, [\mcV_l \trir \mcU_l]$ and $[\mcA' \trir \mcA]$
is admitted, then $\mcN$ is permitted, which contradicts the assumption about $\mcN$.
First suppose that for some $i$, the substitution
$[\mcV_i \trir \mcU_i]$ is not admitted.
By the claim, $[\mcU_i \trir \mcV_i]$ is admitted, so we are done
(remember the first paragraph of the proof).

Now suppose that for every $i$, the substitution $[\mcV_i \trir \mcU_i]$ is admitted,
and consequently $[\mcA' \trir \mcA]$ is not admitted.
By the definition of $\mcA'$, we have $\mcA' \uhrc U = \mcA \uhrc U$ for
every $U \subset A = A'$. 
Hence, we are done if we can show that the substitution $[\mcA \trir \mcA']$ is admitted.
By the definition of $\mcU_i$, $\mcV_i$, $i = 1, \ldots, l$ and $\mcA'$, the
result of the substitution $[\mcB \trir \mcA']$, in any permitted structure,
can be achieved by performing the substitutions
$$[\mcV_1 \trir \mcU_1], \ldots, [\mcV_l \trir \mcU_l]$$
sequentially in the order from left to right.
By assumption, every substitution $[\mcV_i \trir \mcU_i]$ is admitted,
and hence $[\mcB \trir \mcA']$ is admitted.
Since the result of the substitution $[\mcA \trir \mcA']$ can be achieved
by first performing the substitution $[\mcA \trir \mcB]$,
which is admitted by assumption, and then $[\mcB \trir \mcA']$,
it follows that $[\mcA \trir \mcA']$ is admitted.

Now we verify the last statement of the lemma.
When starting with $\mcM$ and then performing
the substitutions $[\mcV_1 \trir \mcU_1]$, $\ldots$, $[\mcV_l \trir \mcU_l]$ and
$[\mcA' \trir \mcA]$ in this order, then, since $\mcN$ is forbidden,
there is a {\em first} structure during this process which is forbidden.
For $\mcM'$ we take the last structure during the process such that it and every
structure before it is permitted.
\hfill $\square$

\begin{theor}\label{upper limit on the probability of an extension axiom}
Assume that $\mcP$, $\mcS_{\mcP}$ and $\mcS_{\mcF}$ are permitted structures such that
$\mcS_{\mcP} \subseteq \mcP$, $|\mcS_{\mcP}| = |\mcS_{\mcF}|$, $\left\| \mcS_{\mcP} \right\| = k$,
$\mcF = \mcP[\mcS_{\mcP} \trir \mcS_{\mcF}]$ is forbidden, but
the substitution $[\mcS_{\mcF} \trir \mcS_{\mcP}]$ is admitted.
Moreover, assume that for
every proper substructure $\mcU \subset \mcS_{\mcP}$, $\mcS_{\mcP} \uhrc |\mcU| = \mcS_{\mcF} \uhrc |\mcU|$.
Let $\alpha$ be the number of different permitted structures with
universe $\{1, \ldots, k\}$ (so $\alpha \geq 2$).\\
(i) For every $n$, the proportion of $\mcM \in \mbK_n$ such that
	\begin{itemize}
	\item[(a)] $\mcM$ contains a copy of $\mcS_{\mcF}$, and
	\item[(b)] the $\mcP/\mcS_{\mcP}$-multiplicity of $\mcM$ is at least 2
	\end{itemize}
	never exceeds $1- 1/(1 + \alpha)$.\\
(ii) Suppose that there exist a permitted structure $\mcC$ and embeddings
$\sigma_1 : \mcP \to \mcC$ and $\sigma_2 : \mcP \to \mcC$ such that
$\sigma_1(|\mcP|) \cap \sigma_2(|\mcP|) = \sigma_1(|\mcS_{\mcP}|)$
and $\sigma_1 \uhrc |\mcS_{\mcP}| = \sigma_2 \uhrc |\mcS_{\mcP}|$.
Then, for every $n$, the proportion of $\mcM \in \mbK_n$ that satisfy
all $(2\left\| \mcP \right\| - k - 1)$-extension axioms never exceeds $1- 1/(1 + \alpha)$.\\
(iii) Suppose that $L$ has no unary relation symbols.
The proportion of $\mcM \in \mbK_n$ such that
\begin{itemize}
	\item[(c)] $\mcM$ satisfies the $\mcS_{\mcF}/\mcU$-extension axiom, where $\mcU \subseteq \mcS_{\mcF}$
	and $\left\|\mcU\right\| = 1$, and
	\item[(d)] the $\mcP/\mcS_{\mcP}$-multiplicity of $\mcM$ is at least 2
\end{itemize}
approaches 0 as $n$ approaches $\infty$.
\end{theor}

\begin{cor}\label{simplification of upper limits theorem}
Suppose that $\mbK$ has the hereditary property and the disjoint amalgamation property.
Also assume that there are permitted structures
$\mcA$, $\mcB$ and $\mcM$ such that $A = B$ and
the substitution $[\mcA \trir \mcB]$ is admitted, but $\mcM[\mcB \trir \mcA]$ is forbidden.\\
Then the proportion of
structures in $\mbK_n$ which satisfy all $(2|M| - 1)$-extension axioms
never exceeds $1 - 1/(1 + \alpha)$, where $\alpha$ is the number of permitted structures
with universe $A$.
If the language has no unary relation symbols then this proportion approaches 0 as $n \to \infty$.
\end{cor}

\noindent
{\em Proof of Corollary~\ref{simplification of upper limits theorem}.}
Assume that $\mbK$ has the hereditary property and disjoint amalgamation
property, and let $\mcA$, $\mcB$ and $\mcM$ satisfy the assumptions of the corollary.
From Lemma~\ref{minimal substitutions} it follows that there are permitted structures
$\mcP$, $\mcS_{\mcP}$ and $\mcS_{\mcF}$ which satisfy the assumptions of 
Theorem~\ref{upper limit on the probability of an extension axiom} and 
$|\mcS_{\mcP}| \subseteq |\mcA|$ and $|\mcP| = |\mcM|$.
Since $\mbK$ has the disjoint amalgamation property, part (ii) of
Theorem~\ref{upper limit on the probability of an extension axiom} 
implies that the proportion of structures in $\mbK_n$ which 
satisfy all $(2\left\|\mcP\right\| - \left\|\mcS_{\mcP}\right\| - 1)$-extension axioms
never exceeds $1 - 1/(1 + \alpha')$,
where $\alpha'$ is the number of permitted structures with universe 
$|\mcS_{\mcP}|$.
Note that if $\alpha$ is the number of permitted structures with universe
$A$, then, since $\left\|\mcS_{\mcP}\right\| \leq |A|$, we have
$1 - 1/(1 + \alpha') \leq 1 - 1(1 + \alpha)$.

Every structure in $\mbK$ which satisfies all 
$(2\left\|\mcP\right\| - \left\|\mcS_{\mcP}\right\| - 1)$-extension axioms
satisfies both (c) and (d) in part (iii) of 
Theorem~\ref{upper limit on the probability of an extension axiom}.
So if the language has no unary relation symbols the proportion
of structures in $\mbK_n$ which satisfy all 
$(2\left\|\mcP\right\| - \left\|\mcS_{\mcP}\right\| - 1)$-extension axioms
must approach 0 as $n \to \infty$.
Since $2|M| - 1 \geq 2|P| - 1 \geq  2\left\|\mcP\right\| - \left\|\mcS_{\mcP}\right\| - 1$
we are done.
\hfill $\square$

\section{Examples}\label{examples concerning notion of admittance}

\noindent
In all examples, $\mbK = \bigcup_{n\in \mbbN} \mbK_n$ 
has the hereditary property.

\begin{exam}\label{example of forbidden weak substructures}{\rm
({\bf Forbidden weak substructures and proof of Theorem~\ref{dichotomy for forbidden weak substructures}
from Theorem~\ref{upper limit on the probability of an extension axiom}.})
Let $L$ have a finite relational vocabulary,
and let $\mbF$ be a set of finite $L$-structures.
For $n \in \mbbN$,
let $\mbK_n$ be the set of all $L$-structures $\mcM$ with universe $\{1, \ldots, n\}$
such that {\em no} $\mcF \in \mbF$ can be weakly embedded into $\mcM$.
Then a structure $\mcA$ is forbidden if and only if some $\mcF \in \mbF$ can
be weakly embedded into $\mcA$.
It follows that there exists (at least) one {\bf \em minimal forbidden} structure 
$\mcF_{min} \in \mbF$ in the sense that every {\em proper} weak substructure of $\mcF_{min}$
is permitted.

If $\mcF_{min}$ does not have any relationship at all, that is, if $\mcF_{min}$ is 
just a finite set of cardinality $m$, say,
then $\mbK_n = \es$ for every $n \geq m$. Since we are only interested in the case
when $\mbK_n \neq \es$ for arbitrarily large $n \in \mbbN$, we now assume that
every minimal forbidden structure has at least one relationship.
From this it follows that $\mbK_n \neq \es$ for every $n \in \mbbN$,
because the assumption ensures that the set $\{1, \ldots, n\}$ without any 
structure belongs to $\mbK_n$.

Let $\mcF_{min}$ be any minimal forbidden structure.
By assumption, for some relation symbol $R$, $R^{\mcF_{min}}$ is nonempty,
so we can remove a relationship $\bar{a}$ from $R^{\mcF_{min}}$ and
call the resulting structure $\mcP$.
Note that $\mcP$ is permitted (since $\mcF_{min}$ is minimal forbidden), and that
$\mcF_{min}$ and $\mcP$ have the same universe which includes $\rng(\bar{a})$.
Let $\mcS_{\mcF} = \mcF_{min} \uhrc \rng(\bar{a})$ and
$\mcS_{\mcP} = \mcP \uhrc \rng(\bar{a})$.
Then $\mcS_{\mcP}$ is permitted, because it is a substructure of $\mcP$, and $\mcP$ is permitted.
If $\rng(\bar{a}) = |\mcF_{min}|$ then $\mcS_{\mcF} = \mcF_{min}$ which is forbidden.
If this holds for {\em every} choice of minimal forbidden $\mcF_{min}$ and
$\mcS_{\mcF}$ as defined above, then ($\ast$) in 
Theorem~\ref{dichotomy for forbidden weak substructures} does not hold, and it is straightforward
to verify that $\mbK$ admits $k$-substitutions for every $k$.
In this case, Theorem~\ref{multiplicity when substitutions are admitted}
implies that, for every extension axiom $\varphi$ of $\mbK$,
the proportion of $\mcM \in \mbK_n$ which satisfy $\varphi$ approaches 1 as $n \to \infty$;
and hence $\mbK$ has a zero-one law for the uniform measure, 
by Remark~\ref{remark about extension axioms and 0-1 laws}.

Now suppose that there is a minimal forbidden $\mcF_{min}$ and $R$ such that for
some $\bar{a} \in R^{\mcF_{min}}$, $\rng(\bar{a})$ is a proper subset of  $|\mcF_{min}|$.
Then ($\ast$) in 
Theorem~\ref{dichotomy for forbidden weak substructures} holds and
$\mcS_{\mcF} = \mcF_{min} \uhrc \rng(\bar{a})$ is a {\em proper} substructure of
$\mcF_{min}$, and since the latter is minimal forbidden,
$\mcS_{\mcF}$ is permitted.
Hence $\mcP$, $\mcS_{\mcP}$ and $\mcS_{\mcF}$ are permitted,
but $\mcF_{min} = \mcP[\mcS_{\mcP} \trir \mcS_{\mcF}]$ is forbidden.
(The notions `admitted' and `weakly admitted' coincide here because
the notions `permitted' and `represented' coincide in this example.)
But since the removal of a relationship from a permitted structure will never (in the present context)
produce a forbidden structure, the substitution
$[\mcS_{\mcF} \trir \mcS_{\mcP}]$ is admitted.
Moreover, by the definition of $\mcS_{\mcF}$ and $\mcS_{\mcP}$,
they agree on all proper subsets of their common universe.
Thus, Theorem~\ref{upper limit on the probability of an extension axiom}
is applicable.
By part (i) of  Theorem~\ref{upper limit on the probability of an extension axiom},
the proportion of $\mcM \in \mbK_n$ such that
\begin{itemize}
	\item[(a)] $\mcM$ contains a copy of $\mcS_{\mcF}$, and
	\item[(b)] the $\mcP/\mcS_{\mcP}$-multiplicity of $\mcM$ is at least 2
\end{itemize}
never exceeds $1 - 1(1 + \alpha)$,
where $\alpha$ is the number of permitted structures with 
universe $\{1, \ldots, |\rng(\bar{a})|\}$.
If $\mbK$ has the disjoint amalgamation property (which is assumed in
Theorem~\ref{dichotomy for forbidden weak substructures}),
then part (ii) of 
Theorem~\ref{upper limit on the probability of an extension axiom}
is applicable, and it follows that the proportion of structures in $\mbK$
which satisfy all $(2|P| - |\rng(\bar{a})| - 1)$-extension axioms
never exceeds $1 - 1/(1 + \alpha)$.
And if the language has no unary relation symbols and $\mbK$ has
the disjoint amalgamation property, then this proportion approaches 0 as
$n \to \infty$, by part (iii) of 
Theorem~\ref{upper limit on the probability of an extension axiom}.
Note that $\left\|\mcF_{min}\right\| = |P|$, so 
Theorem~\ref{dichotomy for forbidden weak substructures} is proved.
}\end{exam}

\begin{exam}\label{example of necessity of no unary symbols}{\rm
This example shows that when, 
in Theorem~\ref{upper limit on the probability of an extension axiom},
it is assumed that the language has
no unary relation symbols,
then this assumption is necessary.
(The author does not have a corresponding example if one adds the assumption
that $\mbK$ has the disjoint amalgamation property, as in 
Theorem~\ref{dichotomy for forbidden weak substructures} 
and Corollary~\ref{simplification of upper limits theorem}.)

Let $P_1$ and $P_2$ be unary relation symbols and let
$L$ be a language the vocabulary of which is finite, relational and contains $P_1$ and $P_2$.
For $n \in \mbbN$, let $\mbK_n$ consist of all $L$-structures $\mcM$ with 
universe $\{1, \ldots, n\}$ such that 
\begin{align*}
&\text{at most one element in $M$ satisfies $P_1(x)$,}\\
&\text{at most one element in $M$ satisfies $P_2(x)$, and}\\
&\mcM \models \neg\exists x, y \big(P_1(x) \wedge P_2(y)\big).
\end{align*}
$\mbK_n$ can also be described in the following way,
by forbidden weak substructures.
Let $\mcA$, $\mcB$ and $\mcC$ have universe $\{1,2\}$ and the following
interpretations: $(P_1)^{\mcA} = \{1\}$, 
$(P_2)^{\mcA} = \{2\}$,
$(P_1)^{\mcB} = \{1, 2\}$, $(P_2)^{\mcB} = \es$,
$(P_1)^{\mcC} = \es$, $(P_2)^{\mcC} = \{1,2\}$ 
and $R^{\mcA} = R^{\mcB} = R^{\mcC} = \es$ for every other relation symbol $R$
in the vocabulary. 
Then $\mbK_n$ can also be described as the set of all $L$-structures $\mcM$ such that
no $\mcF \in \mbF = \{\mcA, \mcB, \mcC\}$ is weakly embeddable in $\mcM$.
Note that $\mbF$ satisfies the condition labelled ($\ast$)
in Theorem~\ref{dichotomy for forbidden weak substructures}, so
if $\alpha$ is the number of permitted structures with universe $\{1\}$,
then the proportion of $\mcM \in \mbK_n$ which satisfy all
$2$-extension axioms never exceeds $1 - 1/(1+\alpha)$.
We have $\alpha \geq 3$, and if the only unary relation symbols of $L$
are $P_1$ and $P_2$ then $\alpha = 3$.

Next, we show that there exists a $0$-extension axiom 
(i.e. an $\mcN/\es$-extension axiom with $N$ a singleton set)
such that the proportion of $\mcM \in \mbK_n$ which satisfy it never exceeds $1/2$.
For every $n$, $\mbK_n$ can be partitioned into three 
parts: one part, $\mbX_n$, consisting of all $\mcM \in \mbK_n$ which satisfy $\exists x P_1(x)$;
another part, $\mbY_n$, consisting of all $\mcM \in \mbK_n$ which satisfy $\exists x P_2(x)$;
and a third part, $\mbZ_n$, consisting of all $\mcM \in \mbK_n$
which do not satisfy either of $\exists x P_1(x)$ or $\exists x P_2(x)$.
The definition of $\mbK_n$ implies that, for each $n$, 
$|\mbX_n| = |\mbY_n| = n|\mbZ_n|$.
Let $\mcA' = \mcA \uhrc \{1\}$.
Then $\mcA'$ is permitted and $\mcA' \models P_1(1)$.
Moreover, for every $n$, the $\mcA'/\es$-extension axiom holds exactly for those
$\mcM \in \mbK_n$ which belong to $\mbX_n$, and we have
\begin{equation*}
\frac{|\mbX_n|}{|\mbK_n|} = \frac{|\mbX_n|}{|\mbX_n| + |\mbY_n| + |\mbZ_n|} =
\frac{n|\mbZ_n|}{n|\mbZ_n| + n|\mbZ_n| + |\mbZ_n|} = \frac{1}{2 + 1/n} \ \to \ \frac{1}{2} \
\text{ as } n \to \infty.
\end{equation*}
}\end{exam}

\medskip

\noindent
Examples~\ref{graphs with a unary predicate}
and~\ref{partially coloured binary relations}
show how Theorem~\ref{upper limit on the probability of an extension axiom}
can be applied. 
They also provide contrast to 
Examples~\ref{example of conditional measure 1}
and~\ref{example of conditional measure 2}, where a
different probability measure is considered.

\begin{exam}\label{graphs with a unary predicate}{\rm
({\bf Graph with a restricted unary predicate})
Let the vocabulary of $L$ consist of a unary relation symbol $Q$ and
a binary relation symbol $R$.
Let $\mbK_n$ be the set of $L$-structures $\mcM$ with universe $\{1, \ldots, n\}$
such that $R^{\mcM}$ is irreflexive and symmetric (i.e. an undirected graph) and
$$\mcM \models \forall x, y \big(R(x,y) \rightarrow (\neg Q(x) \wedge \neg Q(y))\big).$$
We use notation which suggests how Theorem~\ref{upper limit on the probability of an extension axiom}
will be used.
Define $\mcS_{\mcP}$, $\mcS_{\mcF}$, $\mcP$ and $\mcF$ as follows: 
let $|\mcS_{\mcP}| = |\mcS_{\mcF}| = \{a\}$; $Q^{\mcS_{\mcP}} = R^{\mcS_{\mcP}} = \es$;
$Q^{\mcS_{\mcF}} = \{a\}$, $R^{\mcS_{\mcF}} = \es$;
$|\mcP| = \{a,b\}$, $Q^{\mcP} = \es$, $R^{\mcP} = \{(a,b), (b,a)\}$;
and $\mcF = \mcP[\mcS_{\mcP} \trir \mcS_{\mcF}]$.
Then $\mcS_{\mcP}$, $\mcS_{\mcF}$ and $\mcP$ are permitted, 
but $\mcF$ is forbidden, since $\mcF \models Q(a) \wedge R(a,b)$.
Hence, the substitution $[\mcS_{\mcP} \trir \mcS_{\mcF}]$ is not admitted,
but the reverse substitution $[\mcS_{\mcF} \trir \mcS_{\mcP}]$ is admitted,
because we can always remove a $Q$-relationship without producing a forbidden structure.

By Theorem~\ref{upper limit on the probability of an extension axiom} (i),
the proportion of $\mcM \in \mbK_n$ which contain a copy of $\mcS_{\mcF}$
and whose $\mcP / \mcS_{\mcP}$-multiplicity is at least two is not larger than 
$1 - 1/(1+2) = 2/3$.
In this example we can do much better, asymptotically speaking,
and show that the proportion of $\mcM \in \mbK_n$ which
contain a copy of $\mcS_{\mcF}$, or equivalently, which satisfy $\exists x Q(x)$,
approaches 0 as $n \to \infty$.
We can argue as follows to see this.
First let
\begin{align*}
\mbX_n &= \big\{\mcM \in \mbK_n : \text{ $Q(x)$ is satisfied by at least two elements in $|\mcM|$} \big\},\\
\mbY_n &= \big\{\mcM \in \mbK_n : \text{ $Q(x)$ is satisfied by a unique element in $|\mcM|$}\big\}.
\end{align*}
Since 
$$\frac{|\mbY_n|}{|\mbK_n|} \leq \frac{n2^{\binom{n-1}{2}}}{2^{\binom{n}{2}}} = \frac{n}{2^{n-1}} \to 0,$$
as $n \to \infty$,
it is sufficient to show that $|\mbX_n| \leq |\mbY_n|$. 
For $\mcM \in \mbX_n$, let $a \in |\mcM| = \{1, \ldots, n\}$ be minimal
such that $\mcM \models Q(a)$, and let $\mcM'$ be defined as follows:
$|\mcM'| = \{1, \ldots, n\}$, $Q^{\mcM'} = \{a\}$ and let
$R^{\mcM'}$ be the symmetric closure of
$$R^{\mcM} \cup \big\{(b,c) : b \in Q^{\mcM} - \{a\}, \ \ c \in \{1, \ldots, n\} - Q^{\mcM} \big\}.$$
Note that $\mcM' \in \mbY_n$.
It is now easy to verify that the map 
$\mcM \mapsto \mcM'$ from $\mbX_n$ to $\mbY_n$ is injective; thus $|\mbX_n| \leq |\mbY_n|$.

Since $\{\mcM \in \mbK_n : \mcM \models \neg\exists x Q(x)\}$ is  the set of all (undirected)
graphs with vertices $1, \ldots, n$ it follows that, with the uniform probability measure, 
the almost sure theory of $\mbK$ is 
identical to the almost sure theory of all undirected graphs, and consequently
$\mbK$ has a zero-one law for the uniform probability measure.
Since the complete theory of the Fra\"{i}ss\'{e}-limit of $\mbK$ contains the sentence
$\exists x Q(x)$ it is different from the almost sure theory of $\mbK$, for
the uniform measure.
As we will see later, for the `dimension conditional probability measure'
(where dimension equals cardinality in this example),
the almost sure theory of $\mbK$ is identical to the complete theory of the
Fra\"{i}ss\'{e}-limit of $\mbK$.
}\end{exam}

\begin{exam}\label{partially coloured binary relations}{\rm
({\bf Partially coloured binary relation.})
Let the vocabulary of $L$ consist of one binary relation symbol $R$ 
and two unary relation symbols $P_1$, $P_2$.
Let $\mbK_n$ consist of all $L$-structures $\mcM$ with universe $\{1, \ldots, n\}$ such that
\begin{align*}
&\mcM \models \forall x \neg\big(P_1(x) \wedge P_2(x)\big), \ \text{ and}\\
&\mcM \models \forall x, y \big(R(x,y) \ \rightarrow \ \big[\neg\big(P_1(x) \wedge P_1(y)\big)  
\ \wedge \ \neg\big(P_2(x) \wedge P_2(y)\big)\big]\big).
\end{align*}
We can think of $P_i$ as representing the colour `$i$'.
Before using Theorem~\ref{upper limit on the probability of an extension axiom}
to get some information about $\mbK_n$ we consider the proportion
of $\mcM \in \mbK_n$ which satisfy $\exists x P_i(x)$.
Let 
$$\mbX_n = \big\{\mcM \in \mbK_n : (P_1)^{\mcM} = \es \big\}.$$
For every $\mcM \in \mbX_n$, let 
$\mbY_n(\mcM)$ be the set of $\mcN \in \mbK_n$ which satisfy that
$R^{\mcN} = R^{\mcM}$ and either 
\begin{itemize}
\item[\textbullet] \ $(P_1)^{\mcN} = \{a\}$ and $(P_2)^{\mcN} = (P_2)^{\mcM}$ for some $a \notin (P_2)^{\mcM}$, or
\item[\textbullet] \ $(P_1)^{\mcN} = \{a\}$ and $(P_2)^{\mcN} = (P_2)^{\mcM} - \{a\}$ for some $a \in (P_2)^{\mcM}$.
\end{itemize}
It is straightforward to verify that, for every $\mcM \in \mbX_n$, $|\mbY_n(\mcM)| \geq n$,
and, for every $\mcN \in \mbK_n$, the number of $\mcM \in \mbX_n$ such that
$\mcN \in \mbY_n(\mcM)$ is at most 2.
It follows that
\begin{equation*}
\frac{|\mbX_n|}{|\mbK_n|} \ \leq \ \frac{|\mbX_n|}{|\bigcup_{\mcM \in \mbX_n} \mbY_n(\mcM)|} \ \leq \
\frac{|\mbX_n|}{\frac{1}{2}\sum_{\mcM \in \mbX_n}|\mbY_n(\mcM)|} \ \leq \
\frac{|\mbX_n|}{\frac{1}{2}|\mbX_n|n} \ = \  \frac{2}{n} \ \to \ 0,
\end{equation*}
as $n \to \infty$, so the proportion of $\mcM \in \mbK_n$ such that $\mcM \models \exists x P_1(x)$
approaches 1 as $n \to \infty$.
The same argument works for $P_2$.

For an $L$-structure $\mcM$ and $a \in |\mcM|$, let us say that $a$ is {\em blank}
or {\em uncoloured} (in $\mcM$) if $\mcM \models \neg P_1(a) \wedge \neg P_2(a)$.
Let $\mcS_{\mcP}$ have universe $\{a\}$ where $a$ is blank in $\mcS_{\mcP}$ and $R^{\mcS_{\mcP}} = \es$.
Let $\mcS_{\mcF}$ also have universe $\{a\}$ where $a$ has colour 1 in $\mcS_{\mcF}$
(i.e. $\mcS_{\mcF} \models P_1(a)$) and $R^{\mcS_{\mcF}} = \es$.
Then $\mcS_{\mcP}$ and $\mcS_{\mcF}$ are permitted and it is easily seen that
the substitution $[\mcS_{\mcF} \trir \mcS_{\mcP}]$ is admitted, because making a point
blank never violates the conditions for being permitted (with respect to $\mbK$).
But if one point in an $R$-relationship is coloured by $i$, then colouring the other
point in the same $R$-relationship by the same colour $i$ produces a forbidden structure;
so the substitution $[\mcS_{\mcP} \trir \mcS_{\mcF}]$ is not admitted. 
Now we apply Theorem~\ref{upper limit on the probability of an extension axiom}.
Let $|\mcP| = \{a,b\}$, $(P_1)^{\mcP} = \{b\}$ and $R^{\mcP} = \{(a,b)\}$.
Since we know that the proportion of $\mcM \in \mbK_n$ which contain a copy of 
$\mcS_{\mcF}$ (i.e satisfy $\exists x P_1(x)$) approaches 1 as $n \to \infty$,
it follows that for arbitrarily small $\varepsilon > 0$ and all sufficiently large $n$,
the proportion of $\mcM \in \mbK_n$
such that the $\mcP/\mcS_{\mcP}$-multiplicity of $\mcM$ is at least 2 never exceeds 
$\big(1 - 1/(1+3)\big) + \varepsilon = 3/4 + \varepsilon$.
Observe that the $\mcP/\mcS_{\mcP}$-multiplicity of $\mcM$ is at least 2 
if and only if $\mcM$ satisfies the extension axiom
$$\varphi = \forall x \exists y, z \big( [\neg P_1(x) \wedge \neg P_2(x)] \ \rightarrow \
[y \neq z \wedge R(x,y) \wedge R(x,z) \wedge P_1(y) \wedge P_1(z)] \big),$$
so the probability, with the uniform probability measure, that this extension axiom is true never 
exceeds $3/4 + \varepsilon$.
}\end{exam}

\begin{exam}\label{coloured binary relations}{\rm
({\bf Coloured binary relation.})
Let $\mbK_n$ be defined as in Example~\ref{partially coloured binary relations}
{\em except} that we add the condition that there are {\em no} blank elements, that is,
every $\mcM \in \mbK_n$ satisfies $\forall x \big(P_1(x) \vee P_2(x)\big)$.
By Theorem~\ref{transfer of 0-1 laws to uniform measure} and 
Remark~\ref{remark concerning transfer theorem},
for every extension axiom $\varphi$ of $\mbK$, the proportion of $\mcM \in \mbK_n$
which satisfies $\varphi$ approaches 1 as $n \to \infty$.
Since $\mbK$ has the hereditary property and the disjoint amalgamation property, 
Lemma~\ref{minimal substitutions} and
Theorem~\ref{upper limit on the probability of an extension axiom} (part (ii))
implies that there does {\em not} exist permitted $\mcS_{\mcP}$ and $\mcS_{\mcF}$ such that
the substitution $[\mcS_{\mcF} \trir \mcS_{\mcP}]$ is admitted and
$[\mcS_{\mcP} \trir \mcS_{\mcF}]$ is not admitted.
However, since changing one colour to another in a permitted structure may produce
a forbidden structure, there are permitted $\mcA$ and $\mcA'$ (with singleton universes)
such that {\em none} of the substitutions $[\mcA \trir \mcA']$ and $[\mcA' \trir \mcA]$
is admitted.
}\end{exam}

\noindent
Examples~\ref{complete bipartite graphs}
and~\ref{equivalence relations}
(as well as Example~\ref{coloured binary relations}) 
show that if $\mbK$ neither satisfies the conditions of 
Theorem~\ref{multiplicity when substitutions are admitted}, 
nor the conditions of Theorem~\ref{upper limit on the probability of an extension axiom} 
(or Corollary~\ref{simplification of upper limits theorem}),
then it may, or may not, be the case that for every extension axiom
$\varphi$ of $\mbK$ the proportion of $\mcM \in \mbK_n$ which satisfy
$\varphi$ approaches 1 as $n \to \infty$.
In contrast to Examples~\ref{example of necessity of no unary symbols} 
--~\ref{coloured binary relations}, the last two examples of this section do
not have any unary relations.

\begin{exam}\label{complete bipartite graphs}{\rm
({\bf Complete bipartite graph.})
For all $r,s \in \mbbN$, let $\mcK_{r,s}$ denote the undirected graph 
with vertices $a_1, \ldots, a_r, b_1, \ldots, b_s$ and
an edge connecting $a_i$ and $b_j$ for all $i \in \{1, \ldots, r\}$
and $j \in \{1, \ldots, s\}$, and no other edges.
$\mcK_{0,s}$ and $\mcK_{r,0}$ are independent sets (no edges at all)
with $s$ and $r$ vertices, respectively.

For every $n \in \mbbN$, let $\mbK_n$ be the set of all graphs
with vertices $1, \ldots, n$ which are isomorphic to $\mcK_{r,s}$ for some $r,s$.
Clearly, by adding an edge to any represented $\mcM$ with at 
least 3 vertices, we create a forbidden graph.
Also, by removing an edge from any $\mcK_{r,s}$ such that $r+s\geq 3$ and
$\min(r,s) \geq 1$, we create a forbidden graph. 

It is easy to see that if $s, r \geq k+1$, then
$\mcK_{r,s}$ satisfies all $k$-extension axioms of $\mbK = \bigcup_{n\in \mbbN}\mbK_n$.
Also, the proportion of $\mcM \in \mbK_n$ which are isomorphic to some
$\mcK_{r,s}$ with $r,s \geq k+1$ approaches 1 as $n \to \infty$.
It follows that, for every extension axiom $\varphi$ of $\mbK$, the proportion of $\mcM \in \mbK_n$ which
satisfy $\varphi$ approaches 1 as $n \to \infty$.
It is straightforward to verify that the class of represented structures is 
closed under taking substructures (so `permitted' is the same as `represented')
and has the disjoint amalgamation property.
By Corollary~\ref{simplification of upper limits theorem}, there does
not exist any permitted $\mcA$ and $\mcB$ with $A = B$ such that $[\mcA \trir \mcB]$
is admitted and  $[\mcB \trir \mcA]$ is not admitted.
}\end{exam}

\begin{exam}\label{equivalence relations}{\rm
({\bf Equivalence relations})
Here we define $\mbK$ such that (as in the previous example) there are no permitted $\mcA$
and $\mcB$ such that $[\mcA \trir \mcB]$ is admitted and $[\mcB \trir \mcA]$ is not admitted.
In this example, $\mbK$ has an extension axiom $\varphi$  such that the proportion of
$\mcM \in \mbK_n$ in which $\varphi$ is true approaches 0 as $n \to \infty$,
but nevertheless $\mbK$ has a zero-one law.

We represent an equivalence relation on a set $M$ as 
an undirected graph (without loops) with  vertex set $M$ such that if 
$a$ is adjacent to $b$ and $b$ is adjacent to $c \neq a$,
then $a$ and $c$ are adjacent.
Clearly, such a graph, which we call an {\em equivalence graph}, is a disjoint
union of complete graphs.
Let $\mbK_n$ consist of all equivalence graphs with vertices $1, \ldots, n$. 
Equivalently, we could have defined $\mbK_n$ by saying that it consists of all 
undirected graphs with vertices $1, \ldots, n$ in which $\mcV$ is not embeddable, 
where $\mcV$ denotes the graph
with distinct vertices $1, 2, 3$ where $1$ is adjacent with $2$ and $2$ is adjacent
with $3$, but $1$ is not adjacent with $3$.
It is easily seen that $\mbK$ has the disjoint amalgamation property.
By Lemma~\ref{minimal substitutions}, if there would be 
$\mcA$, $\mcB$ with same universe such that $[\mcA \trir \mcB]$ is admitted,
but not $[\mcB \trir \mcA]$, then, because we only have a binary relation symbol,
we could assume that the common universe of $\mcA$ and $\mcB$ has cardinality 2,
and that $[\mcA \trir \mcB]$ means either to remove an edge, or to add an edge.
But it is clear that both the removal of an edge, as well as the addition of an edge,
may produce a forbidden structure, so {\em none} of $[\mcA \trir \mcB]$ and $[\mcB \trir \mcA]$
can be admitted, contradicting the assumption.
Hence, there does not exist $\mcA$ and $\mcB$ such that 
$[\mcA \trir \mcB]$ is admitted, but not $[\mcB \trir \mcA]$.

We now show that if $\mcA$ is the graph having only one vertex $a$ and
$\mcB$ has vertex set $\{a, b\}$ where $a$ is adjacent to $b$ in $\mcB$,
then the probability that $\mcM \in \mbK_n$ satisfies the $\mcB/\mcA$--extension
axiom approaches 0 as $n \to \infty$.
This contrasts the previous example.
Let $\mbX_n$ be the set of $\mcM \in \mbK_n$ which do {\em not} contain any
connected component which is a singleton, and let $\mbX = \bigcup_{n > 1}\mbX_n$.
Note that the class of represented structures with respect to $\mbX$
is closed under taking disjoint unions and extracting connected components;
thus, the class of represented structures with respect to $\mbX$ is {\em adequate}
in the sense of \cite{BY}, which we will use.
For every $n > 1$, $\mbX_n$ contains exactly one connected graph (the complete graph
with vertices $1, \ldots, n$). Therefore, Theorem~7 in \cite{BY} implies that
$$\frac{n|\mbX_{n-1}|}{|\mbX_n|} \to \infty \ \text{ as } \ n \to \infty.$$
Let $\mbY_n$ be the set of $\mcM \in \mbK_n$ that contain at least one
connected component which is a singleton, and let 
$\mbY'_n$ be the set of $\mcM \in \mbK_n$ that contain exactly one
connected component which is a singleton.
Observe that 
$$\mbX_n = \mbK_n - \mbY_n \ \text{ and } \ |\mbY'_n| = n|\mbX_{n-1}|.$$
It follows that 
$$\frac{|\mbX_n|}{|\mbK_n|} \leq \frac{|\mbX_n|}{|\mbY'_n|} =
\frac{|\mbX_n|}{n|\mbX_{n-1}|} \to 0 \ \text{ as } \ n \to \infty.$$
In other words, the proportion of $\mcM \in \mbK_n$ which contain at least
one connected component which is a singleton approaches 1 as $n \to \infty$.
For every such $\mcM$, the $\mcB/\mcA$--extension axiom fails.
Nevertheless, $\mbK$ has a zero-one law for the uniform probability measure,
which follows from Theorem~7 in \cite{BY} and the above observed fact that,
for every $n$, there is a unique connected graph in $\mbK_n$.
}\end{exam}

\section{Proof of Theorem~\ref{upper limit on the probability of an extension axiom}}
\label{proving upper limit on the probability of an extension axiom}

\noindent
Let $L$ have a finite relational vocabulary and let
$\mbK = \bigcup_{n \in \mbbN}\mbK_n$, where every $\mbK_n$ is a set
of $L$-structures with universe $\{1, \ldots, m_n\}$ and $\lim_{n\to\infty}m_n = \infty$.
Suppose that $\mcP$, $\mcS_{\mcP}$ and $\mcS_{\mcF}$ are permitted structures such that
$\mcS_{\mcP} \subseteq \mcP$, $|\mcS_{\mcP}| = |\mcS_{\mcF}|$, $\left\| \mcS_{\mcP} \right\| = k$,
$\mcF = \mcP[\mcS_{\mcP} \trir \mcS_{\mcF}]$ is forbidden, but
the substitution $[\mcS_{\mcF} \trir \mcS_{\mcP}]$ is admitted.
Morover, assume that for
every proper substructure $\mcU \subset \mcS_{\mcP}$, $\mcS_{\mcP} \uhrc |\mcU| = \mcS_{\mcF} \uhrc |\mcU|$.
Let $\alpha$ be the number of different permitted structures with
universe $\{1, \ldots, k\}$ (so $\alpha \geq 2$).

We use the following terminology:

\begin{defin}\label{definition of coexisting pair}{\rm
(i) A pair of structures $(\mcA, \mcB)$ is called a {\bf \em coexisting pair}
if $\mcA$ and $\mcB$ have the same universe.\\
(ii) We say that two coexisting pairs $(\mcA, \mcB)$ and $(\mcA',\mcB')$ are {\bf \em isomorphic}
if there is a bijection $\sigma : |\mcA| \to |\mcA'|$ which is 
an isomorphism from $\mcA$ to $\mcA'$ as well as from $\mcB$ to $\mcB'$.\\
(iii) If $(\mcA, \mcB)$ and $(\mcA',\mcB')$ are isomorphic coexisting pairs then we may
say that $(\mcA',\mcB')$ is a {\bf \em copy} of $(\mcA,\mcB)$.
}\end{defin}

\begin{lem}\label{reducing the multiplicity to 0}
Suppose that $\mcS_{\mcP}$ is a proper substructure of $\mcP$ and that $\mcM$ is represented.
If $(\mcS_{\mcP}^{\mcM}, \mcS_{\mcF}^{\mcM})$ is a copy of the coexisting pair
$(\mcS_{\mcP}, \mcS_{\mcF})$ and $\mcS_{\mcF}^{\mcM} \subseteq \mcM$, then
the $\mcP/\mcS_{\mcP}$-multiplicity of $\mcM[\mcS_{\mcF}^{\mcM} \trir \mcS_{\mcP}^{\mcM}]$ is 0.
\end{lem}

\noindent
{\em Proof.}
Without loss of generality (by just renaming elements)
we may assume that $\mcS_{\mcF} = \mcS_{\mcF}^{\mcM} \subseteq \mcM$
and that $\mcS_{\mcP} = \mcS_{\mcP}^{\mcM}$.
Then $\mcS_{\mcF}$ ($= \mcS_{\mcF}^{\mcM}$) is a substructure
of $\mcM$, and $\mcS_{\mcF}$ and $\mcS_{\mcP}$ ($= \mcS_{\mcP}^{\mcM}$)
have the same universe which is a subset of $|\mcM|$.
By the assumption that 
$\mcP$ is permitted, but $\mcF = \mcP[\mcS_{\mcP} \trir \mcS_{\mcF}]$ is forbidden 
(see before Definition~\ref{definition of coexisting pair})
we have
$\mcS_{\mcF} \neq \mcS_{\mcP}$, 
and, as $\mcS_{\mcF} \subseteq \mcM$ and $|\mcS_{\mcF}| = |\mcS_{\mcP}|$,
$\mcS_{\mcP}$ is not a substructure of $\mcM$.
But $\mcS_{\mcP}$ is a substructure of $\mcM[\mcS_{\mcF} \trir \mcS_{\mcP}]$.

We show that the $\mcP/\mcS_{\mcP}$-multiplicity of $\mcM[\mcS_{\mcF} \trir \mcS_{\mcP}]$ is 0.
Suppose for a contradiction that it is at least 1.
Without loss of generality, we may assume that
$\mcP = \mcF[\mcS_{\mcF} \trir \mcS_{\mcP}]$ is a substructure of $\mcM[\mcS_{\mcF} \trir \mcS_{\mcP}]$,
so in particular, the common universe of $\mcF$  and $\mcP = \mcF[\mcS_{\mcF} \trir \mcS_{\mcP}]$ is a subset of 
the universe of $\mcM[\mcS_{\mcF} \trir \mcS_{\mcP}]$ and of $\mcM$.
For each relation symbol $R$, of arity $r$ say, we
consider the interpretation of $R$ in $\mcM \uhrc |\mcF|$.
If $\bar{a} \in |\mcS_{\mcF}|^r$, then
$$\bar{a} \in R^{\mcM \uhr F} \ \Longleftrightarrow \ \bar{a} \in R^{\mcS_{\mcF}}
\ \Longleftrightarrow \ \bar{a} \in R^{\mcF} \ \text{ (since $\mcS_{\mcF} \subset \mcF$)}.$$
If $\bar{a} \in |\mcF|^r - |\mcS_{\mcF}|^r$, then we use the definition of substitutions 
(Definition~\ref{definition of substitutions}) and get
\begin{align*}
\bar{a} \in R^{\mcM \uhr F} \ &\Longleftrightarrow \ \bar{a} \in R^{\mcM} \ 
\Longleftrightarrow \ \bar{a} \in R^{\mcM[\mcS_{\mcF} \trir \mcS_{\mcP}]} \\
&\Longleftrightarrow \ \bar{a} \in R^{\mcM[\mcS_{\mcF} \trir \mcS_{\mcP}] \uhr F} \ \Longleftrightarrow \
\bar{a} \in R^{\mcF[\mcS_{\mcF} \trir \mcS_{\mcP}]} \ \Longleftrightarrow \ \bar{a} \in R^{\mcF}.
\end{align*}
So whenever $\bar{a} \in |\mcF|^r$ we have $\bar{a} \in R^{\mcM}$ 
if and only if $\bar{a} \in R^{\mcF}$.
Since the argument holds for every relation symbol $R$ it follows that the forbidden structure
$\mcF$ is a substructure of $\mcM$, which contradicts that $\mcM$ is represented.
\hfill $\square$

\begin{defin}{\rm
Let the expression `$\mult(\mcA/\mcB; \mcM) \geq n$' mean 
`the $\mcA/\mcB$-multiplicity of $\mcM$ is at least $n$'.
}\end{defin}

\begin{lem}\label{different mutations}
Suppose that $\mcM, \mcN \in \mbK_n$ are different and that
$\mult(\mcP/\mcS_{\mcP}; \mcM) \geq 2$ and 
$\mult(\mcP/\mcS_{\mcP}; \mcN) \geq 2$.
Let $(\mcS_{\mcP}^{\mcM}, \mcS_{\mcF}^{\mcM})$ and
$(\mcS_{\mcP}^{\mcN}, \mcS_{\mcF}^{\mcN})$
be copies of the coexisting pair $(\mcS_{\mcP}, \mcS_{\mcF})$ such that
$\mcS_{\mcF}^{\mcM} \subseteq \mcM$ and $\mcS_{\mcF}^{\mcN} \subseteq \mcN$.
If $\mcM[\mcS_{\mcF}^{\mcM} \trir \mcS_{\mcP}^{\mcM}] =
\mcN[\mcS_{\mcF}^{\mcN} \trir \mcS_{\mcP}^{\mcN}]$
then $\mcS_{\mcF}^{\mcM}$ and $\mcS_{\mcF}^{\mcN}$ have the same universe
$U$ and $\mcM$ and $\mcN$ are different only on $U$ 
(that is, for every relation symbol $R$, if $\bar{a}$ belongs to exactly one of the relations
$R^{\mcM}$ and $R^{\mcN}$, then $\bar{a} \in U$.)
\end{lem}

\noindent
{\em Proof.}
Let $(\mcS_{\mcP}^{\mcM}, \mcS_{\mcF}^{\mcM})$ and
$(\mcS_{\mcP}^{\mcN}, \mcS_{\mcF}^{\mcN})$
be copies of the coexisting pair $(\mcS_{\mcP}, \mcS_{\mcF})$ such that
$\mcS_{\mcF}^{\mcM} \subseteq \mcM$ and $\mcS_{\mcF}^{\mcN} \subseteq \mcN$.
Then there are maps $\sigma_{\mcM} : |\mcS_{\mcF}^{\mcM}| \to |\mcS_{\mcF}|$ and 
$\sigma_{\mcN}: |\mcS_{\mcF}^{\mcN}| \to |\mcS_{\mcF}|$ such that:
\begin{itemize}
	\item[\textbullet] \ $\sigma_{\mcM}$ is an isomorphism from $\mcS_{\mcF}^{\mcM}$ to $\mcS_{\mcF}$ and 
	from  $\mcS_{\mcP}^{\mcM}$ to $\mcS_{\mcP}$, and 
	\item[\textbullet] \  $\sigma_{\mcN}$ is an isomorphism from $\mcS_{\mcF}^{\mcN}$ to $\mcS_{\mcF}$ and 
	from  $\mcS_{\mcP}^{\mcN}$ to $\mcS_{\mcP}$.
\end{itemize}
Let $\{a_1, \ldots, a_k\}$ be the universe of $\mcS_{\mcF}^{\mcM}$ (and of $\mcS_{\mcP}^{\mcM}$)
and let $\{b_1, \ldots, b_k\}$ be the universe of $\mcS_{\mcF}^{\mcN}$ (and of $\mcS_{\mcP}^{\mcN}$).

Suppose, for a contradiction, that 
\begin{itemize}
\item[(I)] $\mcM[\mcS_{\mcF}^{\mcM} \trir \mcS_{\mcP}^{\mcM}] = 
\mcH = \mcN[\mcS_{\mcF}^{\mcN} \trir \mcS_{\mcP}^{\mcN}]$
and that 
\item[(II)] $\{a_1, \ldots, a_k\} \neq \{b_1, \ldots, b_k\}$.
\end{itemize}
Then
\begin{equation*}
\mcM = \mcH[\mcS_{\mcP}^{\mcM} \trir \mcS_{\mcF}^{\mcM}] \ \text{ and } \
\mcN = \mcH[\mcS_{\mcP}^{\mcN} \trir \mcS_{\mcF}^{\mcN}]. \tag{1}
\end{equation*}
Recall the assumption that for every proper substructure $\mcU \subset \mcS_{\mcP}$,
$\mcS_{\mcP} \uhrc |\mcU| = \mcS_{\mcF} \uhrc |\mcU|$.
Since $(\mcS_{\mcP}^{\mcM}, \mcS_{\mcF}^{\mcM})$ and
$(\mcS_{\mcP}^{\mcN}, \mcS_{\mcF}^{\mcN})$
are copies of $(\mcS_{\mcP}, \mcS_{\mcF})$,
it follows that $\mcS_{\mcP}^{\mcM}$ and $\mcS_{\mcF}^{\mcM}$ agree on
all proper subsets of their common universe; and the same with $\mcM$ replaced by $\mcN$.
From (1) it follows that
\begin{equation*}
\text{if $U \subseteq \{1, \ldots, m_n\}$ and $|U| < k$, then 
$\mcM \uhrc U = \mcH \uhrc U = \mcN \uhrc U$.} \tag{2}
\end{equation*}
Since $\mcH \uhrc \{b_1, \ldots, b_k \} = \mcS_{\mcP}^{\mcN}$
and $\mcM$ is obtained from $\mcH$ by the substitution $\mcM = \mcH[\mcS_{\mcP}^{\mcM} \trir \mcS_{\mcF}^{\mcM}]$,
which {\em only} affects the interpretations of relation symbols on
$\{a_1, \ldots, a_k\}$, 
assumption (II) together with (2) implies that
$$\mcM \uhrc \{b_1, \ldots, b_k \} = \mcS_{\mcP}^{\mcN}.$$
Since the $\mcP/\mcS_{\mcP}$-multiplicity of $\mcM$ is at least 2,
there are $\mcP_i \subseteq \mcM$ and isomorphisms $\sigma_i : \mcP_i \to \mcP$
such that
$\mcS_{\mcP}^{\mcN} \subset \mcP_i$,
$\sigma_i \uhrc |\mcS_{\mcP}^{\mcN}| = \sigma_{\mcN}$, for $i = 1,2$,
and $|\mcP_1| \cap |\mcP_2| = \big\{b_1, \ldots, b_k \big\} 
= |\mcS_{\mcP}^{\mcN}|$.
By assumption (I), $\mcH$ is obtained from $\mcM$ by the substitution 
$\mcH = \mcM[\mcS_{\mcF}^{\mcM} \trir \mcS_{\mcP}^{\mcM}]$ which only affects
the interpretations of relation symbols on $\{a_1, \ldots, a_k\}$.
This together with (II), (2) and the choice of $\mcP_1$ and $\mcP_2$ so that
$|\mcP_1| \cap |\mcP_2| = \big\{b_1, \ldots, b_k \big\}$
implies that for $i = 1$ or $i = 2$, $\mcH \uhrc |\mcP_i| = \mcP_i$.
Choose $i$ so that 
\begin{equation*}
\mcH \uhrc |\mcP_i| = \mcP_i. \tag{3}
\end{equation*}
Since $\mcS_{\mcP}^{\mcN} \subset \mcP_i$ and $\sigma_i : \mcP_i \to \mcP$
is an isomorphism such that $\sigma_i \uhrc |\mcS_{\mcP}^{\mcN}| = \sigma_{\mcN}$, the substitution 
$[\mcS_{\mcP}^{\mcN} \trir \mcS_{\mcF}^{\mcN}]$
changes $\mcP_i$ to a structure which is isomorphic with $\mcF$,
that is, $\mcP_i[\mcS_{\mcP}^{\mcN} \trir \mcS_{\mcF}^{\mcN}] \cong \mcF$,
via the isomorphism $\sigma_i$.
By applying (1) and (3) we get
$$\mcN \uhrc |\mcP_i| = (\mcH[\mcS_{\mcP}^{\mcN} \trir \mcS_{\mcF}^{\mcN}]) \uhrc |\mcP_i| \cong \mcF.$$
Hence the substructure of $\mcN$ with universe $|\mcP_i|$ is isomorphic to the forbidden structure $\mcF$.
Therefore $\mcN$ is not represented, which contradicts that $\mcN \in \mbK_n$.

So if (I) holds then (II) is false and hence all the structures $\mcS_{\mcF}^{\mcM}$, 
$\mcS_{\mcP}^{\mcM}$, $\mcS_{\mcF}^{\mcN}$ and $\mcS_{\mcP}^{\mcN}$ have the same
universe, say $U$. Consequently, from the assumption (I), if $R$ is a relation symbol of arity $r$, say,
and $\bar{a} \in \{1, \ldots, m_n\}^r$ belongs to exactly one of $R^{\mcM}$ and $R^{\mcN}$,
then $\bar{a} \in U$.
\hfill $\square$

\begin{defin}{\rm
(i) For every $L$-structure $\mcM$, let 
$\mathbf{\Sigma}(\mcM; \mcS_{\mcF} \trir \mcS_{\mcP})$ denote
the set of all structures of the form 
$\mcM[\mcS^{\mcM}_{\mcF} \trir \mcS^{\mcM}_{\mcP}]$
where $(\mcS^{\mcM}_{\mcP}, \mcS^{\mcM}_{\mcF})$ is a copy of the coexisting pair
$(\mcS_{\mcP}, \mcS_{\mcF})$ and $\mcS^{\mcM}_{\mcF} \subseteq \mcM$.
(If $\mcM$ contains no copy of $\mcS_{\mcF}$ then $\mathbf{\Sigma}(\mcM; \mcS_{\mcF} \trir \mcS_{\mcP}) = \es$)\\
(ii) For every $n$, let $\mathbf{\Omega}_n$ denote the set of all $\mcM \in \mbK_n$ such
that $\mult(\mcP/\mcS_{\mcP}; \mcM) \geq 2$.\\
(iii) Recall that $\alpha$ denotes the number of different permitted $L$-structures
with universe $\{1, \ldots, k\}$.
}\end{defin}

\begin{lem}\label{almost disjoint mutation sets}
If $\mcM_1, \ldots, \mcM_{\alpha + 1} \in \mathbf{\Omega}_n$ and $\mcM_i \neq \mcM_j$ whenever $i \neq j$,
then
$$\bigcap_{1 \leq i \leq \alpha + 1}\mathbf{\Sigma}(\mcM_i; \mcS_{\mcF} \trir \mcS_{\mcP}) = \es.$$
In other words, for every structure $\mcN$, it can belong to 
$\mathbf{\Sigma}(\mcM; \mcS_{\mcF} \trir \mcS_{\mcP})$
for at most $\alpha$ distinct $\mcM \in \mathbf{\Omega}_n$.
\end{lem}

\noindent
{\em Proof.}
Suppose for a contradiction that $\mcM_1, \ldots, \mcM_{\alpha + 1} \in \Omega_n$ are distinct and
that $\mcN \in \mathbf{\Sigma}(\mcM_i; \mcS_{\mcF} \trir \mcS_{\mcP})$ for every $i \in \{1, \ldots, \alpha + 1\}$.
Then there are copies $(\mcS_{\mcP}^{\mcM_i}, \mcS_{\mcF}^{\mcM_i})$
of $(\mcS_{\mcP}, \mcS_{\mcF})$ such that
$\mcS_{\mcF}^{\mcM_i} \subseteq \mcM_i$ and 
$\mcN = \mcM_i[\mcS_{\mcF}^{\mcM_i} \trir \mcS_{\mcP}^{\mcM_i}]$ for every $i \in \{1, \ldots, \alpha + 1\}$.
By Lemma~\ref{different mutations}, 
all $\mcS_{\mcF}^{\mcM_i}$, $i \in \{1, \ldots, \alpha + 1\}$, have the same universe,
which we denote by $U$, and for every pair $i, j \in \{1, \ldots, \alpha + 1\}$ of distinct numbers,
$\mcM_i$ and $\mcM_j$ are different only on $U$. The assumption that 
$\mcM_i \neq \mcM_j$ if $i \neq j$ now implies that for all distinct $i, j \in \{1, \ldots, \alpha + 1\}$,
$\mcM_i \uhrc U \neq \mcM_j \uhrc U$. 
Since $|U| = k$, this contradicts the choice of $\alpha$, being the number of all different 
permitted $L$-structures with universe $\{1, \ldots, k\}$. 
\hfill $\square$
\\

\noindent
Now we have the tools for proving part (i) of 
Theorem~\ref{upper limit on the probability of an extension axiom},
and then the other parts of the theorem.
Let $\mathbf{\Omega}^*_n$ be the set of all $\mcM \in \mbK_n$ such that
\begin{itemize}
	\item[(a)] $\mcM$ contains a copy of $\mcS_{\mcF}$, and
	\item[(b)] the $\mcP/\mcS_{\mcP}$-multiplicity of $\mcM$ is at least 2.
\end{itemize}
By (b) and the definition of $\mathbf{\Omega}_n$, $\mathbf{\Omega}^*_n \subseteq \mathbf{\Omega}_n$.
Since every $\mcM \in \mathbf{\Omega}^*_n$ contains a copy of $\mcS_{\mcF}$,
it follows that for every $\mcM \in \mathbf{\Omega}^*_n$,
$\mathbf{\Sigma}(\mcM; \mcS_{\mcF} \trir \mcS_{\mcP}) \neq \es$.
Since the substitution $[\mcS_{\mcF} \trir \mcS_{\mcP}]$ is admitted,
$\mathbf{\Sigma}(\mcM; \mcS_{\mcF}\trir\mcS_{\mcP}) \subseteq \mbK_n$
for every $\mcM \in \mbK_n$.
By Lemma~\ref{reducing the multiplicity to 0},
for every $\mcM \in \mathbf{\Omega}^*_n$, 
$\mathbf{\Sigma}(\mcM; \mcS_{\mcF}\trir\mcS_{\mcP}) \subseteq \mbK_n - \mathbf{\Omega}^*_n$.
Lemma~\ref{almost disjoint mutation sets} now implies that 
$$\big|\mbK_n - \mathbf{\Omega}^*_n\big|  \geq
\Big| \bigcup_{\mcM \in \mathbf{\Omega}^*_n} \mathbf{\Sigma}(\mcM; \mcS_{\mcF} \trir \mcS_{\mcP}) \Big| 
\geq \frac{\big| \mathbf{\Omega}^*_n \big|}{\alpha}$$
$$\text{and hence } \quad \alpha\big|\mbK_n - \mathbf{\Omega}^*_n\big| \geq \big|\mathbf{\Omega}^*_n\big|.$$
From this we get
$$\frac{\big|\mbK_n - \mathbf{\Omega}^*_n\big|}{\big|\mbK_n\big|} =
\frac{\big|\mbK_n - \mathbf{\Omega}^*_n\big|}{\big|\mathbf{\Omega}^*_n\big| + 
\big|\mbK_n - \mathbf{\Omega}^*_n\big|} \geq
\frac{\big|\mbK_n - \mathbf{\Omega}^*_n\big|}{
\alpha\big|\mbK_n - \mathbf{\Omega}^*_n\big| + \big|\mbK_n - \mathbf{\Omega}^*_n\big|} =
\frac{1}{\alpha + 1}.$$
Thus, the proportion of $\mcM \in \mbK_n$ {\em not} satisfying both (a) and (b)
is at least $1/(1 + \alpha)$. This concludes the proof of part (i) of 
Theorem~\ref{upper limit on the probability of an extension axiom}.

Part (ii) of Theorem~\ref{upper limit on the probability of an extension axiom}
is a straightforward consequence of part (i).
For if there exist a permitted structure $\mcC$ and embeddings
$\sigma_1 : \mcP \to \mcC$ and $\sigma_2 : \mcP \to \mcC$ such that
$\sigma_1(|\mcP|) \cap \sigma_2(|\mcP|) = \sigma_1(|\mcS_{\mcP}|)$,
$\sigma_1 \uhrc |\mcS_{\mcP}| = \sigma_2 \uhrc |\mcS_{\mcP}|$
and $\mcM \in \mbK_n$ satisfies  all $(2\left\| \mcP \right\| - k - 1)$-extension axioms,
then conditions (a) and (b) in part (i) of  
Theorem~\ref{upper limit on the probability of an extension axiom}
are satisfied.

Now we prove part (iii) of Theorem~\ref{upper limit on the probability of an extension axiom}.
Here we have added the assumption that $L$ has no unary relation symbols, so there is a
unique (up to isomorphism) permitted structure with a singleton universe.
(In fact this is sufficient for what we want to prove.)
Let $\mcU \subset \mcS_{\mcF}$ be such that $\left\|\mcU\right\| = 1$.
Note that since $\mcS_{\mcF} \neq \mcS_{\mcP}$ (and $|\mcS_{\mcF}| = |\mcS_{\mcP}|$)
we have $\left\|\mcS_{\mcF}\right\| > 1$.
Suppose that $\mcM \in \mbK_n$ is such that
\begin{itemize}
	\item[(c)] $\mcM$ satisfies the $\mcS_{\mcF}/\mcU$-extension axiom, and
	\item[(d)] the $\mcP/\mcS_{\mcP}$-multiplicity of $\mcM$ is at least 2.
\end{itemize}
Since $\left\|\mcM\right\| = m_n$, there are $m_n$ distinct copies of $\mcU$ in $\mcM$.
Each one of these copies of $\mcU$ is, by (c), included in a copy of $\mcS_{\mcF}$,
so we get at least $m_n/k$ distinct copies of $\mcS_{\mcF}$ in $\mcM$.
Recall that our assumptions imply that $\mcS_{\mcP} \neq \mcS_{\mcF}$ and
$\mcS_{\mcP} \uhrc |\mcV| = \mcS_{\mcF} \uhrc |\mcV|$ for every proper substructure $\mcV \subset \mcS_{\mcP}$.
Since $\mathbf{\Sigma}(\mcM; \mcS_{\mcF} \trir \mcS_{\mcP})$ contains all
$\mcN$ which can be obtained from $\mcM$ by replacing one copy of $\mcS_{\mcF}$
by a copy of $\mcS_{\mcP}$, we have
$\big|\mathbf{\Sigma}(\mcM; \mcS_{\mcF} \trir \mcS_{\mcP})\big| \geq m_n/k$.
By Lemma~\ref{reducing the multiplicity to 0}, {\em no}
$\mcN \in \mathbf{\Sigma}(\mcM; \mcS_{\mcF} \trir \mcS_{\mcP})$
satisfies (d).
Hence, if $\mbE_n$ is the set of all $\mcM \in \mbK_n$ which satisfy
both (c) and (d), then, by Lemma~\ref{almost disjoint mutation sets},
$$\big|\mbK_n - \mbE_n\big| \geq \frac{m_n |\mbE_n\big|}{k \alpha} \quad
\text{ and hence } \quad \frac{\big|\mbE_n\big|}{\big|\mbK_n\big|} 
\leq \frac{|\mbE_n|}{|\mbK_n - \mbE_n|} 
\leq \frac{k\alpha}{m_n}.$$
As $\lim_{n \to \infty}m_n = \infty$, the proportion of $\mcM \in \mbK_n$ which satisfy both
(c) and (d) approaches 0 as $n$ approaches $\infty$.
This concludes the proof of Theorem~\ref{upper limit on the probability of an extension axiom}.

\section{Conditional probability measures}
\label{conditional measures}

\noindent
In Sections~\ref{forbidden structures}~--~\ref{proving upper limit on the probability of an extension axiom}
we saw that a condition that ensures that every extension axiom is true in almost
all sufficiently large structures is that every substitution involving (only) permitted structures
is admitted. And if this condition does not hold it may happen that some extension axiom
is false in almost all sufficiently large structures.
In this section we start to develop a theory of conditional probability measures on
finite sets of structures.
When using this measure we can include more examples of sets of finite structures for
which any extension axiom is almost surely true in all sufficiently large structures under consideration.
Such examples include Examples~\ref{graphs with a unary predicate},
\ref{partially coloured binary relations} and~\ref{coloured binary relations}, 
and  more generally, coloured structures and partially coloured structures
(as in examples~\ref{example of coloured structures of first kind} --~\ref{partial colourings}).
But there are other examples, such as $\mcK_l$-free graphs ($l \geq 3$) which are not included;
that is, also with the conditional measures considered here there is an extension axiom which almost
surely fails for sufficiently large $\mcK_l$-free graphs.

Although the uniform probability measure
is conceptually simple, it does not necessarily correspond to the
probability measure associated with a method for
randomly generating a structure of some specified kind.  
The conditional measures to be considered are more closely related 
to probability measures associated with random generation of
structures of a given kind.
This is the first point that will be stressed below, after
the next two definitions.

\begin{defin}\label{definition of reduct conditional measure}{\rm
Let $\mbC_0$ and $\mbC_1$ be finite sets of $L$-structures and let $\mbbP_0$ be a probability measure
on $\mbC_0$. Suppose that 
\begin{itemize}
\item[(1)] for every $\mcA \in \mbC_0$ there is at least one $\mcB \in \mbC_1$ such that $\mcA \subseteq_w \mcB$, and 
\item[(2)] for every $\mcB \in \mbC_1$ there is a {\em unique}
$\mcA \in \mbC_0$ such that $\mcA \subseteq_w \mcB$
and whenever $\mcA' \in \mbC_0$ and $\mcA' \subseteq_w \mcB$, then $\mcA' \subseteq_w \mcA$.
We denote such $\mcA$ by $\mcB \uhrc 0$.
\end{itemize}
Then we define the {\bf \em uniformly $\mbbP_0$-conditional probability measure}
$\mbbP_1$ on $\mbC_1$ as follows:
\begin{align*}
&\text{For every $\mcB \in \mbC_1$, the probability of $\mcB$ in $\mbC_1$ is} \\
&\mbbP_1(\mcB) = \frac{1}{\big| \big\{ \mcB' \in \mbC_1 : \mcB' \uhrc 0 = \mcB \uhrc 0 \big\} \big|} 
\cdot \mbbP_0(\mcB \uhrc 0),\\
&\text{and for $\mbX = \{\mcB_1, \ldots, \mcB_n\} \subseteq \mbC_1$ 
(where $\mbX$ is enumerated without repetition)}\\
&\mbbP_1(\mbX) = \sum_{i = 1}^n \mbbP_1(\mcB_i).
\end{align*}
}\end{defin}

\begin{defin}\label{definition of transitive closure of reduct conditional measure}{\rm
More generally, assume that $\mbC_0, \ldots, \mbC_r$ are finite sets of $L$-structures
such that, for every $i = 0, \ldots, r-1$, (1) and (2) in Definition~\ref{definition of reduct conditional measure}
hold if $\mbC_0$ and $\mbC_1$ are replaced by $\mbC_i$ and $\mbC_{i+1}$, respectively.
Let $\mbbP_0$ denote the uniform probability measure on $\mbC_0$
(i.e. all elements of $\mbC_0$ have the same probability $1 / |\mbC_0|$).
By induction, define $\mbbP_{i+1}$ to be the uniformly $\mbbP_i$-conditional probability measure,
for $i = 0, \ldots, r-1$.
We call the probability measure $\mbbP_r$ on $\mbC_r$, thus obtained, the 
{\bf \em uniformly $(\mbC_0, \ldots, \mbC_{r-1})$-conditional probability measure}.
}\end{defin}

\begin{exam}\label{example of conditional measure 1}{\rm
Let us first illustrate the definitions by considering Example~\ref{graphs with a unary predicate},
where $\mbK_n$ is the set of undirected graphs with vertices $1, \ldots, n$
(with edge relation represented by $R$)
and a unary relation symbol $P$ subject to the condition:
$R(a,b)$ $\Longrightarrow$ $\neg P(a)$ and $\neg P(b)$.
We have proved (see Example~\ref{graphs with a unary predicate})
that with the uniform probability measure, the probability
of $\exists x P(x)$ holding in $\mcM \in \mbK_n$ approaches 0 as $n \to \infty$.
Next we show that with a naturally chosen conditional measure,
the probability that $\exists x P(x)$ holds in $\mcM \in \mbK_n$
approaches 1 as $n \to \infty$.

Let $L$ denote the language considered in Example~\ref{graphs with a unary predicate},
with one binary relation symbol $R$ and one unary relation symbol $P$, and let 
$L_0$ be the sublanguage of $L$ whose vocabulary contains only $P$.
For every $n$, let $\mbK_n \uhrc L_0 = \{\mcM \uhrc L_0 : \mcM \in \mbK_n\}$.
Recall from the definition of weak substructure
(Section~\ref{languages, structures and embeddings}), and the discussion after it,
that, for every $\mcM \in \mbK_n$, $\mcM \uhrc L_0$ may also be viewed as 
an $L$-structure (in which the interpretation of $R$ is empty) and it follows that
$\mcM \uhrc L_0 \subseteq_w \mcM$.
It is easy to verify that, for every $n$, if $\mbC_0 = \mbK_n \uhrc L_0$ and $\mbC_1 = \mbK_n$,
then conditions (1) and (2) in Definition~\ref{definition of reduct conditional measure} hold.
Hence, for every $n$, the uniformly $(\mbK_n \uhrc L_0)$-conditional
probability measure on $\mbK_n$ is well-defined.
Now, the claim that the probability, with this measure, that 
(the extension axiom) $\exists x P(x)$ holds in  $\mcM \in \mbK_n$ approaches 1 as $n \to \infty$,
is a consequence of Theorem~\ref{0-1 law for pregeometries}.
But for this simple example 
it suffices to observe that the probability of $\mcM \in \mbK_n$, with
the uniformly $(\mbK_n \uhrc L_0)$-conditional measure,
is the probability of obtaining $\mcM$ by the following generating procedure:
First go through every $i \in \{1, \ldots, n\}$ and with probability $1/2$ 
let it satisfy $P(x)$; then take the set $\{i_1, \ldots, i_m\}$ of all vertices
which do not satisfy $P(x)$,
and for each unordered pair $\{i, j\}$ of elements from 
$\{i_1, \ldots, i_m\}$ assign an edge to it with
probability $1/2$. So the probability that no $i \in \{1, \ldots, n\}$
satisfies $P(x)$ is $1/2^n$, which approaches 0 as $n \to \infty$.
}\end{exam}

\begin{exam}\label{example of conditional measure 2}{\rm
Let us now consider Example~\ref{partially coloured binary relations}
(partially coloured binary relation),
where the vocabulary of $L$ is $\{R, P_1, P_2\}$,  
$R$ is binary and $P_i$, $i = 1, 2$, are unary, and thought of as
``colours''.
$\mbK_n$ consists of all structures with universe $\{1, \ldots, n\}$
such that the universe is partially coloured with respect to the relation $R$,
that is, every element has at most one colour (1 or 2), and it may be uncoloured
(or ``blank''), and whenever $R(a,b)$ holds, then $a$ and $b$ cannot be coloured
with the same colour.

How can we, for any given $k$, design a procedure that generates -- by possibly making some
random assignments on the way -- $\mcM \in \mbK_n$ in such a way that the
probability of ending up with an $\mcM \in \mbK_n$ with exactly $k$ elements
with colour 1 is the same as the proportion of $\mcM \in \mbK_n$ which
have exactly $k$ elements with colour 1?
The author does not know, and the point is that, in general, it may not be
easy to conceive of a generating procedure, of structures from a given set, such that
the probability measure associated with the generating procedure is identical to
the uniform probability measure on the given set of structures. 

Recall, from Example~\ref{partially coloured binary relations}, that there is an 
extension axiom $\varphi$ such that the probability, with the uniform measure,
that $\varphi$ holds in $\mcM \in \mbK_n$ approaches 0 as $n \to \infty$.
But if we apply the following generating procedure of $\mcM \in \mbK_n$,
then, for {\em every} extension axiom $\varphi$, the probability of
ending up with an $\mcM \in \mbK_n$ which satisfies $\varphi$ approaches 1 as $n \to \infty$.
For every $i \in \{1, \ldots, n\}$, with probability $1/3$ let it have
colour 1, colour 2 or be blank; then go through all pairs $(i,j)$ such 
that $i$ and $j$ are not coloured with the same colour and let 
$(i,j) \in R^{\mcM}$  with probability $1/2$. 
The probability of obtaining, in this way, a structure $\mcM \in \mbK_n$ is the same as the probability
of $\mcM$ with the uniformly $(\mbK_n \uhrc L_0)$-conditional
measure on $\mbK_n$, where $L_0$ is the sublanguage of $L$ whose vocabulary is $\{P_1, P_2\}$
and $\mbK_n \uhrc L_0 = \{\mcM \uhrc L_0 : \mcM \in \mbK_n\}$.
By letting the underlying geometry of every structure in $\mbK = \bigcup_{n \in \mbbN} \mbK_n$
be trivial (see Remark~\ref{remarks on pregeometries})
and applying Theorem~\ref{0-1 law for pregeometries} it follows that, 
for every extension axiom $\varphi$ of  $\mbK$, the probability, with the 
uniformly $(\mbK_n \uhrc L_0)$-conditional measure, that $\varphi$ holds in $\mbK_n$
approaches 1 as $n \to \infty$; and by Theorem~\ref{third part of 0-1 law for pregeometries},
$\mbK$ has a zero-one law.
We have in particular shown that the asymptotic probability,
with the uniform probability measure, of a first order definable
property in $\mbK$ may be different from the asymptotic probability
of the same property when the $(\mbK_n \uhrc L_0)$-conditional measure is used.
}\end{exam}

\noindent
Before taking underlying pregeometries into account, we collect a technical lemma
which will be used later.

\begin{lem}\label{probability in C_r can be computed as probabilty in C_r+1}
Suppose that $\mbC_0, \ldots, \mbC_k$ are finite sets of structures
such that, for every $i = 0, \ldots, k-1$, (1) and (2) in Definition~\ref{definition of reduct conditional measure}
hold if $\mbC_0$ and $\mbC_1$ are replaced by $\mbC_i$ and $\mbC_{i+1}$, respectively.
For $r = 1, \ldots, k$, let $\mbbP_r$ denote the uniformly $(\mbC_0, \ldots, \mbC_{r-1})$-conditional
probability measure on $\mbC_r$. 
If $1 \leq r \leq s \leq k$ and $\mcA \subseteq \mbC_r$, then
\begin{align*}
\mbbP_r(\mcA) &= \mbbP_{r+1}\big( \{\mcB \in \mbC_{r+1} : \mcA \subseteq_w \mcB \} \big)\\
\text{and } \mbbP_r(\mcA) &= P_s\big( \{\mcB \in \mbC_s : \mcA \subseteq_w \mcB \} \big).
\end{align*}
\end{lem}

\noindent
{\em Proof.}
The second identity follows from the first by induction, and
the first identity is a straightforward consequence of 
Definitions~\ref{definition of transitive closure of reduct conditional measure}
and~\ref{definition of reduct conditional measure}.
\hfill $\square$

\section{Underlying pregeometries}\label{pregeometries}

\begin{defin}\label{definition of a structure being a pregeometry}{\rm 
(i) We call an $L$-structure $\mcA$ a {\bf \em pregeometry} if 
\begin{itemize}
	\item[(1)] there is a closure operation $\cl_{\mcA}$ on $A$ such that $(A, \cl_{\mcA})$ is a pregeometry,
  \item[(2)] for all $n \in \mbbN$ there is a formula $\theta_n(x_1, \ldots, x_{n+1}) \in L$ such that
	for all $a_1, \ldots, a_{n+1} \in A$, $a_{n+1} \in \cl_{\mcA}(a_1, \ldots, a_n)$ if and only
	if $\mcA \models \theta_n(a_1, \ldots, a_{n+1})$, and
	\item[(3)] if $X \subseteq A$ is closed with respect to $\cl_{\mcA}$ (i.e. $\cl_{\mcA}(X) = X$),
	then $X$
	is closed under interpretations of (eventual) function symbols and constant symbols;
	so $X$ is the universe of a substructure of $\mcA$.
\end{itemize}
(ii) Let $\mbK$ be a class of $L$-structures. 
We call $\mbK$ a {\bf \em pregeometry} if every $\mcA \in \mbK$ is a 
pregeometry and for every $n \in \mbbN$ there is a
formula $\theta_n(x_1, \ldots, x_{n+1}) \in L$ such that
for every $\mcA \in \mbK$ and all $a_1, \ldots, a_{n+1} \in A$, 
$a_{n+1} \in \cl_{\mcA}(a_1, \ldots, a_n)$ if and only if $\mcA \models \theta_n(a_1, \ldots, a_{n+1})$.
}\end{defin}

\begin{rem}\label{remarks on pregeometries}{\rm
For every structure $\mcA$, if $\cl_{\mcA}(X) = X$ for every $X \subseteq A$,
then $\mcA$ is a pregeometry in the sense of 
Definition~\ref{definition of a structure being a pregeometry} (i).
This pregeometry is often called {\bf \em trivial} or {\bf \em degenerate}.
It may happen that for a structure $\mcA$ there is more than one way to define
a pregeometry on $A$. As noted, we always have a trivial pregeometry on $A$.
But if, for example, $\mcA$ is a vector space over some finite field 
(formalized as a first-order structure in a suitable way),
then we can also let $\cl_{\mcA}(X)$ be the linear span of $X$,
and then $\cl_{\mcA}$ becomes a pregeometry on $A$.
When saying that a structure $\mcA$ is a pregeometry we assume that some
particular pregeometry on $A$ (in the sense of 
Definition~\ref{definition of a structure being a pregeometry} (i))
is fixed, and if we say that a class of $L$-structures $\mbK$ is a pregeometry we 
assume that, for every $\mcA \in \mbK$, some pregeometry $\cl_{\mcA}$ is fixed on $A$ and that
the condition in Definition~\ref{definition of a structure being a pregeometry} (ii) holds.
}\end{rem}

\begin{assump}\label{first assumptions on pregeometries}{\rm
For the rest of this section we assume that $\mbK$ is a class
of $L$-structures which is a pregeometry, and that the formulas
$\theta_n(x_1, \ldots, x_{n+1})$ define the pregeometry in the sense
of Definition~\ref{definition of a structure being a pregeometry} (ii).
(Later, in Assumption~\ref{assumptions about pregeometries},
we will add some more assumptions.)
}\end{assump}

\begin{defin}\label{definition concerning permitted structure in pregeometries}{\rm
(i) As in Sections~\ref{forbidden structures} --~\ref{conditional measures}, 
we say that structure 
$\mcA$ is {\bf \em represented (with respect to $\mbK$)} if $\mcA$ is
isomorphic to some structure in $\mbK$. 
We say that $\mcA$ is {\bf \em permitted (with respect to $\mbK$)} 
if it can be embedded into some structure in $\mbK$;
or equivalently, if it is a substructure of some represented structure.
And a structure which is not permitted (with respect to $\mbK$) is
{\bf \em forbidden (with respect to $\mbK$)}.
Note that every represented structure is a pregeometry on which the closure operator is
defined by $\theta_n(x_1, \ldots, x_{n+1})$, $n \in \mbbN$.
This is what we mean when speaking about a pregeometry and closure on a
represented structure.\\
(ii) If $\mcM$ is a pregeometry, then the notation $\mcA \subseteq_{cl} \mcM$
means that $\mcA$ is a substructure of $\mcM$ and $\cl_{\mcM}(A) = A$. 
In words, we express `$\mcA \subseteq_{cl} \mcM$' by saying that
{\bf \em $\mcA$ is a closed substructure of $\mcM$}.
}\end{defin}

\begin{defin}\label{definition of multiplicity for pregeometries}{\rm
The notion of $\mcB / \mcA$-multiplicity is defined as before 
(Definition~\ref{definition of multiplicity}),
{\em except} that we require that $\mcA$ and $\mcB$ are closed in some superstructure.
More precisely:
Suppose that there is 
a represented $\mcN$ such that $\mcA \subset \mcB \subseteq \mcN$ and
both $A$ and $B$ are closed in $\mcN$.
We say that the {\bf \em $\mcB / \mcA$-multiplicity of a (represented) structure $\mcM$ is at least $m$}
if  the following holds:
\begin{itemize}
\item[] whenever $\mcA' \subseteq_{cl} \mcM$ 
and $\sigma : \mcA' \to \mcA$ is an isomorphism, 
then there are
$\mcB'_i \subseteq_{cl} \mcM$ and isomorphisms $\sigma_i : \mcB'_i \to \mcB$, for $i = 1, \ldots, m$,
such that
$\mcA' \subseteq \mcB'_i$, $\sigma_i \uhrc A' = \sigma$
and $B'_i \cap B'_j = A'$ whenever $i \neq j$. 
\end{itemize}
The {\bf \em $\mcB/\mcA$-multiplicity is $m$} if it is at least $m$ but not at least $m+1$.
}\end{defin}

\begin{rem}\label{remark about uniform definability of closure}{\rm
Observe that we can express, in first-order logic, that sets are closed (or not)
in a uniform way. For if $\gamma_n(x_1, \ldots, x_n)$ denotes the formula
$$\neg \exists x_{n+1} \Big( \theta_n(x_1, \ldots, x_n, x_{n+1}) \ \wedge \
\bigwedge_{i=1}^n x_i \neq x_{n+1} \Big),$$
then for every $\mcM \in \mbK$ and all $a_1, \ldots, a_n \in M$,
$\mcM \models \gamma_n(a_1, \ldots, a_n)$ if and only if $\{a_1, \ldots, a_n\}$ is closed
in $\mcM$.
It follows that whenever $\mcM$ is represented and $\mcA \subseteq \mcB \subseteq \mcM$ are
closed substructures of $\mcM$,
then, for every $m \in \mbbN$, there is a sentence $\varphi_m$ such that for every represented $\mcN$,
$\mcN \models \varphi_m$ if and only if the $\mcB / \mcA$-multiplicity of $\mcN$
is at least $m$.
}\end{rem}

\begin{defin}\label{definition of extension axiom for pregeometries}{\rm
For represented $\mcM$ and closed substructures $\mcA \subset \mcB \subseteq \mcM$,
the {\bf \em $\mcB / \mcA$-extension axiom} is the statement expressing that
the $\mcB / \mcA$-multiplicity is at least 1.
As noted in Remark~\ref{remark about uniform definability of closure}, this
statement is expressible with a first-order sentence.
}\end{defin}

\noindent
Note that if the closure operator of (structures in) $\mbK$
is trivial, then the definitions of extension axioms and multiplicity 
coincide with those given earlier; so the earlier setting is a
special case of the current setting.

\begin{defin}\label{definition of polynomial k-saturation for classes}{\rm
Let $\mbK$ be a class of $L$-structures and let $(\mcM_n : n \in \mbbN)$ be 
a sequence of structures from $\mbK$.\\
(i) We say that the sequence $(\mcM_n : n \in \mbbN)$ is {\bf \em polynomially $k$-saturated}
	if there are a sequence of numbers $(\lambda_n : n \in \mbbN)$ with
	$\lim_{n \to \infty}\lambda_n = \infty$ and a polynomial
	$P(x)$ such that for every $n \in \mbbN$:
	\begin{itemize}
		\item[(1)] $\lambda_n \leq |M_n| \leq P(\lambda_n)$, and
		\item[(2)] whenever $\mcN$ is represented and $\mcA \subset \mcB \subseteq \mcN$
		are closed (in $\mcN$) and $\dim_{\mcN}(A) + 1 = \dim_{\mcN}(B) \leq k$, 
		then the $\mcB / \mcA$-multiplicity of $\mcM_n$ is at least $\lambda_n$.
	\end{itemize}
(ii) We say that $\mbK$ is {\bf \em polynomially $k$-saturated} if there are
$\mcM_n \in  \mbK$, for $n \in \mbbN$, such that the
sequence $(\mcM_n : n \in \mbbN)$ is {\bf \em polynomially $k$-saturated}.
}\end{defin}

\begin{exam}\label{examples of polynomially k-saturated pregeometries}{\rm
While it is possible to construct many different $\mbK$ which are polynomially $k$-saturated 
(by application of Theorem~\ref{0-1 law for pregeometries}) the 
kind of pregeometries that are present in examples that the author can construct
are rather limited.
So let us look at examples of $\mbK$ which are  polynomially $k$-saturated
for every $k \in \mbbN$ and which do not have any more structure than what
is necessary for defining the pregeometry.
The cases known are on the one hand the
trivial pregeometry and on the other hand
(possibly projective or affine variants of) linear spaces over a fixed,
but arbitrary, finite field.

If $L$ has empty vocabulary and $\mcE_n$ is the unique $L$-structure 
with universe $\{1, \ldots, n\}$ (with trivial closure operator), 
then it is straightforward to check that $(\mcE_n : n \in \mbbN)$ is
polynomially $k$-saturated for every $k \in \mbbN$.

Now suppose that $\mcG_n$ is a vector space with dimension $n$ with
universe $\{1, \ldots, p^n\}$ over a finite field $F$ of order $p$.
Let $\cl_{\mcG_n}$ be linear span.
To view $\mcG_n$ as a first order structure we can let scalar multiplication
be represented by unary function symbols (one for every element in $F$),
vector addition by a binary function symbol, and let there be a constant symbol for the zero vector.
Then $\{\mcG_n : n \in \mbbN\}$ is a pregeometry in the sense of 
Definition~\ref{definition of a structure being a pregeometry}.
The proof of Lemma~3.5 in \cite{Dj06a} shows that $(\mcG_n : n \in \mbbN)$
is polynomially $k$-saturated, for every $k \in \mbbN$. 
In \cite{Dj06a} it is explained how one can ``transform'' $\mcG_n$ into a
first-order structure, which represents a projective space $\mcP_n$ or affine
space $\mcA_n$ over $F$ of dimension $n$.
By the argument leading to Proposition~3.4 in \cite{Dj06a} it follows that
$(\mcP_n : n \in \mbbN)$ and $(\mcA_n : n \in \mbbN)$ are polynomially $k$-saturated,
for every $k \in \mbbN$. 

There are other ``linear geometries'' (see \cite{CH}) which involve
quadratic forms. These may be candidates for other polynomially $k$-saturated
sequences of pregeometries; but for reasons explained in Problem~3.8 in \cite{Dj06a}, the author
has not been able to prove or disprove it.
}\end{exam}

\noindent
From now on we work within the following context, in addition to the
assumptions already made (see Assumption~\ref{first assumptions on pregeometries}).

\begin{assump}\label{assumptions about pregeometries}{\rm From now on we assume the following:
\begin{itemize}
\item[(1)]$L_0 \subseteq L$ are first-order languages with vocabularies
$V_0$ and $V$, respectively, such that $V - V_0$ is finite and relational.

\item[(2)] $\mbG = \{\mcG_n : n \in \mbbN \}$ is a set of $L_0$-structures
which is a pregeometry, in the sense of 
Definition~\ref{definition of a structure being a pregeometry}.
Moreover, assume that the formulas $\theta_n(x_1, \ldots, x_{n+1}) \in L_0$, for $n \in \mbbN$,
define the pregeometry
in the sense of Definition~\ref{definition of a structure being a pregeometry}.

\item[(3)] For $n \in \mbbN$, $\mbK_n = \mbK(\mcG_n)$ is a set of expansions to $L$ of $\mcG_n$, and
$\mbK = \bigcup_{n \in \mbbN} \mbK_n$.
For each $\mcA \in \mbK$, $\cl_{\mcA}$ is, by definition, the
same as $\cl_{\mcA \uhr L_0}$, where the latter is 
the same as $\cl_{\mcG_n}$ for some $n$, because $\mcA \uhrc L_0 = \mcG_n$ for some $n$.

\item[(4)] Whenever $\mcM$ is represented and $\mcA \subseteq_{cl} \mcM$,
then $\mcA$ is represented.
\end{itemize}
}\end{assump}

\begin{rem}\label{remarks about second set of assumptions}{\rm
(i) Note that point (4) in Assumption~\ref{assumptions about pregeometries}
says that the class of represented structures is closed under closed substructures.\\
(ii) If the closure is trivial, then (4) is equivalent to the hereditary property (for $\mbK$).\\
(iii) Analogues of the main theorems of this section can be stated and proved without
assumption (4), but then, to get such results, the notion of `acceptance of substitutions' 
(Definition~\ref{k-independence hypothesis for classes of structures})
must be modified, and becomes more complicated.
The author opted, in this case, for simplicity rather than some more generality.
}\end{rem}

\begin{defin}\label{definition of dimension reduct}{\rm
Let $\mcA \in \mbK$ and let $d$ be a natural number.\\
(i) The {\bf \em $d$-dimensional reduct} of $\mcA$, denoted $\mcA \uhrc d$,
  is the weak substructure of $\mcA$
  which is defined as follows:
  \begin{itemize}
    \item[(a)] $\mcA \uhrc d$ has the same universe as $\mcA$.
    \item[(b)] Every symbol in the vocabulary of $L_0$ is interpreted in the
    same way in $\mcA \uhrc d$ as in $\mcA$.
    \item[(c)] For every relation symbol $R$ which belongs to the vocabulary
    of $L$ but {\em not} to the vocabulary of $L_0$, and for every tuple $\bar{a}$
    from the universe of $\mcA$,
    $$\bar{a} \in R^{\mcA \uhr d} \Longleftrightarrow \dim_{\mcA}(\bar{a}) \leq d \text{ and }
    \bar{a} \in R^{\mcA}.$$
  \end{itemize}
(ii) $\mbK \uhrc d = \big\{ \mcA \uhrc d : \mcA \in \mbK \big\}$.\\
(iii) $\mbK_n \uhrc d = \big\{ \mcA \uhrc d : \mcA \in \mbK_n \big\}$.
}\end{defin}

\begin{rem}\label{remark about maximal arity and dimension reducts}{\rm
(i) Observe that if there is no relation symbol whose arity is greater than $d$,
then for every $\mcA \in \mbK$, $\mcA \uhrc d = \mcA$; hence $\mbK \uhrc d = \mbK$ and $\mbK_n \uhrc d = \mbK_n$ for every $n$.\\
(ii) By Definition~\ref{definition of dimension reduct}, for every $n \in \mbbN$ and every positive $r \in \mbbN$, the sequence
$\mbK_n \uhrc 0, \mbK_n \uhrc 1, \ldots, \mbK_n \uhrc r$ satisfies the conditions
for $\mbC_0, \ldots, \mbC_r$  in 
Definition~\ref{definition of transitive closure of reduct conditional measure}.
Hence, for every $n \in \mbbN$ and every positive $r \in \mbbN$, the uniformly 
$(\mbK_n \uhrc 0, \ldots, \mbK_n \uhrc r-1)$-conditional
measure is well-defined on $\mbK_n \uhrc r$. 
}\end{rem}

\begin{rem}{\rm
Note that for an $L$-structure $\mcM \in \mbK$ we have different kinds of ``reducts'',
and the same symbol `$\uhr$' is used in all contexts, but the symbol following `$\uhr$' is a key,
besides the context, to what is meant. 
For a sublanguage $L' \subseteq L$, $\mcM \uhrc L'$ is the reduct of $\mcM$ to 
$L'$ in the usual ``language wise'' sense.
For a subset $X \subseteq M$, $\mcM \uhrc X$ denotes the substructure
of $\mcM$ which is generated by $X$.
And for a natural number $d$, $\mcM \uhrc d$ denotes the $d$-dimensional reduct 
of $\mcM$, which is a {\em weak} substructure of $\mcM$, but not necessarily a
substructure.
}\end{rem}

\begin{defin}\label{definition of dimension conditional measure}{\rm
(i) Let $\rho$ be equal to the largest arity of a relation symbol in the vocabulary of $L$. 
Note that if $r \geq \rho$ then for every $\mcA \in \mbK$, $\mcA \uhrc r = \mcA$;
hence $\mbK \uhrc r = \mbK$ and $\mbK_n \uhrc r = \mbK_n$ for every $n$.  \\
(ii) For every $n \in \mbbN$, let $\mbbP_{n,0}$ denote the uniform probability measure on $\mbK_n \uhrc 0$.
For every $n \in \mbbN$ and every positive $r \in \mbbN$, 
let $\mbbP_{n,r}$ denote the uniformly $(\mbK_n \uhrc 0, \ldots, \mbK_n \uhrc r-1)$-conditional
measure on $\mbK_n \uhrc r$.\\
(iii) The uniformly $(\mbK_n \uhrc 0, \ldots, \mbK_n \uhrc \rho - 1)$-conditional
measure $\mbbP_{n,\rho}$ on $\mbK_n = \mbK_n \uhrc \rho$ is also denoted by $\delta_n$ and called
the {\bf \em dimension conditional measure} on $\mbK_n$.
}\end{defin}

\begin{exam}\label{simplest examples of dimension conditional measure}{\rm
Suppose that $L$ and $\mbK_n$ are defined as in any of Examples~\ref{graphs with a unary predicate} 
--~\ref{coloured binary relations}, let $L_0$ be the language with empty vocabulary,
and let the underlying pregeometry be trivial.
If $L_0$ is defined as in the corresponding example, then the dimension conditional
measure on $\mbK_n$ is, by definition, the same as the uniformly 
$(\mbK_n \uhrc 0, \mbK_n \uhrc 1)$-conditional
measure on $\mbK_n$, which in turn is identical to the uniformly
$(\mbK_n \uhrc L_0)$-conditional measure, considered in the mentioned examples;
this follows straightforwardly from the definitions.
Examples with nontrivial underlying pregeometry will appear later.
}\end{exam}

\begin{defin}\label{definition of uniformly bounded}{\rm
We say that the pregeometry $\mbG = \{\mcG_n : n \in \mbbN\}$ is 
{\bf \em uniformly bounded} if there is a function $u : \mbbN \to \mbbN$
such that
for every $n \in \mbbN$ and every $X \subseteq |\mcG_n|$,
$\big|\cl_{\mcG_n}(X)\big| \leq u\big(\dim_{\mcG_n}(X)\big)$.
}\end{defin}

\begin{rem}\label{remark about uniformly bounded pregeometries}{\rm
The trivial pregeometries and the pregeometries obtained from vector spaces
over finite fields are uniformly bounded.
More examples of uniformly bounded pregeometries 
can be obtained by applying the variants of the amalgamation construction first
developed by E. Hrushovski which produce
countably categorical supersimple limit structures 
with rank 1 \cite{Hru, Eva}.
However, the cases of such constructions known to the author do {\em not}
produce pregeometries which are polynomially $k$-saturated for all $k$;
this can be seen by considering the arguments in Section~2 of \cite{Dj04}.
The author does not know an example of
a pregeometry (in the sense of this paper)
$\mbG = \{\mcG_n : n \in \mbbN\}$ which is not uniformly bounded
and such that each $\mcG_n$ is {\em finite}, as we always assume here.
}\end{rem}

\begin{termin}\label{some terminology}{\rm
When saying that two represented structures $\mcA$ and $\mcA'$ 
{\bf \em agree on $L_0$ and on closed proper substructures}
we mean that $\mcA \uhrc L_0 = \mcA' \uhrc L_0$ 
(so in particular, $\cl_{\mcA} = \cl_{\mcA'}$)
and whenever $\mcU \subseteq_{cl} \mcA$ and $\dim_{\mcA}(U) < \dim_{\mcA}(A)$,
then $\mcA \uhrc U = \mcA' \uhrc U$.
}\end{termin}

\noindent
The next definition generalizes the notion of `admitting substitutions' 
from Section~\ref{forbidden structures}
to the context of this section.

\begin{defin}\label{k-independence hypothesis for classes of structures}{\rm
Let $\mcA$ and $\mcA'$ be represented structures.
Note that, in part (i) and (ii) of this definition, the
property defined can only hold if $\mcA$ and $\mcA'$ 
agree on $L_0$ and on closed proper substructures;
so that is the situation which is of interest.
\\
(i) We say that $\mbK$ {\bf \em accepts the substitution $[\mcA \trir \mcA']$ over $L_0$} if 
whenever $\mcM$ is represented and $\mcA \subseteq_{cl} \mcM$,
then there is a represented $\mcN$ such that
$\mcN \uhrc L_0 = \mcM \uhrc L_0$,
$\mcN \uhrc |\mcA| = \mcA'$ and 
if $\mcU \subseteq_{cl} \mcN$,
$\dim_{\mcN}(U) \leq \dim_{\mcN}(A')$
and $\mcU \neq \mcA'$, then $\mcN \uhrc U = \mcM \uhrc U$.
\\
(ii) We say that $\mbK$ {\bf \em accepts $k$-substitutions over $L_0$}
if whenever 
$\mcA$ and  $\mcA'$ are represented structures
which agree on $L_0$ and on closed proper substructures, and 
$\dim_{\mcA}(A) = \dim_{\mcA'}(A') \leq k$,
then $\mbK$ accepts the substitution $[\mcA \trir \mcA']$ over $L_0$.
}\end{defin}

\begin{rem}\label{remark on supremum of arities and acceptance}{\rm
(i) It is easy to see the following: If there is, up to isomorphism, a 
unique represented structure with dimension 0, then $\mbK$ accepts $0$-substitutions
over $L_0$.\\
(ii) Let $\rho$ be the supremum of the arities of all relation symbols that
belong to the vocabulary of $L$ but {\em not} to the vocabulary of $L_0$.
It is straightforward to verify that if $\mbK$  accepts $\rho$-substitutions over $L_0$,
then, for every $k \in \mbbN$, $\mbK$  accepts $k$-substitutions over $L_0$.
}\end{rem}

\noindent
We now give examples of $\mbK$ which accept $k$-substitutions for all $k \in \mbbN$.
After Theorem~\ref{0-law for extension axioms in pregeometries}, which is about
$\mbK$ which do not satisfy this condition, we give more examples, which, for some $k$, do not not
satisfy $k$-substitutions.

\begin{exam}\label{example of coloured structures of first kind}{\rm
({\bf Coloured structures.})
For the sake of having a uniform terminology in this example, and the next,
let us have the following convention.
For $F = \{1\}$ let $L_F$ be the language with empty vocabulary $V_F$ and
let $\mbG^F = \{\mcG_n : n \in \mbbN\}$, where $\mcG_n$ is the unique
$L_F$-structure with universe $\{1, \ldots, n\}$.
In this case call $\mbG^F$ the {\em vector space pregeometry over $\{1\}$}.

For any finite field $F$, the
{\em vector space pregeometry over $F$} refers to the 
pregeometry $\mbG^F = \{\mcG_n : n \in \mbbN\}$ defined in 
Example~\ref{examples of polynomially k-saturated pregeometries};
so $\mcG_n$ is a vector space over $F$ of dimension $n$,
and $L_F$ and $V_F$ is the language and vocabulary, respectively, of $\mcG_n$.

Let $\mbG^F = \{\mcG_n : n \in \mbbN\}$ be the vector space pregeometry over $F$,
where $F$ is a finite field or $\{1\}$.
Then let $l \geq 2$ and assume that $L_{col} \supset L_F$, ``the colour language'' is the 
language with vocabulary $V_{col} = V_F \cup \{P_1, \ldots, P_l\}$ where all $P_i$
are unary relation symbols, representing {\em colours}.
Also assume that $L_{rel} \supset L_F$, ``the language of relations'', has a 
vocabulary $V_{rel}$ such that $V_{rel} - V_F$ contains only finitely many relation symbols, of any arity.
Let $L$ be the language with vocabulary $V = V_{col} \cup V_{rel}$.
For every positive $n \in \mbbN$ define $\mbK_n = \mbK(\mcG_n)$ to be 
set of expansions $\mcM$ of $\mcG_n$ to $L$ 
that satisfy the following three {\em $l$-colouring conditions}:
\begin{itemize}
  \item[(1)] $\mcM \models \forall x \big( P_1(x) \vee \ldots \vee P_l(x)\big)$.
  \item[(2)] For all distinct $i, j \in \{1, \ldots, l\}$, and all $a, b \in M - \cl_{\mcM}(\es)$ 
  such that $a \in \cl_{\mcM}(b)$,
  $\mcM \models \neg\big(P_i(a) \wedge P_j(b)\big)$. (In other words: any two linearly
  dependent non-zero elements must have the same colour.)
  \item[(3)] If $R \in V_{rel}$ has arity $m \geq 2$ and 
  $\mcM \models R(a_1, \ldots, a_m)$, 
    then there are $b, c \in \cl_{\mcM}(a_1, \ldots, a_m)$ such that 
    for every $k \in \{1, \ldots, l\}$, $\mcM \models \neg \big(P_k(b) \wedge P_k(c)\big)$;
    that is, at least two elements in $\cl_{\mcM}(a_1, \ldots, a_m)$ have different colours.
\end{itemize}
It is now straightforward to verify that, for every $F$ considered, 
$\mbK = \bigcup_{n \in \mbbN} \mbK_n$
accepts $k$-substitutions over $L_F$, for every $k \in \mbbN$.
And as mentioned in Example~\ref{examples of polynomially k-saturated pregeometries},
$(\mcG_n : n \in \mbbN)$ is polynomially $k$-saturated for every $k \in \mbbN$.
It is also uniformly bounded. 
Thus, with this setup of $(\mcG_n : n \in \mbbN)$ and $\mbK$ the premises of 
Theorems~\ref{0-1 law for pregeometries} and~\ref{third part of 0-1 law for pregeometries}
(below) are satisfied.
This example and the next will be studied more in
Sections~\ref{l-colourable structures} and~\ref{the uniform measure and l-colourable structures}.
}\end{exam}

\begin{exam}\label{example of coloured structures of second kind}{\rm
({\bf Strongly coloured structures.})
The colourings considered in the previous example are the convention 
within hypergraph theory \cite{Ber, JT}, 
but we would also like to consider another sort of colourings, called
{\em strong colourings} in the hypergraph context \cite{AH},
and we adopt the same terminology.
Here, $\mbG^F$, $L_F$, $L_{col}$, $L_{rel}$ and $L$ are defined as
in Example~\ref{example of coloured structures of first kind}.
Let $\mbK = \bigcup_{n\in\mbbN}\mbK_n$, where
$\mbK_n$ consists of those $L$-expansions $\mcM$ of $\mcG_n$ which satisfy 
(1) and (2) from the previous example and the following {\em strong $l$-colouring condition}:
\begin{itemize}
  \item[(3')] If $R \in V_{rel}$ has arity $m \geq 2$, 
  $\mcM \models R(a_1, \ldots, a_m)$, $b, c \in \cl_{\mcM}(a_1, \ldots, a_m)$ and $b$
  is independent from $c$ (i.e. $b \notin \cl_{\mcM}(c)$), 
   then for every $k \in \{1, \ldots, l\}$, $\mcM \models \neg \big(P_k(b) \wedge P_k(c)\big)$;
   that is, {\em every} pair of mutually independent elements
   $b$ and $c$ in the closure of $a_1, \ldots, a_m$ 
   have different colours.
\end{itemize}
Again, it is straightforward to verify that, for every $F$ considered, 
$\mbK$ accepts $k$-substitutions over $L_F$, for every $k \in \mbbN$.
}\end{exam}

\begin{exam}\label{partial colourings}{\rm ({\bf Other variations of coloured structures})
In the previous two examples, it is also possible to consider 
projective or affine spaces over a finite field, instead of a vector space.
And, by dropping condition (1), one can consider partial colorings 
or strong partial colourings.
For all these variations, $\mbK$ accepts $k$-substitutions for every $k \in \mbbN$.
}\end{exam}

\begin{exam}\label{random relations}{\rm
({\bf Random relations on a vector space})
Let $F$ be a finite field and let $L_F$ and $\mbG^F = \{\mcG_n : n \in \mbbN\}$
be as in Example~\ref{examples of polynomially k-saturated pregeometries}.
Let $L$ be the language whose vocabulary consists of the symbols in the vocabulary
of $L_F$ and, in addition, relation symbols $R_1, \ldots, R_{\rho}$, of any arity.
For every $n$, let $\mbK_n = \mbK(\mcG_n)$ be the set of all $L$-structures $\mcM$ 
such that $\mcM \uhrc L_0 = \mcG_n$.
It is straightforward to verify that, for every $k \in \mbbN$,
$\mbK = \bigcup_{n\in\mbbN}\mbK_n$ accepts $k$-substitutions over $L_F$.
A similar example can be constructed over projective or affine
spaces over $F$.
}\end{exam}

\begin{rem}\label{algebraic approach to random edges}{\rm
({\bf An algebraic approach to adding ``pseudo-random'' edges})
Here we sketch an algebraic approach to expanding $F$-vector spaces by
a binary irreflexive symmetric relation.
The graph structure itself will, in the limit, be the same as the one obtained
in the previous example when only one relation symbol $R_1 = R$ is considered and 
always interpreted as an irreflexive and symmetric relation.
But in the algebraic approach it is not sufficiently clear to the author
how the vector space structure interacts with the graph structure
and therefore the question whether $\mbK$ defined below accepts 2-substitutions over the vector
space language is left open.

Let $F = \mathbb{F}_p$ be the finite field of order $p$, where $p$ is a prime
which is congruent to 1 modulo~4.
As in the previous example, let $L_F$ 
be as in Example~\ref{examples of polynomially k-saturated pregeometries}.
Every field of order $p^n$, denoted $\mathbb{F}_{p^n}$, 
can be viewed as a vector space over $F = \mathbb{F}_p$, and this 
vector space (of dimension $n$), formalised as an $L_F$-structure, is denoted $\mcV_n$.
Let the vocabulary of the ``graph language'', $L_g$, contain only one binary relation symbol $R$,
let $\mbK_n^g$ be the set of undirected graphs (as $L_g$-structures) with vertices $1, \ldots, n$
and let $\mbK^g = \bigcup_{n \in \mbbN} \mbK_n^g$. Then let $L$ be the
language whose vocabulary is the union of the vocabularies of $L_F$ and $L_g$.
Every $\mcV_n$ can be expanded to an $L$-structure, denoted $\mcV_n^g$, so that 
$\mcV_n^g \uhrc L_g$ is
an undirected graph, by letting $\mcV_n^g \models R(a,b)$ if and only if 
$a-b$ is a square in the field $\mathbb{F}_{p^n}$; so $\mcV_n^g \uhrc L_g$ is a Paley graph.
By results about Paley graphs (see Chapter~13 of \cite{Boll})
it follows that, for every extension axiom $\varphi$ 
of $\mbK^g$, $\varphi$ is true in $\mcV_n^g$ for all sufficiently large $n$.
By compactness there is an infinite $L$-structure $\mcV$ such that $\mcV \uhrc L_F$
is a vector space over $F$ and $\mcV \uhrc L_g$ is an undirected graph which satisfies
every extension axiom of $\mbK^g$.
Now we can let $\mcG_n$ be an $n$-dimensional vector space over $F$, viewed as an $L_F$-structure,
and let $\mbK_n = \mbK(\mcG_n)$ be the set of expansions, $\mcM$, to $L$ of $\mcG_n$
such that $\mcM$ is isomorphic with some substructure of $\mcV$.
We may now ask whether it is true that, for every $k$, $\mbK = \bigcup_{n \in \mbbN} \mbK_n$
accepts $k$-substitutions over $L_F$ and/or is polynomially $k$-saturated.
Since $\mcV \uhrc L_g$ satisfies every extension axiom of $\mbK^g$ (and possibly using
more information about Paley graphs) one may be tempted to guess that the answers are yes
in both cases.
However, when dealing with the question of whether $\mbK$ accepts 2-substitutions over $L_F$
we need to understand (it seems) what graphs can appear as $\mcH = \mcM \uhrc L_g$ where
$\mcM$ is a substructure of $\mcV$, so in particular, $M$ is a linearly closed subset of $V$.
This seems to involve deeper understanding of the interaction between
the vector space structure of $\mcV_n$ and the multiplicative structure of $\mathbb{F}_{p^n}$
for all sufficiently large $n \in \mbbN$.
}\end{rem}

\begin{exam}\label{3-hypergraph definable from random graph}{\rm
({\bf Hypergraph and random graph on a vector space})
Let $F$ be a finite field and let $L_F$ and $\mbG^F = \{\mcG_n : n \in \mbbN\}$
be as in Example~\ref{examples of polynomially k-saturated pregeometries}.
Let $L$ be the language whose vocabulary consists of the symbols in the vocabulary
of $L_F$ and, in addition, relation symbols $E$ and $R$ where $E$ is binary
and $R$ is ternary.
For every $n \in \mbbN$, let $\mbK_n = \mbK(\mcG_n)$ be the set of $L$-structures $\mcM$ such
that $\mcM \uhrc L_F = \mcG_n$, $E$ is interpreted as an irreflexive and symmetric relation,
so we call $E$-relationships edges, and, for all $a, b, c \in M$, $(a,b,c) \in R^{\mcM}$
if and only if the subspace spanned by $a, b$ and $c$ contains an odd number of edges.
We show that $\mbK = \bigcup_{n\in\mbbN}\mbK_n$ accepts $k$-substitutions
for all $k \in \mbbN$. As was mentioned in Remark~\ref{remark on supremum of arities and acceptance},
it suffices to show that $\mbK$ accepts $3$-substitutions.
Let $L'$ be the sublanguage of $L$ in which the symbol $R$ has been removed, but all other symbols
have been kept.
Observe that, for every $\mcM \in \mbK$ and for all $a,b,c \in M$, whether $\mcM \models R(a,b,c)$,
or not, is determined by the substructure of $\mcM \uhrc L'$
whose universe is the linear span of $a$, $b$ and $c$.
This implies that it suffices to show that $\mbK$ accepts $2$-substitutions.
Since the only restrictions on $E$ is that it is interpreted as an irreflexive
and symmetric relation, it follows that in whichever way we expand $\mcG_n$ with edges,
we get $\mcM \uhrc L'$ for some $\mcM \in \mbK_n$.
This implies that $\mbK$ accepts $2$-substitutions.

Suppose that $L^*$ is the sublanguage of $L$ where the symbol $E$ has been removed,
but all other symbols have been kept, and let $\mbK^* = \bigcup_{n\in\mbbN}\mbK^*_n$,
where $\mbK^*_n = \{\mcM \uhrc L^* : \mcM \in \mbK_n\}$. 
It is, when writing this, not clear to the author if $\mbK^*$ accepts $3$-substitutions, or not.
}\end{exam}

\begin{exam}\label{more excentric example}{\rm
In this example, pairs of elements as well as elements can be coloured
and some restrictions are imposed.
Suppose that, for every $n$, $\mcG_n$ is a projective space over
the 2-element field and let $L_0$ be the language of $\mcG_n$. 
Let $L \supset L_0$ contain, besides the symbols of $L_0$,
three unary relation symbols $P_1$, $P_2$, $P_3$, three binary 
relation symbols $R_1, R_2, R_3$ and one ternary relation symbol $S$.
We can think of the $P_i$ as colours 
of elements, and the $R_i$ as colours of pairs.
For every $n$, $\mbK_n = \mbK(\mcG_n)$ consists of all expansions 
$\mcM$ of $\mcG_n$ to $L$ which satisfy the following conditions:
\begin{itemize}
	\item[(a)] For every 2-dimensional subspace $X \subseteq M$,
	if no pair $(a,b) \in X^2$ is coloured, then at least one point in $X$ is coloured.
	\item[(b)] For every two dimensional subspace $X \subseteq M$,
	if some pair $(a,b) \in X^2$ is coloured, then there are {\em not} two different
	points in $X$ with the same colour (but two different points may be uncoloured).
	\item[(c)] If $\mcM \models S(a,b,c)$, then $\{a, b, c\}$ is independent 
	and if $(d_1, d_2), (e_1, e_2) \in \cl_{\mcM}(a,b,c)$, then $(d_1, d_2)$ and
	$(e_1, e_2)$ do not have the same colour (but both may be uncoloured).
\end{itemize}
We show that $\mbK$ accepts $3$-substitutions over $L_0$.
Since no relation symbol has arity greater than 3 it follows 
(see Remark~\ref{remark on supremum of arities and acceptance})
that $\mbK$ accepts $k$-substitutions over $L_0$ for every $k \in \mbbN$.

Let $\mcA, \mcA'$ be represented
and assume that $\mcA \uhrc L_0 = \mcA' \uhrc L_0$
and that $\mcA$ and $\mcA'$ agree on all closed proper substructures.
We must show that if $\mcM$ is represented and  $\mcA \subseteq_{cl} \mcM$,
then there exists a represented $\mcN$ such that
$\mcN \uhrc L_0 = \mcM \uhrc L_0$, $\mcN \uhrc A = \mcA'$   and whenever
$\mcU \subseteq_{cl} \mcN$, $\dim_{\mcN}(U) \leq \dim_{\mcN}(A')$, and $U \neq A'$,
then $\mcN \uhrc U = \mcM \uhrc U$.

First suppose that $\dim_{\mcM}(A) = 1$.
Let $\mcM' = \mcM[\mcA \trir \mcA']$, according to Definition~\ref{definition of substitutions}
(Since $\mcA \uhrc L_0 = \mcA' \uhrc L_0$, the substitution involves only interpretations
of relation symbols).
Then go through all $\mcB \subseteq_{cl} \mcM'$ of dimension 2;
whenever we meet such $\mcB$ which is forbidden we can change some binary
relationships ($R_i$, $i = 1,2,3$), but not change any unary relationships ($P_i$, $i = 1,2,3$),
and thus get a permitted substructure.
When this has been done for all 2-dimensional closed substructures, call the
result $\mcM''$; so all 2-dimensional substructures of $\mcM''$ are permitted.
Then we can just remove all $S$-relationships from $\mcM''$
so that in the resulting structure $\mcN$
the interpretation of $S$ is empty. It now follows from the construction of $\mcN$
and (a) -- (c) that $\mcN$ is represented.
And whenever $U \subseteq N$ is 1-dimensional and different from $A'$,
then $\mcN \uhrc U = \mcM \uhrc U$.

Now suppose that $\dim_{\mcM}(A) = 2$.
Let $\mcM' = \mcM[\mcA \trir \mcA']$.
Then $\mcM'$ and $\mcM$ agree on all closed 1- or 2-dimensional subsets
which are different from $A'$.
By removing all $S$-relationships from $\mcM'$ we get $\mcN$ which is represented
and such that $\mcN$ and $\mcM$ agree on all closed 1- or 2-dimensional subsets
which are different from $A'$. Moreover, $\mcN \uhrc A' = \mcA'$.

Finally, suppose that $\dim_{\mcA}(A) = 3$.
Both $\mcA$ and $\mcA'$ satisfy (a) -- (c) (because they are permitted)
and $\mcA$ and $\mcA'$ agree, by assumption, on substructures
of dimension 2. Hence $\mcN = \mcM[\mcA \trir \mcA']$ and $\mcM$
agree on subsets of dimension 2 and on closed subsets of dimension 3
which are different from $A'$. Since $\mcA'$ is represented, and hence satisfies (a) -- (c),
$\mcN$ is represented.
}\end{exam}

\noindent
The next lemma tells that the notion of `accepting $k$-substitutions over $L_0$'
is indeed a generalization of the notion of `admitting $k$-substitutions'.

\begin{lem}\label{admittance of substitutions implies acceptance}
Let $L_0$ be the language with empty vocabulary and let $\mcG_n$
be the unique $L_0$-structure with universe $\{1, \ldots, m_n\}$ (with the trivial pregeometry)
where $\lim_{n\to\infty} m_n = \infty$.
Let $L$ be any language with finite relational vocabulary.
Suppose that, for every $n$, $\mbK_n$ is a set of $L$-structures with universe
$\{1, \ldots, m_n\}$; in other words, $\mbK_n = \mbK(\mcG_n)$ is a set of expansions
of $\mcG_n$ to $L$; and let $\mbK = \bigcup_{n \in \mbbN}\mbK_n$.
For every $k \in \mbbN$, if $\mbK$ admits $k$-substitutions 
(in the sense of Definition~\ref{definition of admitting substitutions}),
then $\mbK$ accepts $k$-substitutions over $L_0$.
\end{lem}

\noindent
{\em Proof.} One just checks that, under the assumptions, $\mbK$
does indeed accept $k$-substitutions over $L_0$, according to
Definition~\ref{k-independence hypothesis for classes of structures}.
\hfill $\square$
\\

\noindent
Recall Assumptions~\ref{assumptions about pregeometries}
and Definition~\ref{definition of dimension conditional measure}~(iii).

\begin{defin}{\rm
For every $n \in \mbbN$ and every $L$-sentence $\varphi$, 
\begin{equation*}
\text{let $\delta_n(\varphi)$ be an abbreviation for
$\delta_n\big(\{\mcM \in \mbK_n : \mcM \models \varphi\}\big)$.}
\end{equation*}
}\end{defin}

\begin{theor}\label{0-1 law for pregeometries}
Let $k > 0$.
Suppose that $(\mcG_n : n \in \mbbN)$ is uniformly bounded, polynomially $k$-saturated and that 
$\mbK = \bigcup_{n \in \mbbN}\mbK(\mcG_n)$ accepts $k$-substitutions over $L_0$.
Then:
\begin{itemize}
	\item[(i)] For every $(k-1)$-extension axiom $\varphi$ of $\mbK$, 
	$\lim_{n \to \infty} \delta_n(\varphi) = 1$.
	\item[(ii)] $\mbK$ is polynomially $k$-saturated.
\end{itemize}
\end{theor}

\begin{theor}\label{third part of 0-1 law for pregeometries}
Suppose that $(\mcG_n : n \in \mbbN)$ is uniformly bounded
and polynomially $k$-saturated for every $k \in \mbbN$.
Also assume that $\mbK = \bigcup_{n \in \mbbN}\mbK(\mcG_n)$ accepts $k$-substitutions over $L_0$
for every $k \in \mbbN$.
Then, for every $L$-sentence $\varphi$, either $\lim_{n \to \infty} \delta_n(\varphi) = 0$ or 
$\lim_{n \to \infty} \delta_n(\varphi) = 1$.
\end{theor}

\noindent
For the last theorem of this section we need a definition.

\begin{defin}\label{definition of independent amalgamation property}{\rm
We say that $\mbK$ has the {\bf \em independent amalgamation property} if
the following holds:
Whenever $\mcA$, $\mcB_1$, $\mcB_2$ are represented,
$\mcA \subseteq_{cl} \mcB_i$, for $i = 1,2$, and
$B_1 \cap B_2 = A$, then there is a represented $\mcC$ such that
$\mcB_i \subseteq_{cl} \mcC$ for $i = 1,2$. 
}\end{defin}

\begin{theor}\label{0-law for extension axioms in pregeometries}
Suppose that $(\mcG_n : n \in \mbbN)$ is uniformly bounded
and polynomially $k$-saturated for every $k \in \mbbN$.
Assume that, up to isomorphism, there is a unique represented structure,
with respect to $\mbK = \bigcup_{n \in \mbbN} \mbK(\mcG_n)$, with dimension 0
(a particular case of this is when $\cl(\es) = \es$).
Let $k \in \mbbN$ be minimal such that $\mbK$ does not accept
$k$-substitutions over $L_0$ and suppose that $\mcA$ and $\mcA'$ are represented
structures (with respect to $\mbK$) such that
$\mcA$ and $\mcA'$ have dimension $k$, agree on $L_0$ and on closed proper substructures,
$\mbK$ accepts the substitution $[\mcA' \trir \mcA]$ over $L_0$,
but does {\rm not} accept the substitution $[\mcA \trir \mcA']$
over $L_0$. Then at least one of the following holds:
\begin{itemize}
	\item[(i)] $\mbK$ does not have the independent amalgamation property.
	\item[(ii)] There are $\beta < 1$ and extension axioms $\varphi$ and $\psi$ such that
	for all sufficiently large $n$, $\delta_n(\varphi \wedge \psi) < \beta$.
	If $k > 1$, then $\lim_{n \to \infty} \delta_n(\varphi \wedge \psi) = 0$. 
\end{itemize}
\end{theor}

\begin{rem}\label{remark about details of amalgamation and phi and psi}{\rm
The proof of Theorem~\ref{0-law for extension axioms in pregeometries}
shows that if the assumptions of the theorem hold and one 
particular instance of the independent amalgamation property
is satisfied, then case (ii) holds;
more information about this instance of independent amalgamation and 
$\varphi$ and $\psi$ is given by the proof.
}\end{rem}

\noindent
The proofs of Theorems~\ref{0-1 law for pregeometries} 
--~\ref{0-law for extension axioms in pregeometries}
are given in the next section.

\begin{exam}\label{Henson examples}{\rm
({\bf Forbidden weak substructures})
We will prove a dichotomy, stated by the corollary below, which is
analogous to Theorem~\ref{dichotomy for forbidden weak substructures}, 
which was proved (using Theorem~\ref{upper limit on the probability of an extension axiom}) 
in Example~\ref{example of forbidden weak substructures}.

Let $\mbG = (\mcG_n : n \in \mbbN)$, where all $\mcG_n$ are $L_0$-structures,
be a pregeometry which satisfies the assumptions of Theorems~\ref{0-1 law for pregeometries} 
--~\ref{0-law for extension axioms in pregeometries}.
We also assume that $\mbG$ has the independent amalgamation property
(in the same sense as in Definition~\ref{definition of independent amalgamation property}
if $\mbK$ is replaced by $\mbG$). 
These assumptions hold for $\mbG = \mbG^F$ as in 
Example~\ref{examples of polynomially k-saturated pregeometries} 
where the members of $\mbG^F$ are vector spaces over the finite field $F$,
as well as for projective and affine versions of these spaces.
Let $L_{rel}$ be a language with relational vocabulary $\{R_1, \ldots, R_s\}$,
and let $L$ be the language whose vocabulary is the union of the vocabularies of
$L_0$ and $L_{rel}$.
Using Henson's terminology in \cite{Hen72}, 
we say that an $L_{rel}$-structure $\mcM$ is {\em decomposable} if there are different
$L_{rel}$-structures $\mcA$ and $\mcB$ such that $M = A \cup B$,
$\mcA \uhrc A \cap B \ = \ \mcB \uhrc A \cap B$ and for every $i = 1, \ldots, s$,
$(R_i)^{\mcM} = (R_i)^{\mcA} \cup (R_i)^{\mcB}$. Otherwise we call $\mcM$ {\em indecomposable}.
Suppose that $\mbF$ is a set of finite indecomposable $L_{rel}$-structures such that if $\mcA, \mcB \in \mbF$
and $\mcA \neq \mcB$, then $\mcA$ is not weakly embeddable into $\mcB$.
Let $\mbK_n = \mbK(\mcG_n)$ be the set of $L$-structures $\mcM$
such that $\mcM \uhrc L_0 = \mcG_n$ and no $\mcF \in \mbF$ can be weakly embedded into $\mcM \uhrc L_{rel}$,
and let $\mbK = \bigcup_{n\in\mbbN}\mbK_n$.
Note that one of the assumptions on $\mbF$ implies that every 
$\mcF \in \mbF$ is {\em minimal} in the sense that if $\mcF'$ is a proper weak substructure
of $\mcF$, then $\mcF'$ can be weakly embedded into
$\mcM \uhrc L_{rel}$ for some $\mcM \in \mbK$.
From the indecomposability of the members of $\mbF$ it follows, 
in essentially the same way as the (straightforward) 
proofs of Lemma~1.1 and Theorem~1.2~(i) in \cite{Hen72},
that $\mbK$ has the independent amalgamation property.

Consider the following statement:

\begin{itemize}
\item[($*$)] There are $\mcF \in \mbF$, a relation symbol $R_i$ 
and $\bar{a} \in (R_i)^{\mcF}$ 
such that $\rng(\bar{a})$ is a proper subset of $F$.
\end{itemize}

\noindent
{\bf Corollary to Theorems~\ref{0-1 law for pregeometries} 
--~\ref{0-law for extension axioms in pregeometries}.}{\em \ 
(i) If ($*$) holds, then there are $\beta < 1$ and extension axioms $\varphi$ and $\psi$ of $\mbK$
such that for all sufficiently large $n$, 
$\delta_n(\varphi \wedge \psi) < \beta$, and if $|\rng(\bar{a})| > 1$,
then $\lim_{n\to\infty} \delta_n(\varphi \wedge \psi) = 0$.\\
(ii) If ($*$) does not hold, then, for every $k \in \mbbN$, 
$\mbK$ accepts $k$-substitutions and is
polynomially $k$-saturated, 
for every extension axiom $\varphi$ of $\mbK$, $\lim_{n \to \infty} \delta_n(\varphi)$ = 1,
and $\mbK$ has a zero-one law with respect to the probability measures $\delta_n$.
\\}

\noindent
{\em Proof.}
We first prove (ii), so suppose that ($*$) does not hold.
We only need to prove that $\mbK$ accepts $k$-substitutions for every $k$,
since the other claims then follow from Theorems~\ref{0-1 law for pregeometries} 
and~\ref{third part of 0-1 law for pregeometries}.
Let $\mcA$ and $\mcA'$ be represented structures, with respect to $\mbK$, that agree on $L_0$ and on 
closed proper substructures, and suppose that $\mcA \subseteq \mcM \in \mbK$.
Moreover, suppose (for a contradiction) that $\mcN = \mcM[\mcA \trir \mcA']$ is forbidden,
so there is $\mcF \subseteq_w \mcN \uhrc L_{rel}$ such that $\mcF$ is isomorphic to some member of $\mbF$.
We may, without loss of generality, assume that for any $i$, if any $R_i$-relationship is removed
from $\mcA'$, giving $\mcA''$, then $\mcM[\mcA \trir \mcA'']$ is represented.
Since $\mcM$ is represented, $F$ must contain some element from
$|\mcA'|$. Since $\mcA'$ is represenetd, $\mcF$ is not a weak substructure of $\mcA' \uhrc L_{rel}$,
so $F$ must also contain some element in $|\mcN| - |\mcA'|$.
As $\mcF \subseteq_w \mcN \uhrc L_{rel} = \mcM[\mcA \trir \mcA'] \uhrc L_{rel}$ 
and $\mcM \in \mbK$, there is some $i$ and
$R_i$-relationship $\bar{a} \in (R_i)^{\mcA'} \subseteq (R_i)^{\mcN}$ such that 
$\rng(\bar{a}) \subseteq |\mcF|$.
But then $\rng(\bar{a})$ is a proper subset of $F$, which contradicts the assumption that
($*$) is false (since we can, if necessary, remove some relationships whose range
includes elements from $|\mcF| - |\mcA'|$, to ``uncover'' $\mcF$ weakly embedded into 
$\mcN \uhrc L_{rel}$).

Now we prove (i), so suppose that ($*$) holds.
Let $\mcF \in \mbF$ and $\bar{a} \in (R_i)^{\mcF}$
be such that $\rng(\bar{a})$ is a proper subset of $F$ and such that
the removal of the $R_i$-relationship $\bar{a}$ produces a structure $\mcP$
which is (weakly) embeddable into $\mcN \uhrc L_{rel}$ for some $\mcN \in \mbK$.
Let $d = |P|$, let $v_1, \ldots, v_d$ be a basis of $\mcG_d$, 
and let $f : P \to \{v_1, \ldots, v_d\}$ be a bijection.
Then let $\mcM$ be the $L$-structure which is obtained by expanding
$\mcG_d$ in such a way that $f : P \to \mcM \uhrc L_{rel}$ becomes an embedding
and if $\bar{b}$ contains an element not in $\{v_1, \ldots, v_d\}$,
then $\bar{b}$ is not a $R_j$-relationship for any $j$.
Then $\mcM \in \mbK$.
To simplify notation, we may assume that $F = P = \{v_1, \ldots, v_d\}$, so $\mcP \subseteq \mcM$.
Let $\mcA = \mcM \uhrc \cl_{\mcM}(\bar{a})$ and let $\mcA'$ be the structure
obtained from $\mcA$ by adding the $R_i$-relationsship $\bar{a}$, but making no other changes.
Then the $L_{rel}$-reduct of $\mcM[\mcA \trir \mcA']$ contains a copy of $\mcF$, so it is forbidden,
and hence $\mbK$ does not accept the substitution $[\mcA \trir \mcA']$ over $L_0$.
But $\mbK$ accepts the substitution $[\mcA' \trir \mcA]$ over $L_0$, because its effect is
only to remove a relationship and this can never create a forbidden structure.

As mentioned before the corollary, $\mbK$ has the independent amalgamation property, so by
Theorem~\ref{0-law for extension axioms in pregeometries}, 
there are $\beta < 1$ and extension axioms $\varphi$ and $\psi$ of $\mbK$
such that for all sufficiently large $n$, 
$\delta_n(\varphi \wedge \psi) < \beta$.
Moreover, if $|\rng(\bar{a})| > 1$ then, as $\rng(\bar{a}) \subseteq P$
and $P = \{v_1, \ldots, v_d\}$ is a basis of $\mcM$,
it follows that $\dim(\mcA) = \dim_{\mcM}(\bar{a}) > 1$, and hence
(by Theorem~\ref{0-law for extension axioms in pregeometries}) 
$\lim_{n\to\infty} \delta_n(\varphi \wedge \psi) = 0$.
\hfill $\square$
}\end{exam}

\begin{exam}\label{example of l-colourable structures}{\rm
({\bf $l$-Colourable structures, and strongly $l$-colourable structures})
Let $F$, $\mbK_n$, $L_{rel}$ and $L$ be as in Example~\ref{example of coloured structures of first kind}
(or as in Example~\ref{example of coloured structures of second kind})
and let 
$\mbC_n = \{\mcM \uhrc L_{rel} : \mcM \in \mbK_n\}$ and 
$\mbC = \bigcup_{n\in \mbbN}\mbC_n$.
Suppose that all relation symbols of $L_{rel}$ have arity at least $2$
and let $R$ be one which has minimal arity, which we denote by $k$. 
Assume that $l \geq k$ (or $l \geq$ maximal arity if we consider {\em strongly} $l$-colourable structures).
Since one can not add arbitrarily many new $R$-relationships
to a sufficiently large independent subset of a structure $\mcM \in \mbC$ without finally getting
forbidden structure, i.e. one that can not be (strongly) $l$-coloured, one can show 
(but we omit the details) that
$\mbC$ does not accept $k$-substitutions over $L_F$.
On the other hand, we can always remove an $R$-relationship from a represented structure
without producing a forbidden one. It follows that
there are represented structures $\mcA$ and $\mcA'$ with dimension $k$ 
which agree on $L_F$ and on closed proper substructures, 
the substitution $[\mcA' \trir \mcA]$ over $L_F$ is accepted, but not the 
substitution $[\mcA \trir \mcA']$ over $L_F$.
It follows from
Theorem~\ref{0-law for extension axioms in pregeometries} and since $k > 1$
that either $\mbC$ does not have the independent amalgamation property {\em or} that
there are extension axioms $\varphi$ and $\psi$ of $\mbC$ such that
$\lim_{n\to\infty} \delta_n(\varphi \wedge \psi) = 0$, where $\delta_n$ is
the dimension conditional measure on $\mbC_n$.

If the only symbol of the vocabulary of $L_{rel}$ which does not belong to 
the vocabulary of $L_F$ is a binary relation symbol $R$ and $l = 2$,
then, by considering a 5-cycle (which cannot be 2-coloured), 
it is easy to see that $\mbC$ does not have
the independent amalgamation property, since that would force a 5-cycle into
some member of $\mbC$. It is also straightforward to see, by considering 5-cycles and 3-cycles,
that if an $L_{rel}$-structure $\mcM$ satisfies all 3-extension axioms of 
$\mbC$, then it is not $2$-colourable. 
In Sections~\ref{l-colourable structures}--\ref{the uniform measure and l-colourable structures}
we will see that, nevertheless, for $F = \{1\}$, i.e. the trivial underlying pregeometry,
and any $L_{rel}$ as in the beginning of the example, $\mbC$ has a zero-one law
for $\delta_n$, as well as for the uniform probability measure.
(The corresponding statement for a finite field $F$, giving a nontrivial underlying pregeometry,
remains open.)
}\end{exam}

\section{Proofs of Theorems~\ref{0-1 law for pregeometries}, 
\ref{third part of 0-1 law for pregeometries} and~\ref{0-law for extension axioms in pregeometries}}
\label{proofs of main results}

\noindent
Remember that Theorems~\ref{0-1 law for pregeometries} --~\ref{0-law for extension axioms in pregeometries}
take place within the setting of Assumptions~\ref{first assumptions on pregeometries}
and~\ref{assumptions about pregeometries}.
Therefore Assumptions~\ref{first assumptions on pregeometries}
and~\ref{assumptions about pregeometries} are active throughout this section.

\subsection{Proof of Theorem~\ref{0-1 law for pregeometries}}\label{proof of 0-1 law for pregeometries}

We are assuming that $\mbG = \{\mcG_n :  n \in \mbbN\}$ is a set of $L_0$-structures
and that $\mbG$ is a pregeometry.
Let $k > 0$.
Suppose that $(\mcG_n : n \in \mbbN)$ is
polynomially $k$-saturated 
and that $\mbK = \bigcup_{n \in \mbbN}\mbK_n$, where $\mbK_n = \mbK(\mcG_n)$, 
accepts $k$-substitutions over $L_0$.
This means that there exists a sequence of numbers $(\lambda_n : n \in \mbbN)$
such that $\lim_{n\to\infty}\lambda_n = \infty$ and a polynomial $P(x)$ such that
for every $n \in \mbbN$:
	\begin{itemize}
		\item[(a)] $\lambda_n \leq |G_n| \leq P(\lambda_n)$, and
		\item[(b)] whenever $\mcA$ and $\mcB$ are represented, $\mcA \subset_{cl} \mcB$
		and $\dim_{\mcB}(A) + 1 = \dim_{\mcB}(B) \leq k$, 
		then the $\mcB / \mcA$-multiplicity of $\mcG_n$ is at least $\lambda_n$.
	\end{itemize}
We must prove the following:
\begin{itemize}
	\item[(i)] For every $(k-1)$-extension axiom $\varphi$ of $\mbK$, 
	$\lim_{n \to \infty} \delta_n(\varphi) = 1$.
	\item[(ii)] $\mbK$ polynomially $k$-saturated.
\end{itemize}
Part (i) will be reduced to the problem of proving that the $\delta_n$-probability
that $\mcM \in \mbK_n$ is sufficiently saturated, 
in the sense of Definition~\ref{definition of (m,k)-saturation} below,
tends to 1 as $n$ tends to infinity.

Recall, from Definition~\ref{definition of dimension conditional measure}~(i),
that $\rho$ is the supremum of the arities of all relation symbols that belong
to the vocabulary of $L$, but not to the vocabulary of $L_0$.
From Assumptions~\ref{first assumptions on pregeometries}
and~\ref{assumptions about pregeometries},
Definition~\ref{definition of dimension reduct}
and Remark~\ref{remark about maximal arity and dimension reducts} it follows that whenever 
$d, n \in  \mbbN$ and $\mcM \in \mbK_n \uhrc d$,
then $\cl_{\mcM}$ coincides with $\cl_{\mcG_n}$ which is the same as $\cl_{\mcM \uhr L_0}$
since $\mcM \uhrc L_0 = \mcG_n$.
Also, if $d \geq \rho$, then for every $\mcM \in \mbK$, $\mcM \uhrc d = \mcM$.

In this proof, and the proofs of Theorems~\ref{third part of 0-1 law for pregeometries}
and~\ref{0-law for extension axioms in pregeometries}, we often work with
$\mbK \uhr d$, for some $d \in \mbbN$,
and consider structures which are represented, permitted, or forbidden,
{\em with respect to $\mbK \uhrc d$}.
Recall, from Definition~\ref{definition of dimension conditional measure}~(iii),
that $\delta_n$ is an abbreviation for $\mbbP_{n,\rho}$.
Essentially, the next definition just repeats point (2) from 
Definition~\ref{definition of polynomial k-saturation for classes} 
in the case of $\mbK \uhrc d$ (instead of $\mbK$), but it will be convenient to use
the terminology defined below.

\begin{defin}\label{definition of (m,k)-saturation}{\rm
(i) Let $d, m \in \mbbN$ and $\mcM \in  \mbK \uhrc d$.
We say that $\mcM$ is {\bf \em $(m,k)$-saturated with respect to $\mbK \uhrc d$}
if the following holds:
\begin{itemize}
	\item[] Whenever $\mcA$ and $\mcB$ are represented {\em with respect to $\mbK \uhrc d$},
	$\mcA \subset_{cl} \mcB$ and
  $\dim_{\mcB}(A) + 1 = \dim_{\mcB}(B) \leq k$, 
		then the $\mcB / \mcA$-multiplicity of $\mcM$ is at least $m$.
\end{itemize}
(i) Since $\mcM \uhrc \rho = \mcM$ for every $\mcM \in \mbK$, we say that $\mcM \in \mbK$ is  
{\bf \em $(m,k)$-saturated with respect to $\mbK$} if $\mcM$ is 
$(m,k)$-saturated with respect to $\mbK \uhrc \rho$.
}\end{defin}

\begin{defin}\label{definition of the sigma-functions}{\rm
For $r \in \mbbN$ we inductively we define functions $\sigma^r : \mbbN \to \mbbN$. 
Let $\sigma^0(x) = x$ for all $x \in \mbbN$.
Let $\sigma^{r+1}(x) = \lfloor \sqrt{\sigma^r(x)}\rfloor$ for all $x \in \mbbN$.
}\end{defin}

\noindent
Note that for every $r \in \mbbN$, $\lim_{n \to \infty} \sigma^r(n) = \infty$.
By assumption, $\lim_{n \to \infty}\lambda_n = \infty$, so for every $r \in \mbbN$,
$\lim_{n \to \infty}\sigma^r(\lambda_n) = \infty$; this will be used later.

Let $\varphi$ be a $(k-1)$-extension axiom.
In order to prove (i) we need to show that
\begin{equation*}
\lim_{n \to \infty} \delta_n\big(\{\mcM \in \mbK_n : \mcM \models \varphi\}\big) = 1. \tag{1}
\end{equation*}
By assumption, $\varphi$ is the $\mcB/\mcA$-extension axiom for some 
$\mcA \subset \mcB \subseteq \mcM$ such that $\mcM$ is represented
with respect to $\mbK = \mbK \uhrc \rho$, both $A$ and $B$ 
are closed in $\mcM$ and
$\dim_{\mcB}(B) \leq k$; in particular $\dim_{\mcB}(A) < \dim_{\mcB}(B)$.
Then, letting $l = \dim_{\mcB}(B) - \dim_{\mcB}(A)$,
there are closed substructures $\mcB_0, \ldots, \mcB_l$ of $\mcM$ such that
$\mcA = \mcB_0 \subset \mcB_1 \subset \ldots \subset \mcB_l = \mcB$ and
$\dim_{\mcB}(B_i) + 1 = \dim_{\mcB}(B_{i+1})$ for $i = 0, \ldots, l-1$.
By Assumption~\ref{assumptions about pregeometries} (4), every $\mcB_i$ is represented.
As noted above, $\lim_{n \to \infty}\sigma^k(\lambda_n) = \infty$.
We now show that if $\mcN$ is represented with respect to $\mbK$ and 
$(\sigma^k(\lambda_n), k)$-saturated, then $\mcN \models \varphi$.
Suppose that $\mcN$ has these properties.
It follows (from Definition~\ref{definition of (m,k)-saturation})
that, for every $i = 0, \ldots, l-1$, the $\mcB_{i+1}/\mcB_i$-multiplicity of 
$\mcN$ is at least $\sigma^k(\lambda_n)$ where $\sigma^k(\lambda_n) \geq 1$ for
all large enough $n$. So if $\mcB'_0 \cong \mcA$  and $\mcB'_0 \subseteq_{\cl} \mcN$,
then there are $\mcB'_i \subseteq_{\cl} \mcN$ such that
$\mcB'_i \cong \mcB_i$ and $\mcB'_{i-1} \subseteq \mcB'_i$ for $i = 1, \ldots, l$.
In particular, $\mcB'_0 \subseteq \mcB'_l \cong \mcB$
and since $\mcB'_0$ was an arbitrary closed copy of $\mcA$ in $\mcN$
it follows that $\mcN$ satisfies the $\mcB/\mcA$-extension axiom, i.e. $\mcN \models \varphi$.
Thus we have shown that in order to prove (1) it is sufficient to show that 
\begin{equation*}
\lim_{n \to \infty} \delta_n\big( \{\mcM \in \mbK_n : 
\mcM \text{ is $(\sigma^k(\lambda_n), k)$-saturated with respect to $\mbK$} \} \big) = 1. \tag{2}
\end{equation*}
For $n \in \mbbN$, let
\begin{align*}
\mbX_n &= \{\mcM \in \mbK_n : \mcM \text{ is $(\sigma^k(\lambda_n), k)$-saturated
with respect to $\mbK$} \},\\
&\text{and for $n,r \in \mbbN$ let}\\
\mbX_{n,r} &= \{\mcM \in \mbK_n \uhrc r : \text{ $\mcM$ is $(\sigma^r(\lambda_n), k)$-saturated
with respect to $\mbK \uhrc r$}\}.
\end{align*}
By Lemma~\ref{transfer of probabilities} below, in order to prove (2) it is sufficient to prove that 
\begin{equation*}
\lim_{n \to \infty} \mbbP_{n,k}( \mbX_{n,k} ) = 1, \tag{3}
\end{equation*}
\noindent

\begin{lem}\label{transfer of probabilities}
For every $n \in \mbbN$, $\delta_n(\mbX_n) = \mbbP_{n,\rho}(\mbX_n) = \mbbP_{n,k}(\mbX_{n,k})$.
\end{lem}

\noindent
For the proof of Lemma~\ref{transfer of probabilities} we need the following:

\begin{lem}\label{transfer of (i,k)-saturation}
Let $i \in \mbbN$.
For every $\mcM \in \mbK$, $\mcM$ is $(i,k)$-saturated with respect to $\mbK$ if and only if
$\mcM \uhrc k$ is $(i,k)$-saturated with respect to $\mbK \uhrc k$.
\end{lem}

\noindent
{\em Proof.}
Observe that for every $\mcM \in \mbK$ and every 
$A \subseteq M$ with $\dim_{\mcM}(A) \leq k$ the following holds:
for any relation symbol $R$, of arity $r$, say, and every $\bar{b} \in A^r$,
$$\bar{b} \in R^{\mcM} \Longleftrightarrow \bar{b} \in R^{\mcM \uhr k}.$$
In other words, $\mcM$ and $\mcM \uhrc k$ agree on all subsets $A$ of dimension at most $k$.
It follows, in particular, that for every $L$-structure $\mcA$ such that $\mcA \uhrc L_0 \in \mbG$
and $\mcA \uhrc L_0$ has dimension at most $k$, $\mcA$ is represented with respect to $\mbK$
if and only if $\mcA$ is represented with respect to $\mbK \uhrc k$.
The lemma is now an immediate consequence of the definition of $(i,k)$-saturation.
\hfill $\square$
\\

\noindent
{\em Proof of Lemma~\ref{transfer of probabilities}.}
Recall that $\rho$ is the supremum of the arities of relation symbols which
belong to the vocabulary of $L$ but not to the vocabulary of $L_0$.
First suppose that $\rho \leq k$.
Let 
$$\mbY_n = \{\mcN \in \mbK_n \uhrc k : \mcM \subseteq_w \mcN \text{ for some } \mcM \in \mbX_n\}.$$
By Lemma~\ref{probability in C_r can be computed as probabilty in C_r+1},
$\mbbP_{n,\rho}(\mbX_n) = \mbbP_{n,k}(\mbY_n)$.
But $\rho \leq k$ implies that, for every $\mcM \in \mbK$, 
$\mcM \uhrc k = \mcM \uhrc \rho = \mcM$. 
Hence,
$\mbX_{n,k} = \mbX_n = \mbY_n$, so 
$\delta_n(\mbX_n) = \mbbP_{n,\rho}(\mbX_n) = \mbbP_{n,k}(\mbX_{n,k})$.

Now suppose that $k < \rho$.
From Lemma~\ref{transfer of (i,k)-saturation} it follows that
$$\mbX_n = \{\mcN \in \mbK_n \uhrc \rho : \mcM \subseteq_w \mcN \text{ for some } \mcM \in \mbX_{n,k}\}$$
By Lemma~\ref{probability in C_r can be computed as probabilty in C_r+1},
$\mbbP_{n,k}(\mbX_{n,k}) = \mbbP_{n,\rho}(\mbX_n) = \delta_n(\mbX_n)$.
\hfill $\square$
\\

\noindent
Thus, it remains to prove (3), i.e. that $\lim_{n \to \infty}\mbbP_{n,k}(\mbX_{n,k}) = 1$.
This will be done by proving, by induction on $r$, that for every $r = 0, \ldots, k$,
$\lim_{n \to \infty}\mbbP_{n,r}(\mbX_{n,r}) = 1$.
In Definition~\ref{definition of substitutions} the notion of 
a substitution $\mcM[\mcA \trir \mcB]$ of $\mcA$ for $\mcB$ inside $\mcM$ was
defined. There it was assumed that the vocabulary of $L$ is relational.
However, eventual function or constant symbols in the vocabulary
of $L$ already belong to the vocabulary of $L_0 \subseteq L$, and, 
in what follows, we only consider substitutions when
$\mcA$ and $\mcB$ agree on $L_0$ and on proper closed substructures
(in the sense of Terminology~\ref{some terminology}).
So in this context, substitutions $\mcM[\mcA \trir \mcB]$, according to 
Definition~\ref{definition of substitutions}, make sense;
and we will use them.

\begin{lem}\label{changes of relations in one step}
Let $0 \leq r < k$, $\mcM \in \mbK_n \uhrc r+1$ and suppose that
$\mcA \subseteq_{cl} \mcM$ and
$\dim_{\mcM}(A) = r+1$. 
Also assume that $\mcB$ is a represented structure with 
respect to $\mbK \uhrc r+1$ such that
$\mcB$ and $\mcA$ agree on $L_0$ and on closed proper substructures.
Then $\mcM[\mcA \trir \mcB] \in \mbK_n \uhrc r+1$.
\end{lem}

\noindent
{\em Proof.}
Let $r$, $\mcM$, $\mcA$ and $\mcB$  satisfy the assumptions of the lemma,
so in particular $\mcA \uhrc L_0 = \mcB \uhrc L_0$.
Note that since $\mcA$ and $\mcB$ have dimension $r+1$ it follows
that $\mcA, \mcB \in \mbK$, because for every $\mcC \in \mbK$ with dimension at most $r+1$
we have $\mcC \uhrc r+1 = \mcC$.
By assumption, $\mcA$ and $\mcB$ agree on $L_0$ and on closed
proper substructures.
The assumption that $\mbK$ accepts $k$-substitutions over $L_0$ implies that
there exists $\mcN \in \mbK_n$ such that 
$\mcN \uhrc L_0 = \mcM \uhrc L_0$, $\mcN \uhrc B = \mcB$
and for every $\mcU \subseteq_{cl} \mcN$ such that $\dim_{\mcN}(U) \leq r+1$ and
$U \neq B$, we have $\mcN \uhrc U = \mcM \uhrc U$.
In particular, $\mcN \uhrc U = \mcM \uhrc U$ for every $U$ with dimension at most $r$.

Since $\mcN \uhrc r+1 \in \mbK_n \uhrc r+1$ it suffices to show that
$\mcM[\mcA \trir \mcB] = \mcN \uhrc r+1$.
For this it is enough to show that for every closed substructure 
$\mcC \subseteq_{cl} \mcM[\mcA \trir \mcB]$
with dimension $r+1$,
\begin{equation*}
\mcN \uhrc C = \mcC. \tag{$*$}
\end{equation*}
Suppose that $\mcC \subseteq_{cl} \mcM[\mcA \trir \mcB]$.
If $\mcC = \mcB$ then, by the choice of $\mcN$, we have $\mcN \uhrc C = \mcN \uhrc B = \mcB$.
If $\mcC \neq \mcB$ then, by the choice of $\mcN$, we have 
$\mcN \uhrc C = \mcM \uhrc C  = \mcC$, 
where the last identity follows because $\mcM = \mcM \uhrc r+1$ and $\mcC$ has dimension $r+1$;
thus ($*$) also holds in case when $\mcC \neq \mcB$.
\hfill $\square$

\begin{lem}\label{changes of relations}
Let $0 \leq r < k$, $\mcM \in \mbK_n \uhrc r+1$ and suppose that
$\mcA \subseteq_{cl} \mcM$ and
$r < \dim_{\mcM}(A) \leq k$. 
Also assume that $\mcB$ is a represented structure with 
respect to $\mbK \uhrc r+1$ such that
$\mcB \uhrc L_0 = \mcA \uhrc L_0$ and for every closed $U \subseteq A = B$
with dimension $r$, $\mcA \uhrc U = \mcB \uhrc U$.
Then $\mcM[\mcA \trir \mcB] \in \mbK_n \uhrc r+1$.
\end{lem}

\noindent
{\em Proof.}
Let $r$, $\mcM$, $\mcA$ and $\mcB$  satisfy the assumptions of the lemma.
By definition of $\mbK \uhrc r+1$, for every $\mcN \in \mbK \uhrc r+1$
and every relation symbol $R$ which does
not belong to the vocabulary of $L_0$, there is {\em no}
$R$-relationship $\bar{a} \in R^{\mcN}$ with dimension greater than $r+1$.
Consequently, the structure $\mcM[\mcA \trir \mcB]$ can be created
by a finite number of substitutions of the kind considered
in Lemma~\ref{changes of relations in one step}.
More precisely: 
There are $\mcN_0, \ldots, \mcN_s \in \mbK_n \uhrc r+1$
and $\mcC_0, \ldots, \mcC_{2s}$ which dimension $r+1$
such that 
\begin{align*}
&\mcM = \mcN_0, \ \mcM[\mcA \trir \mcB] = \mcN_s,\\
&\mcN_{i+1} = \mcN_i[\mcC_{2i} \trir \mcC_{2i+1}], \text{ for } i = 1, \ldots, s, \text{ and}\\
&\text{$\mcC_{2i}$ and $\mcC_{2i+1}$ agree on $L_0$ and on closed proper substructures}.
\end{align*}
By Lemma~\ref{changes of relations in one step}, $\mcN_i \in \mbK_n \uhrc r+1$,
for $i = 0, \ldots, s$, so we are done.
\hfill $\square$

\begin{lem}\label{expansions}
If $0 \leq r < k$ then for
every $\mcM \in \mbK_n \uhrc r$ there is  $\mcM' \in \mbK_n \uhrc r+1$
such that $\mcM' \uhrc r = \mcM$.
\end{lem}

\noindent
{\em Proof.}
If $\mcM \in \mbK_n \uhrc r$ then $\mcM = \mcN \uhrc r$ for some $\mcN \in \mbK_n$.
Take $\mcM' = \mcN \uhrc r+1$. 
Then $\mcM' \in \mbK_n \uhrc r+1$ and $\mcM' \uhrc r = \mcN \uhrc r = \mcM$.
\hfill $\square$

\begin{lem}\label{base case}
For every $n$ and every $\mcM \in \mbK_n \uhrc 0$,
$\mcM$ is $\big(\lambda_n , k \big)$-saturated.
\end{lem}

\noindent
{\em Proof.}
First observe that from 
Definition~\ref{definition of dimension reduct} 
it follows that whenever $\mcM$ is permitted (or, equivalently, in the present
context, represented) with respect to $\mbK \uhrc 0$,
then $\mcM$ is an expansion of $\mcM \uhrc L_0$ ($\cong \mcG_n$ for some $n$) 
obtained by possibly adding
some new relationship(s) involving {\em only} elements in 
$\cl_{\mcM}(\es)$; and whenever $\mcA \subseteq_{\cl} \mcM$ then $\cl_{\mcM}(\es) \subseteq A$

Let  $\mcM \in \mbK_n \uhrc 0$ and let
$\mcA \subseteq_{cl} \mcB$ be permitted structures with respect to $\mbK \uhrc 0$
such that $\dim_{\mcB}(A) + 1 = \dim_{\mcB}(B) \leq k$.
Suppose that $\mcA' \subseteq_{cl} \mcM$ is a copy of $\mcA$ and that 
$\tau : \mcA' \to \mcA$ is an isomorphism.
We must show that there are $\mcB'_i \subseteq_{cl} \mcM$ and isomorphisms 
$\tau_i : \mcB'_i \to \mcB$, for $i = 1, \ldots, \lambda_n$,
such that $\mcA' \subseteq_{cl} \mcB'_i$, 
$\tau_i \uhrc A' = \tau$ and $B'_i \cap B'_j = A'$ whenever $i \neq j$.
As noted in the beginning of the proof, every relationship of $\mcB$ (or of $\mcM$)
which involves some element(s) from $B-A$ (or from $M-A'$) is an $R$-relationship for
some relation symbol $R$ of $L_0$.
Observe that $\tau : A' \to A$ can also be viewed as an isomorphism from
$\mcA' \uhrc L_0$ to $\mcA \uhrc L_0$.
By (b) in the beginning of the proof of Theorem~\ref{0-1 law for pregeometries},
there are $\mcB_i \subseteq_{cl} \mcM \uhrc L_0 = \mcG_n$ 
and isomorphisms $\tau_i : \mcB_i \to \mcB \uhrc L_0$, for $i = 1, \ldots, \lambda_n$,
such that $\mcA' \uhrc L_0 \subseteq_{cl} \mcB_i$, 
$\tau_i \uhrc A' = \tau$ and $B_i \cap B_j = A'$ whenever $i \neq j$.
For $i = 1, \ldots, \lambda_n$, 
let $\mcB'_i \subseteq_{\cl} \mcM$ be such that $\mcB'_i \uhrc L_0 = \mcB_i$.
Then $\mcA' \subseteq_{\cl} \mcB'_i$ for each $i$, and
since, as observed above, every relationship which involves some element(s) from
$M-A'$, or from $B-A$, is an $R$-relationship for some relation symbol $R$ of $L_0$,
it follows that every $\tau_i$ is in fact an isomorphism from $\mcB'_i$ to $\mcB$.
\hfill $\square$

\begin{lem}\label{induction step for pregeometries}
Suppose that $0 \leq r < k$.
For every real $\varepsilon > 0$ there is $n_{\varepsilon} \in \mbbN$ such that
if $n \geq n_{\varepsilon}$, $\mcM \in \mbK_n \uhrc r$ is 
$\big(\sigma^r(\lambda_n), k \big)$-saturated and
$$\mbE_{r+1}(\mcM) = \big\{\mcN \in \mbK_n \uhrc r+1 : \mcN \uhrc r = \mcM \big\},$$
then the proportion of $\mcN \in \mbE_{r+1}(\mcM)$ which are
$\big(\sigma^{r+1}(\lambda_n), k \big)$-saturated
is at least $1 - \varepsilon$.
\end{lem}

\noindent
{\em Proof.}
Let $0 \leq r < k$.
We are assuming that $\{\mcG_n : n \in \mbbN\}$ is a uniformly bounded
pregeometry. Hence there is $\alpha \in \mbbN$ such that if
$\mcA$ is permitted with respect to $\mbK \uhrc r+1$ and has 
dimension at most $k$, then $|A| \leq \alpha$.
Suppose that $\mcM \in \mbK_n \uhrc r$ is 
$\big(\sigma^r(\lambda_n), k \big)$-saturated and let 
$\mbE_{r+1}(\mcM) = \big\{\mcN \in \mbK_n \uhrc r+1 : \mcN \uhrc r = \mcM \big\}$. 
We start by proving that, with the uniform probability measure on $\mbE_{r+1}(\mcM)$,
the probability that a randomly chosen $\mcN \in \mbE_{r+1}(\mcM)$
is $\big(\sigma^{r+1}(\lambda_n), k \big)$-saturated approaches 1 as $n$ tends to $\infty$.
We do this by finding an upper bound (depending on $n$) for the probability that
a randomly chosen
$\mcN  \in \mbE_{r+1}(\mcM)$ is {\em not}
$\big(\sigma^{r+1}(\lambda_n), k \big)$-saturated; 
and then observe that this upper bound approaches 0 as $n$ tends to infinity.
Finally we note that the argument does not depend
on which $\big(\sigma^r(\lambda_n), k \big)$-saturated $\mcM \in \mbK_n \uhrc r$ we consider;
so given $\varepsilon > 0$ there is $n_{\varepsilon}$ which such that for every $n \geq n_{\varepsilon}$
and every $\big(\sigma^r(\lambda_n), k \big)$-saturated $\mcM \in \mbK_n \uhrc r$,
the proportion of  $\mcN \in \mbE_{r+1}(\mcM)$ which are {\em not}
$\big(\sigma^{r+1}(\lambda_n), k \big)$-saturated is at most $\varepsilon$.

Let $\mcN \in \mbE_{r+1}(\mcM)$ and let  
$\mcA \subset_{cl} \mcB$ be represented structures with respect to $\mbK \uhrc r+1$
such that $\dim_{\mcB}(A) + 1 = \dim_{\mcB}(B) \leq k$.
Suppose that $\mcA' \subseteq_{cl} \mcN$ is a copy of $\mcA$ and that 
$\tau : \mcA' \to \mcA$ is an isomorphism.
Let $l_n =  \lfloor\sqrt{\sigma^r(\lambda_n)}\rfloor  = \sigma^{r+1}(\lambda_n)$.
First we find an upper bound for the probability that there does not exist
$\mcB_i \subseteq_{cl} \mcN$ and isomorphisms 
$\tau_i : \mcB_i \to \mcB$, for $i = 1, \ldots, l_n$, such that
$\mcA' \subseteq_{cl} \mcB_i$, $\tau_i \uhrc A' = \tau$, and $B_i \cap B_j = A'$ whenever $i \neq j$.

Let $l'_n = \sigma^r(\lambda_n)$.
Since $\mcM$ is $\big(\sigma^r(\lambda_n), k \big)$-saturated 
there are $\mcB^-_i \subseteq_{cl} \mcM$, $i = 1, \ldots, l'_n$
and isomorphisms $\tau_i : \mcB^-_i \to \mcB \uhrc r$, such
that $\mcA' \uhrc r \subseteq_{cl} \mcB^-_i$,
$\tau_i \uhrc A' = \tau$ and 
$B^-_i \cap B^-_j = A'$ whenever $i \neq j$.
Let $\beta$ be the number of represented structures
with respect $\mbK \uhrc r+1$ with universe included in $\{1, \ldots, \alpha\}$.
Lemma~\ref{changes of relations}
implies that
the probability that the map $\tau_i : B^-_i \to B$
is an isomorphism from $\mcN \uhrc B^-_i$ to $\mcB$ is at least $1/\beta$, 
independently of whether this holds for $j \neq i$.
Let $s$ be a natural number such that $0 \leq s < l_n$.
The probability that for {\em every} $i \in \{sl_n + i, \ldots, (s+1)l_n \}$,
$\tau_i : B^-_i \to B$
is {\em not} an isomorphism from $\mcN \uhrc B^-_i$ to $\mcB$
is at most 
$$\big( 1 - 1/\beta \big)^{l_n}.$$
Let $m_n = |G_n| = |N|$. 
By (a) $\lambda_n \leq m_n \leq P(\lambda_n)$ for all $n \in \mbbN$, where
$P$ is a polynomial.
Since, by assumption, $\lim_{n\to\infty}\lambda_n = \infty$,
we have $\lim_{n\to\infty}m_n = \infty$.
From the definition of $l_n$ as $l_n = \sigma^{r+1}(\lambda_n)$ and the definition of $\sigma^{r+1}$
it follows that there is a polynomial $Q$ such that 
$m_n \leq Q(l_n)$.
The number of ways in which we can choose $\mcA$, $\mcB$, $\mcA'$ and $s$ as above
is not larger than
$$\beta^2 \cdot (m_n)^{\alpha} \cdot l_n \quad \leq \quad 
\beta^2 \cdot (Q(l_n))^{\alpha} \cdot l_n.$$
Moreover, for every choice of such $\mcA$, $\mcB$, $\mcA'$ and $s$, there exist,
for $i = 1, \ldots, l'_n$,
$\mcB^-_i \subseteq_{cl} \mcM$ and
isomorphisms $\tau_i : \mcB^-_i \to \mcB$, with the properties described above.
So if $\mcN$ is not $\big(\sigma^{r+1}(\lambda_n), k \big)$-saturated,
then there exist $\mcA$, $\mcB$, $\mcA'$, $\mcB^-_i$, $\tau_i$, for $i = 1, \ldots, l'_n$, and $s$
as above such that for every $i \in \{sl_n + 1, \ldots, (s+1)l_n \}$,
$\tau_i$ is not an isomorphism from $\mcN \uhrc B^-_i$ to $\mcB$.
Hence, the probability that a randomly chosen $\mcN \in \mbE_{r+1}(\mcM)$ is not 
$\big(\sigma^{r+1}(\lambda_n), k \big)$-saturated 
does {\em not} exceed
$$f_n = \beta^2 \cdot (Q(l_n))^{\alpha} \cdot l_n \cdot \big( 1 - 1/\beta \big)^{l_n}.$$
Since $\lambda_n \to \infty$ as $n \to \infty$ we also have $l_n \to \infty$ as $n \to \infty$.
Because $\beta^2 \cdot (Q(l_n))^{\alpha} \cdot l_n$ is a polynomial in $l_n$ it follows that
$f_n \to 0$ as $n \to \infty$.

Observe that the same expression for $f_n$ works for every 
$\big(\sigma^r(\lambda_n), k \big)$-saturated $\mcM \in \mbK_n \uhrc r$.
So for every $\varepsilon > 0$ there is $n_{\varepsilon}$ such that for every $n \geq n_{\varepsilon}$
and every $\big(\sigma^r(\lambda_n), k \big)$-saturated $\mcM \in \mbK_n \uhrc r$,
the proportion of $\mcN \in \mbE_{r+1}(\mcM)$ which are $\big(\sigma^{r+1}(\lambda_n), k \big)$-saturated
is at least $1 - \varepsilon$.
\hfill $\square$
\\

\noindent
Recall that, for $r = 0, 1, \ldots, k$,
$$\mbX_{n,r} =  \{\mcM \in \mbK_n \uhrc r : \mcM \text{ is $(\sigma^r(\lambda_n), k)$-saturated} \}.$$
From Lemma~\ref{induction step for pregeometries} we can easily derive the following:

\begin{lem}\label{inclusion of the Xes}
For every $r = 0, 1, \ldots, k-1$ and all sufficiently large $n$
(take $0 < \varepsilon < 1/2$, $n_{\varepsilon}$ and $n > n_{\varepsilon}$ 
so that the conclusion of Lemma~\ref{induction step for pregeometries} holds),
$$\mbX_{n,r} \subseteq \{\mcN \uhrc r : \mcN \in \mbX_{n,r+1}\}.$$
\end{lem}

\noindent
{\em Proof.}
Suppose that $\mcM \in \mbX_{n,r}$, so $\mcM$ is $\big(\sigma^r(\lambda_n), k \big)$-saturated.
By Lemma~\ref{induction step for pregeometries}, for all sufficiently large $n$,
$\mbE_{r+1}(\mcM)$ will contain a stucture $\mcN$ which is 
$\big(\sigma^{r+1}(\lambda_n), k \big)$-saturated; hence $\mcN \in \mbX_{n,r+1}$
and $\mcN \uhrc r = \mcM$. 
\hfill $\square$
\\

\noindent
Now we can finish the proof of part (i) of Theorem~\ref{0-1 law for pregeometries} by proving (3),
in other words, that $\lim_{n \to \infty}\mbbP_{n,k}(\mbX_{n,k}) = 1$.
Let $\varepsilon > 0$.
Choose $\varepsilon' > 0$ so that $(1 - \varepsilon')^k \geq 1 - \varepsilon$.
By Lemma~\ref{induction step for pregeometries}, we can choose $n_{\varepsilon'}$ such that 
if $0 \leq r < k$, $n > n_{\varepsilon'}$ and $\mcM \in \mbK_n \uhrc r$ is 
$\big(\sigma^r(\lambda_n), k \big)$-saturated, 
then the proportion of $\mcN \in \mbE_{r+1}(\mcM)$ which are
$\big(\sigma^{r+1}(\lambda_n), k \big)$-saturated
is at least $1 - \varepsilon'$.
By induction we show that, for
$r = 0,1, \ldots, k$ and $n > n_{\varepsilon'}$, 
$$\mbbP_{n,r}(\mbX_{n,r}) \geq (1 - \varepsilon')^r \geq 1 - \varepsilon.$$
The base case $r = 0$ is given by Lemma~\ref{base case}, so assume that $0 < r \leq k$ and
that $\mbbP_{n, r-1}(\mbX_{n,r-1}) \geq (1 - \varepsilon')^{r-1}$.
Let $\mcM_1, \ldots, \mcM_s$ be an enumeration, without repetition, of $\mbX_{n,r}$.
Then let $\mcM'_1, \ldots, \mcM'_t$ be an enumeration, without repetition, of the
set $\{\mcM_1 \uhrc r-1, \ldots, \mcM_s \uhrc r-1\}$.
By the definition of $\mbbP_{n,r}$, the following holds for every $n > n_{\varepsilon'}$:
\begin{align*}
\mbbP_{n,r}(\mbX_{n,r}) &= \mbbP_{n,r}\big(\{\mcM_1, \ldots, \mcM_s\}\big) = 
\sum_{i = 1}^s \mbbP_{n,r}(\mcM_i)\\  
&= \sum_{i = 1}^s \frac{1}{\big|\{\mcN \in \mbK_n \uhrc r : \mcN \uhrc r-1 = \mcM_i \uhrc r-1\}\big|} \cdot
\mbbP_{n,r-1}(\mcM_i \uhrc r-1)\\
&= \sum_{i = 1}^t \frac{\big|\{\mcN \in \mbX_{n,r} : \mcN \uhrc r-1 = \mcM'_i\}\big|}
{\big|\{\mcN \in \mbK_n \uhrc r : \mcN \uhrc r-1 = \mcM'_i\}\big|} \cdot \mbbP_{n,r-1}(\mcM'_i)\\
&= \sum_{i = 1}^t \frac{\big|\{\mcN \in \mbX_{n,r} : \mcN \uhrc r-1 = \mcM'_i\}\big|}
{\big|\mbE_r(\mcM'_i)\big|} \cdot \mbbP_{n,r-1}(\mcM'_i)\\
&\geq (1 - \varepsilon') \sum_{i = 1}^t \mbbP_{n,r-1}(\mcM'_i) 
\quad \quad \text{ (by the choice of $n_{\varepsilon'}$)}\\
&= (1 - \varepsilon')\mbbP_{n,r-1}\big(\{\mcM'_1, \ldots, \mcM'_t\}\big)\\
&\geq (1 - \varepsilon')\mbbP_{n,r-1}(\mbX_{n,r-1}) 
\quad \quad \text{ (by Lemma~\ref{inclusion of the Xes})}\\
&\geq (1 - \varepsilon')(1 - \varepsilon')^{r-1} = (1 - \varepsilon')^r 
\quad \quad \text{ (by the induction hypothesis).}
\end{align*}
Thus (3) is proved, and hence also part (i) of Theorem~\ref{0-1 law for pregeometries}.

Now we prove part (ii) of Theorem~\ref{0-1 law for pregeometries}.
Note that we have proved (2) above, because (3) together with 
Lemma~\ref{transfer of probabilities} implies (2).
By (2), there are, for all $n \in \mbbN$, $\mcM_n \in \mbK_n$ such
that $\mcM_n$ is $(\sigma^k(\lambda_n), k)$-saturated.
Let $\mu_n = \sigma^k(\lambda_n)$, so $\mcM_n$ is $(\mu_n, k)$-saturated,
where $\lim_{n\to\infty}\mu_n = \infty$.
From (a) and the definition of $\sigma^k$ it follows that there is a polynomial $Q$ such that 
$\mu_n \leq |M_n| \leq Q(\mu_n)$ for all $n$.
Since $\mcM_n$ is $(\mu_n, k)$-saturated, the following holds:
If $\mcA \subset_{cl} \mcB$ are represented structures
such that $\dim_{\mcB}(B) \leq k$,
then the $\mcB / \mcA$-multiplicity of $\mcM_n$ is at least $\mu_n$.
From Assumption~\ref{assumptions about pregeometries} (4),
it follows that the sequence $(\mcM_n : n \in \mbbN)$
is polynomially $k$-saturated; and hence $\mbK$ is polynomially $k$-saturated.
This concludes the proof of part (ii), and hence of Theorem~\ref{0-1 law for pregeometries}.

\subsection{Proof of Theorem~\ref{third part of 0-1 law for pregeometries}}
\label{proof of third part of 0-1 law for pregeometries}

We still assume that, for every $k > 0$, $(\mcG_n : n \in \mbbN)$ is polynomially $k$-saturated and
$\mbK = \bigcup_{n \in \mbbN}\mbK_n$, where $\mbK_n = \mbK(\mcG_n)$ ,
accepts $k$-substitutions over $L_0$.
We want to prove that for every $L$-sentence $\varphi$, either
$\lim_{n\to\infty}\delta_n(\varphi) = 0$ or $\lim_{n\to\infty}\delta_n(\varphi) = 1$.
The general idea of the proof follows a well-known pattern: 
we collect into a theory $T_{\mbK}$ all extension axioms of $\mbK$ together 
with sentences which express the pregeometry conditions and 
describe the possible isomorphism types of 
closed substructures of members of $\mbK$.
By part (i) of Theorem~\ref{0-1 law for pregeometries}, $T_{\mbK}$ is consistent.
Then we show that $T_{\mbK}$ is complete by showing that it is countably categorical.
From the completeness, it follows that for every $L$-sentence $\varphi$, either
$T_{\mbK} \models \varphi$ or $T_{\mbK} \models \neg\varphi$. In the first case there is
finite $\Delta \subset T_{\mbK}$ such that $\Delta \models \varphi$ and in the second
case there is finite $\Delta' \subseteq T_{\mbK}$ such that $\Delta' \models \neg\varphi$.
In the first case part (i) of Theorem~\ref{0-1 law for pregeometries} implies that 
$$\lim_{n\to\infty}\delta_n(\{\mcM \in \mbK_n : \mcM \models \Delta\}) = 1,$$
and therefore $\lim_{n\to\infty}\delta_n(\varphi) = 1$.
In the second case we get, in a similar way, that $\lim_{n\to\infty}\delta_n(\neg\varphi) = 1$,
so $\lim_{n\to\infty}\delta_n(\varphi) = 0$.

Now to the details.
We are assuming that $\mbG = \{\mcG_n : n \in \mbbN\}$ is a pregeometry
where the closure operator of every member of $\mbG$ is defined 
by the $L_0$-formulas $\theta_n(x_1, \ldots, x_{n+1})$, $n \in \mbbN$,
according to Definition~\ref{definition of a structure being a pregeometry} and
Assumption~\ref{assumptions about pregeometries}.
In other words, for all $m, n$ and all 
$a_1, \ldots, a_{n+1} \subseteq G_m$, 
\begin{equation*}
a_{n+1} \in \cl_{\mcG_m}(a_1, \ldots, a_n) \text{ if and only if }
\mcG_m \models \theta(a_1, \ldots, a_{n+1}). \tag{4}
\end{equation*}
Also (by Assumption~\ref{assumptions about pregeometries}), 
for every $m$ and every $\mcM \in \mbK_m = \mbK(\mcG_m)$, 
$\cl_{\mcM}$ coincides with $\cl_{\mcG_m}$.
Moreover, the pregeometry $\mbG$ is assumed to be {\em uniformly} locally finite,
so there is $u : \mbbN \to \mbbN$ such that for every $\mcM \in \mbK$ and
every $X \subseteq M$, $|\cl_{\mcM}(X)| \leq u(\dim_{\mcM}(X))$.
We may also assume that for every $k \in \mbbN$ the value $u(k)$ is minimal so that
this holds.

By the {\bf \em finiteness property}, for a pregeometry $(A, \cl)$, we mean
the property that for all $a \in A$ and $X \subseteq A$, $a \in \cl(X)$ if and
only if $a \in \cl(Y)$ for some {\em finite} $Y \subseteq X$.
Besides the finiteness property, all other properties of a pregeometry can, when (4) holds,
be expressed for finite subsets of $A$ by using the formulas 
$\theta_n(x_1, \ldots, x_{n+1})$, $n \in \mbbN$.
Let $T_{preg}$ be the set of sentences which express all properties of a pregeometry
(for finite subsets) except the finiteness property.
Then every $\mcM \in \mbK$ is a model of $T_{preg}$.

Note that, for every $\mcM \in \mbK$ and all $a_1, \ldots, a_n \in M$,
the statement ``$\{a_1, \ldots, a_n\}$ is a closed set (in $\mcM$)'' is 
uniformly expressed by the first-order formula
$$\neg\exists x_{n+1} \Big( \bigwedge_{i=1}^n x_{n+1} \neq x_i \ \wedge \
\theta_n(x_1, \ldots, x_{n+1})\Big),$$
which we denote by $\gamma_n(x_1, \ldots, x_n)$. 
For every positive $m \in \mbbN$, let $s(m)$ be the
the number of nonisomorphic structures of cardinality at most $m$
which occur as closed substructures of members of $\mbK$, and let
$\mcM_{m,1}, \ldots, \mcM_{m,s(m)}$ be an enumeration of all isomorphism types
of such structures.
For $1 \leq i \leq s(m)$, let $\chi_{m,i}(x_1, \ldots, x_m)$ 
describe the isomorphism type of $\mcM_{m,i}$ in such a way that
we require that all variables $x_1, \ldots, x_m$ actually occur in $\chi_{m,i}$.
It means that if $\left\| \mcM_{m,i} \right\| < m$, then 
$\chi_{m,i}(x_1, \ldots, x_m)$ must express that
some variables refer to the same element, by saying `$x_k = x_l$' for some $k \neq l$.
For every $k \in \mbbN$ let $\psi_k$ denote the sentence
$$\forall x_1, \ldots, x_k \exists x_{k+1}, \ldots, x_{u(k)} 
\Big( \gamma_{u(k)}(x_1, \ldots, x_{u(k)}) \ \wedge \ \bigvee_{i=1}^{s(u(k))}
\bigvee_{\pi} \chi_{u(k),i}(x_{\pi(1)}, \ldots, x_{\pi(u(k))}) \Big),$$
where the second disjunction ranges over all permutations $\pi$ of $\{1, \ldots, u(k)\}$.
If $k = 0$ and $u(k) > 0$, then the universal quantifiers do not occur so 
$\psi_0$ is an existential formula.
If $u(0) = 0$, then, by convention, $\psi_0$ is $\forall x (x = x)$.
If $u(k) = k$, then the existential quantifiers do not occur and
$\psi_k$ is a universal formula.
Note that for every $k \in \mbbN$ and every $\mcM \in \mbK$, $\mcM \models \psi_k$.
Let $T_{iso} = \{\psi_k : k \in \mbbN\}$ so every $\mcM \in \mbK$ is a model of $T_{iso}$.

Finally, let $T_{ext}$ consist (exactly) of all extension axioms of $\mbK$ and let
$$T_{\mbK} = T_{preg} \cup T_{iso} \cup T_{ext}.$$
By Theorem~\ref{0-1 law for pregeometries}
and compactness, $T_{\mbK}$ is consistent.
Note that every model of $T_{\mbK}$ is infinite, because we assume that
$(\mcG_n : n \in \mbbN)$ is polynomially $k$-saturated (for every $k > 0$),
which implies that for some sequence $(\lambda_n : n \in \mbbN)$ which tends to infinity
as $n \to \infty$, $\mcG_n$ contains at least $\lambda_n$ different elements. 

\begin{lem}\label{the almost sure theory is a pregeometry}
Suppose that $\mcM \models T_{\mbK}$ and define $\cl_{\mcM}$ as follows:
\begin{itemize}
	\item[(a)] for all $n \in \mbbN$ and all $a_1, \ldots, a_{n+1} \in M$, 
	$a_{n+1} \in \cl_{\mcM}(a_1, \ldots, a_n)$ $\Longleftrightarrow$ 
	$\mcM \models \theta_n(a_1, \ldots, a_{n+1})$.
	\item[(b)] for all $X \subseteq M$ and all $a \in M$, $a \in \cl_{\mcM}(X)$
	$\Longleftrightarrow$ for some finite $Y \subseteq X$, $a \in \cl_{\mcM}(Y)$.
\end{itemize}
Then $(M, \cl_{\mcM})$ is a pregeometry such that for every finite
$X \subseteq M$, $|\cl_{\mcM}(X)| \leq u(\dim_{\mcM}(X))$.
\end{lem}

\noindent
{\em Proof.}
Suppose that $\mcM \models T_{\mbK}$.
Since $T_{preg} \subseteq T_{\mbK}$, it follows from part (a) 
that $\cl_{\mcM}$ satisfies all properties of a pregeometry on finite subsets of $M$.
But (b) guarantees that $\cl_{\mcM}$ has the finiteness property, and then
all other properties follow for all subsets of $M$.
So $(M, \cl_{\mcM})$ is a pregeometry.
Since $T_{iso} \subset T_{\mbK}$ it follows that, for every
$X \subseteq M$, $|\cl_{\mcM}(X)| \leq u(\dim_{\mcM}(X))$.
\hfill $\square$ 
\\

\noindent
To complete the proof of Theorem~\ref{third part of 0-1 law for pregeometries}
we only need to prove:

\begin{lem}\label{countable categoricity for pregeometries}
$T_{\mbK}$ is countably categorical and hence complete.
\end{lem}

\noindent
{\em Proof.}
Let $\mcM$ and $\mcN$ be countable models of $T_{\mbK}$.
We show that $\mcM \cong \mcN$, by a back-and-forth argument.
By symmetry it is sufficient to show the following:
\begin{itemize}
	\item[] Suppose that $\mcA$ is a closed finite substructure of $\mcM$ (or $A = \es$), that $\mcB$
	is a closed finite substructure of $\mcN$ (or $B = \es$), that
	$f : \mcA \to \mcB$ is an isomorphism (if $A$ and $B$ are nonempty) and that $a \in M - A$. 
	Then there are a closed $\mcB' \subseteq M$ such that $B \subset B'$ and an isomorphism 
	$g : \cl_{\mcM}(A \cup \{a\}) \to \mcB'$ which extends $f$.
\end{itemize}
So suppose that $\mcA$ is a closed finite substructure of $\mcM$, that $\mcB$
is a closed finite substructure of $\mcN$, that
$f : \mcA \to \mcB$ is an isomorphism and that $a \in M - A$.
Since $\mcM \models T_{\mbK} \supset T_{iso}$, $\mcA$, $\mcB$ and $\cl_{\mcM}(A \cup \{a\})$ are
isomorphic with closed substructures of members of $\mbK$.
Since $\mcN \models T \supset T_{ext}$, it follows that $\mcN$ satisfies
the $\cl_{\mcM}(A \cup \{a\}) / \mcA$-extension axiom, and as $\mcB \cong \mcA$ there
is a closed $\mcB' \subset \mcN$ such that $\mcB \subset \mcB'$ and an isomorphism
$g : \cl_{\mcM}(A \cup \{a\}) \to \mcB'$ which extends $f$.
Recall the convention that for every
structure $\mcP$ which is isomorphic with a closed substructure of a member of $\mbK$, 
the statement ``there exists a closed copy of $\mcP$'' is 
an extension axiom, called the $\mcP/\es$-extension axiom;
this takes care of the case $A = B = \es$.
\hfill $\square$

\subsection{Proof of Theorem~\ref{0-law for extension axioms in pregeometries}}

Let $\mbG = \{\mcG_n : n \in \mbbN\}$ be a set of $L_0$-structures which
form a uniformly bounded pregeometry, and suppose that
$(\mcG_n : n \in \mbbN)$ is polynomially $k$-saturated for every $k \in \mbbN$.
Assume that there is, up to isomorphism, a unique represented structure with dimension 0;
hence $\mbK$ accepts 0-substitutions over $L_0$.
Suppose that $k$ is {\em minimal} such that 
$\mbK$ does {\em not} accept $k$-substitutions over $L_0$; hence $k > 0$ and
$\mbK$ accepts $(k-1)$-substitutions over $L_0$.
Moreover assume that there are 
represented structures, with respect to $\mbK$, $\mcA$ and $\mcA'$ such that
\begin{itemize}
\item[\textbullet] \ $\mcA$ and $\mcA'$ have dimension $k$, 
\item[\textbullet] \ $\mcA$ and $\mcA'$ agree on $L_0$ and on closed proper substructures,
\item[\textbullet] $\mbK$ accepts the substitution $[\mcA' \trir \mcA]$ over $L_0$, but
\item[\textbullet] $\mbK$ does {\em not} accept the substitution $[\mcA \trir \mcA']$ over $L_0$.
\end{itemize}
Let $\rho$ be the supremum of the arities of all relation symbols which
belong to the vocabulary of $L$ but not to the vocabulary of $L_0$.
By Remark~\ref{remark on supremum of arities and acceptance},
$0 < k \leq \rho$.

In order to prove Theorem~\ref{0-law for extension axioms in pregeometries},
we assume that $\mbK$ has the independent amalgamation property and show
that there are extension axioms $\varphi$ and $\psi$ such that
$\lim_{n \to \infty} \delta_n(\varphi \wedge \psi) = 0$. 
We start with the following, which is straightforward to verify:

\begin{observation}\label{observation about dimension reducts}
For every $L$-structure $\mcM$ and $d \in \mbbN$, $\mcM$ is
represented with respect to $\mbK \uhrc d$ if and only if
there is $\mcM'$ such that $\mcM'$ is represented with respect to $\mbK$
and $\mcM = \mcM' \uhrc d$. 
\end{observation}

\noindent
Note that the notion of `acceptance of $l$-substitutions over $L_0$', 
which was defined for $\mbK$, can equally well be defined for $\mbK \uhrc r$ for
any $r$; the only difference is that the notion
`represented' is in this case with respect to $\mbK \uhrc r$.
By assumption, $\mbK$ accepts $(k-1)$-substitutions over $L_0$.
From Observation~\ref{observation about dimension reducts} it follows
that $\mbK \uhrc (k-1)$ accepts $(k-1)$-substitutions over $L_0$.
Note that for every $\mcM \in \mbK \uhrc k-1$ and every 
relation symbol $R$ in the vocabulary of $L$ but not in the vocabulary
of $L_0$, $\mcM$ does not have any $R$-relationship with dimension greater than $k-1$.
From this and the assumption that $\mbK$ accepts $(k-1)$-substitutions
over $L_0$ it follows that 
\begin{itemize}
	\item[(5)] $\mbK \uhrc (k-1)$ accepts $l$-substitutions over $L_0$ for {\em every} $l \in \mbbN$.
\end{itemize}
By assumption, $\mbK$ accepts the substitution $[\mcA' \trir \mcA]$,
and by Observation~\ref{observation about dimension reducts} it follows
that $\mbK \uhrc k$ accepts the substitution $[\mcA' \trir \mcA]$.

Since $\mcA$ and $\mcA'$ agree on $L_0$ it makes sense to speak about the 
substitution $\mcM[\mcA \trir \mcA']$ if $\mcA \subseteq_{cl} \mcM$,
or $\mcM[\mcA' \trir \mcA]$ if $\mcA' \subseteq_{cl} \mcM$, as was explained
in the paragraph before Lemma~\ref{changes of relations in one step}.
Since $\mcA$ and $\mcA'$ have dimension $k$ and $\mbK$ accepts the substitution
$[\mcA' \trir \mcA]$ over $L_0$, it follows that if $\mcM$ is represented
with respect to $\mbK \uhrc k$, and  $\mcA' \subseteq_{cl} \mcM$, then
$\mcM[\mcA' \trir \mcA]$ is represented
with respect to $\mbK \uhrc k$.
In other words, $\mbK \uhrc k$ {\em admits} the substitution $[\mcA' \trir \mcA]$.
By assumption, $\mbK$ does not accept the substitution $[\mcA \trir \mcA']$.
Therefore we can argue similarly as we just did for the substitution $[\mcA \trir \mcA']$
to conclude that there is $\mcP$ such that $\mcP$ is represented with respect to $\mbK \uhrc k$,
$\mcA \subset_{cl} \mcP$ and
$\mcP[\mcA \trir \mcA']$ is forbidden with respect to $\mbK \uhrc k$.

Since the core of the argument 
(the proof of Lemma~\ref{almost neverness for the k reduct} below)
is an adaptation
of the proof of Theorem~\ref{upper limit on the probability of an extension axiom} 
to the present context,
we introduce the same notation as in 
Section~\ref{proving upper limit on the probability of an extension axiom}.
We rename $\mcA$ and $\mcA'$ with $\mcS_{\mcP}$ and $\mcS_{\mcF}$,
so in particular $\mcS_{\mcP}$ and $\mcS_{\mcF}$ have dimension $k$.
As concluded above, $\mbK \uhrc k$ {\em admits} the substitution $[\mcS_{\mcF} \trir \mcS_{\mcP}]$,
in the sense that whenever $\mcM$ is represented with respect to $\mbK \uhrc k$,
then $\mcM[\mcS_{\mcF} \trir \mcS_{\mcP}]$ is represented with respect to $\mbK \uhrc k$.
Moreover, there is $\mcP$ such that $\mcP$ is represented with respect to $\mbK \uhrc k$,
$\mcS_{\mcP} \subseteq_{cl} \mcP$ and $\mcF = \mcP[\mcS_{\mcP} \trir \mcS_{\mcF}]$
is {\em forbidden} with respect to $\mbK \uhrc k$.
This implies that the dimension of $\mcP$ is strictly larger than the dimension
of $\mcS_{\mcP}$ which is $k$.

By Observation~\ref{observation about dimension reducts}, there is $\widehat{\mcP}$
which is represented with respect to $\mbK$ and such that $\widehat{\mcP} \uhrc k = \mcP$.
We are assuming that $\mbK$ has the independent amalgamation property.
Hence, there are a represented $\mcC$, with respect to $\mbK$,
and embeddings $\tau_i : \widehat{\mcP} \to \mcC$, for $i = 1,2$,
such that $\tau_1 \uhrc |\mcS_{\mcP}| = \tau_2 \uhrc |\mcS_{\mcP}|$ and
$|\mcS_{\mcP}| = \tau_1(|\mcP|) \cap \tau_2(|\mcP|)$; so in particular $\mcS_{\mcP} \subset_{cl} \mcC$.
By replacing $\mcC$ with the closure of 
$\tau_1(|\widehat{\mcP}|) \cup \tau_2(|\widehat{\mcP}|)$ in $\mcC$, we may assume that
$\dim_{\mcC}(|\mcC|) = 2\dim_{\mcP}(|\mcP|) - \dim_{\mcS_{\mcP}}(|\mcS_{\mcP}|) =  2\dim_{\mcP}(|\mcP|) - k$.
Let $c = \dim_{\mcC}(|\mcC|)$.
Since $\dim_{\mcP}(|\mcP|) > k > 0$ (as noted above), we have $c >  \dim_{\mcP}(|\mcP|) > k > 0$, so $c \geq 3$.

If $k = 1$ then let $\mcU$
be the unique closed proper substructure
of $\mcS_{\mcF}$ with dimension 0.
If $k > 1$ then let $\mcU$ be any closed proper substructure of $\mcS_{\mcF}$
with dimension 1.
In both cases $\mcU$ is represented with respect to $\mbK$,
with respect to $\mbK \uhrc k$, and with respect to $\mbK \uhrc k-1$.

Let $\varphi$ denote the $\mcS_{\mcF}/\mcU$-extension axiom and let $\psi$ denote the 
$\mcC/\mcS_{\mcP}$-extension axiom.
We prove that $\lim_{n \to \infty}\delta_n(\varphi \wedge \psi) = 0$.
Let $\mcC' = \mcC \uhrc k$, so $\mcC'$ is represented with respect to $\mbK \uhrc k$,
and note that since the dimension of $\mcS_{\mcP}$ and of $\mcS_{\mcF}$ is $k$ and 
$\mcU \subset_{cl} \mcS_{\mcF}$ we have
$\mcU \uhrc k = \mcU$,
$\mcS_{\mcP} \uhrc k = \mcS_{\mcP}$ and
$\mcS_{\mcF} \uhrc k = \mcS_{\mcF}$.
The next lemma shows that instead of working with $\mbK$, 
$\varphi$ and $\psi$ we can work
with $\mbK \uhrc k$, the $\mcS_{\mcF}/\mcU$-extension axiom and the $\mcC'/\mcS_{\mcP}$-extension axiom.

\begin{lem}\label{reduction to k-reducts}
Let $p$ be the probability, with the measure $\delta_n$, that a structure in $\mbK_n$
satisfies both the $\mcS_{\mcF}/\mcU$-extension axiom ( = $\varphi$) and the 
$\mcC/\mcS_{\mcP}$-extension axiom ( = $\psi$).
Let $q$ be the probability, with the measure $\mbbP_{n,k}$,
that a structure in $\mbK_n \uhrc k$ satisfies both the $\mcS_{\mcF}/\mcU$-extension axiom and the $\mcC'/\mcS_{\mcP}$-extension axiom.
Then $p \leq q$.
\end{lem}

\noindent
{\em Proof.}
Recall that $k \leq \rho$.
By the definitions of $\mbb{P}_{n,k}$ and $\delta_n$, for every $\mcM \in \mbK_n \uhrc k$,
$$\mbb{P}_{n,k}(\mcM) = \delta_n\big(\{\mcN \in \mbK_n : \mcN \uhrc k = \mcM\}\big).$$
As mentioned above, $\mcS_{\mcP} \uhrc k = \mcS_{\mcP}$ and $\mcS_{\mcF} \uhrc k = \mcS_{\mcF}$.
So whenever $\mcN \in \mbK_n$ satisfies the $\mcS_{\mcF}/\mcU$-extension axiom, then
$\mcN \uhrc k$ satisfies the $\mcS_{\mcF}/\mcU$-extension axiom.
And whenever $\mcN \in \mbK_n$ satisfies the $\mcC/\mcS_{\mcP}$-extension axiom, then
$\mcN \uhrc k$ satisfies the $\mcC'/\mcS_{\mcP}$-extension axiom.
Therefore $p$ cannot exceed $q$.
\hfill $\square$
\\

\noindent
By Lemma~\ref{reduction to k-reducts} it suffices to prove that
\begin{itemize}
  \item[(6)] there is $\beta < 1$ such that for all sufficiently large $n$ the probability, 
  with the measure $\mbb{P}_{n,k}$, that a structure in $\mbK_n \uhrc k$
    satisfies both the $\mcS_{\mcF}/\mcU$-extension axiom and the $\mcC'/\mcS_{\mcP}$-extension axiom 
    does not exceed $\beta$; and if $k > 1$, then this probability tends to 0 as $n \to \infty$.
\end{itemize}

\noindent
The claim (6) follows from the next two lemmas and the definition of the measures 
$\mbb{P}_{n,r}$, $r \in \mbbN$.
Remember that $c$ is the dimension of $\mcC$ (and of $\mcC'$).

\begin{lem}\label{almost sureness for the k-1 reduct}
The probability, with the measure $\mbb{P}_{n,k-1}$, that a structure in $\mbK_n \uhrc k-1$
is $(\sigma^c(\lambda_n), c)$-saturated, with respect to $\mbK \uhrc k-1$,
tends to 1 as $n \to \infty$.
\end{lem}

\begin{lem}\label{almost neverness for the k reduct}
Let $\alpha$ be the number of represented structures with 
universe $|\mcS_{\mcF}|$.
Suppose that $\mcM \in \mbK_n \uhrc k-1$ is  
$(\sigma^c(\lambda_n), c)$-saturated with respect to $\mbK \uhrc k-1$ and let
$$\mbE_k(\mcM) = \big\{ \mcN \in \mbK \uhrc k : \mcN \uhrc k-1 = \mcM \big\}.$$
(i) The proportion of structures in $\mbE_k(\mcM)$
which satisfy both the $\mcS_{\mcF}/\mcU$-extension axiom and the 
$\mcC'/\mcS_{\mcP}$-extension axiom never exceeds $1 - 1/(1 + \alpha)$.\\
(ii) If $k > 1$ then the proportion of structures in $\mbE_k(\mcM)$
which satisfy both the $\mcS_{\mcF}/\mcU$-extension axiom and the 
$\mcC'/\mcS_{\mcP}$-extension axiom never exceeds 
$\alpha\left\|\mcS_{\mcF}\right\| \big/ \sigma^c(\lambda_n)$.
Note that this expression does not depend on $\mcM$ and approaches 0 as $n \to \infty$.
\end{lem}

\noindent
{\em Proof of Lemma~\ref{almost sureness for the k-1 reduct}}
Note that when saying that $\mbK \uhrc k-1$ accepts $r$-substitutions over $L_0$ 
we only consider substitutions of the form 
$[\mcA \trir \mcA']$ where $\mcA$ and $\mcA'$ are represented with respect to $\mbK \uhrc k-1$.

Let $\mbK'_n = \mbK_n \uhrc k-1$ and $\mbK' = \mbK \uhrc k-1$.
Let $\mbb{P}'_{n,0}$ be the uniform measure on $\mbK'_n \uhrc 0$ ($= \mbK_n \uhrc 0$)
and for positive $r \in \mbbN$, 
let $\mbb{P}'_{n,r}$ be the $(\mbK'_n \uhrc 0, \ldots, \mbK'_n \uhrc r-1)$-conditional measure
on $\mbK'_n \uhrc r$.
Observe that we have the following:
\begin{align*}
&\text{For $r \leq k-1$,  $\mbK'_n \uhrc r = \mbK_n \uhrc r$ and $\mbb{P}'_{n,r}$ coincides with $\mbb{P}_{n,r}$}\\
&\text{For $r \geq k-1$, $\mbK'_n \uhrc r = \mbK'_n = \mbK_n \uhrc k-1$ and $\mbb{P}'_{n,r}$ coincides with  $\mbb{P}'_{n,k-1}$}
\end{align*}
As $c > k-1$, we in particular have
\begin{align*}
&\mbK'_n \uhrc c = \mbK'_n = \mbK_n \uhrc k-1\\
&\text{and $\mbb{P}'_{n,c}$ coincides with $\mbb{P}'_{n,k-1}$ which in turn coincides with $\mbb{P}_{n,k-1}$.}
\end{align*}
So $\mbb{P}'_{n,c}$ and $\mbb{P}_{n,k-1}$ are the same measure on $\mbK'_n \uhrc c = \mbK_n \uhrc k-1$.
Thus, in order to prove Lemma~\ref{almost sureness for the k-1 reduct} it suffices to show that
the probability, with the measure $\mbb{P}'_{n,c}$,
that a structure in $\mbK'_n \uhrc c$ is $(\sigma^c(\lambda_n), c)$-saturated, with respect to $\mbK' \uhrc c$,
tends to 1 as $n \to \infty$.
If, for $n,r \in \mbbN$, we let
$$\mbX'_{n,r} = \{\mcM \in \mbK'_n \uhrc r : \mcM \text{ is $(\sigma^r(\lambda_n), c)$-saturated}\},$$
then the claim of Lemma~\ref{almost sureness for the k-1 reduct}
is that 
\begin{equation*}\label{conditional probability for the prime-measure}
\lim_{n\to \infty} \mbb{P}'_{n,c}(\mbX'_{n,c}) = 1. \tag{7}
\end{equation*}
By assumption, $(\mcG_n : n \in \mbbN)$ is polynomially $c$-saturated,
and, as mentioned in the beginning of the proof, 
$\mbK'$ (= $\mbK \uhrc k-1$) accepts $r$-substitutions over $L_0$
for every $r \in \mbbN$, so in particular for $r = c$.
In other words, $\mbK'$ satisfies the same assumptions, with respect to $(\mcG_n : n \in \mbbN)$ and
$L_0$, as $\mbK$ did in the
proof of Theorem~\ref{0-1 law for pregeometries}, 
and $\mbP'_{n,c}$ is the $(\mbK'_n \uhrc 0, \ldots, \mbK'_n \uhrc r-1)$-conditional measure
on $\mbK'_n \uhrc r$, where $\mbK'_n \uhrc 0 = \mbK_n \uhrc 0$.
Therefore, the statement of~(7) (and its underlying assumptions) is the same as the statement of~(3)
(and its underlying assumptions)
if we replace $\mbK$, $\mbb{P}_{n,k}$ and $\mbX_{n,k}$ by
$\mbK'$, $\mbb{P}'_{n,c}$ and $\mbX'_{n,c}$, respectively.
Hence,~(7) is proved in exactly the same way as~(3), by just replacing
$\mbK$, $\mbb{P}_{n,r}$ and $\mbX_{n,r}$ with 
$\mbK'$, $\mbb{P}'_{n,r}$ and $\mbX'_{n,r}$, respectively, for $n,r \in \mbbN$.
\hfill $\square$
\\

\noindent
{\em Proof of Lemma~\ref{almost neverness for the k reduct}.}
Suppose that  $\mcM \in \mbK_n \uhrc k-1$ is $(\sigma^c(\lambda_n), c)$-saturated, with respect to $\mbK \uhrc k-1$,
and let 
$$\mbE_k(\mcM) = \big\{ \mcN \in \mbK \uhrc k : \mcN \uhrc k-1 = \mcM \big\}.$$
Let $\alpha$ be the number of represented structures, with respect to $\mbK \uhrc k$,
with universe $|\mcS_{\mcP}|$.
It suffices to show that the proportion of structures in $\mbE_k(\mcM)$ which satisfy both
the $\mcS_{\mcF}/\mcU$-extension axiom and the $\mcC'/\mcS_{\mcP}$-extension axiom
does not exceed $1 - 1/(1 + \alpha)$;
and if $k > 1$ then this proportion approaches 0 as $n\to \infty$.
We will consider the cases $k = 0$ and $k > 0$ one by one.

First assume that $k = 1$.
Then, by the choice of $\mcU$, $\mcU$ has dimension 0 and is represented,
since it is a closed substructure of a represented structure.
By assumption there is a unique, up to isomorphism, represented
structure of dimension 0.
Hence, every represented structure (with respect to $\mbK$, $\mbK \uhrc k$ or $\mbK \uhrc k-1$)
contains a copy of $\mcU$.
Therefore every $\mcM \in \mbK \uhrc k$ which satisfies the $\mcS_{\mcF}/\mcU$-extension axiom
contains a copy of $\mcS_{\mcF}$.
Note that if $\mcN \in \mbK \uhrc k$ satisfies the $\mcC'/\mcS_{\mcP}$-extension axiom,
then the $\mcP/\mcS_{\mcP}$-multiplicity of $\mcN$ is at least 2.
Now we can argue as in Section~\ref{proving upper limit on the probability of an extension axiom}.
More precisely, the proofs of Lemmas~\ref{reducing the multiplicity to 0}, \ref{different mutations} and
\ref{almost disjoint mutation sets} as well as the proof
of part (i) of Theorem~\ref{upper limit on the probability of an extension axiom} carry over to
the present context if we have the following in mind:
The structures $\mcS_{\mcP}$, $\mcS_{\mcF}$, $\mcP$ and $\mcF$ play the same roles in the present context as
in Section~\ref{proving upper limit on the probability of an extension axiom};
in the present context `closed substructures' play the role of `substructures' in 
Section~\ref{proving upper limit on the probability of an extension axiom};
dimension plays the role here that cardinality had in that section;
and $\mbE_k(\mcM)$ plays the role here that `$\mbK_n$' had in that section.
In this way we can conclude that the proportion of 
$\mcN \in \mbE_k(\mcM)$ which contain a copy of $\mcS_{\mcF}$ and
satisfy the $\mcC'/\mcS_{\mcP}$-extension axiom never exceeds $1 - 1/(1 + \alpha)$.

Now suppose that $k > 1$.
Again, the reasoning from Section~\ref{proving upper limit on the probability of an extension axiom}
carries over to the present context.
Since we assume $k > 1$, $\mcU$ has dimension 1 and $\mcU \subset_{cl} \mcS_{\mcF}$.
As noted earlier, $c > k >1$.
Since  $\mcM$ is $(\sigma^c(\lambda_n), c)$-saturated, with respect to $\mbK \uhrc k-1$,
$\mcM$ contains at least $\sigma^c(\lambda_n)$ distinct copies of $\mcU$.
Since $\mcM$ and every $\mcN \in \mbE_k(\mcM)$ agree on all substructures of dimension at most $k-1 \geq 1$,
it follows that every $\mcN \in \mbE_k(\mcM)$ contains at least $\sigma^c(\lambda_n)$ distinct copies of $\mcU$.
Suppose that $\mcN \in  \mbE_k(\mcM)$ satisfies both the $\mcS_{\mcF}/\mcU$-extension axiom and the
$\mcC'/\mcS_{\mcP}$-extension axiom.
First we notice that the satisfaction of the $\mcS_{\mcF}/\mcU$-extension axiom implies that
$\mcN$ contains at least  $\sigma^c(\lambda_n)/\left\|\mcS_{\mcF}\right\|$ 
distinct copies of $\mcS_{\mcF}$ (the copies may {\em partially} overlap, but this poses no problem).
Secondly, the satisfaction of the $\mcC'/\mcS_{\mcP}$-extension axiom implies that the 
$\mcP/\mcS_{\mcP}$-multiplicity of $\mcN$ is at least 2.

As in the previous case (when $k=1$) the proofs of
lemmas~\ref{reducing the multiplicity to 0}, \ref{different mutations} and
\ref{almost disjoint mutation sets} carry over -- with the already mentioned 
provisos -- to this context. 
But we are now able to continue the argument similarly as in the proof
of part (iii) of Theorem~\ref{upper limit on the probability of an extension axiom}.
The number $\sigma^c(\lambda_n)$ plays the same role here as the number
`$m_n$' did in the proof of 
part (iii) of Theorem~\ref{upper limit on the probability of an extension axiom}.
In a similar way as in that proof we can now derive that 
`$\alpha\left\|\mcS_{\mcF}\right\| \big/ \sigma^c(\lambda_n)$'
(instead of `$k\alpha/m_n$' as in the proof of part (iii) of 
Theorem~\ref{upper limit on the probability of an extension axiom})
is an upper bound for the proportion of $\mcN \in \mbE_k(\mcM)$ such that 
$\mcN$ satisfies the $\mcS_{\mcF}/\mcU$-extension axiom and
the $\mcP/\mcS_{\mcP}$-multiplicity of $\mcN$ is at least 2.
\hfill $\square$

\section{Random $l$-colourable structures}\label{l-colourable structures}

\noindent
In this section and the next we consider $l$-colourable,
as well as {\em strongly} $l$-colourable, relational structures 
and zero-one laws for these,
with the uniform probability measure and with a measure which is derived from the dimension conditional
measure with trivial underlying pregeometry. In all cases we have a zero-one law,
and we get the same almost sure theory whether we work with the uniform probability measure
or with the probability measure derived from the dimension conditional measure.
(The notions `zero-one law' and `almost sure theory' are explained 
in Section~\ref{preliminaries about zero-one laws}.)
In the case when one considers the probability measure derived from the dimension conditional
measure the proof only uses methods of formal logic, while in the case when one
considers the uniform probability measure the proof uses, in addition, 
results about the typical distribution of colours, which are proved by combinatorial arguments.
Therefore, we start, in this section, by considering the probability measure derived
from the dimension conditional measure.
In Section~\ref{the uniform measure and l-colourable structures} we state the corresponding
results for the uniform probability measure and complete their proofs.

In this section, $l \geq 2$ is a fixed integer and for each $n \in \mbbN$, $\mbK_n$ is defined as in
Example~\ref{example of coloured structures of first kind} for $F = \{1\}$, 
and $\mbSK_n$ is defined as $\mbK_n$ in
Example~\ref{example of coloured structures of second kind} for $F = \{1\}$.
Note that `$F = \{1\}$' means that the universe of every $\mcM \in \mbK_n$
is $\{1, \ldots, n\}$ and that the pregeometry is trivial
(i.e. $\cl_{\mcM}(X) = X$ for every $\mcM \in \mbK_n$ and every $X \subseteq M$).
As usual, let $\mbK = \bigcup_{n \in \mbbN} \mbK_n$
and $\mbSK = \bigcup_{n \in \mbbN} \mbSK_n$.
The notation $L_{col}$ (the language of the $l$ colours), $L_{rel}$
(the language of relations) and $L$ mean the same as in the mentioned examples.
{\em But we add the assumption that all relation symbols of the vocabulary of $L_{rel}$ have
arity at least 2.} (Colouring unary relations is not so interesting.)
{\em When working with strong $l$-colourings, that is, with $\mbSK$, we also assume that
$l$ is at least as great as the arity of every relation symbol in the vocabulary of $L_{rel}$};
for otherwise the interpretations of some relation symbol(s) will be empty for all 
$l$-coloured structures, and then there is no point in having this (or these) relation symbol(s).
Observe that if there are no relation symbols of arity greater than 2, then $\mbK = \mbSK$,
as the pregeometry is trivial.
A structure which is isomorphic with one in $\mbK$ is called {\bf \em $l$-coloured}.
A structure which is isomorphic with one in $\mbSK$ is called {\bf \em strongly $l$-coloured}.
Note that being $l$-coloured (strongly $l$-coloured) is equivalent to being represented
with respect to $\mbK$ ($\mbSK$).

For each $n$, let 
\begin{align*}
&\mbC_n = \big\{\mcM \uhrc L_{rel} : \mcM \in \mbK_n\big\}, \quad \quad \mbC = \bigcup_{n \in \mbbN}\mbC_n,\\
&\mbS_n = \big\{\mcM \uhrc L_{rel} : \mcM \in \mbSK_n\big\} \quad 
\text{ and }\quad \mbS = \bigcup_{n \in \mbbN}\mbS_n.
\end{align*}
A structure which is isomorphic to one in $\mbC$ (i.e. represented with respect
to $\mbC$) will be called {\bf \em $l$-colourable}.
A structures which is isomorphic to one in $\mbS$ will be called {\bf \em strongly $l$-colourable}.
It is clear that an $L_{rel}$-structure $\mcM$ is (strongly) $l$-colourable if and only if
there is a function $f : M \to \{1, \ldots, l\}$, called an {\bf \em (strong) $l$-colouring}, 
such that the expansion $\mcM'$ of $\mcM$ to $L$,
defined by $\mcM' \models P_i(a)$ if and only if $f(a) = i$, is isomorphic with a member
of $\mbK$ ($\mbSK$).
Therefore we can, when convenient, use (strong) $l$-colouring functions instead of the relation symbols 
$P_1, \ldots, P_l$ to represent (strong) $l$-colourings.

In this section, $\delta_n^{\mbK}$ denotes the dimension conditional measure on $\mbK_n$
and $\delta_n^{\mbSK}$ denotes the dimension conditional measure on $\mbSK_n$
(see Definition~\ref{definition of dimension conditional measure}).
For each $n$, we consider the measures, $\delta_n^{\mbC}$ on $\mbC_n$ 
and $\delta_n^{\mbS}$ on $\mbS_n$
which are inherited from $\mbK_n$ and $\mbSK_n$, respectively, in the following sense:
\begin{align*}
&\text{For every } \mbX \subseteq \mbC_n, \ \ \delta_n^{\mbC}(\mbX) = 
\delta_n^{\mbK}\big(\{\mcM \in \mbK_n : \mcM \uhrc L_{rel} \in \mbX\}\big). \\
&\text{For every } \mbX \subseteq \mbS_n, \ \ \delta_n^{\mbS}(\mbX) = 
\delta_n^{\mbSK}\big(\{\mcM \in \mbSK_n : \mcM \uhrc L_{rel} \in \mbX\}\big).
\end{align*}
For every $L_{rel}$-sentence $\varphi$, let 
$\delta_n^{\mbC}(\varphi) = \delta_n^{\mbC}\big(\{\mcM \in \mbC_n : \mcM \models \varphi\}\big)$
and\\
$\delta_n^{\mbS}(\varphi) = \delta_n^{\mbS}\big(\{\mcM \in \mbS_n : \mcM \models \varphi\}\big)$.

\begin{theor}\label{0-1 law for l-colourable structures} 
For every sentence $\varphi \in L_{rel}$, 
\begin{itemize}
\item[(i)] $\lim_{n \to \infty}\delta_n^{\mbC}(\varphi) = 0$ or 
$\lim_{n \to \infty}\delta_n^{\mbC}(\varphi) = 1$, and
\item[(ii)] $\lim_{n \to \infty}\delta_n^{\mbS}(\varphi) = 0$ or 
$\lim_{n \to \infty}\delta_n^{\mbS}(\varphi) = 1$.
\end{itemize}
\end{theor}

\noindent
Theorem~\ref{0-1 law for l-colourable structures} will be proved in 
Section~\ref{proof of 0-1 law for l-colourable structures}.
We also state the corresponding theorem for the uniform probability measure,
although it will be restated, 
with more detail as Theorems~\ref{main theorem} and~\ref{main theorem for strong colourings},
in Section~\ref{the uniform measure and l-colourable structures}
where its proof will be completed.

\begin{theor}\label{0-1 law for l-colourable structures with the uniform probability measure}
For every sentence $\varphi \in L_{rel}$ the following holds:
\begin{itemize}
\item[(i)] The proportion of $\mcM \in \mbC_n$ in which $\varphi$ is true 
approaches either 0 or 1, as $n$ approaches infinity.
\item[(ii)] The proportion of $\mcM \in \mbS_n$ in which $\varphi$ is true 
approaches either 0 or 1, as $n$ approaches infinity.
\end{itemize}
\end{theor}

\begin{rem}\label{remark about the 0-1 laws}{\rm
(i) Let the relation symbols of $L_{rel}$ be $R_1, \ldots, R_{\rho}$ and
let $I \subseteq \{1, \ldots, \rho\}$.
If we add the restriction that for every $i \in I$, $R_i$ is always interpreted as an
irreflexive and symmetric relation 
(see Remark~\ref{remark that the theorems generalize to symmetric structures}),
then Theorems~\ref{0-1 law for l-colourable structures},
~\ref{0-1 law for l-colourable structures with the uniform probability measure}
and Proposition~\ref{definability of colourings under richness condition}
still hold. 
The proofs in this section are exactly the same even if we add this extra assumption.
But the combinatorial arguments in Section~\ref{the uniform measure and l-colourable structures}, 
needed to complete the proof of 
Theorem~\ref{0-1 law for l-colourable structures with the uniform probability measure}
are sensitive to whether a relation symbol is always interpreted as an irreflexive and symmetric
relation, or not.
For this reason the notation in Section~\ref{the uniform measure and l-colourable structures}
(but not in this section)
specifies which relation symbols are always interpreted as irreflexive and symmetric relations.

(ii) It is open whether Theorems~\ref{0-1 law for l-colourable structures} 
and~\ref{0-1 law for l-colourable structures with the uniform probability measure} 
still hold if $F$ is allowed to be a (fixed) finite field,
thus giving a nontrivial underlying pregeometry,
and $\mbK_n$, $\mbSK_n$, $\mbC_n$ and $\mbS_n$ are, apart from this difference, defined as before.
}\end{rem}

\subsection{Proof of Theorem~\ref{0-1 law for l-colourable structures}}
\label{proof of 0-1 law for l-colourable structures}

\noindent
The proof depends on Theorem~\ref{0-1 law for pregeometries}
which is used in the proof of Lemma~\ref{probability of X approaches 1} below.
Apart from Lemmas~\ref{finding the structure S}
and~\ref{finding the structure U} below, the proof is the same, except for obvious 
changes of notation, in the case of $\mbS$ (strongly $l$-colourable structures)
as in the case of $\mbC$ ($l$-colourable structures).
For this reason, and to avoid cluttering notation and language,
we prove Theorem~\ref{0-1 law for l-colourable structures}
by speaking of $\mbK_n$, $\mbC_n$, $l$-coloured structures and $l$-colourable structures.
Only when proving Lemmas~\ref{finding the structure S}
and~\ref{finding the structure U} will we separate the two cases explicitly.

The general pattern of the proof is a familiar one.
We collect into a theory $T_{\mbC}$ a certain type of extension axioms
(to be called `$l$-colour compatible extension axioms') together with
sentences which describe all possible isomorphism types of structures in $\mbC$.
Then we show that for every $\psi \in T_{\mbC}$,
$\lim_{n\to\infty}\delta_n^{\mbC}(\psi) = 1$, which implies (via compactness) 
that $T_{\mbC}$ is consistent. 
After this we show that $T_{\mbC}$ is complete by showing that it is countably categorical.
The zero-one law is now a straightforward consequence of the previously proven facts,
together with compactness.

\begin{rem}{\rm
We can {\em not} expect that for every extension axiom $\varphi$ of $\mbC$,
$\lim_{n \to \infty} \delta_n^{\mbC}(\varphi) = 1$.
For example, suppose that the vocabulary of $L_{rel}$ contains only one relation symbol
which is binary (which implies $\mbS = \mbC$), and that $l = 2$. 
Then there is {\em no} $2$-colourable $L_{rel}$-structure which satisfies
all $3$-extension axioms of $\mbC$. For if $\mcM$ would be such a structure,
then it is easy to see that $\mcM$ would contain a 3-cycle or a 5-cycle (it does not matter if it
is directed or not) which contradicts that $\mcM$ is 2-colourable.
}\end{rem}

\noindent
In order to define the type of extension axioms that are useful in this context,
we need to find a way of expressing, with an $L_{rel}$-formula, that
two elements in an $L$-structure have the same colour. 
In fact, it suffices to find an $L_{rel}$-formula $\xi(y,z)$ such that
with $\delta_n^{\mbK}$-probability approaching 1 as $n \to \infty$:
if $\mcM \in \mbK_n$ and $a, b \in M$, then $\mcM \models \xi(a,b)$ if and only
if $a$ and $b$ have the same colour in $\mcM$.
The following lemma is a first step in that direction:

\begin{lem}\label{finding the structure S}
There is an (strongly) $l$-colourable structure $\mcS$ and distinct $a, b \in S$ such that
the following hold:
\begin{itemize}
	\item[(a)] Whenever $\gamma : S \to \{1, \ldots, l\}$ is an (strong) $l$-colouring of $\mcS$, then
	$\gamma(a) = \gamma(b)$; in other words, whenever $\mcS$ is (strongly) $l$-coloured then
	$a$ and $b$ get the same colour.
	\item[(b)] For every $i \in \{1, \ldots, l\}$, there is an (strong) $l$-colouring
	$\gamma_i : S \to \{1, \ldots, l\}$ of $\mcS$ such that $\gamma(a) = \gamma(b) = i$.
\end{itemize}
\end{lem}

\noindent
{\em Proof.}
We must treat the case of $\mbC$, i.e. $l$-colourable structures, and the case of
$\mbS$, i.e. strongly $l$-colourable structures, separately.
We start with the case of $\mbC$.
By assumption all relation symbols in the vocabulary of $L_{rel}$ have arity at least 2.
Let $r$ be the minimum of the arities of relation symbols in the vocabulary of $L_{rel}$,
so $r \geq 2$, and let $R$ be a relation symbol in the vocabulary of $L_{rel}$
which has arity $r$.
Let $S = \{0, 1, \ldots, (r-1)l\}$ and 
let $R^{\mcS}$ consist exactly of all tuples $(s_1, \ldots, s_r)$ of {\em distinct}
elements from $S$ such that 
$$\{s_1, \ldots, s_r\} \subseteq S - \{0\} \quad \text{ or } \quad
\{s_1, \ldots, s_r\} \subseteq S - \{1\}.$$
For all other relation symbols $Q$ of the vocabulary of $L_{rel}$, let $Q^{\mcS} = \es$.
Note that there is no relationship in $\mcS$ which contains both $0$ and $1$.

We first show that there is a colouring $\gamma : S \to \{1, \ldots, l\}$ of $\mcS$ such
that $\gamma(0) = \gamma(1) = 1$. This will prove (b), because any permutation
of the colours of an $l$-colouring gives a new $l$-colouring.
Let both $0$ and $1$ be assigned the colour 1.
Then assign the colour $1$ to exactly $r-2$ elements
$s_1, \ldots, s_{r-2} \in S - \{0, 1\}$.
So exactly $r$ elements of $S = \{0,1, \ldots, (r-1)l\}$ 
have been assigned the colour $1$; and these elements are $0, 1, s_1, \ldots, s_{r-2}$.
Hence 
$$\big|\{S - \{0, 1, s_1, \ldots, s_{r-2}\}\big| = 
(r - 1)l + 1 - r = (r-1)(l-1),$$
so $S - \{0, 1, s_1, \ldots, s_{r-2}\}$ can be partitioned into $l-1$
parts each of which contains exactly $r-1$ elements.
Consequently, we can, for each colour $i \in \{2, \ldots, l\}$,
assign the colour $i$ to exactly $r-1$ elements in 
$S - \{0, 1, s_1, \ldots, s_{r-2}\}$.
Since no colour other than $1$ has been assigned to more that $r-1$ elements, the result
is an $l$-colouring of $\mcS$.

We now prove (a).
Assume that $\gamma : S \to \{1, \ldots, l\}$ is a colouring of $\mcS$.
Note that $\big|S - \{0\}\big| = \big|S - \{1\}\big| = (r-1)l$.
By the definition of $\mcS$, every $r$-tuple of distinct elements
$(s_1, \ldots, s_r) \in (S - \{0\})^r$
is an $R$-relationship.
Hence, for every colour $i \in \{1, \ldots, l\}$,
we must have $\big|\gamma^{-1}(i) \cap (S - \{0\})\big| = r-1$.
Suppose that $\gamma(1) = 1$. (If $\gamma(1) \in \{2, \ldots, l\}$ the argument is analogous.)
Assume, for a contradiction, that $\gamma(0) = i \neq 1$.
Above we concluded that $\big|\gamma^{-1}(i) \cap (S - \{0\})\big| = r-1$.
Since $\gamma(0) = i$ we get $\big|\gamma^{-1}(i)\big| = r$, and as
$\gamma(1) \neq i$, we get $\gamma^{-1}(i) \subseteq S - \{1\}$.
Hence, there are distinct $s_1, \ldots, s_r \in \gamma^{-1}(i) \subseteq S - \{1\}$.
By the definition of $\mcS$, $(s_1, \ldots, s_r) \in R^{\mcS}$.
Since $\gamma$ assigns all elements $s_1, \ldots, s_r$ the colour $i$,
this contradicts that $\gamma$ is a colouring of $\mcS$.
So if we take $a = 0$ and $b = 1$, then the lemma holds for this $\mcS$
in the case of (not necessarily strong) $l$-colourings.

Now we prove the lemma in the case of strong $l$-colourings.
Let $S = \{0,1, \ldots, l\}$. 
Let $R$ be any symbol from the vocabulary of $L_{rel}$, so the arity $r$ of $R$ is at least $2$.
By assumption, since we work with strong $l$-colourings now, $2 \leq r \leq l$.
Let $R^{\mcS}$ consist exactly of all tuples $(s_1, \ldots, s_r)$ of {\em distinct}
elements from $S$ such that 
$$\{s_1, \ldots, s_r\} \subseteq S - \{0\} \quad \text{ or } \quad
\{s_1, \ldots, s_r\} \subseteq S - \{1\}.$$
For all other relation symbols $Q$ of the vocabulary of $L_{rel}$, let $Q^{\mcS} = \es$.
Note that there is no relationship in $\mcS$ which contains both $0$ and $1$.
Therefore any assignment of the same colour $i \in \{1, \ldots, l\}$ to $0$ and $1$
can be extended to a strong $l$-colouring of $\mcS$.
Also note that every strong $l$-colouring of $\mcS$ must give all elements in
$S - \{0\}$ different colours; and it must give all elements in $S - \{1\}$ different
colours. Since $|S| = l-1$ there is no other choice but giving $0$ and $1$ the same colour.
Hence the lemma, in the case of strong $l$-colourings, holds for this $\mcS$ with $a = 0$ and $b = 1$.
\hfill $\square$

\begin{notation}\label{notation for S and xi}{\rm
(i) Let $\mcS$ be an $l$-colourable structure  and $a, b \in S$ distinct elements such that 
Lemma~\ref{finding the structure S} is satisfied.
Note that we must have $|S| \geq 3$.
Without loss of generality we assume that $|\mcS| = S = \{1, \ldots, s\}$ for some $s \geq 3$
and that $a = s-1$ and $b = s$.
Hence every assignment of the same colour to $s-1$ and $s$ can be extended to
an $l$-colouring of $\mcS$, and every $l$-colouring of $\mcS$ gives $s-1$ and $s$ the same colour.\\
(ii) Let $\chi_{\mcS}(x_1, \ldots, x_s)$ be a quantifier-free $L_{rel}$-formula which expresses
the $L_{rel}$-iso\-mor\-phism type of $\mcS$; more precisely,
for every $L_{rel}$-structure $\mcM$ and all $a_1, \ldots, a_s \in M$,
$\mcM \models \chi_{\mcS}(a_1, \ldots, a_s)$ if and only if
the map $a_i \mapsto i$ is an isomorphism from $\mcM \uhrc \{a_1, \ldots, a_s\}$ to $\mcS$.\\
(iii) Let $\xi(y,z)$ be the formula 
$$y = z \ \vee \ \exists u_1, \ldots, u_{s-2} \chi_{\mcS}(u_1, \ldots, u_{s-2}, y, z).$$
(iv) For $n, k \in \mbbN$ let $\mbX_{n,k} \subseteq \mbK_n$ be the set of all
$\mcM \in \mbK_n$ which satisfy all $k$-extension axioms with respect to $\mbK$.
}\end{notation}

\begin{lem}\label{probability of X approaches 1}
For every $k \in \mbbN$, $\lim_{n\to\infty}\delta_n^{\mbK}(\mbX_{n,k}) = 1$.
\end{lem}

\noindent
{\em Proof.}
As mentioned in Examples~\ref{examples of polynomially k-saturated pregeometries}
and~\ref{example of coloured structures of first kind},
for every $k \in \mbbN$,
the trivial pregeometry is polynomially $k$-saturated 
and $\mbK$ accepts $k$-substitutions over the language with empty vocabulary.
By Theorem~\ref{0-1 law for pregeometries}~(i),
for every extension axiom $\varphi$ of $\mbK$, $\lim_{n\to\infty}\delta_n^{\mbK}(\varphi) = 1$.
The lemma follows since there are only finitely many $k$-extension axioms.
(In the case of strongly $l$-colourable structures 
we look back at Example~\ref{example of coloured structures of second kind}
instead of Example~\ref{example of coloured structures of first kind}.)
\hfill $\square$

\begin{lem}\label{the property of xi}
Let $\mcM \in \mbK$ and $a, b \in M$.\\
(i) If $\mcM \models \xi(a,b)$ then $a$ and $b$ have the same colour in $\mcM$, 
i.e. for some $i \in \{1, \ldots, l\}$, $\mcM \models P_i(a) \wedge P_i(b)$.\\
(ii) If $k \geq \left\|\mcS\right\|$ and $\mcM \in \mbX_{n,k}$, then
$\mcM \models \xi(a,b)$ if and only if 
$a$ and $b$ have the same colour in $\mcM$.
\end{lem}

\noindent
{\em Proof.}
(i) Suppose that $\mcM \in \mbK$ and $\mcM \models \xi(a,b)$.
If $a=b$ then $a$ and $b$ have the same colour, so suppose that $a \neq b$.
Then there are $m_1, \ldots, m_{s-2} \in M$
such that 
$$\mcM \models \chi_{\mcS}(m_1, \ldots, m_{s-2}, a, b).$$
It follows that the $L_{rel}$-reduct of $\mcM \uhrc \{m_1, \ldots, m_{s-2}, a, b\}$ is
isomorphic with $\mcS$ via the $L_{rel}$-isomorphism $m_i \mapsto i$, for $i = 1, \ldots, s-2$,
$a \mapsto s-1$ and $b \mapsto s$.
Then we get an $l$-colouring of $\mcS$ by letting
$i$ get the same colour as $m_i$, for $i = 1, \ldots, s-2$,
letting $s-1$ get the same colour as $a$, and letting $s$ get the same colour as $b$.
From Lemma~\ref{finding the structure S} it follows that $s-1$ and $s$ must have 
the same colour in $\mcS$; hence $a$ and $b$ must
have the same colour in $\mcM$.

(ii) Let $k \geq \left\|\mcS\right\|$ and $\mcM \in \mbX_{n,k}$.
If $a = b$ then immediately from the definition of $\xi(y,z)$ we get $\mcM \models \xi(a,b)$.
Suppose that $a, b \in M$ are distinct elements which have the same colour in $\mcM$, 
that is, 
for some colour $i \in \{1, \ldots, l\}$, $\mcM \models P_i(a) \wedge P_i(b)$.
By Lemma~\ref{finding the structure S} and Notation~\ref{notation for S and xi}, 
there is an $l$-coloured structure $\mcS_i$ such that
$\mcS_i \uhrc L_{rel} = \mcS$ and $\mcS_i \models P_i(s-1) \wedge P_i(s)$.
Let $\mcS'_i = \mcS_i \uhrc \{s-1, s\}$.
Since $a$ and $b$ have the same colour in $\mcM$, there is no binary relationship of $\mcM$
which includes both $a$ and $b$. Hence $\mcM \uhrc \{a, b\}$ has no other 
relationships than the colour of $a$ and of $b$ which is $i$ in both cases.
By the properties of $\mcS$ (given by Lemma~\ref{finding the structure S}
and Notation~\ref{notation for S and xi}), $\mcS'_i$ has no other relationships
than the colour of $s-1$ and of $s$ which is $i$ in both cases.
Hence, any bijection between $\{s-1, s\}$ and $\{a,b\}$ is an isomorphism
between $\mcS'_i$ and $\mcM \uhrc \{a,b\}$.
Since $\mcM \in \mbX_{n,k}$ and $k \geq \left\|\mcS\right\|$, it follows that
$\mcM$ satisfies the $\mcS_i/\mcS'_i$-extension axiom.
This implies that there are $m_1, \ldots, m_{s-2} \in M$ such that
the map $m_i \mapsto i$, for $i = 1, \ldots, s-2$, $a \mapsto s-1$ and
$b \mapsto s$, is an isomorphism from $\mcM \uhrc \{m_1, \ldots, m_{s-2}, a, b\}$
to $\mcS_i$. Since $\mcS_i \uhrc L_{rel} = \mcS$ we get
$\mcM \models \chi_{\mcS}(m_1, \ldots, m_{s-2}, a, b)$, 
so $\mcM \models \xi(a,b)$.
\hfill $\square$
\\

\noindent
Besides being able to express (with high probability)
with the $L_{rel}$-formula $\xi(y,z)$ that two elements have the same colour,
we also need to be able to represent colours by elements (having those colours)
in a structure, and we must be able to define such elements with
an $L_{rel}$-formula. This is taken care of by
Lemma~\ref{finding the structure U}, Notation~\ref{notation for U etcetera}
and Lemma~\ref{auxiliary properties}, below.
In some more detail, the structure $\mcU$ in the next lemma will help us to
define an $L_{rel}$-formula $\zeta(x_1, \ldots, x_u)$, in Notation~\ref{notation for U etcetera},
where $u \geq l$,
such that if $\mcM \in \mbX_{n,k}$, then 
$\mcM \models \exists x_1, \ldots, x_u \zeta(x_1, \ldots, x_u)$ and if $\mcM \models \zeta(a_1, \ldots, a_u)$,
then the first $l$ elements $a_1, \ldots, a_l$ have different colours in $\mcM$.
The formula $\zeta$ will be used (before Lemma~\ref{lemma on colour compatible extension axioms})
when we define a restricted version of extension axioms for $\mbC$, the
`$l$-colour compatible extension axioms'.

\begin{lem}\label{finding the structure U}
There is an (strongly) $l$-colourable structure $\mcU$ which is not (strongly) $(l-1)$-colourable and such that
$\left\|\mcU\right\|$ is divisible by $l$ and every partition of $|\mcU|$ into $l$
parts of equal size gives rise to an (strong) $l$-colouring of $\mcU$.
\end{lem}

\noindent
{\em Proof.}
We deal with the cases of $l$-colourings and strong $l$-colourings separately and
begin with the case of $l$-colourings.
Let $R$ be a relation symbol from the vocabulary of $L_{rel}$, so the arity, call it $r$,
of $R$ is at least 2. Recall that $l \geq 2$.
Let $\mcU$ be the $L_{rel}$-structure with universe $U = \{1, \ldots, l(r - 1)\}$,
where 
$$R^{\mcU} = \big\{(u_1, \ldots, u_r) \in U^r : i \neq j \Rightarrow u_i \neq u_j \big\},$$
and the interpretation of every other relation symbol is empty.
Then $U$ can be partitioned into $l$ parts, each part with exactly $r - 1$ elements.
Hence every tuple $(u_1, \ldots, u_r) \in U^r$ of distinct elements must contain
$u_i$ and $u_j$ from different parts of the partition. 
Consequently, $\mcU$ is $l$-colourable.
However, if $U$ is partitioned into $l-1$ parts, then at least one part must contain
$r$ distinct elements $u_1, \ldots, u_r$, and since
$(u_1, \ldots, u_r) \in R^{\mcU}$, the partition does not represent an $(l-1)$-colouring of $\mcU$.
Thus, $\mcU$ is not $(l-1)$-colourable.

The case of strong $l$-colourings is even simpler.
Again we take any relation symbol $R$ from the vocabulary of $L_{rel}$. Its arity, say $r$,
is by assumption at least 2.
By the extra assumption when dealing with strongly $l$-colourable structures we in fact have
$2 \leq r \leq l$.
We then let $U = \{1, \ldots, l\}$ and define the interpretations in $\mcU$ as above.
It is clear that $\mcU$ is strongly $l$-colourable, but not strongly $(l-1)$-colourable.
\hfill $\square$

\begin{notation}\label{notation for U etcetera}{\rm
(i) Let, according to Lemma~\ref{finding the structure U},
$\mcU$ be an $l$-colourable, but not $(l-1)$-colourable, structure such that
$\left\|\mcU\right\|$ is divisible by $l$ and every partition of $|\mcU|$ into $l$
parts of equal size gives rise to an $l$-colouring of $\mcU$.
Let the universe of $\mcU$ be $U = \{1, \ldots, u\}$, so $u \geq l$.\\
(ii) Let $\chi_{\mcU}(x_1, \ldots, x_u)$ be a quantifier-free $L_{rel}$-formula which expresses 
the isomorphism type of $\mcU$.\\
(iii) Let $\widehat{\mcU} \in \mbK$ be an expansion of $\mcU$, that is, $\widehat{\mcU}$
is an $l$-colouring of $\mcU$. Without loss of generality we may assume that
the elements $1, \ldots, l \in U$ have different colours in $\widehat{\mcU}$.\\
(iv) Let $I$ be the set of all unordered pairs $\{i,j\} \subseteq U$ such that 
$i$ and $j$ have the same colour in $\widehat{\mcU}$, and let
$\zeta(x_1, \ldots, x_u)$ denote the formula
$$\chi_{\mcU}(x_1, \ldots, x_u) \ \wedge \
\bigwedge_{\{i,j\} \in I}\xi(x_i, x_j) \ \wedge \ 
\bigwedge_{\{i,j\} \notin I} \neg \, \xi(x_i, x_j).$$
}\end{notation}

\begin{lem}\label{auxiliary properties}
(i) Suppose that $k \geq \max(\left\|\mcS\right\|, \left\|\mcU\right\|)$ and $\mcM \in \mbX_{n,k}$. 
Then
\begin{align*}
&\mcM \models \exists x_1, \ldots, x_u \zeta(x_1, \ldots, x_u), \\
&\mcM \models \forall x_1, \ldots, x_u \Big(\zeta(x_1, \ldots, x_u) \rightarrow
\bigwedge_{i < j \leq l}\neg \, \xi(x_i, x_j)\Big),\\
&\mcM \models \forall y, x_1, \ldots, x_u \Big( \zeta(x_1, \ldots, x_u) \rightarrow
\bigvee_{i=1}^l \xi(x_i,y) \Big), \\
&\mcM \models \forall y \xi(y, y) \ \wedge \
\forall y_1, y_2 \big( \xi(y_1, y_2) \ \rightarrow \ \xi(y_2, y_1)\big), \text{ and}\\
&\mcM \models \forall y_1, y_2, y_3 \big( [\xi(y_1, y_2) \ \wedge \ \xi(y_2, y_3)] \ \rightarrow \
\xi(y_1, y_3)\big).
\end{align*}
(ii) If $\psi$ is any one of the sentences in part (i), 
then $\lim_{n \to \infty} \delta_n^{\mbC}(\psi) = 1.$
\end{lem}

\noindent
{\em Proof.}
(i) Recall Notation~\ref{notation for U etcetera} (iii).
Since $k \geq \left\|\mcU\right\| = \left\|\widehat{\mcU}\right\|$ and $\mcM \in \mbX_{n,k}$,
the $\widehat{\mcU}/\es$-extension axiom is satisfied in $\mcM$, so
there are $m_1, \ldots, m_u \in M$ such that the map $m_i \mapsto i$ is an isomorphism
from $\mcM \uhrc \{m_1, \ldots, m_u\}$ to $\widehat{\mcU}$, so in particular it
preserves the colours. As $\mcM \in \mbX_{n,k}$ and $k \geq \left\|\mcS\right\|$,
Lemma~\ref{the property of xi} implies that 
$\mcM \models \zeta(m_1, \ldots, m_u)$.

Now suppose that $b, a_1, \ldots, a_u \in M$ and $\mcM \models \zeta(a_1, \ldots, a_u)$,
that is,
$$\mcM \models \chi_{\mcU}(a_1, \ldots, a_u) \ \wedge \
\bigwedge_{\{i,j\} \in I}\xi(a_i, a_j) \ \wedge \ 
\bigwedge_{\{i,j\} \notin I} \neg \, \xi(a_i, a_j).$$
Together with the definition of $\mcU$ and $I$ (Notation~\ref{notation for U etcetera}),
this implies that if $i < j \leq l$, then $\mcM \models \neg \, \xi(a_i, a_j)$.
Since $\mcM \in \mbX_{n,k}$ and $k \geq \left\|\mcS\right\|$, Lemma~\ref{the property of xi}
implies that if $i < j \leq l$, then $a_i$ and $a_j$ have different colours.
Since there are only $l$ colours, there is $i \leq l$ such that $b$ has the same colour
as $a_i$ in $\mcM$. By Lemma~\ref{the property of xi} again, $\mcM \models \xi(a_i, b)$.
So we have proved that
\begin{align*}
&\mcM \models \forall x_1, \ldots, x_u \Big(\zeta(x_1, \ldots, x_u) \rightarrow
\bigwedge_{i < j \leq l}\neg \, \xi(x_i, x_j)\Big), \ \text{ and}\\
&\mcM \models \forall y, x_1, \ldots, x_u \Big( \zeta(x_1, \ldots, x_u) \rightarrow
\bigvee_{i=1}^l \xi(x_i,y) \Big).
\end{align*}
The relation `$y$ has the same colour as $z$' is an equivalence relation
which under the given conditions is defined by $\xi(y,z)$ (by Lemma~\ref{the property of xi}).
This immediately implies the rest of part (i).

(ii) Let $\psi$ be any one of the sentences in part (i). Since $\psi \in L_{rel}$ we have
$$\{\mcM \in \mbC_n : \mcM \models \psi \} = \{\mcN \uhrc L_{rel} : \mcN \in \mbK_n \text{ and }
\mcN \models \psi\},$$
so by the definition of $\delta_n^{\mbC}$ we get
$\delta_n^{\mbC}(\psi) = \delta_n^{\mbK}(\psi)$, for every $n$.
Therefore it suffices to show that
$\lim_{n \to \infty} \delta_n^{\mbK}(\psi) = 1$.
Take $k \geq \max(\left\|\mcS\right\|, \left\|\mcU\right\|)$.
By Lemma~\ref{probability of X approaches 1},
$\lim_{n\to\infty}\delta^{\mbK}_n(\mbX_{n,k}) = 1$.
By part (i), for every $n$ and every 
$\mcM \in \mbX_{n,k}$, $\mcM$ satisfies $\psi$, 
so $\delta_n^{\mbK}(\psi) \geq \delta_n^{\mbK}(\mbX_{n,k}) \to 1$, as $n \to \infty$.
\hfill $\square$
\\

\noindent
Next, we define `$l$-colour compatible extension axioms'.
Suppose that $\mcB$ is $l$-colourable (and finite) and let $\mcA \subset \mcB$. 
Without loss of generality we assume that  $A = \{1, \ldots, \alpha\}$ and $B = \{1, \ldots, \beta\}$,
so $\alpha < \beta$.
Let $\chi_{\mcA}(x_1, \ldots, x_{\alpha})$ and $\chi_{\mcB}(x_1, \ldots, x_{\beta})$ 
be quantifier-free $L_{rel}$-formulas which express the isomorphism types of $\mcA$ and $\mcB$, respectively;
so for any $L_{rel}$-structure $\mcM$, $\mcM \models \chi_{\mcA}(m_1, \ldots, m_{\alpha})$
if and only if the map $m_i \mapsto i$ is an isomorphism from $\mcM \uhrc \{m_1, \ldots, m_{\alpha}\}$
to $\mcA$; and similarly for $\chi_{\mcB}$.
Let us say that an $l$-colouring $\gamma : \{1, \ldots, \alpha\} \to \{1, \ldots, l\}$
of $\mcA$ is a {\bf \em $\mcB$-good colouring} if it can be extended
to an $l$-colouring $\gamma' : \{1, \ldots, \beta\} \to \{1, \ldots, l\}$ of $\mcB$
(i.e. $\gamma' \uhrc A = \gamma$).
Let $\gamma : \{1, \ldots, \alpha\} \to \{1, \ldots, l\}$ be a $\mcB$-good colouring of $\mcA$
and let $\gamma': \{1, \ldots, \beta\} \to \{1, \ldots, l\}$ be any colouring of $\mcB$ that 
extends $\gamma$.
Let $\tau$ be any permutation of $\{1, \ldots, l\}$.
The idea in what follows is that, for $j \in \{1, \ldots, \alpha\}$,
the colour of $j$ is associated with the colour
of the element which will be substituted for the variable $x_{\tau\gamma(j)}$,
where $\tau\gamma(j) = \tau(\gamma(j))$.
Let $\theta_{\gamma,\tau}(x_1, \ldots, x_l, y_1, \ldots, y_{\alpha})$ be the conjunction of all
$\xi(x_{\tau\gamma(j)}, y_j)$ where $j \in \{1, \ldots, \alpha\}$.
Similarly, let $\theta_{\gamma',\tau}(x_1, \ldots, x_l, y_1, \ldots, y_{\beta})$ be the conjunction of
all $\xi(x_{\tau\gamma'(j)}, y_j)$ where $j \in \{1, \ldots, \beta\}$.
We call the following sentence an 
{\bf \em instance of the $l$-colour compatible $\mcB/\mcA$-extension axiom}:
\begin{align*}
\forall x_1, \ldots, x_u, y_1, \ldots, y_{\alpha} &\exists y_{\alpha + 1}, \ldots, y_{\beta} \big(\\
        \big[\zeta(x_1, \ldots, x_u) \ \wedge \ &\chi_{\mcA}(y_1, \ldots, y_{\alpha}) 
        \ \wedge \ \theta_{\gamma,\tau}(x_1, \ldots, x_l, y_1, \ldots, y_{\alpha})\big] \ \longrightarrow\\
&\big[\chi_{\mcB}(y_1, \ldots, y_{\beta}) \ \wedge \ 
\theta_{\gamma',\tau}(x_1, \ldots, x_l, y_1, \ldots, y_{\beta})\big] \big).
\end{align*}
In the special case that $A = \es$ and $\gamma'$ is an arbitrary $l$-colouring of $\mcB$,
the above formula should be interpreted as 
$$\forall x_1, \ldots, x_u \exists y_1, \ldots, y_{\beta} \big(
\zeta(x_1, \ldots, x_u) \ \longrightarrow \
\chi_{\mcB}(y_1, \ldots, y_{\beta}) \wedge \theta_{\gamma',\tau}(x_1, \ldots, x_l, y_1, \ldots, y_{\beta}) \big).$$
Since there are only finitely many $l$-colourings of any finite structure, there are only finitely
many instances of the $l$-colour compatible $\mcB/\mcA$-extension axiom.
The {\bf \em $l$-colour compatible $\mcB/\mcA$-extension axiom} is, by definition,
the conjunction of all instances of the $l$-colour compatible $\mcB/\mcA$-extension axiom.
If $|\mcB| \leq k+1$ then the $l$-colour compatible $\mcB/\mcA$-extension axiom is also called an 
{\bf \em $l$-colour compatible $k$-extension axiom}.

\begin{lem}\label{lemma on colour compatible extension axioms}
Suppose that $\mcB$ is $l$-colourable (and finite) and let $\mcA \subset \mcB$.
Let $\varphi$ denote the $l$-colour compatible $\mcB/\mcA$-extension axiom.
If $k \geq \max(\left\|\mcS\right\|, \left\|\mcB\right\|)$ and
$\mcM \in \mbX_{n,k}$, 
then $\mcM \models \varphi$.
\end{lem}

\noindent
{\em Proof.}
Let $\mcA$, $\mcB$, $\varphi$, $k$ and $\mcM$ satisfy the premisses of the lemma,
so in particular $\mcM \in \mbX_{n,k} \subseteq \mbK_n$.
We consider only the case when $\left\|\mcA\right\| \geq 1$, since the
case when $\left\|\mcA\right\| = 0$ is analogous.
Without loss of generality we assume that $A = \{1, \ldots, \alpha\}$ and
$B = \{1, \ldots, \beta\}$ where $\alpha < \beta$.
It suffices to prove that every instance of the $l$-colour compatible $\mcB/\mcA$-extension
axiom is true in $\mcM$.

Let $\gamma : \{1, \ldots, \alpha\} \to \{1, \ldots, l\}$ be a $\mcB$-good $l$-colouring of $\mcA$ and
let $\gamma' : \{1, \ldots, \beta\} \to \{1, \ldots, l\}$ be an $l$-colouring of $\mcB$
which extends $\gamma$. 
Also, let $\tau$ be a permutation of $\{1, \ldots, l\}$.
We prove that $\mcM$ satisfies the following instance of the 
$l$-colour compatible $\mcB/\mcA$-extension axiom, where $\theta_{\gamma,\tau}$
is the conjunction of all $\xi(x_{\tau\gamma(j)}, y_j)$
where $j \in \{1, \ldots, \alpha\}$, and
$\theta_{\gamma',\tau}$ is the conjunction of all $\xi(x_{\tau\gamma(j)}, y_j)$
where $j \in \{1, \ldots, \beta\}$:
\begin{align*}
\forall x_1, \ldots, x_u, y_1, \ldots, y_{\alpha} &\exists y_{\alpha + 1}, \ldots, y_{\beta} \big(\\
        \big[\zeta(x_1, \ldots, x_u) \ \wedge \ &\chi_{\mcA}(y_1, \ldots, y_{\alpha}) 
        \ \wedge \ \theta_{\gamma,\tau}(x_1, \ldots, x_l, y_1, \ldots, y_{\alpha})\big] \ \longrightarrow\\
&\big[\chi_{\mcB}(y_1, \ldots, y_{\beta}) \ \wedge \ 
\theta_{\gamma',\tau}(x_1, \ldots, x_l, y_1, \ldots, y_{\beta})\big] \big).
\end{align*}
Note that since $\mcM \in \mbX_{n,k}$ and 
$k \geq \max(\left\|\mcS\right\|, \left\|\mcB\right\|)$, 
we can, and will repeatedly, use Lemma~\ref{the property of xi} which implies that
for all $a, b \in M$, $\mcM \models \xi(a,b)$ if and only if $a$ and $b$ have the
same colour in $\mcM$.

Suppose that 
$$\mcM \models \zeta(m_1, \ldots, m_u) \ \wedge \ \chi_{\mcA}(a_1, \ldots, a_{\alpha}) 
\ \wedge \ \theta_{\gamma,\tau}(m_1, \ldots, m_l, a_1, \ldots, a_{\alpha}).$$
By the definition of $\zeta$ (Notation~\ref{notation for U etcetera} (iii), (iv)),
if $i,j \leq l$ and $i \neq j$, then $m_i$ and $m_j$ have different colours.
Hence, there is a permutation $\pi$ of $\{1, \ldots, l\}$ such that, for every $i \in \{1, \ldots, l\}$,
$m_i$ has colour $\pi(i)$, i.e. $\mcM \models P_{\pi(i)}(m_i)$.
Let $\widehat{\mcB} \in \mbK$ be the expansion of $\mcB$ such that, 
$$\text{for every } j \in \{1, \ldots, \beta\}, \ 
\widehat{\mcB} \models P_{\pi\tau\gamma'(j)}(j).$$
In other words, $j \in \{1, \ldots, \beta\}$ gets the same colour
in $\widehat{\mcB}$ as $m_{\tau\gamma'(j)}$ in $\mcM$, and this colour is 
$\pi\tau\gamma'(j)$.
In particular, this holds whenever $j \leq \alpha$ and $\gamma'$ is replaced by $\gamma$.
Since we assume that
$$\mcM \models \chi_{\mcA}(a_1, \ldots, a_{\alpha}) 
\ \wedge \ \theta_{\gamma,\tau}(m_1, \ldots, m_l, a_1, \ldots, a_{\alpha})$$
it follows that the map $j \mapsto a_j$, for $j \in \{1, \ldots, \alpha\}$,
is an isomorphism from $\widehat{\mcA} = \widehat{\mcB} \uhrc \{1, \ldots, \alpha\}$ 
to $\mcM \uhrc \{a_1, \ldots, a_{\alpha}\}$.
Since $\mcM \in \mbX_{n,k}$ and $k$ is sufficiently large, 
$\mcM$ satisfies the $\widehat{\mcB}/\widehat{\mcA}$-extension axiom.
Hence, there are $a_{\alpha + 1}, \ldots, a_{\beta} \in M$ such that
the map $j \mapsto a_j$, for $j \in \{1, \ldots, \beta\}$, is an
isomorphism from $\widehat{\mcB}$ to $\mcM \uhrc \{a_1, \ldots, a_{\beta}\}$.
This implies that $\mcM \models \chi_{\mcB}(a_1, \ldots, a_{\beta})$ and that, 
for all $j \in \{1, \ldots, \beta\}$,
$\mcM \models P_{\pi\tau\gamma'(a_j)}(a_j)$,
which means that $a_j$ has the same colour as $m_{\tau\gamma'(j)}$, so
$\mcM \models \xi(m_{\tau\gamma'(j)}, a_j)$.
Hence 
$$\mcM \models \theta_{\gamma', \tau}(m_1, \ldots, m_l, a_1, \ldots, a_{\beta}),$$
and we are done.
\hfill $\square$

\begin{cor}\label{extension axioms almost surely hold in l-colourable structures}
For every $l$-colour compatible extension axiom $\varphi$, 
$\lim_{n \to \infty} \delta_n^{\mbC}(\varphi) = 1.$
\end{cor}

\noindent
{\em Proof.}
Let $\varphi$ be an $l$-colour compatible extension axiom.
Since $\varphi \in L_{rel}$ we have
$$\{\mcM \in \mbC_n : \mcM \models \varphi \} = \{\mcN \uhrc L_{rel} : \mcN \in \mbK_n \text{ and }
\mcN \models \varphi\},$$
so by the definition of $\delta_n^{\mcC}$ we get
$\delta_n^{\mbC}(\varphi) = \delta_n^{\mbK}(\varphi)$, for every $n$.
Therefore it suffices to show that
$\lim_{n \to \infty} \delta_n^{\mbK}(\varphi) = 1$.
For some $l$-colourable $L_{rel}$-structures $\mcA \subset \mcB$,
$\varphi$ is the $l$-colour compatible $\mcB/\mcA$-extension axiom.
Take $k \geq \max(\left\|\mcS\right\|, \left\|\mcB\right\|)$.
By Lemma~\ref{probability of X approaches 1},
$\lim_{n\to\infty}\delta^{\mbK}_n(\mbX_{n,k}) = 1$.
By Lemma~\ref{lemma on colour compatible extension axioms}, for every $n$ and every 
$\mcM \in \mbX_{n,k}$, $\mcM$ satisfies $\varphi$, 
so $\delta_n^{\mbK}(\varphi) \geq \delta_n^{\mbK}(\mbX_{n,k}) \to 1$, as $n \to \infty$.
\hfill $\square$
\\

\noindent
For every integer $n > 0$ let $\mcM_{(n,1)}, \ldots, \mcM_{(n,m_n)}$ be an enumeration of all isomorphism types
of $l$-colourable structures of cardinality at most $n$.
Let $\chi^n_i(x_1, \ldots, x_n)$ describe the isomorphism type of $\mcM_{(n,i)}$ in such a way that
we require that all variables $x_1, \ldots, x_n$ actually occur in $\chi^n_i$.
It means that if $\left\| \mcM_{(n,i)} \right\| < n$, then 
$\chi^n_i(x_1, \ldots, x_n)$ must express that
some variables refer to the same element, by saying `$x_k = x_l$' for some $k \neq l$.
For every $n \in \mbbN$ let $\psi_n$ denote the sentence
$$\forall x_1, \ldots, x_n 
\bigvee_{i=1}^{m_n}
\bigvee_{\pi} \chi^n_i(x_{\pi(1)}, \ldots, x_{\pi(n)}),$$
where the second disjunction ranges over all permutations $\pi$ of $\{1, \ldots, n\}$.
Then let $T_{iso} = \{\psi_n : n \in \mbbN, \ n > 0\}$
and note that every $\psi_n$ is true in every $l$-colourable structure.
Let $T_{ext}$ consist of all $l$-colour compatible extension axioms and let 
$T_{col}$ consist of the sentences appearing in part (i) of 
Lemma~\ref{auxiliary properties}.
Finally, let $T_{\mbC} = T_{iso} \cup T_{ext} \cup T_{col}$.
By part (ii) of Lemma~\ref{auxiliary properties},
Corollary~\ref{extension axioms almost surely hold in l-colourable structures}
and compactness, $T_{\mbC}$ is consistent.
Since $T_{ext} \subset T_{\mbC}$, every model of $T_{\mbC}$ is infinite.
In order to prove Theorem~\ref{0-1 law for l-colourable structures}
it is enough to prove that $T_{\mbC}$ is complete.

\begin{lem}\label{countable categoricity}
$T_{\mbC}$ is countably categorical and therefore complete.
\end{lem}

\noindent
{\em Proof.}
Suppose that the $L_{rel}$-structures $\mcM$ and $\mcM'$ are countable models of $T_{\mbC}$.
We will prove that $\mcM \cong \mcM'$ by a back and forth argument, but first
we need some preparation.
Recall that $T_{col} \subset T_{\mbC}$ and that 
$T_{col}$ contains the formulas that appear in part (i) of 
Lemma~\ref{auxiliary properties}.
Therefore there are $m_1, \ldots, m_u \in M$ and 
$m'_1, \ldots, m'_u \in M'$ such that 
$\mcM \models \zeta(m_1, \ldots, m_u)$ and 
$\mcM' \models \zeta(m'_1, \ldots, m_u)$.
Moreover, because of the sentences in $T_{col}$, the following hold:
\begin{itemize}
	\item[\textbullet] $\xi(y,z)$ defines an equivalence relation $R_M$ on $M$ and an
	equivalence relation $R_{M'}$ on $M'$.
	\item[\textbullet] The elements $m_1, \ldots, m_l$ belong to different equivalence classes;
	the elements $m'_1, \ldots, m'_l$ belong to different equivalence classes.
	\item[\textbullet] Every element in $M$ is equivalent to one of $m_1, \ldots, m_l$, 
	so $R_M$ has exactly $l$ equivalence classes; and the same is true for 
	$m'_1, \ldots, m'_l$ and $R_{M'}$.
\end{itemize}
We prove that $\mcM \cong \mcM'$ by a back and forth argument
in which partial isomorphisms between $\mcM$ and $\mcM'$ are extended step by step.
It suffices to prove the following:
\\

\noindent
{\bf Claim.}
{\em 
Suppose that $\mcA$ and $\mcA'$ are finite substructures of 
$\mcM$ and $\mcM'$, respectively, and that $f$ is an isomorphism from 
$\mcA$ to $\mcA'$ such that for all $a \in A$ and all $i \in \{1, \ldots, l\}$, 
$\mcM \models \xi(m_i, a)$ $\Longleftrightarrow$ $\mcM' \models \xi(m'_i, f(a))$.
For every $b \in M - A$ (or $b' \in M' - A'$), 
there are $b' \in M' - A'$ (or $b \in M - A$)
and an isomorphism $g : \mcM \uhrc A \cup \{b\} \to \mcM' \uhrc A' \cup \{b'\}$
such that $g$ extends $f$ (so $g(b) = b'$) and, for every $i \in \{1, \ldots, l\}$,
$\mcM \models \xi(m_i, b)$ $\Longleftrightarrow$ $\mcM' \models \xi(m'_i, b')$.
}
\\

\noindent
Suppose that $\mcA$ and $\mcA'$ are finite substructures of 
$\mcM$ and $\mcM'$, respectively, and that $f$ is an isomorphism from 
$\mcA$ to $\mcA'$ such that for all $a \in A$ and all $i \in \{1, \ldots, l\}$,
$\mcM \models \xi(m_i, a)$ $\Longleftrightarrow$ $\mcM' \models \xi(m'_i, f(a))$.
Let $A = \{a_1, \ldots, a_{\alpha}\}$,
let $b = a_{\alpha + 1} \in M - A$ and let $\mcB = \mcM \uhrc \{a_1, \ldots, a_{\alpha + 1}\}$.
We will find the required $b' \in M' - A'$ by defining a suitable
instance of the $l$-colour compatible $\mcB/\mcA$-extension axiom and then 
use the assumption that $\mcM' \models T_{\mbC} \supset T_{ext}$.

Let 
$$ X = \{m_1, \ldots, m_u, a_1, \ldots, a_{\alpha + 1}\}.$$
Since $\xi(y,z)$ is an existential formula (see Notation~\ref{notation for S and xi} (iii)),
there is a {\em finite} substructure $\mcN \subset \mcM$ such that
$$X \subseteq N, \ \text{ and whenever } \ c, d \in X, \
\text{ then } \ \mcM \models \xi(c,d) \ \Longleftrightarrow \ \mcN \models \xi(c,d).$$
Since $N$ is finite and $\mcN \subset \mcM \models T_{\mbC} \supset T_{iso}$,
$\mcN$ is $l$-colourable.
Let $\gamma^* : N \to \{1, \ldots, l\}$ be an $l$-colouring of $\mcN$,
and define an equivalence relation $\sim^*$ on $N$ by:
$$c \sim^* d \Longleftrightarrow \gamma^*(c) = \gamma^*(d).$$
By the choice of $\mcN$ and Lemma~\ref{the property of xi}~(i), 
for all $c, d \in X$,
$$R_M(c,d) \Longleftrightarrow \mcM \models \xi(c,d)
\Longleftrightarrow \mcN \models \xi(c,d) \Longrightarrow \gamma^*(c) = \gamma^*(d)
\Longleftrightarrow c \sim^* d.$$
This means that the restriction of $R_M$ to $X$ is a refinement of the restriction of $\sim^*$ to $X$.
We have already observed that the restriction of $R_M$ to $X$ has exactly $l$ equivalence classes,
because all $m_1, \ldots, m_l$ belong to different classes. 
Moreover, since $\mcM \models \zeta(m_1, \ldots, m_u)$ we get,
by the definition of $\zeta$, $\mcM \models \chi_{\mcU}(m_1, \ldots, m_u)$,
and since $\mcU$ is $l$-colourable, but not $(l-1)$-colourable, it follows
that $\sim^*$ has exactly $l$ equivalence classes.
It follows that the restriction of $R_M$ to $X$ is the {\em same} relation as the restriction
of $\sim^*$ to $X$. Hence,
$$\text{for all } \ c, d \in X, \ \mcM \models \xi(c,d) \ \Longleftrightarrow \ 
\gamma^*(c) = \gamma^*(d).$$
Therefore, there is a permutation $\tau$ of $\{1, \ldots, l\}$ such that, 
$$\text{for every } \ j \in \{1, \ldots, \alpha + 1\}, \
\mcM \models \xi(m_{\tau\gamma^*(a_j)}, a_j).$$
Let $\gamma' = \gamma^* \uhrc \{a_1, \ldots, a_{\alpha + 1}\}$ and
$\gamma = \gamma^* \uhrc \{a_1, \ldots, a_{\alpha}\}$.
Then let 
$$\theta_{\gamma', \tau}(x_1, \ldots, x_l, y_1, \ldots, y_{\alpha + 1})$$
be the conjunction of all $\xi(x_{\tau\gamma'(a_j)}, y_j)$ 
where $j \in \{1, \ldots, \alpha + 1\}$, and let
$$\theta_{\gamma, \tau}(x_1, \ldots, x_l, y_1, \ldots, y_{\alpha})$$
be the conjunction of all $\xi(x_{\tau\gamma(a_j)}, y_j)$ 
where $j \in \{1, \ldots, \alpha\}$.
Let $\chi_{\mcA}(y_1, \ldots, y_{\alpha})$ and 
$\chi_{\mcB}(y_1, \ldots, y_{\alpha + 1})$ be quantifier-free formulas
which describe the isomorphism types of $\mcA$ and $\mcB$, respectively.
Now the following is an instance of the $l$-colour compatible
$\mcB/\mcA$-extension axiom:
\begin{align*}
\forall x_1, \ldots, x_u, y_1, \ldots, y_{\alpha} &\exists y_{\alpha + 1} \big(\\
        \big[\zeta(x_1, \ldots, x_u) \ \wedge \ &\chi_{\mcA}(y_1, \ldots, y_{\alpha}) 
        \ \wedge \ \theta_{\gamma,\tau}(x_1, \ldots, x_l, y_1, \ldots, y_{\alpha})\big] \ \longrightarrow\\
&\big[\chi_{\mcB}(y_1, \ldots, y_{\alpha + 1}) \ \wedge \ 
\theta_{\gamma',\tau}(x_1, \ldots, x_l, y_1, \ldots, y_{\alpha + 1})\big] \big).
\end{align*}
Since
$$\mcM \models \zeta(m_1, \ldots, m_u) \ \wedge \
\chi_{\mcA}(a_1, \ldots, a_{\alpha}) \ \wedge \
\theta_{\gamma,\tau}(m_1, \ldots, m_l, a_1, \ldots, a_{\alpha}),$$
it follows from the assumptions that
$$\mcM' \models \zeta(m'_1, \ldots, m'_u) \ \wedge \
\chi_{\mcA}(f(a_1), \ldots, f(a_{\alpha})) \ \wedge \
\theta_{\gamma,\tau}(m'_1, \ldots, m'_l, f(a_1), \ldots, f(a_{\alpha})).$$
Since $\mcM'$ satisfies all $l$-colour compatible extension axioms it follows
that there is $b' \in M' - A'$ such that if $g(a_{\alpha+1}) = b'$ and
$g(a) = f(a)$ for all $a \in A$, then
$$\mcM' \models \chi_{\mcB}(g(a_1), \ldots, g(a_{\alpha+1})) \ \wedge \
\theta_{\gamma',\tau}(m'_1, \ldots, m'_l, g(a_1), \ldots, g(a_{\alpha + 1})).$$
It follows that $g$ is an isomorphism from $\mcM \uhrc A \cup \{b\}$
to $\mcM' \uhrc A' \cup \{b'\}$.
Since $\mcM \models \xi(m'_i, b')$ for a unique $i \in \{1, \ldots, l\}$
(by the conclusions in the beginning of the proof),
it also follows that, for every $i \in \{1, \ldots, l\}$,
$\mcM' \models \xi(m'_i, b') \ \Longleftrightarrow \ 
\mcM \models \xi(m_i, b)$, where $b = a_{\alpha+1}$.

Note that the argument also works for the `base case' when $A = A'=\es$ and $f$ is the
empty map; the difference is merely notational. 
If we start out with $b' \in M' - A'$, then we argue symmetrically.
Thus the claim, and hence the lemma, is proved.
\hfill $\square$
\\

\noindent
By the preceeding lemmas, $T_{\mbC}$ is a complete theory such that
whenever $m \in \mbbN$ and $\psi_1, \ldots, \psi_m  \in T_{\mbC}$, 
then $\lim_{n\to\infty}\delta_n^{\mbC}\Big(\bigwedge_{i=1}^m \psi_i\Big) = 1$.
By compactness and completeness it follows that if $T_{\mbC} \models \varphi$, then
$\lim_{n\to\infty}\delta_n^{\mbC}(\varphi) = 1$, and if
$T_{\mbC} \not\models \varphi$, then 
$\lim_{n\to\infty}\delta_n^{\mbC}(\varphi) = 0$.

\subsection{Relationship between the dimension conditional measure and the uniform
measure}\label{relationship between the dimension conditional measure and the uniform measure}

In this section we prove that the results that we have seen for the probability 
measures $\delta_n^{\mbC}$ and $\delta_n^{\mbS}$ transfer to the uniform 
probability measures on $\mbC_n$ and $\mbS_n$, respectively, if one condition about
(strongly) $l$-colourable structures holds. 
In Section~\ref{the uniform measure and l-colourable structures} we prove
that this condition does indeed hold.

\begin{defin}\label{definition of m-rich colouring}{\rm 
Let $m \in \mbbR$.\\
(i) Suppose that $\gamma : S \to \{1, \ldots, l\}$ is a function.
We say that $\gamma$ is {\bf \em $m$-rich} if, for every $i \in \{1, \ldots, l\}$,
$|\gamma^{-1}(i)| \geq m$, that is, at least $m$ members of $S$ are mapped to $i$.\\
(ii) We call  $\mcM \in \mbK$ (or $\mcM \uhrc L_{col}$)  
{\bf \em $m$-richly $l$-coloured} if for every $i \in \{1, \ldots, l\}$,
$\big|\{ a \in M : \mcM \models P_i(a) \}\big| \geq m$.\\
(iii) We also call $\mcM \in \mbSK$ (or $\mcM \uhrc L_{col}$) 
{\bf \em $m$-richly $l$-coloured} if for every $i \in \{1, \ldots, l\}$,
$\big|\{ a \in M : \mcM \models P_i(a) \}\big| \geq m$.
(If $\mcM \in \mbSK$ then it is understood that we are dealing with strong
$l$-colourings, although it was not explicitly reflected in the terminology defined.)
}\end{defin}

\noindent
Recall the notion of an {\em $l$-colour compatible extension axiom}, defined in
Section~\ref{proof of 0-1 law for l-colourable structures}, before 
Lemma~\ref{extension axioms almost surely hold in l-colourable structures}.
These axioms are essentially the same whether we consider (not necessarily strong)
$l$-colourings, or strong $l$-colourings. The only difference is that the 
structures $\mcS$ and $\mcU$ which are implicitly refered to (via the formulas 
$\zeta$ and $\xi$) are different in the two cases.

\begin{theor}\label{transfer of 0-1 laws to uniform measure}
Let $f: \mbbN \to \mbb{R}$ be such that $f(n) / \ln n \to \infty$ as $n \to \infty$.\\
(i) Suppose that the proportion of $\mcM \in \mbK_n$ which are $f(n)$-richly $l$-coloured
approaches 1 as $n \to \infty$.
Then, for every extension axiom $\varphi$ of $\mbK$, the proportion of $\mcM \in \mbK_n$
which satisfy $\varphi$ approaches 1 as $n \to \infty$.
Consequently, $\mbK$ has a 0-1 law for the uniform probability measure.\\
(ii) Suppose that the proportion of $\mcM \in \mbC_n$ which have an 
$f(n)$-rich $l$-colouring approaches 1 as $n \to \infty$.
Then, for every $l$-colour compatible extension axiom $\varphi$ of $\mbC$,
the proportion of $\mcM \in \mbC_n$ which satisfy $\varphi$ approaches 1 as $n \to \infty$.
Moreover, $\mbC$ has a 0-1 law for the uniform probability measure.\\
(iii) Parts (i) and (ii) hold if $\mbK_n$ and $\mbC_n$ are replaced by
$\mbSK_n$ and $\mbS_n$, repectively, and `strong' is added before 
`$l$-colouring'.
\end{theor}

\begin{rem}\label{remark concerning transfer theorem}{\rm
In Section~\ref{the uniform measure and l-colourable structures}
we will prove (Theorem~\ref{new theorem}) 
that there is a constant $\mu > 0$ such that
the proportion of $\mcM \in \mbC_n$ (or $\mcM \in \mbS_n$)
which have a $\mu n$-rich (strong) $l$-colouring approaches 1 as $n \to \infty$.
It follows (Remark~\ref{remark about how to compute mu etc}) that the proportion
of $\mcM \in \mbK_n$ (or $\mcM \in \mbSK_n$) which are 
$\mu n$-richly $l$-coloured approaches 1 as $n \to \infty$.
}\end{rem}

\subsection{Proof of Theorem~\ref{transfer of 0-1 laws to uniform measure}}

\noindent
The proof is exactly the same whether we consider (not necessarily strongly)
$l$-colourable (or $l$-coloured) structures or {\em strongly} $l$-colourable (or $l$-coloured) structures.
This is because we only need to use 
(in Lemma~\ref{the number of different colourings when extension axioms are satisfied} below) 
the properties of the structure $\mcU$ and the formulas $\zeta$ and $\xi$, from
Lemmas~\ref{finding the structure U} and~\ref{the property of xi}, and not
their precise definitions in the respective case.
Therefore we will speak only about $\mbK_n$, $\mbC_n$, $l$-colourings and 
$l$-colourable (or $l$-coloured) structures; the proof in the case of strong $l$-colourings 
is obtained by making the obvious changes of notation and terminology.
Throughout the proof, $l \geq 2$ is fixed so we may occasionally say `colouring' instead of `$l$-colouring'.

Suppose that $f: \mbbN \to \mbb{R}$ is such that $f(n) / \ln n \to \infty$ as $n \to \infty$.
A straightforward consequence is that for every $k \in \mbbN$ and every $0 < \alpha < 1$,
$\lim_{n \to \infty}n^k \cdot \alpha^{f(n)} = 0$.
For if $\beta = 1/\alpha$, then $\beta > 1$, $\ln \beta > 0$ and
$\ln \frac{\beta^{f(n)}}{n^k} = \big(f(n)\ln \beta - k\ln n \big) \to \infty$ as $n \to \infty$;
which gives $\lim_{n\to\infty} \frac{\beta^{f(n)}}{n^k} = \infty$,
and $n^k \cdot \alpha^{f(n)} = \frac{n^k}{\beta^{f(n)}} \to 0$ as $n \to \infty$.

We first prove (i).
As said in Remark~\ref{remark about extension axioms and 0-1 laws},
the zero-one law for $\mbK$, with the uniform measure, follows if 
we can show that for every extension axiom of $\mbK$, the proportion of structures
in $\mbK_n$ which satisfy it approaches 1 as $n \to \infty$.
Suppose that the proportion of $\mcM \in \mbK_n$ which are $f(n)$-richly coloured
approaches 1 as $n \to \infty$.

Let $\varphi$ be an extension axiom of $\mbK$.
It suffices to consider the case when $\varphi$ has only one existential
quantifier, so let $\varphi$ have the form 
$$\forall x_1, \ldots, x_k \exists x_{k+1} \big( \psi(x_1, \ldots, x_k) \rightarrow 
\psi'(x_1, \ldots, x_k, x_{k+1})\big),$$
where $\psi$ and $\psi'$ are quantifier-free.

For every $L_{col}$-structure $\mcA$ with universe $\{1, \ldots, n\}$
for some $n$ (or equivalently, for every $\mcA \in \mbK \uhrc 1$), let
$$\mbE_L(\mcA) = \{\mcM \in \mbK : \mcM \uhrc L_{col} = \mcA \}.$$
Since we  assume that the proportion of $\mcM \in \mbK_n$ which are $f(n)$-richly coloured
approaches 1 as $n \to \infty$, it is sufficient to prove that
\begin{itemize}
	\item[(a)] for every $\varepsilon > 0$ 
	there is $n_{\varepsilon}$ such that
	for every $n > n_{\varepsilon}$, if $\mcA \in \mbK_n \uhrc 1$ is an $f(n)$-rich colouring,
	then the proportion of $\mcM \in \mbE_L(\mcA)$ which satisfy $\varphi$
	is at least $1 - \varepsilon$.
\end{itemize}

\noindent
The proof of (a) is a slight variant of the well known proof that, with the uniform measure,
the probability that an extension axiom is true in a randomly picked structure
(without any restrictions on its relations, and with at least one relation with arity $>$ 1)
with universe $\{1, \ldots, n\}$ approaches 1 as $n$ tends to infinity
(see \cite{Fag, EF, Hod}).

Suppose that $\mcA \in \mbK_n \uhrc 1$ is an $f(n)$-rich colouring.
Let $\alpha$ be the number of nonequivalent quantifier-free $L$-formulas with 
free variables (exactly) $x_1, \ldots, x_{k+1}$.
We show that, with the uniform measure, the probability that $\mcM \in \mbE_L(\mcA)$ 
does {\em not} satisfy $\varphi$ approaches 0 as $n \to \infty$; moreover, the convergence
is uniform in the sense that it depends only on $n = \left\| \mcA \right\|$.
From this (a) follows.

Note that the only restriction on the interpretations of relation symbols from $L_{rel}$
in structures in $\mbE_L(\mcA)$ is that the interpretations respect the colouring of $\mcA$.
Suppose that $\bar{a} = (a_1, \ldots, a_k) \in |\mcA|^k$, $\mcM \in \mbE_L(\mcA)$ and 
$\mcM \models \psi(\bar{a})$. Let $a_{k+1} \in |\mcA| - \rng(\bar{a})$ be any of the at least 
$f(n) - k$ elements not in $\rng(\bar{a})$ which have the colour, say $i$, which is specified for $x_{k+1}$ by $\psi'(x_1, \ldots, x_{k+1})$.
Then the probability, with the uniform measure, that, for such $a_{k+1}$, 
$\mcM \models \psi'(a_1, \ldots, a_k, a_{k+1})$
is at least $1/\alpha$. So the probability that this is not true is at most $1 - 1/\alpha$;
and the probability that $\mcM \not\models \psi'(a_1, \ldots, a_k, a)$ for every one
of the at least $f(n) - k$ elements $a$ outside of $\rng(\bar{a})$ 
with colour $i$ is at most $(1 - 1/\alpha)^{f(n)-k}$.
There are $n^k$ choices of $\bar{a} \in |\mcA|^k$ for which $\exists x_{k+1}\psi'(\bar{a},x_{k+1})$
could fail to be true in $\mcM$, so the probability that $\mcM \not\models \varphi$
is at most $n^k \cdot (1 - 1/\alpha)^{f(n)-k} \to 0$ as $n \to \infty$; by the assumption
about $f(n)$. Since we get the same expression `$n^k \cdot (1 - 1/\alpha)^{f(n)-k}$' 
for every $f(n)$-rich colouring $\mcA \in \mbK_n \uhrc 1$ we have proved (a), and hence (i).

Let $\mcS$ and $\mcU$ be the $L_{rel}$-structures
from Notation~\ref{notation for S and xi} and~\ref{notation for U etcetera},
and let $m = \max(\left\|\mcS\right\|, \left\|\mcU\right\|)$.
Also, let $\xi(y,z)$ be the $L_{rel}$-formula from Notation~\ref{notation for S and xi}.
Fix an arbitrary $k \geq m$ and define
\begin{align*}
\mbX^{\mbK}_n &= \{\mcM \in \mbK_n : \mcM \text{ satisfies all $k$-extension axioms of $\mbK$}\},\\
\mbX^{\mbC}_n &= \{\mcM \in \mbC_n : \mcM = \mcN \uhrc L_{rel} \text{ for some } 
\mcN \in \mbX^{\mbK}_n\},\\
\mbY^{\mbK}_n &= \{\mcM \in \mbK_n : \mcM \text{ is $f(n)$-richly coloured}\},\\
\mbY^{\mbC}_n &= \{\mcM \in \mbC_n : \mcM \text{ has an $f(n)$-rich colouring}\}.
\end{align*}

\begin{lem}\label{added lemma}
Every $\mcM \in \mbX^{\mbC}_n$ satisfies all $l$-colour compatible $k$-extension axioms.
\end{lem}

\noindent
{\em Proof.} 
The notation $\mbX^{\mbK}_n$, introduced before the lemma, denotes the same set of structures as the
notation $\mbX_{n,k}$ defined in Notation~\ref{notation for S and xi}~(iv).
Therefore Lemma~\ref{lemma on colour compatible extension axioms} tells that
every $\mcM \in \mbX^{\mbK}_n$ satisfies all $l$-colour compatible $k$-extension axioms.
Since all such axioms are $L_{rel}$-sentences it follows that
for every $\mcM \in \mbX^{\mbK}_n$, $\mcM \uhrc L_{rel}$ satisfies all
$l$-colour compatible $k$-extension axioms.
The lemma now follows from the definition of $\mbX^{\mbC}_n$.
\hfill $\square$
\\

\noindent
From Lemma~\ref{added lemma} it follows that in order to prove that
the proportion of $\mcM \in \mbC_n$ which satisfy all $l$-colour compatible 
$k$-extension axioms approaches 1 as $n \to \infty$, 
it suffices to show that $\big|\mbX^{\mbC}_n\big| \big/ |\mbC_n| \to 1$
as $n \to \infty$.

\begin{lem}\label{the number of different colourings when extension axioms are satisfied}
For all $\mcM \in \mbX^{\mbC}_n$ the following hold:\\
(i) For every $l$-colouring $\gamma : M \to \{1, \ldots, l\}$ of $\mcM$,
and all $a, b \in M$, $\mcM \models \xi(a,b) \Longleftrightarrow \gamma(a) = \gamma(b)$.\\
(ii) $\mcM$ has a unique $l$-colouring up to permutation of the colours.
\end{lem}

\noindent
{\em Proof.}
As in the proof of the previous lemma, recall that $\mbX^{\mbK}_n$ means
the same as $\mbX_{n,k}$ in Section~\ref{l-colourable structures}.
By Lemma~\ref{auxiliary properties} and the definition of $\zeta$
(in Notation~\ref{notation for U etcetera}), for every $\mcM \in \mbX^{\mbK}_n$, the
following hold:
\begin{itemize}
	\item[\textbullet] \ $\mcU$ is embeddable into $\mcM$.
	\item[\textbullet] \ $\xi(y,z)$ defines an equivalence relation on $M$ with
	exactly $l$ equivalence classes.
\end{itemize}
Since $\mcU$ is an $L_{rel}$-structure and $\xi$ an $L_{rel}$-formula, 
it follows from the definition of $\mbX^{\mbC}_n$ that
the above two points hold for every $\mcM \in \mbX^{\mbC}_n$ as well. 

Let $\mcM \in \mbX^{\mbC}_n$ and let $\gamma : M \to \{1, \ldots, l\}$ be an $l$-colouring of $\mcM$.
Since $\mcU$ is embeddable into $\mcM$ and $\mcU$ is not $(l-1)$-colourable,
it follows that the equivalence relation $\gamma(a) = \gamma(b)$ has exactly $l$ equivalence
classes. Observe that the colouring $\gamma$ gives rise to a unique expansion
of $\mcM$ that belongs to $\mbK$. Therefore Lemma~\ref{the property of xi}~(i) implies that
if $\mcM \models \xi(a,b)$ then $\gamma(a) = \gamma(b)$.
Hence, the equivalence relation defined by $\xi(y,z)$ is a refinement of the 
equivalence relation $\gamma(y) = \gamma(z)$.
Since both equivalence relations have exactly $l$ equivalence classes they must
be the same. In other words, for all $a, b \in M$, $\mcM \models \xi(a,b)$
if and only if $\gamma(a) = \gamma(b)$. Hence (i) is proved.
Part (ii) is now immediate, for if $\gamma$ and $\gamma'$ are two $l$-colourings of $\mcM \in \mbX^{\mbC}_n$,
then
$$\gamma(a) = \gamma(b) \ \Longleftrightarrow \
\mcM \models \xi(a,b) \ \Longleftrightarrow \
\gamma'(a) = \gamma'(b).$$
\hfill $\square$ 
\\

\noindent
Now we have the tools to complete the proof of part (ii) of the theorem.
Observe that with the notation used in the proof of part (i) we have 
$$\mbY^{\mbK}_n = \bigcup \big\{ \mbE_L(\mcA) : \mcA \in \mbK_n \uhrc 1 \text{ is 
an $f(n)$-rich $l$-colouring}\big\},$$
and (a) implies that 
\begin{equation*}
\lim_{n \to \infty} \frac{\big|\mbX^{\mbK}_n \cap \mbY^{\mbK}_n\big|}{\big|\mbY^{\mbK}_n\big|} = 1. \tag{b}
\end{equation*}
Note that for every $l$-colouring of $\mcM \in \mbC$, the colours can be permuted
in $l!$ ways. Therefore, 
\begin{equation*}
|\mbK_n| \geq l! |\mbC_n| \quad \text{ and } \quad 
\big|\mbY^{\mbK}_n\big| \geq l!\big|\mbY^{\mbC}_n\big|. \tag{c}
\end{equation*}

\noindent
Lemma~\ref{the number of different colourings when extension axioms are satisfied}
implies that
\begin{equation*}
\big|\mbX^{\mbK}_n\big| = l!\big|\mbX^{\mbC}_n\big| \quad \text{ and } \quad
\big|\mbX^{\mbK}_n \cap \mbY^{\mbK}_n\big| = l!\big|\mbX^{\mbC}_n \cap \mbY^{\mbC}_n\big|. \tag{d}
\end{equation*}

\noindent
Assume that the proportion of $\mcM \in \mbC_n$ which have an $f(n)$-rich colouring
approaches 1 as $n \to \infty$. In other words,
\begin{equation*}
\lim_{n \to \infty} \frac{\big|\mbY^{\mbC}_n\big|}{|\mbC_n|} = 1. \tag{e}
\end{equation*}
By (c) and (d),
\begin{equation*}
\frac{\big|\mbX^{\mbK}_n \cap \mbY^{\mbK}_n\big|}{\big|\mbY^{\mbK}_n\big|} \leq 
\frac{l!\big|\mbX^{\mbC}_n \cap \mbY^{\mbC}_n\big|}{l!\big|\mbY^{\mbC}_n\big|} 
= \frac{\big|\mbX^{\mbC}_n \cap \mbY^{\mbC}_n\big|}{|\mbC_n|} 
\cdot \frac{|\mbC_n|}{\big|\mbY^{\mbC}_n\big|} \leq 1. \tag{f}
\end{equation*}
Now (b), (e) and (f) imply that
\begin{equation*}
\lim_{n \to \infty} \frac{\big|\mbX^{\mbC}_n \cap \mbY^{\mbC}_n\big|}{|\mbC_n|} = 1. \tag{g}
\end{equation*}
By Lemma~\ref{added lemma} and (g), the proportion of $\mcM \in \mbC_n$ which satisfy all
$l$-colour compatible $k$-extension axioms of $\mbC$ approaches 1 as $n$ approaches $\infty$.
This has been derived for arbitrary $k \geq m$, 
under the assumption that the proportion of $\mcM \in \mbC_n$ which have an $f(n)$-rich colouring
approaches 1, as $n \to \infty$.
Since every $l$-colour compatible extension axiom is an
$l$-colour compatible $k$-extension axiom for all sufficiently large $k$, we have proved:
If the proportion of $\mcM \in \mbC_n$ which have an
$f(n)$-rich colouring approaches 1 as $n \to \infty$, 
then for every $l$-colour compatible extension axiom $\varphi$,
the proportion of $\mcM \in \mbC_n$ which satisfy $\varphi$ approaches 1 as $n \to \infty$.

Now suppose that the proportion of $\mcM \in \mbC_n$ which have an
$f(n)$-rich colouring approaches 1 as $n \to \infty$.
Define $T_{\mbC} = T_{iso} \cup T_{ext} \cup T_{col}$ 
exactly as in Section~\ref{l-colourable structures},
just before Lemma~\ref{countable categoricity}.
By the definition of $T_{iso}$, every $\varphi \in T_{iso}$ is true in
every $\mcM \in \mbC_n$.
By the last statement of the preceeding paragraph,
for every $\varphi \in T_{ext}$ the proportion of
$\mcM \in \mbC_n$ in which $\varphi$ holds approaches 1 as $n \to \infty$.
Recall that the formulas $\xi$ and $\zeta$ (defined in 
Notation~\ref{notation for S and xi} and~\ref{notation for U etcetera})
are $L_{rel}$-formulas. Lemma~\ref{auxiliary properties} and
the definition of $\mbX^{\mbC}_n$ (and of $\mbX_{n,k}$ in Notation~\ref{notation for S and xi})
imply that for every $\varphi \in T_{col}$, the proportion of $\mcM \in \mbC_n$ in
which $\varphi$ is true approaches 1 as $n \to \infty$.
Hence, for every finite $\Delta \subset T_{\mbC}$, 
the proportion of $\mcM \in \mbC_n$ such that $\mcM \models \Delta$
approaches 1 as $n \to \infty$.
By the completeness of $T_{\mbC}$ (Lemma~\ref{countable categoricity}) and compactness,
$\mbC$ has a zero-one law for the uniform probability measure.
Thus, we have proved part~(ii) of Theorem~\ref{transfer of 0-1 laws to uniform measure} and hence
the proof of that theorem is completed
(since, as explained in the beginning of the proof, the proof of part~(iii)
is the same except for obvious changes in notation and terminology).

Observe that, by Lemma~\ref{the number of different colourings when extension axioms are satisfied}
and (g), we have also proved the following:

\begin{prop}\label{definability of colourings under richness condition}
Let $f: \mbbN \to \mbb{R}$ be such that 
$\lim_{n \to \infty} \frac{f(n)}{\ln n} = \infty$. 
Suppose that the proportion of $\mcM \in \mbC_n$ which have an $f(n)$-rich colouring 
approaches 1 as $n \to \infty$.
Then the proportion of $\mcM \in \mbC_n$ such that every $l$-colouring $\gamma$
of $\mcM$ is definable by $\xi(y,z)$,
in the sense that $\mcM \models \xi(a,b) \Leftrightarrow \gamma(a) = \gamma(b)$,
approaches 1 as $n \to \infty$.
Consequently, the proportion of $\mcM \in \mbC_n$ which have a unique $l$-colouring,
up to permutation of colours, approaches 1 as $n \to \infty$.
The same statements hold if $\mbC_n$ is replaced by $\mbS_n$
(in which case the formula $\xi(y,z)$ may be different).
\end{prop}

\section{The uniform probability measure \\ and the typical distribution of colours}
\label{the uniform measure and l-colourable structures}

\noindent
In \cite{KPR}, Kolaitis, Prömel and Rothschild proved
that almost all $l$-colourable undirected graphs are uniquely $l$-colourable
(Corollary~1.23 \cite{KPR}), and the distribution of colours is relatively even
(Corollaries~1.20 and~1.21 \cite{KPR}). They also proved that the class of
$l$-colourable undirected graphs has a zero-one law, 
with the uniform probability measure, which together with
their first main result -- that almost all $\mcK_{l+1}$-free undirected graphs
are $l$-colourable -- implies the other main result, that the class of
$\mcK_{l+1}$-free graphs has a zero-one law. 
In the above context $l\geq 2$ is a fixed integer.
A further study of $l$-colourable graphs was made by Prömel and Steger in
\cite{PS95}, where $l = l(n)$ was allowed to grow, and the authors
found a threshold function $l = l(n)$ for the property of being uniquely $l$-colourable.
As in the previous section, we will let $l \geq 2$ be a fixed integer and study
random (strongly) $l$-colourable relational structures, but now only for the uniform probability measure.

The main results of this section,
Theorems~\ref{main theorem} and~\ref{main theorem for strong colourings}, 
generalize the zero-one law and (almost always) uniqueness of 
$l$-colouring for random $l$-colourable graphs in \cite{KPR} to random (strongly) $l$-colourable 
$L_{rel}$-structures for any relational language $L_{rel}$ subject to some mild assumptions.
They also tell that, almost always, the partition of the universe induced by 
an (strong) $l$-colouring is $L_{rel}$-definable without parameters.
Because of Theorem~\ref{transfer of 0-1 laws to uniform measure}~(ii) 
and Proposition~\ref{definability of colourings under richness condition},
in order to prove these things we only need to show that, for some function
$f : \mbbN \to \mbbR$ such that $\lim_{n\to\infty}f(n)/\ln n = \infty$, 
the proportion of (strongly) $l$-colourable $L_{rel}$-structures $\mcM$ with universe
$\{1, \ldots, n\}$ which have an $f(n)$-rich (strong) $l$-colouring approaches 1 as $n \to \infty$.
We will show (Theorem~\ref{new theorem}) that there is a constant $\mu > 0$ (depending on $l$, $L_{rel}$ and whether we consider $l$-colourings or {\em strong} $l$-colourings) such that the proportion of 
$L_{rel}$-structures $\mcM$ with universe $\{1, \ldots, n\}$ which have {\em only}
$\mu n$-rich (strong) $l$-colourings approaches 1 as $n \to \infty$.
The proof involves counting and estimating the number of
(strongly) multichromatic $m$-tuples and $m$-sets (Definition~\ref{definition of l-colouring etc})
for $m$ ranging from 2 to the maximum arity of the relation symbols.

As in the previous sections we will allow the possibility that certain relation
symbols are always interpreted as irreflexive and symmetric relations
(see Remark~\ref{remark that the theorems generalize to symmetric structures}).
As the arguments in this section are sensitive to whether this restriction
applies to a given relation symbol, we will (in contrast to previous sections)
be careful to let the notation indicate which relation symbols (if any) are 
always interpreted as irreflexive and symmetric relations.
Note that, apart from making this information visible, the notation
below agrees with that which was introduced in the beginning of 
Section~\ref{l-colourable structures}.

\begin{assump}\label{assumptions in last section}{\rm
We fix an integer $l \geq 2$ and a relational language
$L_{rel}$ with vocabulary $\{R_1, \ldots, R_{\rho}\}$, where $\rho > 0$ 
and each $R_k$ has arity $r_k \geq 2$.
}\end{assump}

\begin{defin}\label{definition of l-colouring etc}{\rm
(i) For positive integers $n$, we use the abbreviation $[n] = \{1, \ldots, n\}$.\\
(ii) Let $A$ be a set and let $\gamma : A \to [l]$.
An $m$-tuple $(a_1, \ldots, a_m) \in A^m$ is called {\bf \em monochromatic with respect to $\gamma$}
if $\gamma(a_1) = \ldots = \gamma(a_m)$. Otherwise we call $(a_1, \ldots, a_m)$
{\bf \em multichromatic with respect to $\gamma$}.
Note that if $m \geq 3$, then a multichromatic $m$-tuple may have
repetitions of elements ($a_i = a_j$ for some $i \neq j$).\\
(iii) An $m$-tuple $(a_1, \ldots, a_m) \in A^m$ is called 
{\bf \em strongly multichromatic with respect to $\gamma$} if $\gamma(a_i) \neq \gamma(a_j)$
whenever $i \neq j$.\\
(iv) Let $\mcM = (M, R_1^{\mcM}, \ldots, R_{\rho}^{\mcM})$ be an $L_{rel}$-structure and let 
$\gamma : M \to [l]$.
We say that $\gamma$ is an {\bf \em (strong) $l$-colouring of $\mcM$} if, for every
$k = 1, \ldots, \rho$, every $(a_1, \ldots, a_{r_k}) \in R_k^{\mcM}$
is (strongly) multichromatic with respect to $\gamma$.\\
(v) An $L_{rel}$-structure $\mcM$ is called {\bf \em (strongly) $l$-colourable} if there is $\gamma : M \to [l]$
which is an (strong) $l$-colouring of $\mcM$.\\
(vi) We say that an $L_{rel}$-structure $\mcM$ is {\bf \em uniquely (strongly) $l$-colourable} if it is 
(strongly) $l$-colourable
and for all (strong) $l$-colourings $\gamma$ and $\gamma'$ of $\mcM$ and all $a,b \in M$,
$\gamma(a) = \gamma(b)$ $\Longleftrightarrow$ $\gamma'(a) = \gamma'(b)$.\\
(vii) For every $I \subseteq [\rho]$, $\mbC^I_n$ denotes the set of $l$-colourable 
$L_{rel}$-structures $\mcM$ with universe $[n] = \{1, \ldots, n\}$ such that for every
$k \in I$, $R_k$ is interpreted as an irreflexive and symmetric relation in $\mcM$.
Let $\mbC^I = \bigcup_{n \in \mbbN}\mbC^I_n$, where $\mbbN$ is the set of positive integers.
(viii) For every $I \subseteq [\rho]$, $\mbS^I_n$ denotes the set of strongly $l$-colourable 
$L_{rel}$-structures $\mcM$ with universe $[n]$ such that for every
$k \in I$, $R_k$ is interpreted as an irreflexive and symmetric relation in $\mcM$.
Let $\mbS^I = \bigcup_{n \in \mbbN}\mbS^I_n$.\\
(iv) For $\alpha \in \mbbR$, a function $\gamma : [n] \to [l]$ is called {\bf \em $\alpha$-rich} if
$|f^{-1}(i)| \geq \alpha$ for every $i \in [l]$.
}\end{defin}

\noindent
As usual when the uniform probability measure is considered, 
the phrase {\bf \em `almost all $\mcM \in \mbC^I$ has property $P$'} means that the
proportion of $\mcM \in \mbC^I_n$ which have property $P$ approaches 1 as $n$ approaches infinity.
The phrase {\bf \em `$\mbC^I$ has a zero-one law'} means that for every $L_{rel}$-sentence $\varphi$, 
either $\varphi$ or its negation, $\neg\varphi$, is satisfied by almost all $\mcM \in \mbC^I$. 
(And similarly for $\mbS^I$ in place of $\mbC^I$.)

\begin{theor}\label{main theorem}
For every $I \subseteq [\rho]$ the following hold:\\
(i) There is an $L$-formula $\xi(x,y)$ such that for almost all $\mathcal{M} \in \mathbf{C}^I$
the following holds: for every $l$-colouring $\gamma : M \to [l]$ of $\mathcal{M}$ and all $a, b \in M$,
$\gamma(a) = \gamma(b)$ if and only if $\mathcal{M} \models \xi(a,b)$.\\
(ii) Almost all $\mathcal{M} \in \mathbf{C}^I$ are uniquely $l$-colourable.\\
(iii) $\mathbf{C}^I$ has a zero-one law.
\end{theor}

\begin{theor}\label{main theorem for strong colourings}
Suppose that every relation symbol has arity $\leq l$.
For every $I \subseteq [\rho]$, all three parts (i), (ii) and (iii) of 
Theorem~\ref{main theorem} hold if $\mbC^I$ is replaced by $\mbS^I$
and `strong' is added before `$l$-colouring/colourable'.
\end{theor}

\noindent
Recall from Remark~\ref{remark about the 0-1 laws} that
if $I \subseteq \{1, \ldots, \rho\}$ and $\mbC_n = \mbC^I_n$ and $\mbS_n = \mbS^I_n$,
then Theorem~\ref{transfer of 0-1 laws to uniform measure}
and Proposition~\ref{definability of colourings under richness condition}
hold. Hence, Theorems~\ref{main theorem} and~\ref{main theorem for strong colourings}
are immediate consequences of Theorem~\ref{new theorem} below and
Theorem~\ref{transfer of 0-1 laws to uniform measure}
and Proposition~\ref{definability of colourings under richness condition}.

\begin{theor}\label{new theorem}
(i) For every $I \subseteq [\rho]$ there are constants $\mu, \lambda > 0$ such that,
for all sufficiently large $n$, the proportion of 
$\mcM \in \mbC^I_n$ which have an $l$-colouring that is {\em not} $\mu n$-rich is at most
$2^{-\lambda n^m}$, where $m$ is the maximum arity of the relation symbols (so $m \geq 2$). 
Consequently, the proportion of $\mcM \in \mbC^I_n$  which have only $\mu n$-rich
$l$-colourings approaches 1 as $n \to \infty$. \\
(ii) If at least one relation symbol has arity $\leq l$, then 
part (i) also holds if $\mbC^I_n$ is replaced by $\mbS^I_n$,
`$l$-colouring' by `strong $l$-colouring', and $m$ is the largest arity $\leq l$ (so $m \geq 2$).\\
\end{theor}

\begin{rem}\label{remark about how to compute mu etc}{\rm
(i) In both parts of Theorem~\ref{new theorem},
the proof shows how to compute $\mu$ from the number of colours, $l$, 
and the arities $r_1, \ldots, r_{\rho}$
of the relation symbols of the language $L_{rel}$.\\
(ii) Theorem~\ref{new theorem} gives a bit more than what has been said above; namely
that the assumption in part~(i) of Theorem~\ref{transfer of 0-1 laws to uniform measure}
is true, which can be seen as follows.
Let $\mbK_n$ and $\mbSK_n$ be defined as in the previous section.
By Theorem~\ref{new theorem}~(i), there are constants $\mu, \lambda > 0$ such that
if $\mbY^{\mbK}_n$ is the set of $\mcM \in \mbK_n$ which are $\mu n$-richly $l$-coloured,
then, for all sufficiently large $n$, 
$$\big| \{\mcM \uhrc L_{rel} : \mcM \in \mbK_n - \mbY^{\mbK}_n\}\big| \ \Big/ \ |\mbC_n| 
\ \leq \
2^{-\lambda n^2}.$$
Since for each $\mcM \in \mbK_n$, $\mcM \uhrc L_{rel}$ can be $l$-coloured, or equivalently, 
expanded to an
$L$-structure (using the notation of the previous section), in at most
$l^n = 2^{\beta n}$ (for some $\beta > 0$) ways, we get
$\big| \mbK_n - \mbY^{\mbK}_n \big| \big/ |\mbK_n| \leq 2^{\beta n - \lambda n^2} \to 0$ as $n \to \infty$.
Therefore the assumption in part~(i) of Theorem~\ref{transfer of 0-1 laws to uniform measure}
holds, and it follows that, for every $k \in \mbbN$, the proportion of $\mcM \in \mbK_n$ which satisfy all
$k$-extension axioms of $\mbK$ approaches 1 as $n \to \infty$.
The same argument can be carried out for $\mbSK_n$, $\mbS_n$ and strong $l$-colourings.
}\end{rem}

\begin{exam}{\rm
Here follows applications of Theorems~\ref{main theorem}
and~\ref{main theorem for strong colourings}.

(i) Let $\mathcal{F}$ be the Fano plane as a 3-hypergraph, that is,
$\mathcal{F}$ has seven vertices and seven 3-hyperedges
(3-subsets of the vertex set) such that every pair of distinct vertices
is contained in a unique 3-hyperedge.
If $\mathbf{K}_n$ is the set of all 3-hypergraphs with
vertices $1, \ldots, n$ in which $\mathcal{F}$ is not embeddable,
and $\mathbf{K} = \bigcup_{n \in \mathbb{N}} \mathbf{K}_n$,
then almost all members of $\mbK$ are 2-colourable \cite{PeSch}.
Since $\mathcal{F}$ cannot be weakly embedded into any
2-colourable 3-hypergraph, it follows from Theorem~\ref{main theorem}
that $\mbK$ has a zero-one law for the uniform
probability measure.

(ii) Let $\mathcal{G}$ be the 3-hypergraph with vertices $1,2,3,4,5$
and 3-hyperedges $\{1,2,3\}$, $\{1,2,4\}$, $\{3,4,5\}$, and let
$\mathbf{K}_n$ be the set of 3-hypergraphs with vertices $1, 2, \ldots, n$
in which $\mathcal{G}$ is not weakly embeddable.
Then almost all members of $\mathbf{K} = \bigcup_{n \in
\mathbb{N}}\mathbf{K}_n$
are strongly 3-colourable \cite{BM}.
(Tripartite in \cite{BM} means the same as strongly
3-colourable here.) Since $\mathcal{G}$ cannot be weakly
embedded into any strongly 3-colourable 3-hypergraph it follows
from Theorem~\ref{main theorem for strong colourings} 
that $\mathbf{K}$ has a zero-one law for the
uniform probability measure.
}\end{exam}

\noindent
In Section~\ref{counting multichromatic tuples} we derive 
an upper bound on the number of multichromatic $m$-tuples
if the $l$-colouring $\gamma : [n] \to [l]$ is {\em not} $\frac{n}{a}$-rich and $a$ is 
sufficiently large. 
Then we show that the number of multichromatic $m$-sets are fairly tightly
controlled by the number of multichromatic $m$-tuples.
These results are used in Section~\ref{proof of new theorem} where we prove
part (i) of Theorem~\ref{new theorem}.
In Section~\ref{counting strongly multichromatic tuples} we consider {\em strongly} multichromatic
$m$-tuples and $m$-sets and derive similar results as in Section~\ref{counting multichromatic tuples}
which are used in
Section~\ref{proof of second part of new theorem} where part (ii) of 
Theorem~\ref{new theorem} is proved.

\subsection{Counting multichromatic tuples and sets}\label{counting multichromatic tuples}

\begin{notation}\label{our notation}{\rm
(i) Let $n, m, l \in \mbbN$ and suppose that $n \geq l \geq 2$ and $m \geq 2$.\\ 
(ii) Let $\gamma : [n] \to [l]$. \\
(iii) Let $\mult(n,\gamma,m)$ denote the number of ordered $m$-tuples $(a_1, \ldots, a_m) \in [n]^m$
which are multichromatic with respect to $\gamma$.\\
(iv) For every $i \in [l]$, let $p(n,\gamma,i) = \big|\gamma^{-1}(i)\big|$, so
$p(n,\gamma,i)$ is the number of elements in~$[n]$ which are assigned the colour $i$ by $\gamma$.
}\end{notation}

\noindent
The number of $(a_1, \ldots, a_m) \in [n]^m$ which are monochromatic
with respect to $\gamma$ is $\sum_{i=1}^l p^m(n,\gamma,i)$, from which it follows that
\begin{equation}\label{expression for mult}
\mult(n,\gamma,m) \ = \ n^m \ - \ \sum_{i=1}^l p^m(n,\gamma,i).
\end{equation}

\begin{rem}\label{maximum of f}{\rm
{\em Let $k \in \mbbN$ and $\alpha \in \mathbb{R}^+$. 
The function 
$f_{k,m}(x_1, \ldots, x_k) \ = \ \sum_{i=1}^k x_i^m$
constrained by
$x_1 + \ldots + x_k = \alpha$ and $x_i \geq 0$, for $i = 1, \ldots, k$,
attains its minimal value in the point $(\alpha/k, \ldots, \alpha/k)$, 
and hence this value is $\alpha^m \big/ k^{m-1}$.}
This fact is easily proved by using the method of Lagrange multipliers \cite{Had}.
Alternatively, one can use a variant of Hölder's inequality:
In the result with number 16 in \cite{HLP} (p. 26), take $r=1$, $s=m$ and $a = (x_1, \ldots, x_k)$
and the claim ``$\mathfrak{M}_r(a) < \mathfrak{M}_s(a)$ unless ...'' becomes
$(x_1 + \ldots + x_k)/k < \big(\frac{1}{k}f_{k,m}(x_1, \ldots, x_k)\big)^{1/m}$
unless all $x_i$ are equal. Since we assume that $x_1 + \ldots + x_k = \alpha$, the claim 
follows by taking the $m$th power on both sides.
}\end{rem}

\begin{lem}\label{upper bound of multichromatic m-tuples for nonrich colourings}
Let $a > 0$.
If  $\gamma : [n] \to [l]$ is {\rm not} $\frac{n}{a}$-rich, then
$$\mult(n,\gamma,m) \ \leq \
\Bigg( 1 \ - \ \bigg[\frac{a - 1}{a}\bigg]^m\frac{1}{(l-1)^{m-1}} \Bigg)n^m.$$
\end{lem}

\noindent
{\em Proof.}
Let $a > 0$.
Suppose that $\gamma : [n] \to [l]$ is not $\frac{n}{a}$-rich, 
which means that, for some $i \in [l]$,
$p(n,\gamma,i) < \frac{n}{a}$. 
For simplicity of notation (and without loss of generality) assume that $i = l$.
Then 
\begin{equation}\label{getting rid of p(n,gamma,1)}
n \ - \ p(n,\gamma,l) \ > n \ - \ \frac{n}{a} \ = \ \frac{a-1}{a}n.
\end{equation}
Now we have
\begin{align*}
&\mult(n,\gamma,m) \ =\\ 
&= \ n^m \ - \ \sum_{i=1}^{l-1} p^m(n,\gamma,i) \ - \ p^m(n,\gamma,l)\\
&\leq \ n^m \ - \ \frac{\big[n - p(n,\gamma,l)\big]^m}{(l-1)^{m-1}} \ - \ p^m(n,\gamma,l)\\
&\hspace{40mm} 
\text{ by Remark~\ref{maximum of f} with $\alpha = n - p(n,\gamma,l)$ and $k = l-1$}\\
&< \ n^m \ - \ \bigg[\frac{a-1}{a}\bigg]^m\frac{n^m}{(l-1)^{m-1}} \ - \ 
p^m(n,\gamma,l)
\quad \text{ by (\ref{getting rid of p(n,gamma,1)})}\\
&\leq \ \Bigg( 1 \ - \ \bigg[\frac{a - 1}{a}\bigg]^m\frac{1}{(l-1)^{m-1}} \Bigg)n^m.
\hspace{20mm} \square
\end{align*}

\begin{notation}\label{notation for multichromatic sets}{\rm
(i) As usual, by a {\em $k$-set} we mean a set of cardinality $k$.\\
(ii) For every integer $k \geq 2$, every $n \in \mbbN$ and every
$\gamma : [n] \to [l]$, let $\overline{\mult}(n,\gamma,k)$ be the 
number of $k$-subsets $\{a_1, \ldots, a_k\} \subseteq [n]$ such that
there are $i,j \in [k]$ with $\gamma(a_i) \neq \gamma(a_j)$.
We call such a $k$-set $\{a_1, \ldots, a_k\}$ {\bf \em multichromatic}.\\
(iii) For integers $1 \leq i \leq k$, let $\perm(i,k)$ be the
number of (ordered) $k$-tuples $(a_1, \ldots, a_k)$ of elements of an $i$-set $A$ 
such that every $a \in A$ occurs at least once in $(a_1, \ldots, a_k)$.
}\end{notation}

\begin{lem}\label{connection between multichromatic tuples and sets}
Let $m_{\max} \geq 2$ be an integer and
suppose that $\sigma_n : [n] \to [l]$ and $\gamma_n : [n] \to [l]$, for $n \in \mbbN$.
Moreover, assume that for all $2 \leq m \leq m_{\max}$
there are constants $c_m, d_m > 0$ such that for all sufficiently large $n$,
$$c_m n^m \ \leq \ \mult(n,\sigma_n,m)  -  \mult(n,\gamma_n,m) \ \leq \ d_m n^m.$$
Then, for all $2 \leq m \leq m_{\max}$, there are constants $c'_m, d'_m > 0$ such that
for all sufficiently large $n$,
$$c'_m n^m \ \leq \ \overline{\mult}(n,\sigma_n,m)  -  \overline{\mult}(n,\gamma_n,m) 
\ \leq \ d'_m n^m.$$
\end{lem}

\noindent
{\em Proof.}
Suppose that for all for all $2 \leq m \leq m_{\max}$ there are
$c_m, d_m > 0$ such that for all sufficiently large $n$,
\begin{equation}\label{the assumption about mult being bounded from both directions}
c_m n^m \ \leq \ \mult(n,\sigma_n,m)  -  \mult(n,\gamma_n,m) \ \leq \ d_m n^m.
\end{equation}
Note that if an $m$-tuple $(a_1, \ldots, a_m) \in [n]^m$ is multichromatic with respect to 
$\gamma : [n] \to [l]$, then $\big|\{a_1, \ldots, a_m\}\big| = i$ for some $2 \leq i \leq m$.
Therefore (with the notation introduced before the lemma), for every $m$ 
and every $\gamma : [n] \to [l]$ we have
\begin{equation}\label{equation connecting mult and overline-mult}
\mult(n,\gamma,m) \ = \ \sum_{i=2}^m \ \overline{\mult}(n,\gamma,i) \cdot \perm(i,m).
\end{equation}
We use induction on $m = 2, \ldots, m_{\max}$.
If $m = 2$ then
(\ref{the assumption about mult being bounded from both directions})
and (\ref{equation connecting mult and overline-mult}) give
$$\frac{c_2}{\perm(2,2)} \ \leq \ 
\overline{\mult}(n,\sigma_n,2)  -  \overline{\mult}(n,\gamma_n,2) \ \leq \
\frac{d_2}{\perm(2,2)},$$
so we can take $c'_2 = c_2 \big/ \perm(2,2)$ and $d'_2 = d_2 \big/ \perm(2,2)$.

As induction hypothesis, suppose that for $m = 2, \ldots, k < m_{\max}$
there are $c'_m, d'_m > 0$ such that for all sufficiently large $n$,
\begin{equation}\label{induction hypothesis}
c'_m n^m \ \leq \ \overline{\mult}(n,\sigma_n,m)  -  \overline{\mult}(n,\gamma_n,m) 
\ \leq \ d'_m n^m.
\end{equation}
By assumption, for all sufficiently large $n$ we have
$$c_{k+1} n^{k+1} \ \leq \ \mult(n,\sigma_n,k+1)  -  \mult(n,\gamma_n,k+1) \ \leq \ d_{k+1} n^{k+1},$$
and by (\ref{equation connecting mult and overline-mult}) with $m=k+1$ it follows
that for all sufficiently large $n$,
\begin{align*}
c_{k+1} n^{k+1} \ &\leq \
\sum_{i=2}^k \ \perm(i,k+1)\Big[\overline{\mult}(n,\sigma_n,i) - \overline{\mult}(n,\gamma_n,i)\Big]  \\ 
+ \ &\perm(k+1,k+1)\Big[\overline{\mult}(n,\sigma_n,k+1)  -  \overline{\mult}(n,\gamma_n,k+1)\Big]
\ \leq \ d_{k+1} n^{k+1}. 
\end{align*}
By the induction hypothesis, (\ref{induction hypothesis}) holds for $m = 2, \ldots, k$
and all sufficiently large $n$. Hence, for all sufficiently large $n$,
\begin{align*}
&c_{k+1} n^{k+1} \ \leq \\
&\perm(k+1,k+1)\Big[\overline{\mult}(n,\sigma_n,k+1)  -  \overline{\mult}(n,\gamma_n,k+1)\Big]
\ + \ O\big(n^k\big) 
\ \leq \ d_{k+1} n^{k+1},
\end{align*}
so there must be $c'_{k+1}, d'_{k+1} > 0$ such that for all sufficiently large $n$,
$$c'_{k+1} n^{k+1} \ \leq \ \overline{\mult}(n,\sigma_n,k+1)  -  \overline{\mult}(n,\gamma_n,k+1) 
\ \leq \ d'_{k+1} n^{k+1}. \quad \quad \square$$

\subsection{Proof of the first part of Theorem~\ref{new theorem}}\label{proof of new theorem}

\noindent
We continue to use the terminology and notation introduced in Notation~\ref{our notation}
and~\ref{notation for multichromatic sets}.
Recall that all arities $r_1, \ldots, r_{\rho}$ of the relation symbols $R_1, \ldots, R_{\rho}$
are at least 2.
Let $m_{\max} = \max(r_1, \ldots, r_{\rho})$.
For all $l,m \geq 2$ we have
\begin{equation*}
\frac{1}{(l-1)^{m-1}} \ > \ \frac{1}{l^{m-1}}.
\end{equation*}
Let $a > l$ be large enough so that, whenever $2 \leq m \leq m_{\max}$,
\begin{equation}\label{choice of a}
\bigg[\frac{a-1}{a}\bigg]^m \frac{1}{(l-1)^{m-1}} \ > \ \frac{1}{l^{m-1}}.
\end{equation}

\noindent
For every $n \in \mbbN$ such that $n \geq l$, fix $\sigma_n : [n] \to [l]$ 
such that, for every $i \in [l]$, 
\begin{equation}\label{property of sigma}
\frac{n}{l} - 1 \ \leq \ p(n,\sigma_n,i) \ \leq \ \frac{n}{l} + 1.
\end{equation}
Then, for all sufficiently large $n$, $\sigma_n$ is $\frac{n}{a}$-rich
(because we chose $a > l$).
Observe that if $\gamma : [n] \to [l]$, then the number of $\mcM \in \mbC^I_n$
for which $\gamma$ is an $l$-colouring is
\begin{equation*}
2^{\sum_{k \in [\rho]-I} \mult(n,\gamma,r_k) \ + \ \sum_{k \in I} \overline{\mult}(n,\gamma,r_k)}.
\end{equation*}
Therefore,
\begin{equation}\label{lower bound for X-n}
\big|\mbC^I_n\big| \ \geq \  
2^{\sum_{k \in [\rho]-I} \mult(n,\sigma_n,r_k) \ + \ \sum_{k \in I} \overline{\mult}(n,\sigma_n,r_k)}.
\end{equation}

\noindent
A lower bound of $\mult(n,\sigma_n,m)$ is obtained as follows:
\begin{align*}
\mult(n,\sigma_n,m) \ = \ n^m \ - \ \sum_{i=1}^l p^m(n,\sigma_n,i) \
&\geq \ n^m \ - \ l\bigg(\frac{n}{l} + 1\bigg)^m \quad \quad \text{ by (\ref{property of sigma})}\\
&= \ n^m \ - \ \frac{n^m}{l^{m-1}} \ \pm \ O\big(n^{m-1}\big),
\end{align*}
so
\begin{equation}\label{lower bound for sigma-n}
\mult(n,\sigma_n,m) \ \geq \ \bigg( 1 - \frac{1}{l^{m-1}}\bigg) n^m \ \pm \ O\big(n^{m-1}\big).
\end{equation}
For every $n \in \mbbN$, choose $\gamma_n : [n] \to [l]$ such that
\begin{equation}\label{gamma-n not rich}
\text{$\gamma_n$ is {\em not} $\frac{n}{a}$-rich, and}
\end{equation}
for every $\gamma : [n] \to [l]$ which is {\em not} $\frac{n}{a}$-rich,
\begin{align}\label{gamma-n is maximal such}
&\sum_{k \in [\rho] - I} \mult(n,\gamma,r_k) \ + \ \sum_{k \in I} \overline{\mult}(n,\gamma,r_k) \\
\leq \ 
&\sum_{k \in [\rho] - I} \mult(n,\gamma_n,r_k) \ + \ \sum_{k \in I} \overline{\mult}(n,\gamma_n,r_k). \nonumber
\end{align}
Let $\mbX_n \subseteq \mbC^I_n$ be the set of all $\mcM \in \mbC^I_n$ which have
an $l$-colouring which is {\em not} $\frac{n}{a}$-rich.
It suffices to prove that $|\mbX_n| \big/ |\mbC^I_n| \to 0$ as $n \to \infty$.
If $\mcM \in \mbX_n$ then there is an $l$-colouring $\gamma : [n] \to [l]$
of $\mcM$ which is not $\frac{n}{a}$-rich, and the number of $\mcN \in \mbC^I_n$ for which
$\gamma$ is an $l$-colouring is at most 
$2^{\sum_{k \in [\rho]-I} \mult(n,\gamma,r_k) \ + \ \sum_{k \in I} \overline{\mult}(n,\gamma,r_k)}$.
Since the number of functions $\gamma : [n] \to [l]$ is $l^n = 2^{\beta n}$, for some $\beta > 0$,
it follows from~(\ref{gamma-n not rich}) and~(\ref{gamma-n is maximal such}) that 
\begin{equation}\label{upper bound of complement of X-n}
\big|\mbX_n\big| \ \leq \ 
2^{\beta n \ + \ 
\sum_{k \in [\rho]-I} \mult(n,\gamma_n,r_k) \ + \ \sum_{k \in I} \overline{\mult}(n,\gamma_n,r_k)}.
\end{equation}

\noindent
From (\ref{gamma-n not rich}) and Lemma~\ref{upper bound of multichromatic m-tuples for nonrich colourings}
it follows that
\begin{equation}\label{upper bound for gamma-n}
\mult(n,\gamma_n,m) \ \leq \
\Bigg( 1 \ - \ \bigg[\frac{a-1}{a}\bigg]^m \frac{1}{(l-1)^{m-1}} \Bigg)n^m.
\end{equation}

\noindent
Note that for all $n, m$ and $\gamma : [n] \to [l]$ we have
$\overline{\mult}(n,\gamma,m) \leq \mult(n,\gamma,m) \leq n^m$.
Therefore,~(\ref{lower bound for sigma-n}) and~(\ref{upper bound for gamma-n})
imply that
\begin{align*}
n^m \ &\geq \ \mult(n,\sigma_n,m) - \mult(n,\gamma_n,m) \\
&\geq \ \Bigg[ 1 - \frac{1}{l^{m-1}} 
\ - \ \bigg( 1 \ - \ \bigg[\frac{a-1}{a}\bigg]^m \frac{1}{(l-1)^{m-1}} \bigg) \Bigg]n^m \ 
\pm \ O\big(n^{m-1}\big)\\
&= \ \Bigg[ \bigg[\frac{a-1}{a}\bigg]^m \frac{1}{(l-1)^{m-1}} \ - \ \frac{1}{l^{m-1}} \Bigg]n^m \
\pm \ O\big(n^{m-1}\big).
\end{align*}
Together with~(\ref{choice of a}) this implies that there is $c > 0$ such that
whenever $2 \leq m \leq m_{\max}$ and $n$ is sufficiently large
\begin{equation}\label{the double bound of mult}
cn^m \ \leq \ \mult(n,\sigma_n,m) \ - \ \mult(n,\gamma_n,m) \ \leq \ n^m.
\end{equation}
Lemma~\ref{connection between multichromatic tuples and sets} now implies that
for all $2 \leq m \leq m_{\max}$ there are $c'_m > 0$ such that for all sufficiently large $n$,
\begin{equation}\label{the double bound of overline-mult}
c'_m n^m \ \leq \ \overline{\mult}(n,\sigma_n,m) \ - \ \overline{\mult}(n,\gamma_n,m).
\end{equation}
By~(\ref{lower bound for X-n}) and~(\ref{upper bound of complement of X-n}) we have
\begin{align}\label{the fraction to be estimated}
&\big|\mbX_n\big| \Big/ \big|\mbC^I_n\big| \leq \\
&\leq 2^{\beta n \ + \ \sum_{k \in [\rho] - I} \big[ \mult(n,\gamma_n,r_k) \ - \ \mult(n,\sigma_n,r_k)\big]
\ + \ \sum_{k \in I} \big[ \overline{\mult}(n,\gamma_n,r_k) \ - \ \overline{\mult}(n,\sigma_n,r_k)\big]}.
\nonumber
\end{align}
From~(\ref{the double bound of mult}) and~(\ref{the double bound of overline-mult})
it follows that for all sufficiently large $n$,
\begin{align}
&k \in [\rho] - I \ \Longrightarrow \ 
\mult(n,\gamma_n,r_k) - \mult(n,\sigma_n,r_k) \ \leq \ -c n^{r_k}, \ \text{ and}\\
\label{final estimates}
&k \in I \ \Longrightarrow \ 
\overline{\mult}(n,\gamma_n,r_k) - \overline{\mult}(n,\sigma_n,r_k) \ \leq \ -c'_{r_k} n^{r_k},
\end{align}
where $c, c'_{r_k} > 0$ for all $k \in [\rho]$.
Since $r_k \geq 2$ for all $k \in [\rho]$ it follows 
from~(\ref{the fraction to be estimated})--(\ref{final estimates}) that,
for $m = m_{\max}$ and some $\lambda > 0$, we have (for all large enough $n$)
\begin{equation*}
\big|\mbX_n\big| \Big/ \big|\mbC^I_n\big| \ \leq \ 
2^{-\lambda n^m} \ \to \ 0 \quad \text{ as } \  n \to \infty.
\end{equation*}
In other words, the proportion of $\mcM \in \mbC^I_n$ which only have $\frac{n}{a}$-rich 
$l$-colourings approaches 1 as $n \to \infty$;
in part~(i) of Theorem~\ref{new theorem} we can take $\mu = \frac{1}{a}$.

\subsection{Counting strongly multichromatic tuples and sets}\label{counting strongly multichromatic tuples}

\begin{notation}\label{notation for strongly multichromatic tuples}{\rm
(i) Let $n, m, l \in \mbbN$ and suppose that $n \geq l \geq m \geq 2$.\\ 
(ii) Let $\gamma : [n] \to [l]$. \\
(iii) Let $\smult(n,\gamma,m)$ denote the number of ordered $m$-tuples $(a_1, \ldots, a_m) \in [n]^m$
which are strongly multichromatic with respect to $\gamma$.\\
(iv) Let $\overline{\smult}(n,\gamma,m)$ be the number of $m$-subsets 
$\{a_1, \ldots, a_m\} \subseteq [n]$ such that $\gamma(a_i) \neq \gamma(a_j)$
whenever $i \neq j$.\\
(v) For every $i \in [l]$, let $p(n,\gamma,i) = \big|\gamma^{-1}(i)\big|$, so
$p(n,\gamma,i)$ is the number of elements in~$[n]$ which are assigned the colour $i$ by $\gamma$.
}\end{notation}

\noindent
Observe that
\begin{equation}\label{connection between smult for sets and tuples}
\smult(n,\gamma,m) \ = \ m! \ \overline{\smult}(n,\gamma,m)
\end{equation}
and
\begin{equation}\label{the expression for smult}
\overline{\smult}(n,\gamma,m) \ = 
\sum_{1 \leq i_1 < \ldots < i_m \leq l} p(n,\gamma,i_1) \ldots p(n,\gamma,i_m).
\end{equation}

\noindent
For $k \geq m$, let 
\begin{equation*}
g_{k,m}(x_1, \ldots, x_k) \ = 
\sum_{1 \leq i_1 < \ldots < i_m \leq k} x_{i_1} \ldots x_{i_m}.
\end{equation*}

\noindent
The next lemma is a special case of an inequality found in \cite{HLP} (p. 52),
see Remark~\ref{remark on an inequality} below,
but here we give a proof based on the better known method of Lagrange multipliers \cite{Had}.

\begin{lem}\label{upper bound of f for strongly multichromatic}
Let $\alpha > 0$.
Subject to the constraints $x_1 + \ldots + x_k = \alpha$ and
$x_i \geq 0$ for all $i \in [l]$, $g_{k,m}$ attains
its maximum in $(\alpha/k, \ldots, \alpha/k)$, and hence the maximum is
$$g_{k,m}(\alpha/k, \ldots, \alpha/k) \ = \ \binom{k}{m}\bigg( \frac{\alpha}{k} \bigg)^m.$$
\end{lem}

\noindent
{\em Proof.}
Suppose that, under the given constraints, $g_{k,m}$ attains its maximum
in $(a_1, \ldots, a_k)$. Then at least one $a_i$ is non-zero, hence positive.
We show that for every $j$, $a_j = a_i$. By the constraint $a_1 + \ldots + a_k = \alpha$
it follows that $a_j = \alpha/k$ for all $j \in [k]$. 
For simplicity of notation, and without loss of generality, assume that $i=1$,
so $a_1 > 0$. Since the following argument works out in the same way for all $j = 2, \ldots, k$,
let's assume that $j = 2$ (simplifying notation again).

Let $h(x_1, x_2) = g_{k,m}(x_1, x_2, a_3, \ldots, a_k)$.
Since we assume that $g_{k,m}$ attains its maximum, under the given constraints,
in $(a_1, \ldots, a_k)$, it follows that $h(x_1, x_2)$ attains its maximum, under the
constraints $x_1 + x_2 = \alpha - a_3 - \ldots - a_k$ 
(where $\alpha - a_3 - \ldots - a_k > 0$ since $a_1 > 0$)
and $x_1, x_2 \geq 0$, in $(a_1, a_2)$.
Observe that
\begin{align*}
h(x_1, x_2) \ = \ &\sum_{3 \leq i_1 < \ldots < i_{m-2} \leq k} x_1 x_2 a_{i_1} \ldots a_{i_{m-2}} \
+ \ \sum_{3 \leq i_1 < \ldots < i_{m-1} \leq k} x_1 a_{i_1} \ldots a_{i_{m-1}}\\
+ \ &\sum_{3 \leq i_1 < \ldots < i_{m-1} \leq k} x_2 a_{i_1} \ldots a_{i_{m-1}} \
+ \ \sum_{3 \leq i_1 < \ldots < i_m \leq k} a_{i_1} \ldots a_{i_m}\\
= \ &\sum_{3 \leq i_1 < \ldots < i_{m-2} \leq k} x_1 x_2 a_{i_1} \ldots a_{i_{m-2}} \
+ \ \sum_{3 \leq i_1 < \ldots < i_{m-1} \leq k} (x_1 + x_2) a_{i_1} \ldots a_{i_{m-1}}\\
+ \ &\sum_{3 \leq i_1 < \ldots < i_m \leq k} a_{i_1} \ldots a_{i_m}.
\end{align*}
Subject to the constraint $x_1 + x_2 = \alpha - a_3 - \ldots - a_k$, the part
$$\sum_{3 \leq i_1 < \ldots < i_{m-1} \leq k} (x_1 + x_2) a_{i_1} \ldots a_{i_{m-1}} \
+ \ \sum_{3 \leq i_1 < \ldots < i_m \leq k} a_{i_1} \ldots a_{i_m}$$
is constant, and hence $h(x_1, x_2)$ attains its maximum in the same point as
$h^*(x_1, x_2) = cx_1 x_2$, where $c > 0$ is a constant.
The reason that we can assume that $c > 0$ is that if $m > 2$, then
there are at least $m-2$ non-zero $a_i$'s with $i > 2$; because otherwise $g_{k,m}$
would be zero in $(a_1, \ldots, a_k)$ and then this point could not be a maximum,
contrary to assumption.
Thus it suffices to show that, for any $\beta > 0$, if 
$h^*(x_1, x_2) = cx_1 x_2$ attains its maximum in $(b_1, b_2)$
under the constraints $x_1 + x_2 = \beta$, $x_1, x_2 \geq 0$,
then $b_1 = b_2$. This is easily proved by (for example) using Lagrange multipliers \cite{Had}.

Given that $a_1 = \ldots = a_k$, the constraints on $g_{k,m}$ imply
that $a_i = \alpha/k$ for all $i$, and insertion of $(x_1, \ldots, x_k) = (a_1, \ldots, a_k)$
in the expression of $g_{k,m}$ shows that its maximum, subject to the constraints, is 
$\binom{k}{m} \big( \alpha / k \big)^m$.
\hfill $\square$

\begin{rem}\label{remark on an inequality}{\rm
Lemma~\ref{upper bound of f for strongly multichromatic} is 
a special case of the result with number 52 in \cite{HLP} (p. 52),
which is attributed to Maclaurin \cite{Mac}.
In the notation of that result, but with the letter $n$ replaced by $k$, 
we have ``$p_1 > (p_2)^{1/2} > \ldots > (p_k)^{1/k}$ unless ...'', so in particular
``$p_1 > (p_m)^{1/m}$ unless ...'', which, with the notation here and because
$x_1 + \ldots + x_k = \alpha$, becomes
$\alpha/k > \big(g_{k,m}(x_1, \ldots, x_k)/\binom{k}{m}\big)^{1/m}$ unless all $x_i$ are equal.
By raising both sides to the $m$th power we get the statement of the lemma.
}\end{rem}

\begin{lem}\label{upper bound for gamma in the strongly l-coloured case}
Let $a > 0$.
If $\gamma : [n] \to [l]$ is {\rm not} $\frac{n}{a}$-rich,
then 
$$\overline{\smult}(n,\gamma,m) \ \leq \
\Bigg[\frac{1}{(l-1)^m} \binom{l-1}{m} \ + \ \frac{1}{a(l-1)^{m-1}} \binom{l-1}{m-1}\Bigg]n^m,$$
where the left term within the large parentheses vanishes if $m = l$.
\end{lem}

\noindent
{\em Proof.}
Suppose that $a > 0$ and that $\gamma : [n] \to [l]$ is not $\frac{n}{a}$-rich. 
Then, for some $i \in [l]$, we have $p(n,\gamma,i) < \frac{n}{a}$.
For simplicity of notation, and without loss of generality, assume that $i = l$. 
Then:
\begin{align*}
\overline{\smult}(n,\gamma,m) \ = \
&\sum_{1 \leq i_1 < \ldots < i_m \leq l} p(n,\gamma,i_1) \cdot \ldots \cdot p(n,\gamma,i_m)
\quad \text{by (\ref{the expression for smult})} \\
= \ &\sum_{1\leq i_1 < \ldots < i_m \leq l-1} p(n,\gamma,i_1) \ldots p(n,\gamma,i_m) \\
+ \ &\sum_{1\leq i_1 < \ldots < i_{m-1} \leq l-1} 
p(n,\gamma,i_1) \ldots p(n,\gamma,i_{m-1})p(n,\gamma,l) \\
&\text{where the first sum vanishes if $m=l$}\\
< \ 
&\sum_{1\leq i_1 < \ldots < i_m \leq l-1} p(n,\gamma,i_1) \ldots p(n,\gamma,i_m)\\
+ \ &\Bigg[\sum_{1\leq i_1 < \ldots < i_{m-1} \leq l-1} p(n,\gamma,i_1) \ldots p(n,\gamma,i_{m-1})\Bigg]
\frac{n}{a} \quad \text{ by assumption}\\
\leq \ &\binom{l-1}{m}\Bigg(\frac{n - p(n,\gamma,l)}{l-1}\Bigg)^m \ + \
\frac{n}{a}\binom{l-1}{m-1} \Bigg(\frac{n - p(n,\gamma,l)}{l-1}\Bigg)^{m-1}\\
&\text{by Lemma~\ref{upper bound of f for strongly multichromatic}, twice,
with $k = l-1$, $\alpha = n - p(n,\gamma,l)$, and}\\ 
&\text{with $m$ in the first application and
$m-1$ in the second application}\\
\leq \ &\binom{l-1}{m} \frac{n^m}{(l-1)^m} \ + \ \binom{l-1}{m-1} \frac{n^m}{a(l-1)^{m-1}}\\
= \ &\Bigg[ \frac{1}{(l-1)^m} \binom{l-1}{m} \ + \ \frac{1}{a(l-1)^{m-1}} \binom{l-1}{m-1}\Bigg]n^m. 
\quad \quad \quad \square
\end{align*}

\subsection{Proof of the second part of Theorem~\ref{new theorem}}\label{proof of second part of new theorem}

\noindent
Let $m_{\max} =  \max(r_1, \ldots, r_{\rho})$, where $r_1, \ldots, r_{\rho} \geq 2$ are
the arities of the relation symbols $R_1, \ldots, R_{\rho}$ of the vocabulary of $L_{rel}$.
Suppose that at least one relation symbol has arity $\leq l$.
We use the notation from the previous section 
(Notation~\ref{notation for strongly multichromatic tuples}).
Since for every $m > l$ and every $\gamma : [n] \to [l]$,
{\em no} $(a_1,\ldots,a_m) \in [n]^m$ is strongly multichromatic
with respect to $\gamma$, 
we may, without loss of generality, assume that $m_{\max} \leq l$.

Let $\mbX_n \subseteq \mbS^I_n$ be the set of all $\mcM \in \mbS^I_n$ which have
a strong $l$-colouring which is {\em not} $\frac{n}{a}$-rich, 
where $a > l$ is a number that will be specified
after we have made some estimates. 
In order to prove part (ii) of Theorem~\ref{new theorem} it is enough to prove
that we can choose $a$ so that
\begin{equation}\label{equation giving second part of theorem}
\lim_{n \to \infty} \frac{\big|\mbX_n \big|}{\big| \mbS^I_n \big|} \ = \ 0.
\end{equation}
For every $n \geq l$, fix $\sigma_n : [n] \to [l]$ 
such that, for every $i \in [l]$, 
\begin{equation}\label{definition of sigma-n for strong colourings}
\frac{n}{l} - 1 \ \leq \ p(n,\sigma_n,i) \ \leq \ \frac{n}{l} + 1.
\end{equation}
Then every $\sigma_n$ is $\frac{n}{a}$-rich (because $a > l$). For all sufficiently large $n$,
\begin{equation}\label{lower bound for X-n for strong colourings}
\big| \mbS^I_n \big| \ \geq \
2^{\sum_{k \in [\rho]-I} \smult(n,\sigma_n,r_k) \ + \ \sum_{k \in I} \overline{\smult}(n,\sigma_n,r_k)}.
\end{equation}

\noindent
For every $n \in \mbbN$, choose $\gamma_n : [n] \to [l]$ such that
\begin{equation}\label{first property of gamma-n for strong colourings}
\text{$\gamma_n$ is {\em not} $\frac{n}{a}$-rich, and}
\end{equation}
for every $\gamma : [n] \to [l]$ which is {\em not} $\frac{n}{a}$-rich,
\begin{align}\label{second property of gamma-n for strong colourings}
&\sum_{k \in [\rho] - I} \smult(n,\gamma,r_k) \ + \ \sum_{k \in I} \overline{\smult}(n,\gamma,r_k) \\
\ \leq \
&\sum_{k \in [\rho] - I} \smult(n,\gamma_n,r_k) \ + \ \sum_{k \in I} \overline{\smult}(n,\gamma_n,r_k). \nonumber
\end{align}
Since, for some $\beta > 0$, there are at most $l^n = 2^{\beta n}$ $l$-colourings $\gamma : [n] \to [l]$,
and every $\mcM \in \mbX_n$ has an $l$-colouring which is not $\frac{n}{a}$-rich, 
it follows
from~(\ref{first property of gamma-n for strong colourings}) 
and~(\ref{second property of gamma-n for strong colourings}) that
\begin{equation}\label{upper bound of complement of X-n for strong colourings}
\big| \mbX_n \big| \ \leq \
2^{\beta n + \sum_{k \in [\rho] - I} \smult(n,\gamma_n,r_k) + 
\sum_{k \in [\rho]} \overline{\smult}(n,\gamma_n,r_k)}.
\end{equation}
From~(\ref{connection between smult for sets and tuples}), (\ref{lower bound for X-n for strong colourings})
and~(\ref{upper bound of complement of X-n for strong colourings}) it follows that
in order to prove (\ref{equation giving second part of theorem}) it suffices to
show that, for every $2 \leq m \leq m_{\max}$
there is a constant $\lambda_m > 0$ such that for all sufficiently large $n$,
\begin{equation}\label{comparing sigma-n and gamma-n for strong colourings}
\overline{\smult}(n,\sigma_n,m) - \overline{\smult}(n,\gamma_n,m) \ \geq \ 
\lambda_m n^m.
\end{equation}

\noindent
We have
\begin{align}\label{lower bound for sigma-n for strong colourings}
\overline{\smult}(n,\sigma_n,m) \ &= \ 
\sum_{1 \leq i_1 < \ldots < i_m \leq l} p(n,\sigma_n,i_1) \ldots p(n,\sigma_n,i_m) 
\quad \text{by (\ref{the expression for smult})}  \\
&\geq \ \binom{l}{m}\bigg( \frac{n}{l} - 1 \bigg)^m \quad \text{by the definition of $\sigma_n$} \nonumber \\
&= \ \binom{l}{m}\frac{n^m}{l^m} \ \pm \ O\big(n^{m-1}\big). \nonumber
\end{align}
By Lemma~\ref{upper bound for gamma in the strongly l-coloured case},
\begin{equation}\label{upper bound of gamma-n for strong colourings}
\overline{\smult}(n,\gamma_n,m) \ \leq \
\Bigg[\frac{1}{(l-1)^m} \binom{l-1}{m} \ + \ \frac{1}{a(l-1)^{m-1}} \binom{l-1}{m-1}\Bigg]n^m,
\end{equation}
where the left term within the parentheses vanishes if $m = l$.
From (\ref{lower bound for sigma-n for strong colourings}) 
and~(\ref{upper bound of gamma-n for strong colourings}) we get
\begin{align*}
&\overline{\smult}(n,\sigma_n,m) \ - \ \overline{\smult}(n,\gamma_n,m) \\
&\geq \ \Bigg[ \frac{1}{l^m}\binom{l}{m} \ - \ 
\frac{1}{(l-1)^m}\binom{l-1}{m} \ - \
\frac{1}{a(l-1)^{m-1}}\binom{l-1}{m-1} \Bigg] n^m \ \pm \
O\big(n^{m-1}\big).
\end{align*}
Observe that the rightmost term within the large parentheses 
can be made arbitrarily small by choosing $a$ large enough.
Thus, to prove (\ref{comparing sigma-n and gamma-n for strong colourings})
it suffices to show that
\begin{equation}\label{comparison of terms not involving a}
\frac{1}{l^m}\binom{l}{m} \ - \ \frac{1}{(l-1)^m}\binom{l-1}{m} \ > \ 0.
\end{equation}

\noindent
But~(\ref{comparison of terms not involving a}) holds because whenever
$2 \leq m \leq l$ we have
\begin{align*}
\frac{1}{l^m} \binom{l}{m} \ &= \ 
\frac{1}{l^m} \cdot \frac{\prod_{i=0}^{m-1}(l-i)}{m!} \ = \ 
\frac{1}{m!} \prod_{i=0}^{m-1} \frac{l-i}{l} \\ 
&= \ \frac{1}{m!} \prod_{i=0}^{m-1} \bigg(1 - \frac{i}{l}\bigg) \ > \ 
\frac{1}{m!} \prod_{i=0}^{m-1} \bigg(1 - \frac{i}{l-1}\bigg) \ = \ 
\frac{1}{(l-1)^m}\binom{l-1}{m}.
\end{align*}


\begin{thebibliography}{99}\label{the bibliography}



\bibitem{AH} G. Agnarsson, M. M. Halld\'{o}rsson, Strong colorings of hypergraphs,
in G. Persiano, R. Solis-Oba (Eds.),
{\em Approximation and Online Algorithms}, Lecture Notes in Computer Science 3351,
Springer Verlag (2005).


\bibitem{Aig} M. Aigner, {\em Combinatorial Theory},
Springer-Verlag (1997).

\bibitem{ABBM} N. Alon, J. Balogh, B. Bollob\'{a}s, R. Morris,
The structure of almost all graphs in a hereditary property,
{\em Journal of Combinatorial Theory, Series B},
Vol. 101 (2011) 85-110.

\bibitem{BBS08} J. Balogh, B. Bollob\'{a}s, M. Simonovits,
The typical structure of graphs without given excluded subgraphs,
{\em Random Structures and Algorithms},
Vol. 34 (2009) 305--318.

\bibitem{BBS10} J. Balogh, B. Bollob\'{a}s, M. Simonovits,
The fine structure of octahedron-free graphs,
{\em Journal of Combinatorial Theory, Series B},
Vol. 101 (2011) 67--84.

\bibitem{BM} J. Balogh, D. Mubayi,
Almost all triangle-free triple systems are tripartite,
preprint.

\bibitem{Ber} C. Berge, {\em Hypergraphs}, Elsevier Science Publ. (1989).

\bibitem{Boll} B. Bollob\'{a}s, {\em Random graphs} (Second Edition), 
Cambridge University Press (2001).

\bibitem{Bur-book} S. N. Burris, {\em Number Theoretic Density and Logical Limit Laws},
Mathematical Surveys and Monographs, Vol. 86,
American Mathematical Society (2001).


\bibitem{BY} S. Burris, K. Yeats, Sufficient conditions for labelled 0-1 laws,
{\em Discrete Mathematics \& Theoretical Computer Science}, Vol. 10 (2008) 147-156. 

\bibitem{Che92} G. Cherlin, Combinatorial Problems Connected with Finite Homogeneity, 
{\em Contemporary Mathematics} Vol. 131, 1992 (Part 3) 3--30.

\bibitem{CH} G. Cherlin, E. Hrushovski, {\em Finite Structures with Few Types}, Annals of Mathematics Studies 152, 
Princeton University Press, (2003).

\bibitem{Dj04} M. Djordjevic, On first-order sentences without finite models,
{\em The Journal of Symbolic Logic}, Vol. 69 (2004), 329-339.

\bibitem{Dj06a} M. Djordjevi\'{c}, The finite submodel property and $\omega$-categorical expansions of pregeometries,
{\em Annals of Pure and Applied Logic}, Vol. 139 (2006), 201-229.


\bibitem{EF} H-D. Ebbinghaus, J. Flum, {\em Finite Model Theory}, Second Edition, 
Springer Verlag, 1999.

\bibitem{Eva} D. M. Evans, $\aleph_0$-categorical structures with a predimension, 
{\em Annals of Pure and Applied Logic} 116 (2002) 157-186.


\bibitem{EKR} P. Erdös, D. J. Kleitman, B. L. Rothschild,
Asymptotic enuemeration of $K_n$-free graphs,
in {\em International Colloquium on Combinatorial Theory, Atti dei Convegni Lincei 17},
Vol. 2, Rome 1976, 19--27.


\bibitem{Fag} R. Fagin, Probabilities on finite models, {\em The Journal of Symbolic Logic}, 
Vol. 41 (1976) 50-58.

\bibitem{Gleb} Y. V Glebskii, D. I. Kogan, M. I. Liogonkii, V. A. Talanov, 
Volume and fraction of satisfiability of formulas of the lower predicate calculus, 
{\em Kibernetyka} Vol. 2 (1969) 17-27.

\bibitem{Had} G. F. Hadley, {\em Nonlinear and dynamic programming}, Addison-Wesley (1964).

\bibitem{HLP} G. H. Hardy, J. E. Littlewood, G. P\'{o}lya,
{\em Inequalities}, Cambridge University Press (1934).

\bibitem{Hen72} W. Henson, Countable homogeneous relational structures and
$\aleph_0$-categorical theories, {\em The Journal of Symbolic Logic}, 
Vol. 37 (1972) 494--500.

\bibitem{Hod} W. Hodges, {\em Model theory}, Cambridge University Press (1993).

\bibitem{Hru} E. Hrushovski, Simplicity and the Lascar group, manuscript (2003).

\bibitem{Imm} N. Immerman, Upper and lower bounds for first-order expressibility,
{\em Journal of Computer and Systems Sciences}, Vol. 25 (1982) 76--98.

\bibitem{JT} T.R. Jensen, B. Toft, {\em Graph Coloring Problems},
Wiley-Interscience (1995).

\bibitem{KPR} Ph. G. Kolaitis, H. J. Prömel, B. L. Rothschild,
$K_{l+1}$-free graphs: asymptotic structure and a 0-1 law,
{\em Transactions of The American Mathematical Society}, Vol. 303 (1987) 637--671. 

\bibitem{Mac} C. Maclaurin, A second letter to Martin Folkes, Esq.;
concerning the roots of equations, with the demonstration of other rules in algebra,
{\em Philosophical Transactions}, Vol. 36 (1729) 59--96.

\bibitem{Ober} W. Oberschelp, Asymptotic 0-1 laws in combinatorics,
In D. Jungnickel (editor), {\em Combinatorial theory, Lecture Notes in Mathematics},
Vol. 969 (1982), 276--292, Springer Verlag.

\bibitem{PeSch} Y. Person, M. Schacht, Almost all hypergraphs without
Fano planes are bipartite, in Claire Mathieu (editor),
{\em Proceedings of the Twentieth Annual ACM-SIAM
Symposium on Discrete Algorithms (SODA 09)},
217--226, ACM Press (2009).

\bibitem{Poi} B. Poizat, Deux ou trois choses que je sais de $L^n$,
{\em The Journal of Symbolic Logic}, Vol. 47 (1982) 641--658.

\bibitem{PS} H. J. Prömel, A. Steger, The asymptotic number of graphs not
containing a fixed color-critical subgraph, {\em Combinatorica},
Vol. 12 (1992) 463--473.

\bibitem{PS95} H. J. Prömel, A. Steger, Random $l$-colourable graphs,
{\em Random Structures and Algorithms}, Vol. 6 (1995) 21--37.

\bibitem{PST} H. J. Prömel, A. Steger, A. Taraz,
Asymptotic enumeration, global structure, and constrained evolution,
{\em Discrete Mathematics}, Vol. 229 (2001) 213--233.

\bibitem{Win} P. Winkler, Random structures and zero-one laws,
N. W. Sauer et al. (editors), {\em Finite and Infinite Combinatorics in Sets and Logic},
399--420, Kluwer Academic Publishers (1993).


\end{thebibliography}
\end{document}